\RequirePackage{fix-cm}
\documentclass[a4paper,twoside,tbtags,leqno]{article}
\usepackage{amsmath}
\usepackage{amsfonts}
\usepackage{amssymb} 
\usepackage{stackrel} 
\usepackage{latexsym}
\usepackage{epsfig}
\usepackage{graphicx}
\usepackage{oldgerm}
\usepackage{amsthm}%\usepackage{theorem}   
\newcommand\DN{\newcommand}

\numberwithin{equation}{section}
\newcounter{Const} 
\DN\Ct{\refstepcounter{Const}c_{\theConst}}%\label{#1}
\DN\cref[1]{c_{\ref{#1}}}	
 \theoremstyle{definition}
 \newtheorem{definition}{Definition}[section]

 \newtheorem{theorem}{Theorem}[section]
 \newtheorem{lemma}[theorem]{Lemma}
 \newtheorem{corollary}[theorem]{Corollary}
 \newtheorem{proposition}[theorem]{Proposition}
\newtheorem{remark}{Remark}[section]

\DN\lref[1]{Lemma~\ref{#1}}\DN\tref[1]{Theorem~\ref{#1}}
\DN\pref[1]{Proposition~\ref{#1}}\DN\sref[1]{Section~\ref{#1}}
\DN\ssref[1]{Subsection~\ref{#1}}\DN\dref[1]{Definition~\ref{#1}} 
\DN\rref[1]{Remark~\ref{#1}} \DN\corref[1]{Corollary~\ref{#1}}
\DN\eref[1]{Example~\ref{#1}}

\DN\bs{\bigskip}\DN\ms{\medskip}
\DN\N{\mathbb{N}}\DN\R{\mathbb{R}}%\DN\n{\mathsf{n}}
\DN\map[3]{#1\!:\!#2\!\to\!#3}\DN\ot{\otimes} \DN\ts{\times }
\DN\PD[2]{\frac{\partial#1}{\partial#2}} \DN\half{\frac{1}{2}}
\DN\Rtwo{\R ^2}
\DN\Rd{\R ^d}\DN\Rtwom{\R ^{2m}}\DN\RtwoN{\R ^{2\N }}
\DN\RdN{\R ^{d\N }} %\DN\RdN{(\R ^d)^{\mathbb{N}}}
\DN\limi[1]{\lim_{#1\to\infty}} \DN\limz[1]{\lim_{#1\to0}}
\DN\limsupi[1]{\limsup_{#1\to\infty}} 	\DN\liminfi[1]{\liminf_{#1\to\infty}} 	
\DN\elaw{\stackrel{\mathrm{law}}{=}}
\DN\As[1]{{\textbf{$ ($#1$ )$}}} %\DN\At[1]{\textbf{(#1)}}
\DN\Section{\section}\DN\Ssection{\subsection}\DN\SSsection{\subsubsection}
\DN\aaaaa{\noindent {\em Proof. }}\DN\bbbbb{ \qed \medskip }%\DN\qedproof{\qed \end{\proof}}

\DN\sinbeta{\mathrm{sin} _{\beta }}
\DN\sinone{\mathrm{sin} _{1 }}
\DN\sintwo{\mathrm{sin} _{2 }}
\DN\sinfour{\mathrm{sin} _{4 }}

\DN\partialsi{\partial_{s_i}}

\DN\HH{\mathsf{H}}

\DN\vR{\varphi _{\RR }}
\DN\vRt{\tilde{\varphi} _{\RR }}
\DN\RR{R}

\DN\simpr{\sim _{\pp ,\rr }}

\DN\uL{}%\DN\uL{}%
\DN\DOlaL{depends only on $ \la \circ \lab ^{-1}$}

\DN\DOMUL{depends only on $ \mu \circ \lab ^{-1}$}

\DN\thetazero{\Theta _{\mm ,\nn }}
\DN\thetaone{\theta _{\pp ,\rr }^1}
\DN\thetatwo{\theta _{\pp ,\rr }^2}
\DN\thetathree{\theta _{\mm ,\nn }^3}
\DN\thetafour{\theta _{\qq }^4}

\DN\Psiz{\Psi }
\DN\existence{existence} \DN\theexistence{the existence}
\DN\uniqueness{uniqueness} \DN\theuniqueness{the uniqueness}
\DN\pathwiseuniqueness{pathwise uniqueness} \DN\thepathwiseuniqueness{the pathwise uniqueness}

\DN\CsC{C_1 * C_2}
\DN\Qb{Q_{\mathbf{b}}}
\DN\Qbt{Q_{\mathbf{b}}^{[0,t]}}
\DN\sigmaconst{Assume that the coefficient $ \sigma^m $ is constant}

\DN\Pwh{\hat{P}_{w}}
\DN\Pw{P_{w}}
\DN\Pwt{P_{w,t}}
\DN\cmu{\chin ^2 d\mug ^{[1]}}
\DN\RA{\Rightarrow }
\DN\DOmu{depends only on $ \mu $}

\DN\fg{\widetilde{b}^{m}}%{f_{\mathrm{Gin}}}

\DN\fgsigma{\widetilde{\sigma}^{m}}%{f_{\mathrm{Gin}}}
\DN\fgfgsigma{\widetilde{\sigma}^{m}}%{f_{\mathrm{Gin}}}

\DN\SBhat{\{ \widetilde{\sigma }^m , \, \widetilde{ b }^m \}}
\DN\SB{(\widetilde{\sigma }^m , \, \widetilde{ b }^m )} 
\DN\nNpqk{$ \nn \in \NNNthree $ with $ \nn = (\pqr )$}
\DN\UV{(\mathbf{u},\mathsf{v})}

\DN\upathA{\upath (\mathbf{A})} 
\DN\lpathA{\lpath ^{-1}(\mathbf{A})} 
\DN\upathw{\{ \upath (\www )\, ;\, \mrT (\www ) < \infty \}}
\DN\lpathw{\{ \ww \, ;\, \mrT (\lpath (\ww ) ) < \infty \}}

\DN\ubx{(\uuu ,\bx ) }
\DN\vbx{(\vvv ,\bx ) }
\DN\bx{\mathbf{b},\mathsf{x}}

\DN\AAEFFl{\As{A1}--\As{A4}, \As{\Bone}, and \As{\CONE}--\As{\CTWO} for an $ \ell $}
\DN\AAEFF{\As{A1}--\As{A4}, \As{\Bone}, and \As{\CONE}--\As{\CTWO}} 
\DN\AAEE{\As{A1}--\As{A4} and \As{\Bone}--\As{\Btwo}}
\DN\AAcEE{\As{A1}--\As{A4}, and \As{\Bone}--\As{\Btwo}}
\DN\AAE{\As{A1}--\As{A4}, and \As{\Bone}}

\DN\AsDD{\As{\Done} and \As{\Dtwo}}

\DN\nAm{\As{\nbj} and \As{\muAC} for $ \mu $}
\DN\nNm{\As{\nbj} and \As{\muAC}}
\DN\NNMM{\As{\nbj}, \As{\muAC} and \As{\muTT} for $ \mu $}
\DN\NNMMg{\As{\nbj}, \As{\muAC} and \As{\muTT} for $ \mug $}
\DN\Inmg{\As{\iFc}, \As{\muAC} for $ \mug $, \As{\Cthree}, and \As{\nbj}}
\DN\Inmgs{\iFcs, \As{\muAC} for $ \mug $, \As{\Cthree}, and \As{\nbj}}
\DN\INAT{\As{\iFc}, \NNMM}

\DN\INA{\As{\iFc}, \As{\nbj}, and \As{\muAC} for $ \mu $}
\DN\INAa{\As{\iFc}, \As{\nbj}, and \As{\muAC} for $ \mut $}
\DN\NaN{\As{\sIn} and \As{\nbj}}

\DN\sIn{$\mathbf{SIN}$}
\DN\iFc{$\mathbf{IFC}$}
\DN\iFcs{\As{$\mathbf{IFC}$}$ _{\mathbf{s}}$}
\DN\nbj{$\mathbf{NBJ}$}
\DN\muAC{$\mathbf{AC}$}
\DN\muTT{$\mathbf{TT}$}

\DN\Cone{\muTT}
\DN\Ctwo{\muAC}
\DN\Cthree{\sIn}
\DN\Cfour{\nbj}
\DN\Conefour{\As{\Cone} and \As{\Ctwo} for $ \mu $, \As{\Cthree}, and \As{\Cfour}}
\DN\CCCC{\As{\iFc}, \As{\Ctwo} for $ \mu $, \As{\Cthree}, and \As{\Cfour}}
\DN\IASN{\As{\iFc}, \As{\Ctwo} for $ \mu $, \As{\Cthree}, and \As{\Cfour}}
\DN\IASNF{\As{\iFc}, \As{\Ctwo} for $ \mu $, \As{\Cthree},  \As{\Cfour}, and \As{$\mathbf{MF}$}}

\DN\IASNs{\iFcs, \As{\Ctwo} for $ \mu $, \As{\Cthree}, and \As{\Cfour}}
\DN\IASNa{\As{\iFc}, \As{\Ctwo} for $ \mut $, \As{\Cthree}, and \As{\Cfour}}
\DN\IASNaa{\As{\iFc}, \As{\Ctwo} for $ \mut $, \As{\Cthree}, and \As{\Cfour}}
\DN\ISNs{\iFcs, \As{\Cthree}, and \As{\Cfour}}
\DN\ISNss{\iFcs, \As{\Cthree}$_{\mathbf{s}}$, and \As{\Cfour}$_\mathbf{s}$}
\DN\ISN{\As{\iFc}, \As{\Cthree}, and \As{\Cfour}}

\DN\Conethree{\As{\Cone}, \As{\Ctwo}, and \As{\Cfour}} 

\DN\Dzero{Assume that $ \Pts $ is \9}

\DN\MST{Make the same assumption as \tref{l:79}}

\DN\ASTpath{\As{$ \mathbf{T}_{\mathrm{path}}1$}$ _{\mathbf{s}}$}
\DN\ASTpatH{\As{$ \mathbf{T}_{\mathrm{path}}2$}$ _{\mathbf{s}}$}

\DN\Done{$ \mathbf{C}_{\mathrm{path}}1$}
\DN\Dtwo{$ \mathbf{C}_{\mathrm{path}}2$}

\DN\Bone{B1}\DN\Btwo{B2}
\DN\CONE{C1}\DN\CTWO{C2}

\DN\XoneX{(\mathbf{X}^1,\mathsf{X}^{1*})}
\DN\YopX{\mathbf{Y}^{m}\oplus \mathbf{X}^{m*}}
\DN\YopXz{\mathbf{Y}_0^{m}\oplus \mathbf{X}_0^{m*}}

\DN\FFF{\Fsm \XzBXm \oplus \mathbf{X}^{m*}}

\DN\barF{\widetilde{F}}

\DN\FXBiopX{\Fsmi \XzBXm \oplus \mathbf{X}^{m*}}

\DN\Ft{\{ \mathcal{F}_t \}}
\DN\FtB{$\{ \mathcal{F}_t \}$-Brownian motion }
\DN\Fm{F^m}
\DN\Fms{F_{\mathbf{s}}^m} %12/11変更 sfからbfへ%\DN\Fms{F_{\labm (\mathsf{s})}^m}
\DN\Fmss{F_{\mathbf{s}}^m}% \rref{r:60z}の直前\As{\iFc}で使用。 
\DN\Fi{F_{\mathbf{s}}^{\infty ,i}}

\DN\lBlhatm{\mathbf{B}^m,\mathbf{X}^{m*}}

\DN\lBl{\lab (\mathsf{X} _0),\mathbf{B}(\mathsf{X}),\lpathX }

\DN\lBlmM{\mathbf{B}^m ,\mathbf{X}^{m*}}

\DN\Stwoq{\sS _{\rr }}%\DN\Stwoq{\Rtwo _{\rr }}
\DN\Stwoqm{\Sm _{\rr } }%\DN\Stwoqm{(\Rtwo _{\rr })^m }
\DN\Rtowm{(\Rtwo )^m}

\DN\bB{\mathbf{b}} %\DR\b{\uL{\mathbf{b}}}
\DN\Wb{W_{\bB }}
\DN\B{\uL{\mathbf{B}}}
\DN\zt{\zeta}
\DN\U{U}\DN\V{V}
\DN\UU{U} \DN\VV{V}
\DN\ol[1]{\overline{#1 }}%\DN\ol[1]{\overline{#1 }^{\Pts }}

\DN\Blm{(\hat{\mathbf{B}}^m(\mathsf{X}),\lpathX )}
\DN\PPPm{\Pt ^m }
\DN\uuu{\mathbf{u}}
\DN\vvv{\mathbf{v}}

\DN\wti{\w _t^i}
\DN\w{\mathit{w}}
\DN\wt{\mathit{w}_t}
\DN\ww{\mathsf{w}}
\DN\wwt{\ww _t}
\DN\www{\mathbf{w}}\DN\wwwt{\www _t}\DN\wtt{\www (t)}
\DN\wwww{\uL{\mathbf{w}}}
\DN\wW{\uL{\ww }}
\DN\la{\lambda }

\DN\laa{\lambda ^{\aa }}
\DN\Plaa{P _{\laa }}

\DN\lal{\la \circ {\lab }^{-1}}
\DN\LaL{\la \circ {\lab }^{-1}}

\DN\Pl{\PP _{\la }}

\DN\Pll{\uL{P}} 

\DN\Plp{\Pl \circ \lpath ^{-1}}

\DN\mul{\mu \circ {\lab }^{-1}}

\DN\Pml{\mathbf{P}_{\mul } }

\DN\zti{ 0 \le t < \infty }\DN\zzti{ 0 < t < \infty }
\DN\bj{big jump\ }

\DN\Vnu{V_{\nu }}
\DN\Eq{Eq. \!\!}
\DN\A{\mathsf{A}}

\DN\y{\mathbf{y}}

\DN\pirc{\pi _r^c}
\DN\pisc{\pi _s^c}
\DN\pitc{\pi _t^c}
\DN\pir{\pi _r}
\DN\pis{\pi _s}

 \DN\St{\sSS ^{\bft}}
 \DN\Su{\sSS ^{\bfu}}
 \DN\Sv{\sSS ^{\bfv}}
 \DN\Sw{\sSS ^{\bfw}}

\DN\bft{\mathbf{t}}
\DN\bfu{\mathbf{u}}
\DN\bfv{\mathbf{v}}
\DN\bfw{(\bfu ,\bfv )} 
\DN\sfx{\mathsf{x}}
\DN\sfy{\mathsf{y}}
\DN\sfz{\mathsf{z}}

\DN\ER{\mathcal{R} }
\DN\NNN{\mathsf{N}}
\DN\NNNthree{\NNN _3 }
\DN\nnNNN{\nn \in \NNN }

\DN\IFC{IFC}
\DN\ISDE{\eqref{:50a}--\eqref{:50c}}
\DN\ISDEb{\text{\eqref{:50a}--\eqref{:50b}}\ }
\DN\ISDEbb{\text{\eqref{:50a}--\eqref{:50b}}}

\DN\ISDEsixB{\eqref{:60a}--\eqref{:60b}}

\DN\ISDEsixb{\ISDEsixc}
\DN\ISDEsixbb{\ISDEsixcc}
\DN\ISDEsixc{\eqref{:60a}--\eqref{:60c}\ }
\DN\ISDEsixcc{\eqref{:60a}--\eqref{:60c}}
\DN\ISDEhij{\eqref{:60h}--\eqref{:60j}\ }
\DN\ISDEhijj{\eqref{:60h}--\eqref{:60j}}

\DN\nn{\mathsf{n}} \DN\mm{\mathsf{m}}
\DN\vq{\varphi _{\rr }}
\DN\vr{\varphi _{\RR }}

\DN\chin{\chi _{\nn }}
\DN\chinn{\chi _{\nn +1}}
\DN\az{a_{0}}
\DN\ak{a_{\qq }}
\DN\akk{a_{\qq +1}}
\DN\Ka{\mathsf{K}[\mathbf{a}]}
\DN\Kaa{\mathsf{K}[\mathbf{a}^+]}
\DN\Kak{\mathsf{K}[\mathit{a}_{\qq }]}
\DN\Kakk{\mathsf{K}[\mathit{a}_{\qq }^+]}
\DN\Ha{\mathsf{H}[\mathbf{a}]}

\DN\Hb{\Ha _{\pqr }}
\DN\Hak{\mathsf{H}[\mathbf{a}]_{\qq }}

\DN\Han{\Ha _{\nn }}
\DN\Hann{\Ha _{\nn +1}}

\DN\HanC{\Han ^{\circ }}
\DN\Hns{\Ha _{\nn , \mathsf{s}}}

\DN\RRpr{\mathbb{R}_{\pp ,\rr }}
\DN\RRprC{\RRpr ^{\circ }} 
\DN\RRprCs{\RRprC (\mathsf{s})}
\DN\RRprs{\RRpr (\mathsf{s})}

\DN\HHns{\mathbb{H}_{\nn }(\mathsf{s})}
\DN\Hnsm{\mathbb{R}_{\pp ,\rr , \mm }}

\DN\Hmn{\Han \cap \IIm }
\DN\Hmnc{\HanC \cap \IIm }
\DN\HI{\HanC \backslash \IIm }

\DN\LLL{\langle}\DN\rRRR{\rangle}
\DN\HImn{\LLL \Hmn \rRRR }
\DN\HImnC{\LLL \Hmnc \rRRR } 
 \DN\HIiC{ \HanC \backslash \LLL \Hmnc \rRRR }

\DN\HmnP{\Pitwo (\Hmn )} \DN\HmnPc{\Pitwo (\Hmnc )}

\DN\II{\mathsf{I}}
\DN\IIm{\II _{\mm }} %\DN\IImn{\II _{\nn }(\mm )}
\DN\IImm{\II _{\mm +1}} %\DN\IImmn{\II _{\nn }(\mm +1)}
\DN\IIi{\II _{\infty }}

\DN\ff{\mathsf{f}} \DN\bm{b^m}

\DN\YmX{(\mathbf{Y}^{m} , \mathbf{X}^{m*})}
\DN\XmX{(\mathbf{X}^{m} , \mathbf{X}^{m*})}
\DN\Xms{ \mathbf{X}^{m*} }
\DN\xm{_{\mathsf{X}}^m}
\DN\xmhat{_{\hat{\mathsf{X}}}^m}

\DN\xxx{\mathbf{x}}
\DN\xxxm{\mathbf{x}^m}

\DN{\RRR}{\As{A1}--\As{\muTT}}
\DN{\RRRRRR}{\As{A1}--\As{\nbj}}
\DN{\RRRRRRs}{\As{A1}$_ \aa $--\As{\nbj}$_ \aa $}

\DN\rhog{\rho _{\mathrm{Gin}}}
\DN\rhogx{\rho _{\mathrm{Gin},x}}
\DN\kg{\mathsf{K}_{\mathrm{Gin}}}
\DN\mugone{\mug ^{[1]}}
\DN\mug{\mu _{\mathrm{Gin}}}
\DN\mugx{\mu _{\mathrm{Gin},x}}
\DN\mugxx{\mu _{\mathrm{Gin},\mathbf{x}}}
\DN\muSinb{\muone _{\sinbeta }}
\DN\muAi{\mu _{\mathrm{Ai}}}

\DN\mub{\mu _{\mathrm{Be},\alpha }}
\DN\muone{\mu ^{[1]}}
\DN\mum{\mu ^{[m]}}

\DN\rairytwo{\rho _{\mathrm{Ai},2}}
\DN\rhohat{\hat{\rho}}

\DN\sS{S} \DN\sSS{\mathsf{S}}
\DN\SSS{\mathbf{S}}

\DN\Sm{\sS ^m} \DN\SmSS{\Sm \ts \sSS }
\DN\Srm{\Sr ^m} \DN\SrmSS{\Srm \ts \sSS }

\DN\SSsdeg{\mathsf{S}_{\mathrm{sde}}}
\DN\SSsde{\mathsf{S}_{\mathrm{sde}}}
\DN\SSSsde{\mathbf{S}_{\mathrm{sde}}}
\DN\SSSsdeone{\mathsf{S}_{\mathrm{sde}}^{[1]}}

\DN\SSsdemt{\mathsf{S}_{\mathrm{sde}}^m(t,\mathsf{X})}
\DN\SSsdemtw{\mathsf{S}_{\mathrm{sde}}^m(t,\ww )}
\DN\SSSsdemt{\mathbf{S}_{\mathrm{sde}}^m(t,\mathsf{X})}
\DN\SSSsdemtg{\mathbf{S}_{\mathrm{sde}}^m(t,\mathsf{X})}
\DN\SSSsdemtw{\mathbf{S}_{\mathrm{sde}}^m(t,\ww )}

\DN\SSrm{\sSS _{r}^m }
\DN\SSrNk{\sSS _{r,N-k}}
\DN\SSk{\sSS ^{[k]}}
\DN\SSksingle{\Ssi ^{k}}
\DN\Ssi{\sSS _{\mathrm{s.i.}}}
\DN\Ss{\sSS _{\mathrm{s}}}
\DN\Ssip{\sSS _{\mathrm{s.i.}}^{+f}}
\DN\SSone{\sSS ^{\mathbf{1}}}
\DN\SoneSSS{\sS \ts \sSS }

\DN\Sk{\sS ^{k}}
\DN\Srr{\sS _{\rr }}

\DN\Srk{\Sr ^{k}}
\DN\Sr{\sS _{r}}

\DN\SN{\sS ^{\mathbb{N}}} 

\DN\HHa{\mathsf{H}_{\mathsf{a}}}
\DN\HHz{\mathsf{H}} 

\DN\dlog{\mathsf{d}}
\DN\dmu{\dlog ^{\mu }}
\DN\dpsi{\dlog ^{\mupsi }}
\DN\dmuone{\dlog ^{\mu ^{1}}}
\DN\dmuhat{\hat{\dlog } ^{\mu }}

\DN\dG{\dlog ^{\mathrm{Gin} }}%\DN\dG{\dlog ^{\mug }}
\DN\dSine{\dlog ^{\mu _{\sinbeta }} }

\DN\sigmaxms{\sigma _{\mathsf{x}}^m }
\DN\bbbxms{\mathit{b} _{\mathsf{x}}^m }
\DN\sigmaXms{\sigma \xm }
\DN\bbbXms{\mathit{b} \xm }
\DN\bbbXmshat{\mathit{b} \xmhat }
\DN\bbb{\mathit{b}}
\DN\aaa{\mathit{a}}
\DN\bb{\mathsf{b}}
\renewcommand\aa{\mathsf{a}}

\DN\dom{\mathcal{D}}
\DN\Damu{\mathcal{D}^{\amu }} 
\DN\di{\dom _{\circ }}

\DN\Lmu{L^2(\sSS ,\mu )}
\DN\Lmg{L^{2}(\sSS ,\mug )}
 \DN\Llocp{L_{\mathrm{loc}}^{p}}
 \DN\Lloctwo{L_{\mathrm{loc}}^{2}}

\DN\LlocG{\Lloctwo (\muone _{\mathrm{Gin}})}

\DN\LlocpG{\Llocp (\muone _{\mathrm{Gin}})}

\DN\LchiG{L^2 (\chin ^2\muone _{\mathrm{Gin}})}
\DN\Lc{L^2 (\chin ^2 \muone _{\mathrm{Gin}})}

\DN\LlocSine{\Lloctwo (\muone _{\sinbeta })}
\DN\LchiSine{L ^2 (\chin \muone _{\sinbeta })}
\DN\LlocpSine{\Llocp (\muone _{\sinbeta })}

\DN\Eamu{\mathcal{E}^{\amu }}
\DN\Emu{\mathcal{E}^{\mu } }
\DN\E{\mathcal{E} }

\DN\amu{a ,\mu }
\DN\amuone{a ,\muone }
\DN\DDD{\mathbb{D}}
\DN\DDDa{\mathbb{D}^{a}}
\DN\DDDr{\DDD _{r}}
\DN\D{\mathbf{D}}
\DN\PP{\mathsf{P}}
\DN\PPa{\PP ^{\aa }}
\DN\PPsa{\PPs ^{\aa }}
\DN\amut{a,\mut }

\DN\ulab{\mathfrak{u} }
\DN\lab{\mathfrak{l} } 
\DN\labm{\lab ^m } 
\DN\labi{\lab _i}

\DN\ulabm{\mathfrak{u} _{[m]}}
\DN\ulabone{\mathfrak{u} _{[1]}}
\DN\upath{\ulab _{\mathrm{path}}}
\DN\lpath{\lab _{\mathrm{path}}} 
\DN\lpathmm{\lpath ^m} 
\DN\lpathmstar{\lpath ^{m* }} 
\DN\lpathmstarM{\lpath ^{(m*) }} 
\DN\lpathi{\lpath ^i} 
\DN\lpathX{\lpath (\mathsf{X}) }

\DN\OFPsF{(\Omega ,\mathcal{F},\{ \PPs \}, \{ \mathcal{F}_t \} )}
\DN\OF{(\Omega , \mathcal{F} )}
\DN\OFF{(\Omega ,\mathcal{F},\{ \mathcal{F}_t \} )}
\DN\OFP{(\Omega ,\mathcal{F}, P )}

\DN\OFPF{(\Omega ,\mathcal{F}, P ,\{ \mathcal{F}_t \} )}
\DN\OFPFF{(\Omega ',\mathcal{F}', \pP ',\{ \mathcal{F}_t' \} )}

\DN\2{(\Omega ,\mathcal{F}, \Ps , \{ \mathcal{F}_t \} )}

\DN\OFPFts{(\Omega ,\mathcal{F}, \Pts ,\{ \mathcal{F}_t \} )}
\DN\OFPFmg{(\Omega ,\mathcal{F}, \PP _{\mug } , \{ \mathcal{F}_t \} )}
\DN\OFPFs{(\Omega ,\mathcal{F}, \PPs , \{ \mathcal{F}_t \} )}

\DN\uPs{under $ \Ps $}

\DN\qQ{Q}
\DN\qQs{Q_{\mathbf{s}}}
\DN\QQ{\mathsf{Q}} 
\DN\QQs{\QQ _{\mathsf{s}}} 
\DN\QQmg{\QQ _{\mug }} 
\DN\OFQFs{(\Omega ,\mathcal{F}, \QQs , \{ \mathcal{F}_t \} )}
\DN\OFQF{(\Omega ,\mathcal{F}, \QQ , \{ \mathcal{F}_t \} )}

\DN\OFqF{(\Omega ,\mathcal{F}, \qQ , \{ \mathcal{F}_t \} )}
\DN\OFqFs{(\Omega ,\mathcal{F}, \qQs , \{ \mathcal{F}_t \} )}
\DN\Upsilont{\widetilde{\Upsilon }}
\DN\PmglB{\Pmg \circ (\lab (\mathsf{X}_0),\mathbf{B}(\mathsf{X}))^{-1}}

\DN\XlP{$ \mathbf{X}=\lpathX $ under $ \Pm $}
\DN\1{\sum_{p=1}^2 (x_{j+1,p}-x_{j,p})\partial_p \partial_k \widetilde{\dG _l } (t (x_{j+1}-x_{j}) + x_{j} ,\mathsf{s}) }

\DN\3{\eqref{:60h}--\eqref{:60j}}
\DN\ginSDE{\eqref{:35c}--\eqref{:35e}}

\DN\4{\sum_{p=1}^2 (x_{j+1,p}-x_{j,p})\partial_p \widetilde{\dG _l } (t (x_{j+1}-x_{j}) + x_{j} ,\mathsf{s}) }

\DN\5{\tTTT _{\mathrm{path}}^{\{1\}}}

\DN\6{Assume that $ \XB $ under $ \Ps $ is \7.}
\DN\7{a weak solution of \ISDEsixb satisfying \iFcs}
\DN\8{weak solutions of \ISDEsixb satisfying \iFcs}
\DN\9{a weak solution of \ISDEsixb satisfying \iFcs\ for $ \PXz $-a.\! s.\! $ \mathbf{s}$}

\DN\mac{\As{$ \mu $-AC}}
\DN\mgac{\As{$ \mug $-AC}}

\DN\PPPP{\PsB \dsbb }
\DN\mB{$ \Upsilon $-a.s.\! $ (\mathbf{s},\bB ) $}

\DN\PPP{\Pts \circ \Pitwo ^{-1}}

\DN\pP{P}
\DN\pPm{\pP _{\mu }}
\DN\pPs{\Pts }
\DN\pPF{\pP _{\{ \Fs \}}}
\DN\pPFa{\pP _{\{ \Fsa \}}}

\DN\Ehat{\mathcal{C}}%\DN\Ehat{\widehat{\mathcal{C}}}
\DN\Ehatm{\Ehat^{m}}
\DN\Ehatmt{\Ehatm _t}

\DN\Pitwo{\Pi }
\DN\Pione{\Pi _1}

\DN\muairy{\mu_{\mathrm{Ai}, 2}}
\DN\muairyone{\mu_{\mathrm{Ai}, 1}}
\DN\muairyfour{\mu_{\mathrm{Ai}, 4}}
\DN\muairybeta{\mu_{\mathrm{Ai}, \beta }}

\DN\muairyk{\mu_{\mathrm{Ai}, 2}^{[k]}}
\DN\muairyonek{\mu_{\mathrm{Ai}, 1}^{[k]}}
\DN\muairyfourk{\mu_{\mathrm{Ai}, 4}^{[k]}}
\DN\muairybetak{\muairybeta ^{[k]}}
\DN\muairybetax{\mu_{\mathrm{Ai}, \beta ,x}}

\DN\musin{\mu_{\sintwo }}
\DN\musinone{\mu_{\sinone }}
\DN\musinfour{\mu_{\sinfour }}
\DN\musinbeta{\mu_{\sinbeta }}

\DN\HHm{\HHgm }
\DN\HHgm{\mathbf{H}^m} %Deleted see memo

\DN\uH{\ulab ^{-1} (\mathsf{H}) }
\DN\lH{\lab (\mathsf{H}) }

\DN\F{\uL{F}}%\DN\F{\mathbb{F}}
\DN\FFFs{{\Fs }}
\DN\Fs{\F _{\mathbf{s}}}
\DN\Fsa{\Fs ^{\aa }}

\DN\Fst{\Fs ^{[0,t]}}
\DN\Fsit{\Fs ^{ \infty , [0,t]}}

\DN\FopX{\Fsmbar \XB }

\DN\Fsmbar{\barF _{\mathbf{s}}^m}
\DN\Fsm{\uL{\Fms }}%\DN\Fsm{\Fs ^{m}} 
\DN\Fsmi{\uL{F _{\mathbf{s}} ^{m,i}}}%\DN\Fsmi{\Fs ^{m,i}} 
\DN\Fsn{\barF _{\mathbf{s}}^{n}} %\DN\Fsn{\Fs ^{n}} 
\DN\Fsi{\uL{F _{\mathbf{s}}^{\infty}}} %\DN\Fsi{\F _{\mathbf{s}}^{\infty}} 
\DN\sB{\sbb }%%\DN\sB{\mathbf{s},\mathbf{B}} %%試験的に
\DN\sbb{\mathbf{s},\bB } 
\DN\dsbb{\uL{(d\mathbf{s}d\bB )}}
\DN\sBc{_{\sbb } }%%\DN\sBc{_{\mathbf{s},\mathbf{B}} }%%試験的に
\DN\FisB{\uL{F _{\sbb }^{\infty}}} %% \DN\FisB{\F _{\sB }^{\infty}} %%試験的に
\DN\Fsbi{F_{\sB }^{\infty}}

\DN\WT{W } 				%%\DN\WT{C_{T}}1/8hennkou %\DN\WT{W_{\mathrm{path}}}
\DN\WS{\WT (\sSS )}			%[] \DN\WS{C([0,\infty);\sSS )}%\DN\WTS{ \WT (\sSS )}
\DN\WSsi{\WT (\Ssi )}
 \DN\WSs{\WT (\Ss )}

\DN\WSN{\WT (\SN )} 			%[] \DN\WSN{C([0,\infty);\SN )}
\DN\WSz{\WT _{\mathbf{0}}(\SN )} 	%[] \DN\WSz{\WSN _{\mathbf{0}}}
\DN\WSNs{\WT ^{\mathbf{s}}(\SN )}
\DN\WSsiNE{\WT _{\mathrm{NE}}(\Ssi )}
\DN\Wsde{\WT (\SSsde )}

\DN\WRdm{\WT (\Rdm )}
\DN\WRNz{\WT _{\mathbf{0}} (\RdN )}
\DN\WRN{\WT (\RdN )} 

\DN\WRtzm{\WT _{\mathbf{0}} (\Rtm )}
\DN\WRtmstar{\WT (\R ^{m*} )}
\DN\WRdzm{\WT _{\mathbf{0}} (\Rdm )}
\DN\WRmstar{\WT (\sS ^{m*} )}

\DN\WWW{\SN \ts \WW }
\DN\WWtm{\WRtzm \ts \WT (\RtN ) }
\DN\WWdm{\WRdzm \ts \WRN }

\DN\Rtm{\R ^{2m}}
\DN\Rtmstar{\R ^{2m*}}
\DN\RtN{\R ^{2\N }}
\DN\Rdm{\R ^{dm}}

\DN\WTRNm{ \WT (\Rtm )} 
\DN\WTzt{ \WT _{\mathbf{0},t}}

\DN\WSm{ \WT (\sS ^m )} 

\DN\Wsol{\mathbf{W} _{\mathrm{sol}}}
\DN\Wsols{\Wsol }
\DN\WWWsol{\mathbb{W} ^{\mathrm{sol}}} 
\DN\WWfix{\mathbf{W} ^{\mathrm{fix}}} 
\DN\WWWfix{\mathbb{W} ^{\mathrm{fix}}} 
\DN\WWfixsB{\mathbf{W} _{\sbb }^{\mathrm{fix}}}

\DN\Bbar{\mathcal{B}(\WSN )_{\sbb }}
\DN\Bt{\mathcal{B}_t }
\DN\Btm{\mathcal{B}_t^m }
\DN\Bot{\mathcal{B}_{\mathbf{o},t}}
\DN\Bw{\widetilde{\mathcal{B} }}

\DN\TTT{\mathrm{Tail}}
\DN\tTTT{\widehat{\mathcal{T}}} 
\DN\Tail{\mathcal{T}} %%6/20変更
\DN\TS{\Tail \, (\sSS )}
\DN\TSsi{\Tail (\Ssi )}

\DN\tTpath{\tTTT _{\mathrm{path}}}

\DN\Tpath{\mathcal{T} _{\mathrm{path}} }
\DN\Tpathone{\Tpath ^{\{1\}} }

\DN\TpathTSN{\Tpath (\SN )}

\DN\TpathTSNone{\Tpathone (\SN }
\DN\tTpathTSN{\tTTT _{\mathrm{path}} (\SN )}
\DN\tTpathTSNm{\tTTT _{\mathrm{path}}^m (\SN )}
\DN\tTpathTSNone{\tTTT _{\mathrm{path}}^{\{1\}} (\SN ;}
\DN\TpathSone{\tTTT _{\mathrm{path}}^{\{1\}} (\sSS ; }

\DN\TpSNP{\overline{\Tpath (\SN )}^{P}}
\DN\TpSNsB{\Tpath (\SN )_{\sbb }} 
\DN\TpSNsBB{\Tpath (\SN )_{\sbb }'}

\DN\TpathS{\tTpath (\sSS )}
\DN\TSone{\Tail ^{\{1\}}(\sSS ; }
\DN\TpathSsi{\tTpath (\Ssi )}

\DN\mupsi{\mu _{\Psiz }}
\DN\mupsitaila{\mu _{\Psiz ,\mathrm{Tail}}^{\mathsf{a}}}
\DN\musl{\mu _{\mathsf{s}}^{\lab }}
\DN\mmusl{(\mu _{\mathsf{s}})^{\lab }}
\DN\mut{\mu _{\mathrm{Tail}}^{\aa }}
\DN\mutl{\mu _{\mathrm{Tail}}^{\aa }^{\lab }}
\DN\mmutl{(\mu _{\mathrm{Tail}}^{\aa })^{\lab }}

\DN\PPPs{\mathbb{P}_{\mathbf{s}}}
\DN\PPPsb{\mathbb{P}_{\mathbf{s},\bB }}
\DN\PX{\pP \circ \mathbf{X}^{-1} }
\DN\PXz{\pP \circ \mathbf{X}_0^{-1} }
\DN\PsX{\Ps \circ \mathbf{X}^{-1} }

\DN\PPs{\PP _{\mathsf{s}}}
\DN\Pt{\widetilde{P}} 
\DN\Ptm{\widetilde{P}_m} 
\DN\Pts{\Pt _{\mathbf{s}}}
\DN\PsB{\Pt _{\sbb }}%%DN\PsB{\Pt _{\sB }}%%試験的に
\DN\Ptsb{\Pt _{\sbb }}
\DN\PsBg{\PsB ^{\mathrm{Gin}}}

\DN\Pss{\mathbf{P} _{\mathbf{s}}}
\DN\Ps{P _{\mathbf{s}}}
\DN\PBr{P _{\mathrm{Br}}^{\infty}} \DN\PBrm{P _{\mathrm{Br}}^{m}}
\DN\Pssf{\PPs }
\DN\Pmulg{\mathbf{P}_{\mug \circ {\lab }^{-1} }} 
\DN\Pm{\mathsf{P}_{\mu }}
\DN\Pmg{\mathsf{P}_{\mug }}
\DN\Ptmg{\Pt _{\mug }}
\DN\Ptmgm{\Pt _{\mug }^m}
\DN\Ptztm{\Ptzt ^m}
\DN\Pmut{\mathsf{P}_{\mut }}

\DN\PmT{\Pl ^{\wW _{\mathbf{t}}}} %保存すること：戻すかもしれない
 \DN\PmU{\Pl ^{\wW _{\mathbf{u}}}}
 \DN\PmV{\Pl ^{\wW _{\mathbf{v}}}}
 \DN\PmW{\Pl ^{\wW _{\mathbf{w}}}}

\DN\hp{h _{\pp }}
\DN\XB{(\mathbf{X},\mathbf{B})}
\DN\Xt{\mathbf{X}^{[0,t]}}
\DN\XBt{(\Xt ,\mathbf{B}^{[0,t]})}

\DN\wb{(\wwww ,\bB )}
\DN\wsb{\wwww (\sbb )}
\DN\sXB{(\mathbf{s},\mathbf{B},\mathbf{X})}
\DN\XzBXm{(\mathbf{B}^m,\mathbf{X}^{m*})}
\DN\BX{(\mathbf{B},\mathbf{X})}
\DN\cB{(\cdot ,\mathbf{B})}
\DN\vw{(\mathbf{w},\bB ) }
\DN\sBX{(\mathbf{s},\mathbf{B},\mathbf{X})} %%試験的に
\DN\sbw{(\mathbf{s},\bB ,\wwww )}
\DN\FsiXB{\Fsi \XB }
\DN\Fsiwb{\Fsi \wb }

\DN\WW{\WSN \times \WRNz } 
\DN\WWm{\WSm \times \WRNz } 
\DN\Swc{\Wsols \times \WRNz }

\DN\mrX{\mathsf{m}_{r,T} (\mathsf{X})}
\DN\mrW{\mathsf{m}_{r,T} (\mathsf{w})}
\DN\mrXX{\mathsf{m}_{r,T} (\mathbf{X})}
\DN\mrXl{\mathsf{m}_{r,T} (\lpathX )}
\DN\mrw{\uL{\mrT (\lpath (\ww ))}}

\DN\mrT{\mathsf{m}_{r,T} }
\DN\MrT{\mathsf{M}_{r,T}}

\DN\q{q}

\DN\pq{\pp ,\rr }\DN\pqr{\pp ,\qq ,\rr } 

\DN\pp{\mathsf{p}}\DN\qq{\mathsf{q}}
\DN\rr{\mathsf{r}}
\DN\ttt{t}

\DN\jjjj{\mathbf{j}}
\DN\ellell{\ell }

\DN\g{\mathsf{g}}\DN\h{\mathsf{h}}
\DN\ijn{_{i,j=1}^{n}}

\DN\Ito{It$ \mathrm{\hat{o}}$}

\begin{document}
\title{\textsf{Infinite-dimensional stochastic differential equations \\and tail $ \sigma$-fields }
\\
{\small \texttt{To appear in Probability Theory and Related Fields}}
}

\author{ Hirofumi Osada$^\dagger $ \and Hideki Tanemura$ ^*$}

\maketitle

\noindent 
$ ^\dagger$\\
Faculty of Mathematics, Kyushu University, Fukuoka 819-0395, Japan, 
\\Tel.: +81-92-802-4489,  Fax: +81-92-802-4405 \\
\texttt{email: osada@math.kyushu-u.ac.jp}

\ms 

\noindent 
$ ^*$\\
Department of Mathematics, Keio university, Kohoku-ku, Yokohama 223-8522, Japan. 
\\
\texttt{tanemura@math.keio.ac.jp}

\begin{abstract} 
We present general theorems solving the long-standing problem of the existence and pathwise uniqueness of strong solutions of infinite-dimensional stochastic differential equations (ISDEs) called interacting Brownian motions. 
These ISDEs describe the dynamics of infinitely-many Brownian particles 
moving in $ \mathbb{R}^d $ with free potential $ \Phi $ and 
mutual interaction potential $ \Psi $. 

We apply the theorems to essentially all interaction potentials of Ruelle's class such as the Lennard-Jones 6-12 potential and Riesz potentials, and to logarithmic potentials appearing in random matrix theory. 
We solve ISDEs of the Ginibre interacting Brownian motion and the sine$ _{\beta}$ interacting Brownian motion with $ \beta = 1,2,4$. 
We also use the theorems in separate papers for the Airy and Bessel interacting Brownian motions. 
One of the critical points for proving the general theorems is to establish a new formulation of solutions of ISDEs in terms of tail $ \sigma $-fields of labeled path spaces consisting of trajectories of infinitely-many particles. 
These formulations are equivalent to the original notions of solutions of ISDEs, and more feasible to treat in infinite dimensions. 
\\
\textbf{keywords}
\DN\and{: }
{interacting Brownian motions \and infinite-dimensional stochastic differential equations \and random matrices \and strong solutions \and pathwise uniqueness}
\end{abstract}

%%[] ? \paragraph{Paragraph headings} Use paragraph headings as needed.

{
\tableofcontents 
}

\section{Introduction.}\label{s:1} 
We study infinite-dimensional stochastic differential equations (ISDEs) of 
\begin{align*}&
\quad \quad \quad \mathbf{X}=(X^i)_{i\in\N } \in 
C([0,\infty);\RdN ) ,
\quad \text{ where } \RdN = (\Rd )^{\N }, 
\end{align*}
describing infinitely-many Brownian particles moving in $ \Rd $ 
with free potential $ \Phi = \Phi (x)$ and interaction potential 
$ \Psi = \Psi (x,y) $. %, where $ \RdN = (\Rd )^{\N }$. 
The ISDEs of $ \mathbf{X}=(X^i)_{i\in\N } $ are given by 
\begin{align}\label{:11a}&
dX_t^i = d B_t^i - \frac{\beta }{2}\nabla_x \Phi (X_t^i) dt - 
 \frac{\beta }{2} \sum_{j\not=i}^{\infty} \nabla_x \Psi (X_t^i,X_t^j) dt 
\quad (i \in \N )
.\end{align}
Here $ \mathbf{B}=(B^i)_{i=1}^{\infty}$, $ \{B^i\}_{i\in\N }$ are independent copies of $ d $-dimensional Brownian motions, 
$ \nabla_x = (\PD {}{x_i} )_{i=1}^d$ is the nabla in $ x $, and 
$ \beta $ is a positive constant called inverse temperature. 
The process $ (\mathbf{X},\mathbf{B})$ is defined on a filtered probability space 
$ (\Omega ,\mathcal{F},P,\{ \mathcal{F}_t \} )$.

The study of ISDEs was initiated by Lang \cite{lang.1,lang.2}, and continued by Shiga \cite{shiga}, 
Fritz \cite{Fr}, and the second author \cite{T2}, and others. 
In their respective work, the free potential $ \Phi $ is assumed to be zero and 
 interaction potentials $\Psi $ are of class $ C_0^3 (\Rd )$ 
or exponentially decay at infinity. 
This restriction on $ \Psi $ excludes polynomial decay and logarithmic growth 
interaction potentials, which are essential from the viewpoint of statistical physics and random matrix theory. 
The following are examples of such ISDEs. 

\noindent {\bf Sine$ _{\beta }$ interacting Brownian motions (\sref{s:V1}): } 
 \\
Let $ d= 1 $, $ \Phi (x) = 0 $, $ \Psi (x,y) = - \log|x-y|$. %, and $ \beta > 0 $. 
We set 
\begin{align}\label{:11b}&
dX_t^i = dB_t^i + \frac{\beta }{2}\limi{r} \sum_{|X_t^i-X_t^j|<r,\ j\not=i } 
\frac{1}{X_t^i-X_t^j} dt \quad (i\in\mathbb{N})
.\end{align}
This ISDE with $ \beta = 2 $ is called the Dyson model in infinite dimensions by Spohn \cite{Spo87}. 
\medskip

\noindent {\bf Airy$_{\beta} $ interacting Brownian motions: }\\
Let $ d=1$, $ \Phi (x) = 0 $, $ \Psi (x,y) = - \log|x-y|$. 
%and $ \beta=1,2,4$. 
We set 
\begin{align}\label{:11c}&
dX^i_t=dB^i_t 
+ \frac{\beta}{2} \lim_{r\to\infty} \Big\{ \big( 
 \sum_{j \not= i, \ |X^j_t |<r}\frac{1}{X^i_t -X^j_t } \big) 
-\int_{|x|<r}\frac{\rhohat (x)}{-x}dx 
\Big \}dt 
\quad (i\in\mathbb{N}) 
.\end{align}
Here we set 
\begin{align}\notag &% \label{:11d}&
\rhohat (x) = \frac{1_{(-\infty , 0)}(x)}{\pi} \sqrt{-x}
.\end{align}
The function 
$\hat{\rho}$ is the scaling limit of the semicircle function
$$
\sigma_{\rm semi}(x) = \frac{1}{2\pi}\sqrt{4-x^2}\mathbf{1}_{(-2,2)}(x)
$$
with scaling $n^{1/3}\sigma_{\rm semi}(xn^{-2/3}+2)$
 associated with soft-edge scaling. 

We solve \eqref{:11c} for $ \beta = 1,2,4$ in \cite{o-t.airy} 
using a result presented in this paper. 
As the solutions $ \mathbf{X}=(X_t^i)_{i\in\N }$ do not collide with each other, 
we label them in decreasing order such that $ X_t^{i+1} < X_t^{i}$ 
for all $ i \in \N $ and $ \zti $. 
Let $ \beta = 2$. The top particle $ X_t^1$ is called the Airy process 
and has been extensively studied by \cite{Joh02,johansson.02,p-spohn} and others. 
In \cite{o-t.sm}, we calculate the space-time correlation functions of 
the unlabeled dynamics $ \mathsf{X}_t=\sum_{i=1}^{\infty}\delta_{X_t^i}$ and 
prove that they are equal to those of \cite{KT07b} and others. 
In particular, our dynamics for $ \beta =2$ are the same as 
the Airy line ensemble constructed in \cite{CH14}. 

\noindent 
{\bf Bessel$_{\alpha,\beta } $ interacting Brownian motions:}\\
Let $ d = 1 $, $ 1 \le \alpha < \infty $, 
$ \Phi (x) = - \frac{\alpha }{2} \log x $ and $ \Psi (x,y) = - \log|x-y|$. 
We set 
\begin{align} \label{:11f} 
& dX_t^i = dB_t^i + \frac{\beta }{2}
\{ \frac{\alpha }{2X_t^i } + \sum _{ j\not = i }^{\infty}
\frac{1}{X_t^i - X_t^j} \} dt \quad (i \in \N )
.\end{align}
Here particles move in $ (0,\infty)$. 

\noindent {\bf Ginibre interacting Brownian motions (\tref{l:35}): } \\
Let $ d=2$, $ \Psi (x,y) = -\log |x-y| $, and $ \beta = 2 $. 
We introduce the ISDE 
\begin{align}\label{:11g}&
dX_t^i = dB_t^i + \limi{r} \sum_{|X_t^i-X_t^j|<r,\ j\not=i } 
\frac{X_t^i-X_t^j}{|X_t^i-X_t^j|^2} dt \quad (i\in\mathbb{N})
\end{align}
and also 
\begin{align}\label{:11h}&
dX_t^i = dB_t^i - X_t^i dt + \limi{r} \sum_{|X_t^j|<r,\ j\not=i } 
\frac{X_t^i-X_t^j}{|X_t^i-X_t^j|^2} dt \quad (i\in\mathbb{N})
.\end{align}
We shall prove that ISDEs \eqref{:11g} and \eqref{:11h} have 
the same strong solutions, reflecting the dynamical rigidity of 
two-dimensional stochastic Coulomb systems 
called the Ginibre random point field (see \sref{s:3}).

All these examples are related to random matrix theory. 
ISDEs \eqref{:11b} and \eqref{:11c} with $ \beta = 1,2,4 $ 
are the dynamical bulk and soft-edge scaling limits of the finite particle systems of 
Gaussian orthogonal/unitary/symplectic ensembles, respectively. 
ISDE \eqref{:11f} with $ \beta = 1,2,4 $ 
is the hard-edge scaling limit of those of the Laguerre ensembles. 
 ISDEs \eqref{:11g} and \eqref{:11h} are dynamical bulk scaling limits 
of the Ginibre ensemble, which is a system of eigen-values of 
non-Hermitian Gaussian random matrices.

\medskip

The next two examples have interaction potentials of 
Ruelle's class \cite{ruelle.2}. 
The previous results also exclude these potentials 
because of the polynomial decay at infinity. 
We shall give a general theorem applicable to substantially all Ruelle's class 
potentials (see \tref{l:V2}). 

\noindent {\bf Lennard-Jones 6-12 potential: }\\
Let $ d= 3$ and 
$\Psi_{6,12}(x) = \{|x|^{-12}-|x|^{-6}\} $. 
The interaction $ \Psi_{6,12}$ is called 
the Lennard-Jones 6-12 potential. The corresponding ISDE is 
\begin{align} \label{:11i}&
dX^i_t = dB^i_t + \frac{\beta }{2} \sum ^{\infty}_{j=1,j\ne i} \{
\frac{12(X ^i_t-X^j_t)}{|X^i_t-X^j_t|^{14}} - 
\frac{6(X ^i_t-X^j_t)}{|X^i_t-X^j_t|^{8}}\, \} dt \quad (i\in \N )
. \end{align}
Since the sum in \eqref{:11i} is absolutely convergent, we do not include 
prefactor $ \limi{r}$ unlike other examples \eqref{:11b}, \eqref{:11c}, \eqref{:11g}, 
and \eqref{:11h}. 

\noindent {\bf Riesz potentials: }\\
Let $ d < a \in \N $ and $\Psi_a(x)=(\beta / a )|x|^{-a}$. 
The corresponding ISDE is 
\begin{align}&\label{:11j}
dX^i_t=dB^i_t + \frac{\beta }{2}
\sum ^{\infty}_{j=1,j\ne i} \frac{X ^i_t-X^j_t}{| X^i_t-X^j_t|^{a+2}}dt 
\quad (i\in \N ) 
.\end{align}
At first glance, ISDE \eqref{:11j} resembles \eqref{:11b} and \eqref{:11g} 
because \eqref{:11j} corresponds to the case $ a =0 $ in \eqref{:11b} and 
\eqref{:11g}. 
However, the sums in the drift terms converge absolutely 
unlike those in \eqref{:11b} and \eqref{:11g}.

\medskip

In the present paper, we introduce a new method of establishing 
the existence of strong solution and the pathwise uniqueness of solution
 of the ISDEs, including the ISDEs 
with long-range interaction potentials. 
Our results apply to a surprisingly wide range of models and, 
in particular, all the examples above (with suitably chosen $ \beta $).

In the previous works, the \Ito\ scheme was used directly in infinite dimensions. This scheme requires the \lq\lq(local) Lipschitz continuity'' of coefficients because it relies on the contraction property of Picard approximations. Hence, a smart choice of stopping times is pivotal but is difficult for long-range potentials in infinite-dimensional spaces.
We {\em do not} apply the \Ito\ scheme to ISDEs {\em directly} but apply it infinitely-many times to an infinite system of finite-dimensional SDEs with consistency (IFC), which we explain in the sequel. 

Our method is based on several novel ideas and consists of three main steps. 
The first step begins by reducing the ISDE to a differential equation of a random point field $ \mu $ satisfying 
\begin{align}\label{:13a}&
 \nabla_x \log \muone (x,\mathsf{s})= -\beta \big\{ 
 \nabla_x \Phi (x) + 
 \limi{r} \sum_{|s_i |< r}^{\infty} \nabla_x \Psi (x, s_i ) \big\} 
.\end{align}
Here $ \mathsf{s}=\sum_i \delta_{s_i}$ is a configuration, 
$ \muone $ is the 1-Campbell measure of $ \mu $ as defined in \eqref{:20i}, and 
$ \nabla_x \log \muone $ is defined in \eqref{:20k}. 
We call $ \nabla_x \log \muone $ the {\em logarithmic derivative} of $ \mu $. 
Equation \eqref{:13a} is given in an informal manner here, and we refer to 
\sref{s:2} for the precise definition. 

The first author proved in \cite{o.isde} with \cite{o.tp,o.rm,o.rm2} that, 
if \eqref{:13a} has a solution $ \mu $ satisfying the assumptions 
\As{A2}--\As{A4} in \sref{s:5c} and \sref{s:E}, then the ISDE \eqref{:11a} 
has a solution $ (\mathbf{X}, \mathbf{B})$ starting at 
$ \mathbf{s}=(s_i)_{i\in\mathbb{N}}$. 
Here a solution means a pair of a stochastic process 
$ \mathbf{X}$ and $ \RdN $-valued Brownian motion $ \mathbf{B}=(B^i)_{i\in\N }$ 
satisfying \eqref{:11a}. We thus see that 
the quartet of papers \cite{o.tp,o.isde,o.rm,o.rm2} achieved 
the first step. 
We note that such solutions in \cite{o.isde} are not strong solutions explained below and 
that the existence of strong solutions and pathwise uniqueness of 
solution of ISDEs \eqref{:11a} were left open in \cite{o.isde}.

In the second step, we introduce the IFC mentioned above 
 using the solution $ (\mathbf{X}, \mathbf{B})$ in the first step. 
That is, we consider a family of the finite-dimensional SDEs of 
$ \mathbf{Y}^m =(Y^{m,i})_{i=1}^m $, $ m \in \N $, given by 
\begin{align}\label{:13b}
dY_t^{m,i} = 
 dB_t^i - \frac{\beta }{2} \nabla_x & \Phi (Y_t^{m,i} ) dt - \frac{\beta }{2} 
 \sum_{j\not=i}^{m} \nabla_x \Psi (Y_t^{m,i}, Y_t^{m,j} ) dt 
 \\ \notag & 
- \frac{\beta }{2}
 \limi{r} \sum_{j= m + 1,\ |X_t^{j}|<r}^{\infty} 
 \nabla_x \Psi (Y_t^{m,i}, X_t^{j} ) dt 
\end{align}
with initial condition 
\begin{align} & \notag 
%\label{:13c} 
\mathbf{Y}_0^m %(\mathbf{s},\mathbf{B}^m,\mathbf{X}^{m*}) 
= \mathbf{s}^m 
.\end{align}
Here, for each $ m \in \N $, we set 
$ \mathbf{s}^m = (s_i)_{i=1}^{m} $ for $ \mathbf{s}=(s_i)_{i=1}^{\infty}$, 
and $ \mathbf{B}^m = (B^i)_{i=1}^{m} $ denotes the $ (\Rd )^m$-valued Brownian motions. 
Note that we regard $ \mathbf{X}$ as ingredients of the coefficients of the SDE \eqref{:13b}. Hence the SDE \eqref{:13b} is time-inhomogeneous although the original ISDE \eqref{:11a} is time-homogeneous. 

 Under suitable assumptions, SDE \eqref{:13b} has 
a pathwise unique, strong solution $ \mathbf{Y}^{m} $. 
Hence, $ \mathbf{Y}^{m}$ is a function of $ \mathbf{s}$, $ \mathbf{B}$, and $ \mathbf{X}$: 
\begin{align} & \notag %\label{:13d}
\mathbf{Y}^{m} = \mathbf{Y}^{m}\sBX = 
\mathbf{Y}_{\mathbf{s},\mathbf{B}}^{m}(\mathbf{X}) 
.\end{align}
As a function of $ (\mathbf{s},\mathbf{B})$, the process $ \mathbf{Y}^{m}$ depends 
only on the first $ m $ components $ (\mathbf{s}^{m},\mathbf{B}^{m})$. 
Since we take $ m $ to go to infinity, the limit, if it exists, 
depends on the whole $ (\mathbf{s},\mathbf{B})$. 
As a function of $ \mathbf{X}$, 
the solution $ \mathbf{Y}^m$ depends only on $\mathbf{X}^{m*}$, where we set 
\begin{align}&\notag %\label{:13e} &
 \mathbf{X}^{m*}=(X^{i})_{i=m+1}^{\infty}
.\end{align}
Hence, 
$ \mathbf{Y}^m(\mathbf{s},\mathbf{B},\cdot)$ is 
$ \sigma [\mathbf{X}^{m*}] $-measurable. We therefore write 
\begin{align} & \notag %\label{:13f}
 \mathbf{Y}^{m} = 
 \mathbf{Y}^{m}(\mathbf{s},\mathbf{B},(\mathbf{0}^m,\mathbf{X}^{m*})) 
= \mathbf{Y}_{\mathbf{s},\mathbf{B}}^{m}((\mathbf{0}^m,\mathbf{X}^{m*})) 
.\end{align}
Here $ \mathbf{0}^m = (0,\ldots,0)$ is the $ (\Rd )^m $-valued constant path. 
The value $ 0 $ does not have any special meaning. Just for the notational reason, 
we put $ \mathbf{0}^m$ here. 
 
With the pathwise uniqueness of the solutions of SDE \eqref{:13b}, 
we see that $ \mathbf{X}^m = (X^i)_{i=1}^{m}$ 
is the unique strong solution of \eqref{:13b}. Hence we deduce that 
\begin{align}\label{:13g}&
\mathbf{X}^m = 
\mathbf{Y}_{\mathbf{s},\mathbf{B}}^{m}((\mathbf{0}^m,\mathbf{X}^{m*})) 
\quad \text{ for all } m \in \N 
.\end{align}
This implies that $ \mathbf{X}^{m}$ becomes a function of $ \mathbf{s} $, 
$ \mathbf{B} $, and $ \mathbf{X}^{m*}$. 
The dependence on $ \mathbf{X}^{m*}$ inherits the coefficients of the SDE \eqref{:13b}.

Relation \eqref{:13g} provides the crucial consistency we use. 
From this we deduce that the maps $ \mathbf{Y}_{\mathbf{s},\mathbf{B}}^{m} $ have a limit 
$ \mathbf{Y}_{\mathbf{s},\mathbf{B}}^{\infty}$ 
as $ m $ goes to infinity at least for $ \mathbf{X}$ in the sense that 
the first $ k $ components $ \mathbf{Y}_{\mathbf{s},\mathbf{B}}^{m,k}$ of 
$ \mathbf{Y}_{\mathbf{s},\mathbf{B}}^m $ coincide with 
$ \mathbf{Y}_{\mathbf{s},\mathbf{B}}^{m',k}$ for all $ k \le m , m'$. 
Furthermore, $ \mathbf{X}$ is a fixed point of the limit map 
$ \mathbf{Y}_{\mathbf{s},\mathbf{B}}^{\infty}$: 
\begin{align}\label{:13h}&
\mathbf{X} = \mathbf{Y}_{\mathbf{s},\mathbf{B}}^{\infty} (\mathbf{X}) 
. \end{align}
Hence, the limit map $ \mathbf{Y}_{\mathbf{s},\mathbf{B}}^{\infty} $ is 
a function of $ (\mathbf{s},\mathbf{B})$ 
and $ \mathbf{X}$ itself through $ \{\mathbf{X}^{m*}\}_{m\in\N }$. 
The point is that, for each fixed $ (\mathbf{s},\mathbf{B})$, the limit 
$ \mathbf{Y}_{\mathbf{s},\mathbf{B}}^{\infty}= 
 \mathbf{Y}^{\infty} (\mathbf{s},\mathbf{B},\cdot)$ 
as a function of $ \mathbf{X} $ is measurable with respect to 
the tail $ \sigma $-field $ \Tpath (\RdN ) $ 
of the labeled path space $\WRN :=C([0,\infty);\RdN )$. 
Here we set $ \Tpath (\RdN ) $ such that 
\begin{align} & \notag %\label{:13i}
\Tpath (\RdN ) = \bigcap_{m=1}^{\infty} \sigma [\mathbf{X}^{m*}] 
.\end{align}

Let $ \Pts $ be the distribution of the solution $ \XB $ of ISDE (1.1) starting at $ \mathbf{s}$. 
By definition $ \Pts $ is a probability measure on $ \WRN \times \WRNz $, where 
$ \WRNz =\{ \mathbf{w} \in \WRN ; \mathbf{w}_0 = \mathbf{0} \} $. 
We write $ (\mathbf{w},\mathbf{b}) \in \WRN \times \WRNz $. 
Let $ \PsB $ be the distribution of the first component under conditioning the second component at $ \bB $, that is, 
$ \PsB $ is the regular conditional distribution of $ \Pts $ such that 
$$ \PsB = \Pts (\mathbf{w} \in \cdot \, | \bB ) 
.$$
By construction $ \PsB $ is the distribution of $ \XB $ starting at $ \mathbf{s}$ 
conditioned at $ \mathbf{B}=\bB $. 

Because $ \mathbf{Y}_{\mathbf{s},\mathbf{B}}^{\infty}$ is 
a $\Tpath (\RdN )$-measurable function in $ \mathbf{X}$ 
for fixed $ (\sB )$,
we deduce that, if $ \Tpath (\RdN )$ is trivial with respect to $ \PsB $,
then $ \mathbf{Y}_{\mathbf{s},\mathbf{B}}^{\infty}$ is
a function of $ (\mathbf{s}$, $ \mathbf{B})$
being independent of $ \mathbf{X}$ for $ \PsB $-a.s. 
Hence, from the identity \eqref{:13h}, the existence of strong solutions and 
the pathwise uniqueness of them are related to $ \PsB $-triviality of $\Tpath (\RdN ) $. 

In \tref{l:61} (First tail theorem), 
we shall give a sufficient condition of the existence of the strong solutions and 
the pathwise uniqueness in terms of the property of $ \PsB $-triviality of $ \Tpath (\RdN ) $. 
This condition is necessary and sufficient as we see in \tref{l:64}. 
The critical point here is that, to some extent, we regard the labeled path tail 
$ \sigma $-field $ \Tpath (\RdN ) $ 
%of the labeled path space $ C([0,T];\RdN ) $ 
as a {\em boundary} of ISDE \eqref{:11a} and pose {\em boundary conditions} in terms of 
 $ \PsB $-triviality of it. 

The formalism regarding a strong solution as a function $ F $ 
of path space is at the heart of the Yamada--Watanabe theory. 
They proved the equivalence between 
$ \{$the existence of a weak solution + the pathwise uniqueness$ \}$ and 
(the existence of) a unique strong solution 
(see Theorem 1.1 \cite[163p]{IW}). 
Our main theorems, Theorems \ref{l:61} and \ref{l:65}, clarify the relation 
between this property of solutions and tail triviality of the labeled path space. 
We shall provide the existence of a strong solution and the pathwise uniqueness of solutions of ISDEs. 
In this sense, this is a counterpart of Yamada--Watanabe's result in infinite dimensions. 
Our argument is however completely different from the Yamada--Watanabe theory. 
Indeed, the existence of a weak solution has been established in the first step, 
and we shall prove the pathwise uniqueness and the existence of strong solutions together 
using the analysis of the tail $ \sigma $-field of the labeled path space.

\smallskip

In the third step, we prove $ \PsB $-triviality of $ \Tpath (\RdN ) $. 
Let $ \sSS $ be the set of configurations on $ \Rd $. 
Then, by definition, $ \sSS $ is the set given by 
\begin{align} & \label{:14a}
\sSS = \{ \mathsf{s} = \sum_{i} \delta _{s_i}\, ;\, 
 \text{ $\mathsf{s} ( K ) < \infty $ for all compact sets $ K \subset \Rd $} \} 
.\end{align}
By convention, we regard the zero measure as an element of $ \sSS $. 
Each element $ \mathsf{s}$ of $ \sSS $ is called a configuration. 
We endow $ \sSS $ with the vague topology, which makes $ \sSS $ a Polish space. 
Let $ \TS $ be the tail $ \sigma $-field 
of the configuration space $ \mathsf{S}$ over $ \Rd $: 
\begin{align}\label{:14b}&
\TS = \bigcap_{r=1}^{\infty} \sigma [\pirc ] 
.\end{align}
Here $ \Sr = \{ |x|\le r \} $ and 
$\pirc $ is the projection $ \map{\pirc }{\sSS }{\sSS }$ such that 
$ \pirc (\mathsf{s}) = \mathsf{s} (\cdot \cap \Sr ^c)$. 

In \tref{l:70} (Second tail theorem), we deduce $ \PsB $-triviality of 
$ \Tpath (\RdN ) $ from $ \mu $-triviality of $ \sSS $, 
where $ \mu $ is the solution of equation \eqref{:13a}. 
Because $ \WRN $ is a huge space and 
the tail $ \sigma $-field $ \Tpath (\RdN ) $ 
is not topologically well behaved, this step is difficult to perform. 

The key point here is the absolute continuity condition \As{AC} given in \sref{s:7}, that is, the condition such that 
the associated unlabeled process $ \mathsf{X}=\{ \mathsf{X}_t \} $, where 
$ \mathsf{X}_t = \sum_{i=1}^{\infty} \delta _{X_t^i}$, 
starting from $ \la $ satisfies 
\begin{align}\label{:14c}&
\Pl \circ \mathsf{X}_t^{-1} \prec \mu \quad \text{ for each } \zzti 
,\end{align}
where $ \Pl $ is the distribution of $ \mathsf{X} $ such that 
$ \Pl \circ \mathsf{X}_0^{-1} = \la $. Here we write 
$ m_1 \prec m_2 $ if $ m_1 $ is absolutely continuous with respect to $ m_2 $. 
This condition is satisfied if $ \lambda = \mu $ and 
$ \mu $ is a stationary probability measure 
of the unlabeled diffusion $ \mathsf{X}$. 
Under this condition, $ \PsB $-triviality of $ \Tpath (\RdN ) $ follows from $ \mu $-triviality of $ \TS $.

One of the key points of the proof of the second tail theorem is \pref{l:73}, 
which proves that the finite-dimensional distributions of 
$ \mathsf{X} $ satisfying \eqref{:14c} 
restricted on the tail $ \sigma $-field are the same as 
the restriction of the product measures of $ \mu $ 
on the tail $ \sigma $-field of the product of the configuration space $ \sSS $ 
(see \eqref{:73a}).

The difficulty in controlling $ \Tpath (\RdN ) $ under the distribution given by 
 the solution of ISDE \eqref{:11a} is that 
the labeled dynamics $ \mathbf{X}=(X^i)_{i\in\N }$ 
have no associated stationary measures 
because they would be an approximately infinite product of Lebesgue measures (if they exist). 
Instead, we study the associated unlabeled dynamics $ \mathsf{X}$ 
as above. 
We shall assume that $ \mathsf{X} $ satisfies the absolute continuity condition. 
Then we immediately see that 
its single time distributions $ \Pl \circ \mathsf{X}_t^{-1}$ 
$ (t \in (0,\infty ))$ starting from $ \la $ are tail trivial on $ \sSS $. 
From this, with some argument, we deduce triviality of the tail 
$ \sigma $-field $ \Tpath (\sSS )$ of 
the unlabeled path space $ C([0,\infty);\sSS )=: \WS $. 
We then further lift triviality of $ \Tpath (\sSS )$ to that of 
$ \Tpath (\RdN ) $. 
To implement this scheme, we employ a rather difficult treatment of the map 
$ \mathbf{Y}_{\mathbf{s},\mathbf{B}}^{\infty}$ in \eqref{:13h}.

We write the second and third steps abstractly. 
This scheme is quite robust and conceptual, 
and can be applied to many other types of ISDEs involving stochastic integrals beyond the \Ito\ type. 
Indeed, we do not use any particular structure of \Ito\ stochastic analysis in these steps. 
Furthermore, even if $ \mu $ is {\em not} tail trivial, 
we can decompose it to the tail trivial random point fields and solve the ISDEs in this case as well.

In \cite{Fr}, Fritz constructed non-equilibrium solutions in the sense that the state space 
 of solutions $ \mathbf{X}$ is explicitly given. In \cite{Fr}, it was assumed that 
potentials are of $ C_0^3(\Rd )$ and the dimension $ d $ is less than or equal to $ 4 $. 
While we were preparing the manuscript, 
Tsai \cite{tsai.14} constructed non-equilibrium solutions 
for the Dyson model with $ 1 \le \beta < \infty $. 
His proof relies on the monotonicity specific to one-dimensional particle systems and uses translation invariance of the solutions. Hence it is difficult to apply his method to 
multi-dimensional models, and 
Bessel$ _{\alpha ,\beta }$ and Airy$ _{\beta }$ interacting Brownian motions even 
if in one dimension. 
Moreover, the reversibility of the unlabeled processes of the solutions 
is unsolved in \cite{tsai.14}.

The present paper is organized as follows. 
In \sref{s:2}, we introduce notation used throughout the paper and recall some notions for random point fields.

From \sref{s:5} to \sref{s:Z} we devote to the general theory concerning on ISDEs. 
In \sref{s:5}, we state the main general theorems (Theorems \ref{l:5A}--\ref{l:5B}). 
We explain the role of the main assumptions \As{\iFc}, \As{TT}, \As{AC}, \As{\Cthree}, and \As{\Cfour} 
in \sref{s:33} and present a list of conditions. 
In \sref{s:6}, clarify the relation between a strong solution and a weak solution 
satisfying \As{\iFc} and triviality of $\Tpath (\RdN ) $ in \tref{l:61}. 
We do this in a general setting beyond interacting Brownian motions. 
We present and prove \tref{l:61} (First tail theorem). 
In \sref{s:7}, 
we derive triviality of $\Tpath (\RdN )  $   from that of $ \TS $. 
We present and prove \tref{l:70} (Second tail theorem). 
In \sref{s:Z}, we prove the main theorems (Theorems \ref{l:5A}--\ref{l:5B}). 

From \sref{s:30} to \sref{s:4} we study the Ginibre interacting Brownian motion, which is one of the most prominent examples of interacting Brownian motions with the logarithmic potential. 
%
%In \sref{s:30}, we study the Ginibre interacting Brownian motion
In \sref{s:3}, we present preliminary results. 
In \sref{s:32}, we state the result for the Ginibre ensemble (\tref{l:35}). 
In \sref{s:4}, we prove \tref{l:35} by employing \tref{l:61} and \tref{l:70}. 

From \sref{s:D} to \sref{s:V} we devote to applying the general theory to 
the class of the ISDEs called interacting Brownian motions. 
We shall prepare feasible sufficient conditions for applications. 
In \sref{s:D}, we quote results on weak solutions of ISDEs 
and the related Dirichlet form theory. 
In \sref{s:E}, we give sufficient conditions of 
assumptions \As{\sIn} and \As{\nbj} used in Theorems \ref{l:5A}--\ref{l:5B}. 
In \sref{s:9}, we give a sufficient condition of assumption \As{\iFc}. 
In \sref{s:Q}, we devote to the sufficient conditions of \NaN\ for 
$ \mu $ with non-trivial tails. 
In \sref{s:V}, we give various examples and prove Theorems \ref{l:V1}--\ref{l:V2}. 
In \sref{s:J} (Appendix), we prove the tail decomposition of random point fields.

We shall explain the main assumptions in the present paper in \sref{s:33} and present a list of assumptions 
in Table \ref{list}. 

\section{Preliminary: logarithmic derivative and quasi-Gibbs measures.}\label{s:2} 
Let $ \sS $ be a closed set in $ \Rd $ such that the interior $ \sS _{\mathrm{int}}$
 is a connected open set satisfying $ \overline{\sS _{\mathrm{int}}} = \sS $ and 
that the boundary $ \partial \sS $ has Lebesgue measure zero. 
Let $ \sSS $ be the configuration space over $ \sS $. 
The set $ \sSS $ is defined by \eqref{:14a} by replacing $ \Rd $ with $ \sS $. 

A symmetric and locally integrable function 
$ \map{\rho ^n }{\sS ^n}{[0,\infty ) } $ is called 
the $ n $-point correlation function of a random point field $ \mu $ 
on $ \sSS $ with respect to the Lebesgue measure if $ \rho ^n $ satisfies 
\begin{align} & \notag %\label{:20g}
\int_{A_1^{k_1}\ts \cdots \ts A_m^{k_m}} 
\rho ^n (x_1,\ldots,x_n) dx_1\cdots dx_n 
 = \int _{\sSS } \prod _{i = 1}^{m} 
\frac{\mathsf{s} (A_i) ! }
{(\mathsf{s} (A_i) - k_i )!} d\mu 
 \end{align}
for any sequence of disjoint bounded measurable sets 
$ A_1,\ldots,A_m \in \mathcal{B}(\sS ) $ and a sequence of natural numbers 
$ k_1,\ldots,k_m $ satisfying $ k_1+\cdots + k_m = n $. 
When $ \mathsf{s} (A_i) - k_i < 0$, according to our interpretation, 
${\mathsf{s} (A_i) ! }/{(\mathsf{s} (A_i) - k_i )!} = 0$ by convention.

Let $ \tilde{\mu }^{[1]}$ be the measure on $(S\times \mathsf{S}, \mathcal{B}(S)\otimes \mathcal{B}(\mathsf{S}))$ determined by
$$
\tilde{\mu }^{[1]}(A\times \mathsf{B})=
\int_{\mathsf{B}} \mathsf{s}(A)\mu(d\mathsf{s}), 
\quad A\in\mathcal{B}(S), \; \mathsf{B}\in \mathcal{B}(\mathsf{S}).
$$
The measure $\tilde{\mu }^{[1]}$ is called the one-Campbell measure of $\mu$. 
In case $\mu$ has one-correlation function $\rho^1$, there exists a regular conditional probability $\tilde{\mu}_x$ of $ \mu $ satisfying
\begin{align*}
&\int_{A}\tilde{\mu}_x (\mathsf{B})\rho^1(x)dx = 
\tilde{\mu }^{[1]} (A\times \mathsf{B}), 
\quad A\in\mathcal{B}(S), \; \mathsf{B}\in \mathcal{B}(\mathsf{S}).
\end{align*}
The measure $\tilde{\mu}_x$ is called the Palm measure of $\mu$ \cite{kal}. 

In this paper, we use the probability measure 
$\mu_x(\cdot)\equiv \tilde{\mu}_x(\cdot -\delta_x)$, 
which is called the reduced Palm measure of $\mu$. Informally, $ \mu_x $ is given by 
\begin{align} & \notag %\label{:20h}
 \mu _{x} = \mu (\cdot - \delta_x | \, \mathsf{s} (\{ x \} ) \ge 1 )
.\end{align} 
We consider the Radon measure $ \muone $ on $ \sS \times \sSS $ such that 
\begin{align}\label{:20i}&
 \muone (dx d\mathsf{s}) = \rho ^1 (x) \mu _{x} (d\mathsf{s}) dx 
.\end{align}
In the present paper, we always use $ \muone $ instead of $\tilde{\mu }^{[1]}$. 
Hence we call $ \muone $ the one-Campbell measure of $ \mu $. 

For a subset $ A $, we set $ \map{\pi_{A}}{\sSS }{\sSS } $ by 
 $ \pi_{A} (\mathsf{s}) = \mathsf{s} (\cdot \cap A )$. 
 We say a function $ f $ on $ \sSS $ is local if 
$ f $ is $ \sigma [\pi_{K}]$-measurable 
for some compact set $ K $ in $ \sS $. 
 For such a local function $ f $ on $ \sSS $, 
we say $ f $ is smooth if 
$ \check{f}=\check{f}_{O }$ is smooth, where 
$ O $ is a relative compact open set in $ \sS $ 
such that $ K \subset O $. Moreover, $ \check{f}_{O } $ 
is a function defined on 
$ \sum_{k=0}^{\infty} O^k $ such that 
$ \check{f}_{O } (x_1,\ldots x_k)$ restricted on $ O^k $ 
is symmetric in $ x_j $ ($ j=1,\ldots,k$) for each $ k $ such that 
$ \check{f}_{O }(x_1,\ldots,x_k) = f (\mathsf{x})$ for 
$ \mathsf{x} = \sum _i \delta _{x_i}$ and that 
$ \check{f}_{O } $ is smooth in $ (x_1,\ldots,x_k)$ for each $ k $. 
Here, for $ k=0$, $ \check{f}_{O } $ is a constant function on $ \{ \mathsf{s};\, \mathsf{s}(O ) = 0 \} $. 
Because $ \mathsf{x}$ is a configuration and $ O $ is 
relatively compact, the cardinality of the particles of $ \mathsf{x}$ 
is finite in $ O $. Note that $ \check{f}_{O } $ has a consistency such that 
\begin{align}\label{:20j}&
\check{f}_{O }(x_1,\ldots,x_k) = 
\check{f}_{O' } (x_1,\ldots,x_k) 
\quad \text{ for all }(x_1,\ldots,x_k) \in (O \cap O')^k 
.\end{align}
Hence $ f (\mathsf{x}) = \check{f}_{O } (x_1,\ldots x_k) $ 
is well defined. 

Let $ \di $ be the set of all bounded, local smooth functions on $ \sSS $. 
We set 
\begin{align*}&
C_{0}^{\infty}(\sS )\otimes \di = 
\{ \sum_{i=1}^N f_i (x)g_i(\mathsf{s})\, ;\, f_i \in C_{0}^{\infty}(\sS ),\ 
g_i \in \di ,\ N \in \N 
\}
.\end{align*}
Let $ \sS _r = \{ s \in \sS \, ;\, | s | \le r \} $. 
We write $ f \in L_{\mathrm{loc}}^p(\sS \ts \sSS , \muone )$ if 
$ f \in L^p(\Sr \ts \sSS , \muone )$ for all $ r \in\N $. 
For simplicity we set 
$ L_{\mathrm{loc}}^p(\muone ) = L_{\mathrm{loc}}^p(\sS \ts \sSS , \muone )$. 

\begin{definition}\label{d:21}
An $ \Rd $-valued function $ \dmu $ is called 
 {\em the logarithmic derivative} of $\mu $ if 
$ \dmu \in L_{\mathrm{loc}}^1(\muone )$ and, 
for all $\varphi \in C_{0}^{\infty}(\sS )\otimes \di $, 
\begin{align}\label{:20k}&
 \int _{\sS \times \sSS } \dmu (x,\mathsf{y})\varphi (x,\mathsf{y}) \muone (dx d\mathsf{y}) = 
 - \int _{\sS \times \sSS } \nabla_x \varphi (x,\mathsf{y}) \muone (dx d\mathsf{y}) 
.\end{align}
\end{definition}

If the boundary $ \partial \sS $ is nonempty and particles 
hit the boundary, then $ \dmu $ would contain a term 
arising from the boundary condition. 
 For example, if the Neumann boundary condition 
is imposed on the boundary, then there would be a local time-type drift. 
In this sense, it would be more natural to suppose that $ \dmu $ 
is a distribution rather than $ \dmu \in L_{\mathrm{loc}}^1(\muone )$. 
Instead, we shall later assume that particles never hit the boundary, 
and the above formulation is thus sufficient in the present situation. 
It would be interesting to generalize the theory, including the case with the boundary condition; however, we do not pursue this here.

A sufficient condition of the explicit expression of 
the logarithmic derivative of random point fields is given in \cite[Theorem 45]{o.isde}. 
Using this, one can obtain the logarithmic derivative of random point fields appearing in random matrix theory such as sine$_{\beta } $, Airy$_{\beta }, $ ($ \beta=1,2,4$), Bessel$_{2,\alpha } $ ($ 1\le \alpha $), 
and the Ginibre random point field (see \lref{l:W2}, \cite{o.rm2,o-t.airy}, \cite{h-o.bes}, \ \lref{l:33}, respectively). 
For canonical Gibbs measures with Ruelle-class interaction potentials, one can easily calculate the logarithmic derivative employing Dobrushin-Lanford-Ruelle (DLR) equation (see \lref{l:W8}). 

\medskip

Let $ \sSS _r^{m} = \{ \mathsf{s} \in \sSS \, ;\, \mathsf{s} (\Sr ) = m \} $. 
 Let $ \Lambda_r $ be the Poisson random point field 
 whose intensity is the Lebesgue measure on $ \Sr $ and set 
 $ \Lambda_r^m = \Lambda_r (\cdot \cap \sSS _r^{m} ) $. 
We set maps $ \map{\pir , \pirc }{\sSS }{\sSS }$ such that 
$ \pir (\mathsf{s}) = \mathsf{s} (\cdot \cap \Sr )$ and 
$ \pirc (\mathsf{s}) = \mathsf{s} (\cdot \cap \Sr ^c)$. 

\begin{definition}\label{d:22}
 A random point field $ \mu $ is called 
 a $ ( \Phi , \Psi )$-quasi Gibbs measure if its regular conditional probabilities 
 $$ 
 \mu _{r,\xi }^{m} = 
 \mu (\, \pir \in \cdot \, | \, \pirc (\mathsf{x}) = 
 \pirc (\mathsf{\xi }),\, \mathsf{x}(\Sr ) = m )
 $$
satisfy, for all $r,m\in \mathbb{N}$ and $ \mu $-a.s.\! $ \xi $, 
\begin{align}\label{:22b}&
\cref{;10Q}^{-1} e^{-\mathcal{H}_r(\mathsf{x}) } \Lambda_r^{m} (d\mathsf{x}) \le 
 \mu _{r,\xi }^{m} (d\mathsf{x}) \le 
\cref{;10Q} e^{-\mathcal{H}_r(\mathsf{x}) } \Lambda_r^{m} (d\mathsf{x}) 
.\end{align}
 Here $ \Ct \label{;10Q} = \cref{;10Q} (r,m,\xi )$ 
 is a positive constant depending on $ r ,\, m ,\, \xi $. 
 For two measures $ \mu , \nu $ on a $ \sigma $-field $ \mathcal{F} $, we write 
$ \mu \le \nu $ if $ \mu (A) \le \nu (A) $ for all $ A \in \mathcal{F} $. 
Moreover, $ \mathcal{H}_r $ is the Hamiltonian on $ \Sr $ defined by 
\begin{align}& \notag %\label{:A2y}
\mathcal{H}_r (\mathsf{x}) = \sum_{x_i\in \Sr } \Phi (x_i) 
+ \sum_{ j < k ,\, x_j, x_k \in \Sr } \Psi (x_j,x_k)
\quad \text{ for } \mathsf{x} = \sum_i \delta_{x_i} 
.\end{align}
\end{definition}
%{}{}{G}

\begin{remark}\label{r:A2} \thetag{1} 
From \eqref{:22b}, we see that for all $r,m\in \mathbb{N}$ and $ \mu $-a.s.\! $ \xi $, 
$ \mu _{r,\xi }^{m} (d\mathsf{x}) $ have (unlabeled) Radon-Nikodym densities 
$ m_{r,\xi }^{m} (d\mathsf{x}) $ with respect to $ \Lambda_r^m $. 
This fact is important when we decompose quasi-Gibbs measures with respect to 
tail $ \sigma $-fields in \lref{l:X2}. 
Clearly, canonical Gibbs measures $ \mu $ with potentials $ (\Phi , \Psi )$ 
are quasi-Gibbs measures, and their densities 
$ m_{r,\xi }^{m} (d\mathsf{x}) $ with respect to $ \Lambda_r^{m} $ 
are given by the DLR equation. % via prescribed potentials $ \Phi $ and $ \Psi $. 
That is, for $ \mu $-a.s.\! $ \xi = \sum_j \delta_{\xi _j} $, 
\begin{align} \notag %\label{:}&
& m_{r,\xi }^{m} (d\mathsf{x}) = \frac{1}{\mathcal{Z}_{r,\xi }^{m} }
\exp \{ - \mathcal{H}_r (\mathsf{x}) - 
\sum_{ x_i \in \Sr ,\, \xi _j \in \Sr ^c } \Psi (x_i, \xi _j)\} 
.\end{align}
Here $ \mathcal{Z}_{r,\xi }^{m} $ is the normalizing constant. 
For random point fields 
appearing in random matrix theory, interaction potentials are logarithmic functions, 
where the DLR equations do not make sense as they are. 
\\\thetag{2}
If $ \mu $ is a $ ( \Phi , \Psi )$-quasi Gibbs measure, 
then $ \mu $ is also a $ ( \Phi + \Phi _{\mathrm{loc,bdd}} , \Psi )$-quasi Gibbs measure 
for any locally bounded measurable function $ \Phi _{\mathrm{loc,bdd}} $. 
In this sense, the notion of \lq\lq quasi-Gibbs'' is robust for the perturbation of free potentials. 
In particular, both the sine$_{\beta } $ and Airy$_{\beta } $ random point fields are 
$ (0, - \beta \log |x-y| )$-quasi Gibbs measures, where $ \beta = 1,2,4$ (see \cite{o.rm2,o-t.airy}). 
\\\thetag{3} 
From \eqref{:22b} we see that $ \Phi $ and $ \Psi $ are locally bounded from below 
with respect to the $ L^{\infty}$-norm by the Lebesgue measure. 
\end{remark}

We collect notation we shall use throughout the paper: 
For a subset $ A $ of a topological space, we shall denote by 
$ \WT (A)=C([0,\infty);A)$ 
the set consisting of $ A $-valued continuous paths on $ [0,\infty) $.

We set $ \| w \|_{C([0,T];\sS )} = \sup_{t\in[0,T]} |w (t)|$. 
We then equip $ \WSN $ with Fr\'{e}chet metric 
$ \mathrm{dist}(\cdot,*) $ given by, for 
$ \mathbf{w}=(w_n)_{n\in\mathbb{N}} $ and 
$ \mathbf{w}'=(w_n')_{n\in\mathbb{N}} $, 
\begin{align}& \notag 
\mathrm{dist}(\mathbf{w},\mathbf{w}') = 
\sum_{T=1}^{\infty} 2^{-T}\Big\{
\sum_{n=1}^{\infty} 2^{-n} \min\{ 1, \|w_n-w_n'\|_{C([0,T];\sS )} \} 
\Big\}
.\end{align}

We introduce unlabeling maps $ \map{\ulabm }{\sS ^{m} \ts \sSS }{\sSS }$ 
 ($ m \in \mathbb{N} $) such that 
\begin{align} & \label{:23a}
 \ulabm ((\mathbf{x},\mathsf{s})) = 
 \sum_{i=1}^{m} \delta _{x_i} + \mathsf{s} \quad \text{ for } 
 \mathbf{x} = (x_i) \in \sS ^{m} , \quad \mathsf{s} \in \sSS 
.\end{align}
By the same symbol $ \ulabm $, we also denote the map 
 $ \map{ \ulabm }{\sS ^{m}}{\sSS } $ such that 
 $ \ulabm (\mathbf{x}) = \sum_{i} \delta_{x_i} $, where 
 $ \mathbf{x} = (x_i)_{i=1}^m $ and $ m \in \mathbb{N} $. 
Let $ \map{\ulab }{\SN }{\sSS }$ such that 
\begin{align}\label{:23b}&
\ulab ((s_i)_{i=1}^{\infty}) = \sum_{i=1} ^{\infty} \delta_{s_i}
.\end{align}
We often write $ \mathbf{s}= (s_i)_{i=1}^{\infty}$ and 
$ \mathsf{s} = \sum_{i=1}^{\infty} \delta _{s_i}$. 
Thus \eqref{:23b} implies $ \ulab (\mathbf{s}) = \mathsf{s}$. 
For $ \www = \{\wwwt \} = \{ (\wti )_{i\in\mathbb{N}} \}_{t\in[0,\infty)} \in \WSN $, 
we set $ \upath (\www ) $, called the unlabeled path of $ \www $, by 
\begin{align}\label{:23c}&
\upath (\www )_t = \ulab (\wwwt ) = \sum_{i} \delta_{\wti }
.\end{align}
We note that $ \upath (\www ) $ is not necessary an element of $ \WS $. 
See \rref{r:35} for example. 
% For $ \www =(\w ^i )$ we write $ \ww = \upath (\www )$. 

Let $ \Ss $ be the subset of $ \sSS $ with no multiple points. 
Let $ \Ssi $ be the subset of $ \Ss $ consisting of an infinite number of points. 
Then by definition $ \Ss $ and $ \Ssi $ are given by 
\begin{align} \label{:23g} &
\Ss = \{ \mathsf{s} \in \sSS \, ;\, \mathsf{s} ( \{ x \} ) \le 1 \text{ for all } x \in \sS \}, \ 
\Ssi = \{ \mathsf{s} \in \Ss \, ;\, \mathsf{s} (\sS ) = \infty \} 
.\end{align} 
A Borel measurable map $ \map{\lab }{\Ssi }{\SN }$ is called a label if 
$ \ulab \circ \lab (\mathsf{s}) = \mathsf{s}$ for all $ \mathsf{s} \in \Ssi $. 
%There exist many labels $ \lab $. 
%We shall fix one throughout the paper. 

Let $ \WSs $ and $ \WSsi $ be the sets consisting of $ \Ss $- and $ \Ssi $-valued continuous path 
on $ [0,\infty)$. 
Each $ \mathsf{w} \in \WSs $ can be written as $ \mathsf{w}_t = \sum_i \delta_{w_t^i}$, 
where $ w^i $ is an $ \sS $-valued continuous path defined on an interval $ I_i $ of the form 
$ [0,b_i)$ or $ (a_i,b_i)$, where $ 0 \le a_i < b_i \le \infty $. 
Taking maximal intervals of this form, we can choose $ [0,b_i)$ and $ (a_i,b_i)$ uniquely up to labeling. 
We remark that 
$ \lim_{t\downarrow a_i} |w_t^i| =\infty $ and 
$ \lim_{t\uparrow b_i} |w_t^i| =\infty $ for $ b_i < \infty $ for all $ i $. 
We call $ w^i$ a tagged path of $ \mathsf{w}$ and $ I_i $ the defining interval of $ w^i $. 
We set 
\begin{align}
\label{:23gg}&
 \WSsiNE = 
\{ \mathsf{w} \in \WSsi \, ;\, I_i = [0,\infty ) \text{ for all $ i $}\} 
.\end{align}
We say tagged path $ w^i $ of $ \mathsf{w}$ does not explode if $ b_i = \infty $, and does not enter if 
$ I_i = [0,b_i)$, where $ b_i $ is the right end of the interval where $ w^i $ is defined. 
Thus $ \WSsiNE $ is the set consisting of non-explosion and non-entering paths. 
%If each tagged path $ w^i =\{ w_t^i \} $ of $ \mathsf{w}_t = \sum_i \delta_{w_t^i}$ does not explode, 
Then we can naturally lift each unlabeled path $ \ww \in \WSsiNE $ to the labeled path 
$ \www =(\w ^i )_{i\in\N } \in \WSN $
using a label $ \lab = (\labi )_{i\in\N }$ such that $ \mathbf{w}_0 = \lab (\mathsf{w}_0)$. 
Indeed, we can do this because each tagged particle can carry the initial label $ i $ forever. 
We write this correspondence by $ \lpath (\ww ) =(\lpathi (\ww ))_{i\in\N }$ and set $ \www $ as 
\begin{align}\label{:23i}&
\www = \lpath (\ww ) \text{ with } \www _0 = \lab (\ww _0) 
.\end{align}
Then $ \w ^i =\lpathi (\ww ) $ by construction. We set 
\begin{align}\notag &% \label{:23j}&
\ww ^{m*} = \sum_{i>m}\delta_{\w ^i}
,\end{align}
where $ \sum_{i>m}\delta_{\w ^i}=
\{ \sum_{i>m}\delta_{\wt ^i} \}_{t\in[0,\infty)} $. 
For an unlabeled path $ \ww $, we call the path 
\begin{align}
\label{:23p}&
\mathbf{w}^{[m]} = (\lpath ^1(\ww ),\ldots,\lpath ^m(\ww ),
\sum_{i>m}\delta_{\w ^i })
\end{align}
the $ m $-labeled path. 
Similarly, for a labeled path $ \www =(\w ^i ) \in \WT (\SN )$ 
we set $ \www ^{[m]} $ by 
\begin{align}
\label{:23q}&
\www ^{[m]} = (w^1,\ldots,w ^m, \sum_{i>m}\delta_{w ^i})
.\end{align}
\begin{remark}\label{r:23}
$ \upath (\www )_t = \ulab (\wt )$ by \eqref{:23c}, whereas 
$ \lpath (\ww )_t \not= \lab (\ww _t)$ in general. 
\end{remark}

\section{The main general theorems: \tref{l:5A}--\tref{l:5B}.} \label{s:5} 
\subsection{ISDE}
Let $ \mathbf{X}= (X ^i)_{i\in\mathbb{N}}$ be an $ \SN $-valued continuous process. 
We write $ \mathbf{X} = \{ \mathbf{X}_t \}_{t\in[0,\infty )} $ and 
$ X^i = \{ X_t^i \}_{t\in[0,\infty )}$. 
For $ \mathbf{X}$ and $ i\in\mathbb{N}$, 
we define the unlabeled processes 
$ \mathsf{X}=\{ \mathsf{X}_t \}_{t\in[0,\infty )} $ and 
$ \mathsf{X}^{i\diamondsuit} = 
\{ \mathsf{X}_t^{i\diamondsuit} \}_{t\in[0,\infty )} $ 
as 
$ \mathsf{X}_t = \sum_{i\in\mathbb{N}} \delta_{X_t^i} 
$ and $
\mathsf{X}_t^{i\diamondsuit} = \sum_{j\in\mathbb{N} ,\ j\not=i } \delta_{X_t^j} 
$.

Let $ \mathsf{H}$ and $ \SSsde $ be Borel subsets of $ \sSS $ such that 
\begin{align}\label{:50x}& 
\mathsf{H} \subset \SSsde \subset \Ssi 
\bigcap \{ \mathsf{s}\,;\, \mathsf{s} (\partial \sS ) = 0 \} 
.\end{align}
Let $ \ulabone $ be as in \eqref{:23a}. 
Define $ \SSSsde \subset \SN $ and $ \SSSsdeone \subset \sS \times \sSS $ by %
\begin{align}\label{:50y}&
\SSSsde = \ulab ^{-1} (\SSsde ) 
,\quad %\\\label{:50z}&
\SSSsdeone = \ulabone ^{-1} (\SSsde ) 
.\end{align}

Let 
$ \map{\sigma }{\SSSsdeone }{\mathbb{R}^{d^2}}$ and 
$ \map{\bbb }{\SSSsdeone }{\Rd }$ be Borel measurable functions, 
where $ d $ is the dimension of the Euclidean space $ \Rd $ that 
includes $ \sS $. 
In infinite dimensions, it is natural to consider coefficients 
$ \sigma $ and $ \bbb $ defined only on a suitable subset 
$ \SSSsdeone $ of $ \sS \ts \sSS $. 
Let $ \map{\lab }{\Ssi }{\SN }$ be the label introduced in \sref{s:2}. 
We consider an ISDE of $ \mathbf{X} =(X^i)_{i\in\mathbb{N}}$ 
starting from $ \lH $ with state space $ \SSSsde $ such that 
% starting at $ \mathbf{s}$: 
\begin{align}\label{:50a}&
dX_t^i = \sigma (X_t^i,\mathsf{X}_t^{i\diamondsuit}) dB_t^i + 
\bbb (X_t^i,\mathsf{X}_t^{i\diamondsuit}) dt \quad (i\in\mathbb{N}) 
,\\\label{:50b}&
\mathbf{X} \in \WT (\SSSsde ) 
,\\\label{:50c}& %[]大変更注意
\mathbf{X}_0 % = \lab (\mathsf{s}) %情報がないので、この項を削る
\in \lH %= \mathbf{s} 
.\end{align}
Here $ \mathbf{B}=(B^i)_{i\in\N }$ is 
an $ \RdN $-valued Brownian motion; 
that is, $ \{B^i\}_{i\in\N }$ are independent copies of 
a $ d$-dimensional Brownian motion starting at the origin.

\begin{remark}\label{r:50}
Note that $ \lH \subset \uH $ and that $ \uH $ is much larger than $ \lH $. 
For example, if we take $ \lab (\mathsf{s}) = (s_i)$ as $ |s_i| \le |s_{i+1}|$ for all $ i \in\N $, 
then $ \mathbf{X}_t $ will soon exit from $ \lH $. This is why we take 
$ \SSSsde = \ulab ^{-1}(\SSsde )$ in \eqref{:50y} rather than $ \SSSsde = \lab (\SSsde )$. 
\end{remark}

From \eqref{:50b} the process $ \mathbf{X}$ moves in the set $ \SSSsde $ 
where the coefficients $ \sigma $ and $ \bbb $ are well defined. 
Moreover, each tagged particle $ X^i $ of $\mathbf{X}= (X^i )_{i\in\N }$ never explodes . 
By \eqref{:50b}, $ \mathsf{X}_t \in \SSsde $ for all $ t \ge 0$, and in particular the initial starting point 
$ \mathbf{s}$ in \eqref{:50c} is supposed to satisfy 
 $ \mathbf{s} \in \lH \subset \SSSsde $, which implies 
$ \ulab (\mathbf{s}) \in \mathsf{H} \subset \SSsde $. 

By \eqref{:50x}, $ \mathsf{H}$ is a subset of $ \SSsde $. 
We shall take $ \mathsf{H}$ in such a way that \eqref{:50a}--\eqref{:50c} 
has a solution for each $ \mathbf{s} \in \lH $. 
To detect a sufficiently large subset $ \mathsf{H} $ satisfying this 
%where \ISDEb hold for all $ \mathbf{s} \in \lH $, 
is an important step to solve the ISDE. 

The meaning of $ \mathsf{H}$ to be large is however a problem at this stage 
because there is no natural measure 
on the infinite product space $ \SN $. 
In practice, we equip $ \sSS $ with a random point field $ \mu $ such that 
$ \mu (\mathsf{H}) = \mu (\SSsde )= \mu (\mathsf{\sSS })=1$. 
We thus realize $ \mathsf{H}$ as a support of $ \mu $. 
We shall later assume \As{A1} in \sref{s:D} % \eqref{:51a} 
to relate $ \mu $ with \eqref{:50a} 
in such a way that the unlabeled dynamics $ \mathsf{X}$ 
of the solution $ \mathbf{X}$ is $ \mu $-reversible. 
% 
%Hence $ \mathbf{X} $ moves in $ \uH $. 
%The role of $ \mu $ is thus to specify the real state space of $ \mathbf{X}$. 
In this sense the random point field $ \mu $ is associated with ISDE \eqref{:50a}. 

We remark that we can extend $ \lH $ in \eqref{:50c} to $ \uH $ by retaking 
other labels $ \lab $. Because the coefficients of ISDE \eqref{:50a} have symmetry, this causes no problems. 

Essentially, following \cite[Chapter IV]{IW} in finite dimension, we present a set of notions 
related to solutions of ISDE. In Definitions \ref{d:51}--\ref{d:58}, $ \OFPF $ is a general probability space. 

%\marginpar{[]}

\begin{definition}[weak solution]\label{d:51} 
By a weak solution of ISDE \eqref{:50a}--\eqref{:50b}, we mean 
 an $ \sS ^{\mathbb{N}} \ts \RdN $-valued stochastic process $ \XB $ 
defined on a probability space $ (\Omega , \mathcal{F}, P )$ 
with a reference family $ \{ \mathcal{F}_t \}_{t \ge 0 } $ such that

\noindent 
\thetag{i} $ \mathbf{X}=(X^i)_{i=1}^{\infty} $ 
is an $ \SSSsde $-valued continuous process. 
Furthermore, $ \mathbf{X}$ is adapted to 
$ \{ \mathcal{F}_t \}_{t \ge 0 }$, 
that is, $ \mathbf{X}_t $ is 
$ \mathcal{F}_t /\Bt $-measurable for each $ \zti $, where 
\begin{align}& \label{:50k}
\Bt = 
\sigma [ \www _s ; 0\le s \le t ,\, \www \in \WSN ] 
.\end{align}
\thetag{ii} $ \mathbf{B} = (B^i)_{i=1}^{\infty}$ is an $ \RdN $-valued 
\FtB 
with $ \mathbf{B}_0 = \mathbf{0}$, 
\\
\thetag{iii} 
the family of measurable $ \{ \mathcal{F}_t \}_{t \ge 0 } $-adapted 
processes $ \Phi $ and $ \Psi $ defined by 
\begin{align}\notag &% & \label{:50l}
\Phi ^i(t,\omega ) = \sigma 
(X_t^i(\omega ),\mathsf{X}_t^{i\diamondsuit}(\omega ))
,\quad % \\ \notag &
\Psi ^i(t,\omega ) = \bbb 
(X_t^i(\omega ),\mathsf{X}_t^{i\diamondsuit}(\omega ))
\end{align}
belong to $ \mathcal{L}^{2} $ and $ \mathcal{L}^1 $, respectively. 
Here $ \mathcal{L}^{p} $ is the set of all measurable $ \{ \mathcal{F}_t \}_{t \ge 0 } $-adapted 
processes $ \alpha $ such that 
$ E[ \int_0^T|\alpha (t,\omega)|^p dt ] < \infty $ for all $ T $. 
Here we can and do take a predictable version of $ \Phi ^i $ and 
$ \Psi ^i$ (see pp 45-46 in \cite{IW}). 
\\
\thetag{iv} with probability one, the process $ \XB $ satisfies for all $ t $ 
\begin{align}\notag &% \label{:50m}&
X_t^i - X_0^i = 
\int_0^t \sigma (X_u^i,\mathsf{X}_u^{i\diamondsuit}) dB_u^i 
 + 
\int_0^t 
\bbb (X_u^i,\mathsf{X}_u^{i\diamondsuit}) du \quad (i\in\mathbb{N}) 
.\end{align}
\end{definition}

\begin{definition}[weak solution on $ \mathbf{A}$]\label{d:20a} 
We say the ISDE \eqref{:50a}--\eqref{:50b} has a weak solution 
on a Borel set $ \mathbf{A} $ if it has a weak solution for arbitrary initial distribution $ \nu $ such that $ \nu (\mathbf{A}) = 1 $. 
\end{definition}

We say $ \mathbf{X}$ is a weak solution if the accompanied Brownian motion $ \mathbf{B}$ 
is obvious or not important. 
A solution $\mathbf{X}$ staring at $ \mathbf{x}$ means $ \mathbf{X}$ is a solution 
such that $ \mathbf{X}_0=\mathbf{x}$ a.s. 
\begin{remark}\label{r:20a}
In \cite[Chap.\, IV]{IW}, the state space and the set of the initial starting points of SDEs are the same and taken to be $ \Rd $. In the present paper, the set of the initial starting points is $ \lH $. $ \lH $ is a subset of $ \SSSsde $. So we introduced the notion \lq\lq weak solution on $ \mathbf{A}$'' in \dref{d:20a}. 
\end{remark}

\begin{definition}[uniqueness in law]\label{d:53}
We say that the uniqueness in law of weak solutions on $ \lH $ for \eqref{:50a}--\eqref{:50b} 
holds if whenever $ \mathbf{X}$ and $ \mathbf{X}'$ are two weak solutions 
whose initial distributions coincide, then the laws of the processes 
$ \mathbf{X}$ and $ \mathbf{X}'$ on the space $ \WSN $ coincide. 
If this uniqueness holds for an initial distribution $ \delta _{\mathbf{s}}$, then we say 
the uniqueness in law of weak solutions for \eqref{:50a}--\eqref{:50b} starting at $ \mathbf{s}$ holds. 
\end{definition}

\begin{remark} 
For each $ \mathbf{s} \in \lH $ take $ \delta_{\mathbf{s}}$ 
as the initial law of the ISDE \eqref{:50a}--\eqref{:50b}. 
Then the uniqueness in \dref{d:53} is equivalent to 
the uniqueness of the law of weak solutions starting at each $ \mathbf{s} \in \lH $. 
We refer to Remark 1.2 in \cite[162 p]{IW} for the corresponding result for finite-dimensional SDEs. 
\end{remark}

\begin{definition}[pathwise uniqueness]\label{d:54}
We say that the pathwise uniqueness of weak solutions of 
\eqref{:50a}--\eqref{:50b} on $ \lH $ holds 
if whenever $ \mathbf{X}$ and $ \mathbf{X}'$ are two weak solutions 
defined on the same probability space $ (\Omega , \mathcal{F} ,P )$ with the same reference family $ \{ \mathcal{F}_t \}_{t \ge 0 }$ and the same $ \RdN $-valued \FtB 
$ \mathbf{B} $ such that $ \mathbf{X}_0=\mathbf{X}_0' \in \lH $ a.s., 
then 
\begin{align}\notag &% \label{:50o}&
P (\text{$ \mathbf{X}_t=\mathbf{X}_t'$ for all $ t \ge 0 $}) = 1 
.\end{align} 
\end{definition}

\begin{definition}[pathwise uniqueness of weak solutions starting at $ \mathbf{s}$]
\label{d:54b}
We say that the pathwise uniqueness of weak solutions of 
\eqref{:50a}--\eqref{:50b} starting at $ \mathbf{s}$ holds 
if the condition of \dref{d:54} holds for 
$ \mathbf{X}_0=\mathbf{X}_0' = \mathbf{s} $ a.s. 
\end{definition}

We now state the definition of strong solution, which is analogous to 
Definition 1.6 in \cite[163 p]{IW}. 
Let $ \PBr $ be the distribution of an $ \RdN $-valued Brownian motion 
$ \mathbf{B} $ with $ \mathbf{B}_0 = \mathbf{0}$. 
Let $ \WRNz = \{ \www \in \WRN \, ; \www _0 =\mathbf{0} \}$. 
Clearly, $ \PBr (\WRNz )= 1 $. 

Let $ \Bt $ be as \eqref{:50k}. 
Let $ \Bt (\PBr )$ be the completion of 
$ \sigma [ \www _s ; 0\le s \le t ,\, \www \in \WRNz ] $ with respect to $ \PBr $. 
Let $ \mathcal{B}(\PBr ) $ be the completion of $ \mathcal{B}(\WRNz ) $ with respect to $ \PBr $. 

\begin{definition}[a strong solution starting at $ \mathbf{s}$] \label{d:55} 
A weak solution $ \mathbf{X}$ of \eqref{:50a}--\eqref{:50b} 
with an $ \RdN $-valued $ \Ft $-Brownian motion 
$\mathbf{B}$ is called a strong solution starting at $ \mathbf{s}$ defined on $ \OFPF $ 
if $ \mathbf{X}_0=\mathbf{s}$ a.s.\, and if there exists a function 
$ \map{\Fs }{\WRNz }{\WSN }$ such that 
$ \Fs $ is $ \mathcal{B}(\PBr ) /\mathcal{B}(\WSN )$-measurable, and that 
$ \Fs $ is $ \Bt (\PBr ) /\Bt $-measurable for each $ t $, and 
that $ \Fs $ satisfies 
\begin{align}\notag &% \label{:50q}&
\mathbf{X} = \Fs (\mathbf{B}) \quad \text{ a.s.}
\end{align}
Also we call $ \Fs $ itself a strong solution starting at $ \mathbf{s}$. 
\end{definition}

\begin{definition}[a unique strong solution starting at $ \mathbf{s}$] \label{d:55b} 
We say \eqref{:50a}--\eqref{:50b} has a unique strong solution starting at $ \mathbf{s}$ if 
% \eqref{:50a}--\eqref{:50b} has a strong solution $ \Fs $ starting at $ \mathbf{s}$ 
there exists a $ \mathcal{B}(\PBr ) /\mathcal{B}(\WSN )$-measurable function 
$ \map{\Fs }{\WRNz }{\WSN }$ 
such that, for any weak solution $ (\hat{\mathbf{X}},\hat{\mathbf{B}})$ 
of \eqref{:50a}--\eqref{:50b} starting at $ \mathbf{s}$, it holds that 
\begin{align}\notag &% \label{:50r}&
\hat{\mathbf{X}}=\Fs (\hat{\mathbf{B}}) \quad \text{ a.s.}
\end{align}
and if, for any $ \RdN $-valued $ \{\mathcal{F}_t\} $-Brownian motion $ \mathbf{B} $ 
defined on $ \OFPF $ with $ \mathbf{B}_0=\mathbf{0}$, the continuous process $ \Fs (\mathbf{B})$ 
is a strong solution of \eqref{:50a}--\eqref{:50b} starting at $ \mathbf{s}$. 
Also we call $ \Fs $ a unique strong solution starting at $ \mathbf{s}$. 
\end{definition}

We next present a variant of the notion of a unique strong solution. 

\begin{definition}[a unique strong solution under constraint] \label{d:58} 
For a condition \As{$ \bullet$}, we say \eqref{:50a}--\eqref{:50b} has a unique strong solution starting at $ \mathbf{s}$ under the constraint \As{$ \bullet$} if 
% \eqref{:50a}--\eqref{:50b} has a strong solution $ \Fs $ starting at $ \mathbf{s}$ 
%
there exists a $ \mathcal{B}(\PBr ) /\mathcal{B}(\WSN )$-measurable function 
$ \map{\Fs }{\WRNz }{\WSN }$ 
%
%satisfying \As{$ \bullet$} % 
such that for any weak solution 
$ (\hat{\mathbf{X}},\hat{\mathbf{B}})$ of \eqref{:50a}--\eqref{:50b} starting at $ \mathbf{s}$ 
satisfying \As{$ \bullet$}, it holds that 
\begin{align}
\label{:50s}&
\hat{\mathbf{X}}=\Fs (\hat{\mathbf{B}}) \quad \text{ a.s.}
\end{align}
and if for any $ \RdN $-valued $ \{\mathcal{F}_t\} $-Brownian motion $ \mathbf{B}$
defined on $ \OFPF $ with $ \mathbf{B}_0=\mathbf{0}$ the continuous process $ \Fs (\mathbf{B})$ 
is a solution of \eqref{:50a}--\eqref{:50b} starting at $ \mathbf{s}$ satisfying \As{$ \bullet$}. 
Also we call $ \Fs $ a unique strong solution starting at $ \mathbf{s}$ 
under the constraint \As{$ \bullet$}. 
\end{definition}

\begin{remark}\label{r:55}
\thetag{1} The meaning of strong solutions is similar to the conventional situation in \cite[pp 159--167]{IW}. The difference is that we consider solutions starting at a point $ \mathbf{s}$. 
In \cite{IW}, initial distributions are taken over all probability measures on the state space. 
\\\thetag{2}
Similarly as \dref{d:58} we can introduce the notions of constrained versions of uniqueness in 
Definitions \ref{d:53}--\ref{d:54b}. 
\end{remark}

\subsection{Main Theorem I (\tref{l:5A}): 
$ \mu $ with trivial tail.} \label{s:5d}

Let $ \XB $ be an $ \sS ^{\mathbb{N}} \ts \RdN $-valued continuous process 
defined on a filtered space $ \OFPF $. 
We assume that $ \OFP $ is a standard probability space. 
Then the regular conditional probability 
\begin{align}\notag &% \label{:53e}&
\Ps = P ( \cdot | \mathbf{X}_0=\mathbf{s})
\end{align}
exists for $ P \circ \mathbf{X}_0 ^{-1}$-a.s.\ $ \mathbf{s}$. 
We investigate $ \XB $ under $ \Ps $, and thus regard $ \XB $ as a stochastic process defined on the filtered space $ \2 $.

For 
$ \mathbf{X} = (X^i)_{i\in\N } $ we set 
$\mathsf{X}_t^{m*} = \sum_{i=m+1}^{\infty} \delta_{X_t^i }$ as before. 
Define %the functions 
\begin{align} \notag %\label{:53a}
& \text{
$ \map{ \sigmaXms } 
{[0,\infty) \ts \sS ^{m} }{\mathbb{R}^{d^2}}$\ and \ 
 $ \map{ \bbbXms }{[0,\infty) \ts \sS ^{m} }{\Rd }$ }
\end{align}
such that, for $ (u,\mathbf{v}) \in \sS ^m$ and 
$ \mathsf{v} = \sum_{i=1}^{m-1} \delta_{v_i} \in \sSS $, where 
$\mathbf{v}=(v_1,\ldots,v_{m-1}) \in \sS ^{m-1} $, 
\begin{align} \label{:53b}
&
 \sigmaXms ( t, (u, \mathbf{v})) = 
 {\sigma } (u , \mathsf{v} + \mathsf{X}_t^{m*}) ,\quad 
\bbbXms ( t, (u, \mathbf{v})) = {\bbb } (u , \mathsf{v} + \mathsf{X}_t^{m*})
.\end{align}
We write $ \lab (\mathsf{s}) =(s_i)_{i\in\N }= \mathbf{s} $ and 
$ \mathsf{s}_m^* = \sum_{i=m+1}^{\infty} \delta_{s_i}$. 
Recall $ \mathbf{X}_0 = \lab (\mathsf{s}) $. 
We then have $ \mathsf{X}_0^{m*} = \mathsf{s}_m^* $ by construction. 
We remark that the coefficients $ \sigmaXms $ and $ \bbbXms $ depend on 
both unlabeled path $ \mathsf{X}^{m*} $ and the label $ \lab $, 
although we suppress $ \lab $ from the notation. 
Let 
\begin{align} \label{:35a} & %\notag &
\SSSsdemtw 
= \{ \mathbf{s}^m = (s_1,\ldots,s_m) \in \sS ^m \,;\, 
\ulab (\mathbf{s}^m) + \ww _t^{m*} \in \SSsde 
 \} 
,\end{align}
where $ \ww _t^{m*}=\sum_{i=m+1}^{\infty}\delta_{w_t^i}$ for 
$ \ww _t = \sum_{i=1}^{\infty} \delta_{w_t^i}$. 
By definition, $ \SSSsdemtw $ 
is a subset of $ \sS ^m $ depending on $ \ww _t^{m*}$. 
In particular, $ \SSSsdemtw $ is a time-dependent domain. 
Let 
\begin{align} & \notag 
 \mathbf{Y}^m=(Y^{m,i})_{i=1}^m ,\quad 
 \mathbf{Y}^{m,i\diamondsuit} = (Y^{m,j})_{j\not=i}^m , 
\quad 
 \mathsf{Y}_t^{m,i\diamondsuit}=\sum_{j\not=i}^m \delta_{Y_t^{m,j}}
.\end{align}

We introduce the SDE with random environment $ \mathsf{X}$ 
defined on $ \2 $ describing $ \mathbf{Y}^m $ given by 
\begin{align}\label{:53f}&
dY_t^{m,i} = 
\sigmaXms (t, (Y_t^{m,i},\mathbf{Y}_t^{m,i\diamondsuit})) dB_t^i + 
 \bbbXms (t, (Y_t^{m,i},\mathbf{Y}_t^{m,i\diamondsuit})) dt
% \quad (i=1,\ldots,m)
,\\\label{:53g}& 
 \mathbf{Y}_t^m \in \SSSsdemt \quad \text{ for all } t 
,\\\label{:53h}&
 \mathbf{Y}_0^{m} = \mathbf{s}^m ,
\quad \text{ where $ \mathbf{s}^m=(s_1,\ldots,s_m) $ for $ \mathbf{s}=(s_i)\in\SN $}
.\end{align}
A triplet of $ \{ \mathcal{F}_t \} $-adapted, continuous process 
$ (\mathbf{Y}^m,\mathbf{B}^m ,\mathbf{X}^{m*} )$ on $ \2 $ satisfying 
\eqref{:53f}--\eqref{:53h} is called a weak solution. 

\begin{remark}\label{r:31} \thetag{1} 
Equation \eqref{:53f} makes sense because $ \mathbf{X}^{m*}$, $ \mathbf{B}^m $, 
and $\mathbf{Y}^m$ are all defined on the same filtered space $ \2 $. 
We remark that \eqref{:53f} depends on $ \mathbf{X}^{m*} $, and that 
$ \mathbf{X}^{m*} $ is regarded as a part of the coefficient of \eqref{:53f}. 
We emphasize $ \XB $ is a priori given in SDE \eqref{:53f}. 
We consider SDE \eqref{:53f} only for $ \mathbf{B}^m$, but 
not for an arbitrary $ \Ft $-Brownian motion $ \hat{\mathbf{B}}^m $. 
\\
\thetag{2} 
The SDE \eqref{:53f} is {\em not} a conventional type because the coefficient depends on $ \mathbf{X}$. We can regard $ \mathbf{X}$ as a random environment, and call this SDE of random environment type. Random environment type SDEs had appeared in homogenization problem (see \cite{o.homo,PV} for example). 
In this case, random environment and Brownian motion in SDEs are usually 
independent of each other. This is not the case in the present situation. 
If $ (\hat{\mathbf{B}}^m,\hat{\mathsf{X}}^{m*})$ is equivalent in law to 
$ (\mathbf{B}^m,\mathsf{X}^{m*}) $, then we can replace $ (\mathbf{B}^m,\mathsf{X}^{m*}) $ by $ (\hat{\mathbf{B}}^m,\hat{\mathsf{X}}^{m*})$ in \eqref{:53f}. 
The new SDE is equivalent to \eqref{:53f} in the sense that the former has a weak solution if and only if the latter has one. 
We emphasize that $ \mathbf{B}^m $ and $ \mathbf{X}^{m*}$ are $ \{ \mathcal{F}_t\} $-adapted and can depend on each other. 
%
% The original Brownian motion $ \mathbf{B}$ is a functional of $ \mathbf{X}$ 
% in all the examples in the present paper. 
% We can nevertheless take the first $ m $-components $ \mathbf{B}^m$ of $ \mathbf{B}$ as 
% the Brownian motion $ \mathbf{B}^m $ in \eqref{:53f} 
% because $ \mathbf{B}^m $ is $ \{ \mathcal{F}_t\} $-adapted. 
\\\thetag{3} 
The triplet $ (\mathbf{X}^m,\mathbf{B}^m,\mathbf{X}^{m*})$ 
made of the original weak solution $ \XB $ of the ISDE \eqref{:50a}--\eqref{:50c} 
is a weak solution of \eqref{:53f}--\eqref{:53h}. 
This yields the crucial identity \eqref{:13g}. 
\end{remark}

We define the notion of strong solutions and a unique strong solution of \eqref{:53f}--\eqref{:53h}. 
Let $ \PPPm $ be the distribution of $ (\lBlhatm ) $ under $ \2 $: 
\begin{align}
\notag & 
\PPPm = \Ps \circ (\lBlhatm )^{-1}
.\end{align}
Let 
$ \WRdzm = \{ \mathbf{w} \in \WRdm \, ; 
\mathbf{w}_0 =\mathbf{0} \}$ as before and set 
\begin{align}\notag & 
\Ehatm = 
\overline{\mathcal{B} (\WWdm ) } ^{\Pt ^m }
,\\\notag & 
\Ehatmt = \overline{\Bt (\WWdm ) }^{\PPPm }
.\end{align}
Here $ \Bt (\WWdm )= \sigma [(\mathbf{v}_s,\mathbf{w}_s) ; 0\le s \le t ,\, 
(\mathbf{v},\mathbf{w})\in \WWdm ] $. 
Let 
$ \Btm = \sigma [\mathbf{w}_s; 0\le s \le t ,\, \mathbf{w} \in \WT (\Rdm )]$. 
%
%We also set $ \mathbf{s}= \lab (\mathsf{s})$ and $ \mathbf{s}^m=\labm (\mathsf{s}) $. 
% 
We state the definition of strong solution. 
\begin{definition}[strong solution for $ \XB $ starting at $ \mathbf{s}^m $] \label{d:41} 
 $\mathbf{Y} ^{m}$ is called a strong solution of 
\eqref{:53f}--\eqref{:53h} for $ \XB $ \uPs\ if 
$ (\mathbf{Y}^{m},\mathbf{B}^m,\mathbf{X}^{m*})$ satisfies \eqref{:53f}--\eqref{:53h} and 
there exists a $ \Ehatm $-measurable function 
\begin{align}
\notag &% \label{:42e}&
\map{\Fms }{\WWdm }{\WT (\Rdm ) }
\end{align}
such that $ \Fms $ is $ \Ehatmt /\Btm $-measurable for each $ t $, and $ \Fms $ satisfies %
\begin{align}\notag &% &\label{:42f} 
\mathbf{Y}^m = \Fms (\lBlhatm ) \quad \text{ $ \Ps $-a.s.}
\end{align}
\end{definition}
\begin{remark}\label{r:42}
Our definition of a strong solution is different from that of Definition 1.6 in \cite[163 p]{IW} 
with the following points: We consider solutions starting at a point $ \mathbf{s}^m $ only. 
The main difference is that both the $ \Ft $-Brownian motion $ \mathbf{B}$ and 
the process $ \mathbf{X}^{m*}$ are a priori given and fixed. 
Hence the solution $ \mathbf{Y}^m $ is a function of not only $ \mathbf{B}$ but also 
$ \mathbf{X}^{m*}$. 
This means, if we put an arbitrary $ \Ft $-Brownian motion $ \mathbf{B}'$ 
into $ \Fms $ as $ \Fms (\mathbf{B}',\mathbf{X}^{m*})$, 
then $ \Fms (\mathbf{B}',\mathbf{X}^{m*})$ is not necessary a solution. 
We call $ \mathbf{X}^{m*}$ an environment processes. 
We note that there is no environment process in the framework of Definition 1.6 in \cite[163 p]{IW}. 
We shall take the limit $ m \to \infty $, and prove that the effect of $ \mathbf{X}^{m*}$ will vanish in the limit. As a result, the limit ISDE becomes conventional. 
Vanishing the effect of $ \mathbf{X}^{m*}$ as $ m \to \infty $ is a key to our argument. 
We will do this by the second main theorem (\tref{l:70}). 
\end{remark}

\begin{definition}
[a unique strong solution for $ \XB $ starting at $ \mathbf{s}^m $]\label{d:42} 
The SDE \eqref{:53f}--\eqref{:53h} is said to have 
a unique strong solution for $ \XB $ \uPs\ 
if there exists a function $ \Fms $ satisfying the conditions in \dref{d:41} and, 
for any weak solution $ (\hat{\mathbf{Y}}^m ,\mathbf{B}^m,\mathbf{X}^{m*})$ of 
\eqref{:53f}--\eqref{:53h} \uPs, 
\begin{align}
\notag &% \label{:42g}& \quad \quad 
\hat{\mathbf{Y}}^m = \Fms (\lBlhatm ) \quad \text{for $ \Ps $-a.s.}
\end{align}
\end{definition}

\smallskip 

The function $ \Fms $ in \dref{d:41} is also called a strong solution starting at $ \mathbf{s}^m $. 
The SDE \eqref{:53f}--\eqref{:53h} is said to have a unique strong solution $ \Fms $ 
if $ \Fms $ satisfies the condition in \dref{d:42}. 
We note that the function $ \Fms $ is unique for $ \PPPm $-a.s.\,in this case. 

We recall that these two notions are different from those of the infinite-dimensional counterparts 
Definitions \ref{d:55} and \dref{d:55b} because the SDE \eqref{:53f}--\eqref{:53h} is of random environment type. 

We introduce the IFC condition of $ \XB $ defined on $ \OFPF $ as follows. 
 
\smallskip 

\noindent \As{\iFc} \ 
The SDE \eqref{:53f}--\eqref{:53h} 
has a unique strong solution $ \Fms (\mathbf{B}^m,\mathbf{X}^{m*})$ 
for $ \XB $ under $ \Ps $ %defined on $ \2 $ 
for $ P\circ \mathbf{X}_0 ^{-1}$-a.s. $ \mathbf{s}$ 
for all $ m \in \N $.

\smallskip 
\noindent 
For convenience we introduce a quenched version of \As{\iFc}: 

\noindent \iFcs\ 
The SDE \eqref{:53f}--\eqref{:53h} 
has a unique strong solution $ \Fms (\mathbf{B}^m,\mathbf{X}^{m*})$ 
for $ \XB $ under $ \Ps $ %defined on $ \2 $ 
for all $ m \in \N $. 

\smallskip 
\noindent 
By definition \As{\iFc} holds if and only if 
\iFcs\ holds for $ P\circ \mathbf{X}_0 ^{-1}$-a.s. $ \mathbf{s}$. 

\smallskip 

The SDE \eqref{:53f}--\eqref{:53h} is time inhomogeneous and 
the state space of the solution $ \mathbf{Y}^m $ given by \eqref{:53g} 
also depends on time $ t $ through $ \mathbf{X}_t^{m*} $. 
Because the SDE \eqref{:53f} is {\em finite}-dimensional, 
one can apply the classical theory of SDEs directly. 
A feasible sufficient condition of \As{\iFc} is given in \sref{s:9}.

% \begin{definition}[IFC solutions]\label{d:59}
% $ \XB $ under $ P $ is called an IFC solution of \eqref{:50a}--\eqref{:50b} 
% if it is a weak solution satisfying the IFC condition. 
% \end{definition}
%\footnote{初期条件に付いて、これでよいか、IFC条件との表現の整合性を考える。$ \Ps $のほうがよいかもしれない、SDEに第3行を付け加え。よくわからん。}
We remark that a continuous process $ \XB $ satisfying \As{\iFc} is not necessary 
a weak solution. See \dref{d:61AIFC} and \rref{r:63}.

%The next assumptions of probability measures $ P $ on $ \Omega $ and $ \mu $ on $ \sSS $
%are related to tail triviality of the labeled path space $\TpathTSN $ with respect to label. 
Let $ \TS $ be the tail $ \sigma$-field of $ \sSS $ defined as \eqref{:14b}. 
A random point field $ \mu $ on $ \sS $ is called tail trivial if 
$ \mu (A) \in \{ 0,1 \} $ for all $ A \in \TS $. 
Let $ \mathbf{X}=(X^i)_{i\in\N }$ be a continuous process defined on $ \OFPF $ 
and $ \mathsf{X}$ be the associated unlabeled process such that %be as \eqref{:20d}. 
$ \mathsf{X}_t=\sum_i \delta_{X_t^i}$. 
Let $ \WSsiNE $ be as \eqref{:23gg}. 
Let $ \map{\mrT }{\WT (\SN ) }{\mathbb{N}\cup\{ \infty \} }$ be such that 
 \begin{align}\label{:53n}&
 \mrT (\mathbf{w}) = \inf \{ m \in \mathbb{N}\, ;\, 
 \min_{t\in[0,T]}|w _t^n| > r \text{ for all } n \in \mathbb{N} \text{ such that } n > m \} 
 .\end{align}
We make assumptions of $ \mu $ and dynamics $ \mathbf{X}$ under $ \pP $. 

\medskip 

\noindent 
\As{\muTT} \ \ \ 
$ \mu $ is tail trivial. 

\noindent 
\As{\muAC} \ \ \ 
$ \pP \circ \mathsf{X}_t^{-1} \prec \mu $ for all $ \zzti $. 

\noindent 
\As{\sIn} \ \!\,
$ \pP (\mathsf{X} \in \WSsiNE ) = 1 $. 

\noindent \As{\nbj} \ 
$ \pP ( \mrXX < \infty ) = 1 $ for each $ r , T \in \mathbb{N} $. 

\medskip

We define the conditions \As{\muAC}, \As{\sIn}, and \As{\nbj} for 
a probability measure $ \widehat{\pP }$ on $ \WRN $ 
by replacing $ \mathsf{X}$ and $ \mathbf{X}$ by $ \mathsf{w}$ and $\mathbf{w}$, respectively. 
We remark here \As{\muAC}, \As{\sIn}, and \As{\nbj} are conditions depend only on the distribution of 
$ \mathbf{X}$.

We remark that, if $ \XB $ under $ \pP $ satisfies $ \As{\sIn}$, then $ \XB $ under $ \Ps $ satisfies 
\As{\sIn} for $ \pP \circ \mathbf{X}_0^{-1} $-a.s.\,$ \mathbf{s}$, where $ \Ps = \pP (\cdot | \mathbf{X}_0=\mathbf{s})$. The same holds for \As{\nbj}. This is however not the case for \As{\muAC}. 
Similarly as  \iFcs, we write  \As{\sIn}$ _{\mathbf{s}}$ and \As{\nbj}$_{\mathbf{s}} $ 
when we emphasize dependence on $ \mathbf{s}$.

It is known that all determinantal random point fields on continuous spaces are tail trivial \cite{bqs,ly.18,o-o.tail}. These results are a generalization of that of determinantal random point fields 
on discrete spaces \cite{BLPS,ST03,ly.03}. 

In \sref{s:7}, we deduce triviality of $ \TpathTSN $ from that of $ \TS $ through 
the tail $ \sigma $-field $ \Tpath (\sSS )$ of the unlabeled path space under these assumptions. 
We shall introduce the scheme carrying the tail $ \sigma $-field of $ \sSS $ 
to the tail $ \sigma $-field of $\WSN $. 

The assumption \As{\nbj} is crucial for the passage from the unlabeled dynamics 
 $ \mathsf{X} $ to the labeled dynamics $ \mathbf{X}$. 
If $ \pP ( \mrXX = \infty ) > 0 $, then we can not use this scheme. 
To catch the image for $ \mrXX = \infty $, 
we shall give an example of path $ \mathbf{w}$ such that $ \mrT (\mathbf{w}) = \infty $ in \rref{r:35}. 
This example indicates the necessity of \As{\Cfour}.

%[]
Let $ \XB $ be a weak solution of \ISDEb  defined on $ \OFPF $. 

If $ \Ps = \pP (\cdot | \mathbf{X}_0 = \mathbf{s})$ is a regular conditional probability, 
then $ \XB $ under $ \Ps $ is a weak solution of \ISDEb starting at $ \mathbf{s}$ 
for $ \pP \circ \mathbf{X}_0^{-1}  $-a.s.\,$ \mathbf{s}$. 

Conversely, suppose that $ \{ \Ps \} $ is a family of probability measures on $ \OFF $, 
given for $ \mathbf{m}$-a.s.\,$ \mathbf{s} \in \SN $, 
such that $ \XB $ under $ \Ps $ is a weak solution 
starting at $ \mathbf{s}$. 
If $ \Ps ( A )$ is $ \overline{\mathcal{B}(\SN )}^{\mathbf{m}} $-measurable in $ \mathbf{s} $ 
for any $ A \in \mathcal{B}(\WSN ) $, 
then $ \XB $ under $ \pP := \int \Ps \mathbf{m}(d\mathbf{s})$ is a weak solution of \ISDEb such that 
$\mathbf{m} =  \pP \circ \mathbf{X}_0^{-1}  $. 

Taking these into account, we introduce the following condition for a family  
of strong solutions $ \{ \Fs  \} $ of \ISDEb\ given for 
$  \pP \circ \mathbf{X}_0^{-1} $-a.s.\,$ \mathbf{s}$. 

\medskip
\noindent 
\As{$\mathbf{MF}$} \quad 
$ \pP (\Fs ( \mathbf{B} ) \in A ) $ is 
 $ \overline{\mathcal{B}(\SN )}^{ \pP \circ \mathbf{X}_0^{-1}} $-measurable in $ \mathbf{s} $ 
for any $ A \in \mathcal{B}(\WSN ) $. 

\medskip
\noindent 
%Thus the condition \As{$\mathbf{MF}$} implies $ \{ \Fs  \} $ is a measurable family of 
%strong solutions. 
For a family of strong solutions $ \{ \Fs \} $ satisfying \As{$\mathbf{MF}$} 
% given  satisfying \As{$\mathbf{MF}$} 
 we set 
\begin{align}\label{:5Az}&
\pPF  = \int \pP ( \Fs (\mathbf{B} ) \in \cdot ) \pP \circ \mathbf{X}_0^{-1} (d\mathbf{s}) 
.\end{align}

Let $ \XB $  be a weak solution under $ \pP $. 
Suppose that $ \XB $ is a unique strong solution  under 
$ \Ps $ for $ \pP \circ \mathbf{X}_0^{-1}$-a.s.\,$ \mathbf{s}$, where 
$ \Ps = \pP (\cdot | \mathbf{X}_0=\mathbf{s})$. 
Let $ \{ \Fs \}$ be the unique strong solutions given by $ \XB $ under $ \Ps $. 
Then \As{$\mathbf{MF}$} is automatically satisfied and 
\begin{align}\label{:5Ay}&
 \pPF  = \pP \circ \mathbf{X}^{-1}
.\end{align}
Indeed, $ \mathbf{B}$ is a Brownian motion under both $ \pP $ and $ \Ps $. 
Then %, and %we have, 
for $  \pP \circ \mathbf{X}_0^{-1} $-a.s.\,$ \mathbf{s}$ % and we have 
%Furthermore, $ \Fs $ is a unique strong solution. Then we have 
\begin{align}\label{:5Aa}&
\pP ( \Fs (\mathbf{B} ) \in \cdot )  = \Ps ( \Fs (\mathbf{B} ) \in \cdot ) = \Ps (\mathbf{X} \in \cdot )
.\end{align}
Hence we deduce \eqref{:5Ay} from \eqref{:5Az} and \eqref{:5Aa}.

\begin{definition}\label{dfn:43}
For a condition \As{$ \bullet $}, 
we say \ISDEb has a family of unique strong solutions $ \{ \Fs \} $ starting at $ \mathbf{s}$  for 
$ \pP \circ \mathbf{X}_0^{-1} $-a.s.\,$ \mathbf{s}$ 
under the constraints of \As{$\mathbf{MF}$} and \As{$ \bullet $} 
if  %$ \{ \Fs \} $ satisfies \As{$\mathbf{MF}$} and 
$ \{ \Fs  \} $ satisfies \As{$\mathbf{MF}$} and $  \pPF  $ satisfies \As{$ \bullet $}. Furthermore, 
\thetag{i} and  \thetag{ii} are satisfied.

\noindent\thetag{i} 
%$ \Fs $ is a strong solution  of \eqref{:50a}--\eqref{:50b} starting at $ \mathbf{s}$ 
%for $ \pP \circ \mathbf{X}_0^{-1}$-a.s.\,$ \mathbf{s}$. 
%\noindent \thetag{ii} 
For any weak solution $ (\hat{\mathbf{X}},\hat{\mathbf{B}})$ under $ \hat{P}$ 
of \eqref{:50a}--\eqref{:50b} with 
$$ 
\hat{P}\circ \hat{\mathbf{X}}_0^{-1} \prec P \circ \mathbf{X}_0^{-1} 
$$
%$ \hat{\mathbf{X}}_0 \elaw \mathbf{X}_0 $ 
 satisfying \As{$ \bullet $}, it holds that, for $\hat{P} \circ \hat{\mathbf{X}}_0^{-1}$-a.s. $  \mathbf{s}$, 
\begin{align}\notag &%\label{:5Ab}&
\hat{\mathbf{X}}=\Fs (\hat{\mathbf{B}}) \quad 
\text{ $ \hat{P}_{\mathbf{s}} $-a.s.}
,\end{align}
where $ \hat{P}_{\mathbf{s}} = \hat{P}(\cdot | \hat{\mathbf{X}}_0={\mathbf{s}})$. 

\noindent 
\thetag{ii} 
For an arbitrary $ \RdN $-valued $ \{\mathcal{F}_t\} $-Brownian motion $ \mathbf{B}$
defined on $ \OFPF $ with $ \mathbf{B}_0=\mathbf{0}$, 
$ \Fs (\mathbf{B})$ is a strong solution of \eqref{:50a}--\eqref{:50b}  satisfying \As{$ \bullet $} 
starting at $ \mathbf{s}$ 
for $ P \circ \mathbf{X}_0^{-1}$-a.s.\,$  \mathbf{s}$. 
\end{definition}

\begin{theorem} \label{l:5A} 
Assume \As{\muTT} for $ \mu$. 
Assume that \ISDEb has a weak solution $\XB $ under $ P $ satisfying \IASN. 
Then \ISDEb has a family of unique strong solutions $ \{ \Fs \} $ starting at $ \mathbf{s}$ 
 for $ \pP \circ \mathbf{X}_0^{-1} $-a.s.\,$ \mathbf{s}$ 
%under the constraints such that $ \{ \Fs  \} $ satisfies \As{$\mathbf{MF}$} and that $ \pPF $ satisfies \IASN. 
under the constraints of \As{$\mathbf{MF}$}, \IASN. 
\end{theorem}

The first corollary is a quench result. 
\begin{corollary}\label{l:53c1}
Under the same assumptions as \tref{l:5A} the following hold. 

\noindent \thetag{1} 
$ \XB $ under $ \Ps := \pP (\cdot | \mathbf{X}_0=\mathbf{s})$ is a strong solution starting at $ \mathbf{s}$ 
for $ \pP \circ \mathbf{X}_0^{-1} $-a.s.\,$ \mathbf{s}$. 

\noindent 
\thetag{2} Let $ (\mathbf{X}',\mathbf{B}') $ be any weak solution of \ISDEb 
defined on $ (\Omega',\mathcal{F}', P' ,\{ \mathcal{F}_t' \} )$ satisfying \IASN. 
Assume that $ \pP '\circ \mathbf{X}_0'^{-1} $ is absolutely continuous with respect to 
$ \pP \circ \mathbf{X}_0^{-1} $. 
Then $ (\mathbf{X}',\mathbf{B}') $ under $ \pP _{\mathbf{s}}' $ 
becomes a strong solution starting at $ \mathbf{s}$ satisfying $ \mathbf{X}'= \Fs (\mathbf{B}')$ 
for $ \pP '\circ \mathbf{X}_0'^{-1} $-a.s.\,$ \mathbf{s}$. Furthermore, 
\begin{align}&\notag %\label{:}&
 \Ps '\circ \mathbf{X}'^{-1}= \Ps \circ \mathbf{X}^{-1}
\end{align}
for $ \pP '\circ \mathbf{X}_0'^{-1} $-a.s.\,$ \mathbf{s}$. 
Here $ \pP _{\mathbf{s}}' = \pP '(\cdot | \mathbf{X}_0' =\mathbf{s})$. 

\noindent \thetag{3} 
 For any Brownian motion $ \mathbf{B}''$, $ (\Fs (\mathbf{B}''),\mathbf{B}'')$ 
becomes a strong solution of \ISDEb  satisfying \ISNss\ starting at $ \mathbf{s}$ 
for $ \pP \circ \mathbf{X}_0^{-1} $-a.s.\,$ \mathbf{s}$.  
\end{corollary}

The second corollary is an anneal result. 
\begin{corollary}\label{l:53c2}
Under the same assumptions as \tref{l:5A} the following hold. 

\noindent \thetag{1} 
The uniqueness in law of weak solutions of \ISDEb holds under the constraints of \IASN. 

\noindent \thetag{2} 
The pathwise uniqueness of weak solutions of \ISDEb holds under the constraints of \IASN. 
\end{corollary}

\begin{remark}\label{r:53z}
Because we exclude $ t = 0 $ in \As{\muAC}, 
\tref{l:5A} is valid even if $ \pP \circ \mathsf{X}_0^{-1} $ is singular to $ \mu $. 
\end{remark}

\begin{remark}\label{r:53a} 
We study ISDEs on $ \SN $. 
It is difficult to solve the ISDEs on $ \SN $ directly. One difficulty in treating $ \SN $-valued ISDEs is that $ \SN $ does not have any good measures. To remedy this situation, we introduce 
the representation of $ \SN $ as an infinite sequence of infinite-dimensional spaces: 
\begin{align}\label{:53l}&
\sSS ,\quad \sS \ts \sSS ,\quad \sS ^2 \ts \sSS ,\quad 
\sS ^3 \ts \sSS ,\quad \sS ^4 \ts \sSS 
 ,\quad \ldots 
. \end{align}
Each space in \eqref{:53l} has a good measure called 
the $ m$-Campbell measure (see \eqref{:90A}). % \cite[(2.27)]{o.isde}. 
Using \eqref{:53b}, we can rewrite \eqref{:53f} as 
\begin{align} \label{:53m}&
dY_t^{m,i} = 
 \sigma (Y_t^{m,i},\mathsf{Y}_t^{m,i\diamondsuit} + \mathsf{X}_t^{m*}) dB_t^i + 
 \bbb (Y_t^{m,i},\mathsf{Y}_t^{m,i\diamondsuit} + \mathsf{X}_t^{m*}) dt
.\end{align}
Thus $ (\mathbf{Y}^m,\mathsf{X}^{m*})$ is 
 an $ \SmSS $-valued process, and the scheme of infinite-dimensional spaces 
$ \{ \SmSS \}_{m=0}^{\infty} $ 
in \eqref{:53l} is useful. 
\end{remark}

\medskip 
\noindent

\begin{remark}\label{r:355}
We call \As{\Cfour} no big jump condition because 
for any path $ \mathbf{w} = (w^i)_{i\in\N } $ such that $\mathsf{m}_{r,T}(\mathbf{w})=\infty$ 
\begin{align}\notag &% 
\sup_{ i \in\N}\sup_{0\le s,t \le T}|w_s^i - w_t^i|
{\bf 1}_{[0,r]}(\min_{0\le u \le T}|w_u^i|)=\infty
\end{align}
and so for any $\delta>0$
\begin{align}\notag &% \label{BJ}
\sup_{ i \in\N}\sup_{\substack{0\le s,t \le T \\ |s-t|\le \delta}}|w_s^i - w_t^i|
{\bf 1}_{[0,r]}(\min_{0\le u \le T}|w_u^i|)=\infty,
\end{align}
which implies the existence of paths which visit $S_r$ 
during [0,T] and have modulus of continuity bigger than $\ell$ for any $\ell\in\N$.
%%We give an example of path such that $ \mrT (\mathbf{w}) = \infty $ in \rref{r:35} below. 
\end{remark}

\begin{remark}\label{r:35} 
An example of a path $ \mathbf{w}=(w^i)_{i\in\N } \in \WT (\RtwoN ) $ such that 
$ \mathsf{m}_{r,T}(\mathbf{w}) = \infty $ is as follows. 
%We take $ \lab (\mathbf{s})$ such that $ |s_i| \le |s_{i+1}|$ for all $ i \in \N $. 
Let $ t_i= \sum_{j=1}^i {2^{-j}}$ and 
\begin{align*}&
w_t^i = 
\begin{cases}
(1,i) & t \in [0, t_i]\cup [t_{i+1} , \infty) 
\\
\mathrm{linear} & [t_i,t_i + 2^{-i-2}]
\\
(0,0) & t= t_i + 2^{-i-2}
\\
\mathrm{linear} & [t_i + 2^{-i-2},t_{i+1}]
.\end{cases}
\end{align*}
All particles sit on the vertical line $ \{ (1,y); y \in \R^+ \} $ 
at time zero. The $ i $-th particle sits at $ (1,i)$, jumps off 
at time $ t_i $ and touch the origin at time $ t_i + 2^{-i-2}$. 
Then it springs up to the original position $ (1,i)$.
We need \As{\Cfour} to exclude this type of \lq\lq \bj '' paths. 
We note that $ \upath (\mathbf{w}) \not\in \WS $ although $ \mathbf{w}\in \WT (\RtwoN ) $. 
We conjecture that, if $ \upath (\mathbf{w}) \in \WSsiNE $, 
then $ \mathsf{m}_{r,T}(\mathbf{w}) < \infty $ is automatically satisfied. 
\end{remark}

\subsection{Main theorem II (\tref{l:5B}): $ \mu $ with non-trivial tail.}\label{s:55}

In this section, we relax \As{\muTT} of $ \mu $ by the tail decomposition of $ \mu $ as follows.

Let $ \mut $ be the regular conditional probability of $ \mu $ conditioned by $ \TS $: 
%the tail $ \sigma $-field $ \TS $: % defined as 
 \begin{align}\label{:54a}&
 \mut = \mu ( \, \cdot \, | \TS ) ( \aa ) 
 .\end{align}
Because $ \sSS $ is a Polish space, such a regular conditional probability exists and satisfies 
%the decomposition 
\begin{align}\label{:54b} &
\mu (\A ) = \int _{\sSS } \mut (\A ) \mu (d \aa )
. \end{align}
By construction, $ \mut (\A ) $ is a $ \TS $-measurable function in $ \aa $ 
for each $ \A \in \mathcal{B}(\sSS ) $.

Let $ \mathsf{H}$ be a subset of $ \SSsde $ in \eqref{:50c} and $ \mu (\mathsf{H}) = 1 $. 
We assume there exists a version of $ \mut $ with a subset of $ \mathsf{H}$, 
denoted by the same symbol $ \HHz $, such that $ \mu (\HHz )= 1$ and that, for all $\mathsf{a}\in \HHz $, 
\begin{align}\label{:54c}& 
 \mut \text{ is tail trivial}
,\\\label{:54d}&
 \mut (\{ \bb \in \sSS ; \mut = \mu _{\mathrm{Tail}}^{\bb } \} ) = 1 
,\\\label{:54e}&
\text{$ \mut $ and $ \mu_{\mathrm{Tail}}^{\mathsf{b}}$ 
are mutually singular on $ \TS $ 
if $ \mut \not=\mu_{\mathrm{Tail}}^{\mathsf{b}}$}
.\end{align}

Let $ \stackrel[\TTT ]{}{\sim} $ be the equivalence relation such that 
$ \mathsf{a} \stackrel[\TTT ]{}{\sim} \mathsf{b} $ if and only if 
\begin{align} \notag &% \label{:54f} &
 \mut = \mu _{\mathrm{Tail}}^{\bb } 
.\end{align}
Let $ \HHa = \{ \mathsf{b}\in \mathsf{H}; \mathsf{a}\stackrel[\TTT ]{}{\sim} \mathsf{b}\} $. 
Then $ \HHz $ can be decomposed as a disjoint sum 
\begin{align}\label{:54g} &
 \HHz = \sum_{ [\aa ] \in \HHz / \stackrel[\TTT ]{}{\sim}} \HHa 
.\end{align}
From \eqref{:54d}, we see that $ \mut ( \HHa ) = 1 $ for all $ \mathsf{a} \in \HHz $. 

For a labeled process $ \mathbf{X}=(X^i) $ on $ \OFPF $ 
we set $ \mathsf{X}_t=\sum_i \delta_{X_t^i}$ as before. 
Let $ \lab $ be a label. 
We assume $ \mathbf{X}_0 = \lab \circ \mathsf{X}_0 $ for $ \pP $-a.s. 
We note that plural labels satisfy the relation 
$ \mathbf{X}_0 = \lab \circ \mathsf{X}_0 $ for $ \pP $-a.s.\, in general. 
%We assume $ \mu $, $ \mut $, and $ \HHz $ satisfy the properties as above. 
For $ \mu $ as above we assume 
\begin{align}
\label{:54h}&
\mu = \pP \circ \mathsf{X}_0^{-1}
.\end{align}
We set $ \pP ^{\aa } = \pP (\, \cdot \, | \mathsf{X}_0^{-1}(\TS )) |_{\mathsf{X}_0 = \aa } $. 
Then by \eqref{:54a} and \eqref{:54b} 
\begin{align}%\notag &%
\label{:54j}&
 \pP ^{\aa } = \int P (\cdot | \mathsf{X}_0 = \mathsf{s}) \mut (d\mathsf{s})
.\end{align}
We can rewrite $ \pP ^{\aa } $ as 
\begin{align}%\notag &% 
\label{:54i}&
 \pP ^{\aa } = \int P (\cdot | \mathbf{X}_0 = \mathbf{s}) \mut \circ \lab ^{-1} (d\mathbf{s})
.\end{align}
From \eqref{:54a} and \eqref{:54h} we easily see $ \mut = \pP ^{\aa }\circ \mathsf{X}_0^{-1}$ and 
\begin{align} \label{:54J}&
\mut \circ \lab ^{-1} = \pP ^{\aa }\circ \mathbf{X}_0^{-1}
.\end{align}
Let $ \Ps ^{\aa } = \pP ^{\aa }(\cdot | \mathbf{X}_0=\mathbf{s})$. 
Then $ \Ps ^{\aa } = P (\cdot | \mathbf{X}_0 = \mathbf{s}) $ 
for $ \mut \circ \lab ^{-1}$-a.s. $ \mathbf{s}$ and for $ \mu $-a.s. $ \aa $. 

\begin{theorem} \label{l:5B} %
Assume $ \mu $, $ \mut $, and $ \HHz $ satisfy %the conditions above: namely, 
 \eqref{:54c}--\eqref{:54e} for all $ \aa \in \HHz \subset \SSsde $, $ \mu (\HHz ) = 1 $, 
and \eqref{:54h}. 
Assume that $\XB $ under $ \pP $ is a weak solution 
of \ISDEb satisfying \As{\iFc}, \As{\sIn}, and $ \As{\nbj}$. 
Assume that, for $ \mu $-a.s.\ $ \aa $, $\XB $ under $ \pP ^{\aa } $ satisfies \As{\muAC} for $ \mut $. 
Then, for $ \mu $-a.s. $ \aa $, \ISDEb has a family of unique strong solutions $ \{ \Fsa \}$ 
starting at $ \mathbf{s}$ for $ \pP ^{\aa } \circ \mathbf{X}_0^{-1} $-a.s.\,$ \mathbf{s}$ 
%under the constraints 
%such that $ \{ \Fsa  \} $ satisfies \As{$\mathbf{MF}$} and that $ \pPFa $ satisfies \IASNa. 
under the constraints of \As{$\mathbf{MF}$}, \IASNa. 
\end{theorem}

\begin{remark}\label{r:5B} \thetag{1} 
We shall prove in \lref{l:X2} that a version $ \{ \mut \} $ satisfying \eqref{:54c}--\eqref{:54e} exists 
if $ \mu $ is a quasi-Gibbs random point field satisfying \As{A2} in \sref{s:5c}. 
All examples in the present paper are such quasi-Gibbs random point fields. % satisfying \As{A2}. 
\\\thetag{2}
The unique strong solution $ \Fsa $ in \tref{l:5B} yields the corollaries similar to \corref{l:53c1} 
and \corref{l:53c2}. 
In particular, for $ \mu $-a.s.\ $ \aa $, $ \XB = (\Fsa (\mathbf{B}),\mathbf{B})$ under $ \Ps ^{\aa }$ 
for $\mut \circ \lab ^{-1} $-a.s.\ $ \mathbf{s}$. 
We note here $ \mut \circ \lab ^{-1} = \pP ^{\aa }\circ \mathbf{X}_0^{-1}$ by \eqref{:54J}. 
\\\thetag{3}
In \tref{l:5B}, the assumption \lq\lq $\XB $ under $ \pP ^{\aa } $ satisfies 
\As{\muAC} for $ \mut $ for $ \mu $-a.s.\ $ \aa $'' is critical, 
and does not hold in general. We shall give sufficient conditions in \tref{l:QA} and \tref{l:QB}. 
From these theorems we deduce all the examples in the present paper satisfy the assumption. 
\end{remark}

%%%%%%%%%%%%%%%%%%%%%%%%%%%%%%%%%%%%%%%%%%%%%%%%%%%%%%%%

\Ssection{The role of the five assumptions: 
\As{\iFc}, \As{TT}, \As{\muAC}, \As{\Cthree}, and \As{\Cfour}.}
\label{s:33}

In \sref{s:5d}, we introduced the five significant assumptions: 
\As{\iFc}, \As{TT}, \As{\muAC}, \As{\Cthree}, and \As{\Cfour}. 
In this subsection, we explain the role of these assumptions in the proof of \tref{l:5A}. 
We also explain the role of other main assumptions used in the present paper. 

To prove \tref{l:5A}, we use the strategy introduced in \sref{s:1}. 
One of the critical points of the strategy is the reduction of ISDE to the infinite system of finite-dimensional stochastic differential equations (SDEs). For this, we use the pathwise unique strong solutions of the finite-dimensional SDEs associated with ISDE. 
The condition \As{\iFc} claims the finite-dimensional SDEs have such a unique strong solution. 
So the \As{\iFc} is pivotal to the reduction of ISDE to the IFC schemes.

Another critical point of the proof is tail triviality of the labeled path space under the distribution of the weak solution.
We shall deduce this from tail triviality of $ \mu $, denoted by \As{TT},
in a general framework as the second tail theorem in \sref{s:7}. 

The key idea for this is the passage from \As{TT} to that of the path space of the labeled dynamics. Because of \As{\muAC}, we deduce triviality of 
the tail $ \sigma $-field of $ \mathsf{S}$ under the single time marginal distributions of $ \mathsf{X}$. 
From this we shall deduce tail triviality of the unlabeled path space 
under the finite-dimensional distributions of the unlabeled dynamics $ \mathbf{X}$, 
where the meaning of tail is spatial (see \eqref{:74y} for definition). 

We next deduce tail triviality of the labeled path space from that of the unlabeled path space. 
The map $ \lpath $ in \eqref{:23i} from the unlabeled path space to the labeled one 
plays an essential role in our argument. 
The assumption \As{\sIn} is necessary for the construction of this map. 
Here \As{\sIn} is an abbreviation of \lq\lq unlabeled path spaces on \textbf{s}ingle, \textbf{i}nfinite configurations with \textbf{n}o explosion of tagged particles''. 
To carry out the passage, 
%from the unlabeled path space to the labeled path space, 
we shall use \As{\Cfour} in addition to \As{\Cthree}. 
%Thus \As{\Cfour} guarantees the passage from the unlabeled dynamics to the labeled dynamics. 
%
%Indeed, we shall use \As{\Cfour} and \As{\muAC} to prove Second tail theorem (\tref{l:70}) in \sref{s:7}. 

In \sref{s:6}, we shall deduce the existence of a unique strong solution of ISDE 
from tail triviality of labeled path space in \tref{l:61}. 
We call \tref{l:61} the first tail theorem. 
The conditions in \tref{l:61} are denoted by \ASTpath\ and \ASTpatH.

The assumptions \As{\iFc}, \As{\Cone}, \As{\Ctwo}, \As{\Cthree}, and \As{\Cfour} are used 
in the proof of \tref{l:70} (the second tail theorem) in \sref{s:7}. 
\As{\Cone} and \As{\Ctwo} are the assumptions of \tref{l:74}, which claims 
 triviality of the labeled path space at the cylindrical level. 
\As{\Cone}, \As{\Ctwo}, \As{\Cthree}, and \As{\Cfour} are the assumptions of \tref{l:77}, 
which proves \AsDD. 
%%[]

Conditions \As{A1}--\As{A4} given in \sref{s:D} and \sref{s:E} are related to Dirichlet forms; these are the conditions for random point fields $ \mu $ from the viewpoint of Dirichlet form theory. 
Recall that \As{\iFc} is the most critical condition for the theory. 
We shall give the feasible sufficient condition of \As{\iFc} in terms of 
the assumptions \As{B1}--\As{B2} and \As{C1}--\As{C2} in \sref{s:9}. 
We shall present them in \tref{l:UA} and \tref{l:UB}. 

We present a table for these conditions in Table \ref{list}.

\begin{table}%[htb]
\caption{List of Conditions} \label{list}
 \begin{tabular}{|l|c|c||l|} \hline
 & Assumption & Place & The main role of the assumption
 \\ \hline %\hline
 & \As{\iFc},\ \iFcs & \sref{s:5d} & 
\tref{l:5A}: To construct a scheme of finite-dimensional SDEs\\ \cline{2-4}
%\\\hline
 & \As{TT} & \sref{s:5d} & \tref{l:5A} and the second tail theorem (\tref{l:70})\\ \cline{2-4}
 & \As{AC} & \sref{s:5d} & \tref{l:5A} and the second tail theorem (\tref{l:70})\\ \cline{2-4}
 & \As{SIN} & \sref{s:5d} & \tref{l:5A} and the second tail theorem (\tref{l:70})\\ \cline{2-4}
 & \As{\Cfour} & \sref{s:5d} & \tref{l:5A} and the second tail theorem (\tref{l:70})\\ \cline{2-4}
 & \As{$\mathbf{MF}$} & \sref{s:5d} & \tref{l:5A} and \dref{dfn:43}: To construct $ \pPF $ from $ \{ \Fs \} $\\ \cline{2-4}
 %%& IFC solution & \dref{d:59} & A weak solution satisfying \As{\iFc}
%Assumption of the first and second tail theorems 
%\\ \cline{2-4}
 % & \As{\iFc} & Def 5.9 & A condition for IFC solutions 
%\\
\hline
 & \ASTpath & \sref{s:61} & The first tail theorem (\tref{l:61}) 
\\ \cline{2-4}
 & \ASTpatH & \sref{s:61} & The first tail theorem (\tref{l:61}) 
%\\ \cline{2-4}
\\\hline
 & \As{\Done} & \sref{s:72} & \tref{l:79}, which yields the second tail theorem \\ \cline{2-4}
 & \As{\Dtwo} & \sref{s:72} & \tref{l:79}, which yields the second tail theorem 
 \\\hline
 & \As{A1}-\As{A3} & Section 9 & To construct weak solutions of (ISDE) via Dirichlet forms
\\ \cline{2-4}
& \As{A4} & Section 10 & To derive \As{SIN} and \As{\Cfour} \\ \cline{2-4}
 & \As{B1}-\As{B2} & Section 11 & To derive \As{\iFc} \\ \cline{2-4}
 & \As{C1}-\As{C2} & Section 11. 3 & To derive \As{B2} \\\hline 
 \end{tabular}
\end{table}

\Section{Solutions and tail $ \sigma $-fields: First tail theorem (\tref{l:61})} 
\label{s:6}

This section proves the existence of a strong solution,  the pathwise uniqueness 
of solutions, and that the ISDE \ISDEsixb has a unique strong solution. 
%This section proves the existence of a strong solution and the pathwise uniqueness 
%of the solution of \ISDEsixb below. 
% We shall present a new formulation of solutions of ISDEs 
% called the IFC solution (see \dref{dfn:61}). Using this formulation, we shall derive 
% the existence and pathwise uniqueness of a strong solution in \tref{l:61}. 
% 
The ISDEs studied in this section are more general than those in \sref{s:1} and \sref{s:5}. 
Naturally, the ISDEs in the previous sections are typical examples 
that our results (\tref{l:61}, \tref{l:64}, and \tref{l:65}) can be applied to. 

Throughout this section, $ \mathbf{X} = (X^i)_{i\in\N }$ is 
an $ \SN $-valued, continuous $ \Ft $-adapted processes defined on $\2 $ 
starting at $ \mathbf{s}$, which is indicated by the subscript $ \mathbf{s}$ 
in $ \Ps $. 
$ \mathbf{B} =(B^i)_{i\in\N }$ is an $ \RdN $-valued, standard 
$ \Ft $-Brownian motion starting at the origin. 
$ \PBr $ is the distribution of $ \mathbf{B}$. 
%the standard Brownian motion on $ \RdN $ starting at the origin $ \mathbf{0} \in \RdN $. 
Thus, 
\begin{align}
\label{:62y}& 
\Ps (\mathbf{X}_0=\mathbf{s}) = 1 ,\quad 
\Ps (\mathbf{B} \in \cdot ) = \PBr 
.\end{align}
We shall fix the initial starting point $ \mathbf{s}$ throughout \sref{s:6}. 
\Ssection{General theorems of the uniqueness and existence of strong solutions of ISDEs.}
 \label{s:61}
In this subsection, we introduce ISDE \ISDEsixb and state one of the main theorems (\tref{l:61}: First tail theorem). 

Let $ \Wsol $ be a Borel subset of $ \WSN $. 
Let $ \mathcal{B} (\Wsol )$ be the Borel $ \sigma $-field of $ \Wsol $. 
Let $ \mathcal{B}_t (\Wsol ) $ be the sub $ \sigma $-field of $ \mathcal{B}(\Wsol ) $ 
such that 
$$ \mathcal{B}_t (\Wsol ) = \sigma [\mathbf{w}_u; 0\le u \le t \, , \, 
\mathbf{w}\in \Wsol ]
.$$
Following \cite{IW} in finite dimensions, we shall introduce SDEs in infinite dimensions. 

\begin{definition}\label{d:60}
$ \mathcal{A}^{d,r} $ is the set of all functions 
$ \map{\alpha }{[0,\infty )\ts \Wsol }{\Rd \otimes \R ^r }$ such that \\
\thetag{1} $ \alpha $ is 
$ \mathcal{B}([0,\infty ))\ts \mathcal{B}(\Wsol ) / \mathcal{B} 
(\Rd \otimes \R ^r )$-measurable, 
\\\thetag{2} 
$ \Wsol \ni \mathbf{w}\mapsto \alpha (t,\mathbf{w}) \in \Rd \otimes \R ^r $ is 
$ \mathcal{B}_t (\Wsol ) / \mathcal{B} 
(\Rd \otimes \R ^r )$-measurable for each $ t \in [0,\infty )$. 
%We write $ \alpha (\cdot ) _t = \alpha (t, \cdot )$. 
\end{definition}

Let $ \map{\sigma ^i }{\Wsol }{\WT (\R ^{d^2}) }$ and $ \map{b^i}{\Wsol }{\WT (\Rd ) }$ such that 
$ \sigma ^i (\mathbf{w})_t \in \mathcal{A}^{d,d} $ and 
$ b^{i} (\mathbf{w})_t \in \mathcal{A}^{d,1} $. 
We assume $ \sigma^{i} \in \mathcal{L}^2 $ and $ b^{i} \in \mathcal{L}^1 $, 
where $ \mathcal{L}^p $ is the same as \dref{d:51}. 
We introduce the ISDE on $ \SN $ of the form 
\begin{align}\label{:60a} & 
dX_t^i = \sigma ^{i}(\mathbf{X})_t dB_t^i + 
b^{i}(\mathbf{X})_t dt \quad (i\in\mathbb{N})
,\\\label{:60b}&
\mathbf{X} \in \Wsol 
,\\\label{:60c} &
\mathbf{X}_0 = \mathbf{s}
.\end{align}
%Here $ \mathbf{X} = (X^i)_{i\in\mathbb{N}}$ and $ \mathbf{B} = (B^{i})_{i\in\mathbb{N}}$. 
Here $ (\mathbf{X},\mathbf{B})$ is defined on $ \2 $, 
%$ \mathbf{X} = (X^i)_{i\in\N }$ is an $ \SN $-valued, continuous $ \Ft $-adapted processes, 
and $ \mathbf{B}=(B^i)_{i\in\mathbb{N}}$ is an $ \RdN $-valued, $ \Ft $-Brownian motion 
as before. 

The definition of a weak solution and a strong solution, and the related notions 
are similar to those of \sref{s:5}.

We remark that, in this section, 
we do not assume $ \Wsol = \upath ^{-1} (\WT (\mathsf{W})) $ 
for some $ \mathsf{W} \subset \sSS $ unlike the previous sections. 
This is because we intend to clarify the relation between the strong and pathwise notions of solutions of ISDE and tail triviality of the labeled path space. 
Indeed, our theorems (\tref{l:61} and \tref{l:64}) 
clarify a general structure of the relation between 
the existence of a strong solution and the pathwise uniqueness of the solutions of 
ISDE \ISDEsixc and 
 triviality of the tail $ \sigma $-field of the labeled path space 
$ \TpathTSN $ defined by \eqref{:61a} below. 

We make a minimal assumption for this structure. As a result, ISDEs in this section are much general than before. In \sref{s:7}, we return to the original situation, and deduce tail triviality of the labeled path space from that of the configuration space.

The correspondence between ISDE {\ISDEb} and \ISDEsixB\ is as follows. 
\begin{align}\notag &% \label{:60d}%\notag %
\Wsol = \{\mathbf{w}\in \WT (\SSSsde ); \mathbf{w}_0 \in \lH \}, 
\\ \notag &% \label{:60e}& 
\sigma^i(\mathbf{X})_t = \sigma (X_t^i,\mathsf{X}_t^{i\diamondsuit}) ,\ 
b ^i(\mathbf{X})_t = \bbb (X_t^i,\mathsf{X}_t^{i\diamondsuit}) 
.\end{align}
Here $ \mathsf{H}$, $ \sigma $, and $ \bbb $ are given in ISDE \eqref{:50a} and \eqref{:50c}. 
Moreover, $ \Wsol $ corresponds to both $ \lH $ and $ \ulab ^{-1}(\SSsde ) $. 

The final form of our general theorems (Theorems \ref{l:5A}--\ref{l:5B}) 
are stated in terms of random point fields. 
We emphasize that there are many interesting random point fields 
satisfying the assumptions, such as the sine, Airy, Bessel, and Ginibre random point fields, and all canonical Gibbs measures with potentials of Ruelle's class (with suitable smoothness of potentials such that the associated ISDEs make sense). 

We take the viewpoint not to pose the explicit conditions of the coefficients 
$ \sigma ^{i}$ and $ b ^{i}$ to solve ISDE \ISDEsixcc, 
but to assume the existence of a weak solution and 
the pathwise uniqueness of solutions of the associated infinite system of 
finite-dimensional SDEs instead. 

For a given $ \XB $ defined on $ \2 $ satisfying \eqref{:62y}, 
we introduce the infinite system of the finite-dimensional SDEs \3. 
\begin{align}\label{:60h}&
dY_t^{m,i} = \sigma ^{i} (\YopX )_t d B_t^i + 
b^{i}(\YopX )_t dt \quad ( i=1,\ldots,m) 
,\\\label{:60i}&
 \YopX \in \Wsols 
,\\\label{:60j}&
 \YopXz = \mathbf{s}
,\end{align}
where $ \mathbf{Y}^{m} = (Y^{m,1},\ldots,Y^{m,m})$ 
is an unknown process defined on $ \2 $, and 
$ \YopX = (Y^{m,1},\ldots,Y^{m,m},X^{m+1},X^{m+2},\ldots ) 
$.

%From \eqref{:60j}, we deduce that $ \YopX _0 = \mathbf{s}$, and hence, 
The process $\mathbf{Y}^{m} $ 
denotes a solution of \ISDEhij starting at $\mathbf{s}^{m} = (s_1,\ldots,s_m)$. 
The notion of a strong solution and a unique strong solution of \ISDEhij 
for $ \XB $ \uPs\ 
is defined as \dref{d:41} and \dref{d:42} with an obvious modification. 
Let $ {\mathbf{B}}^m =({B}^i)_{i=1}^m $ be 
%the $ \Ft $-Brownian motion consisting of 
the first $ m $-components of the $ \Ft $-Brownian motion 
$ \mathbf{B}=(B^i)_{i=1}^{\infty}$. 
%We make the assumption: 
The following assumption corresponds to {\iFcs} in \sref{s:5d}. 

\medskip 
\noindent 
\iFcs\ \ 
SDE \3 has a unique strong solution 
$ \mathbf{Y}^m = \Fmss ({\mathbf{B}}^m,\mathbf{X}^{m*}) $ 
for $ \XB $ under $ \Ps $ for each $ m \in \N $. 
%starting at $ \mathbf{s}$ defined on $ \2 $ 

\begin{remark}\label{r:60z}
 The meaning of SDE \3 is not conventional because the coefficients include additional randomness 
$ \mathbf{X}^{m*}$, that is, 
$ \mathbf{X}^{m*}$ is interpreted as ingredients of the coefficients of the SDE \eqref{:60h}. 
Furthermore, $ \mathbf{B}^m $ and $ \mathbf{X}^{m*}$ can depend on each other. See \rref{r:31}. 
Another interpretation of SDE \3 is that $ (\mathbf{B}^m,\mathbf{X}^{m*})$ is regarded as input to the system rather than the interpretation such that 
$ \mathbf{X}^{m*} $ is regarded as a part of the coefficient. 
Solving SDE \3 then means constructing a function of $ (\mathbf{B}^m,\mathbf{X}^{m*})$. 
Thus a strong solution means a functional of $ (\mathbf{B}^m,\mathbf{X}^{m*})$. 
\end{remark}

We assume $ \XB $ under $ \Ps $ is \7. 
Then we obtain the crucial identity:
\begin{align}
\label{:60m}&
\mathbf{Y}^m = \mathbf{X}^m
.\end{align}
As we saw in \sref{s:1}, \eqref{:60m} plays an important role in the whole theory.

We set $ \mathbf{w}^{m*}=(w^i)_{i=m+1}^{\infty}$ for 
$ \mathbf{w}=(w^i)_{i=1}^{\infty}\in \WSN $. Let 
\begin{align}\label{:61a}&
\TpathTSN = \bigcap_{m=1}^{\infty} \sigma [\mathbf{w}^{m*}] 
.\end{align}
By definition, $ \TpathTSN $ is the tail $ \sigma $-field of $ \WSN $ with respect to the label. For a probability measure $ P $ on $ \WSN $, we set 
\begin{align}\label{:61b}&
 \TpathTSNone ;P ) = \{ \mathbf{A} \in \TpathTSN \, ;\, P (\mathbf{A}) = 1 \} 
.\end{align}

For the continuous process $ \XB $ given at the beginning of this section, we set 
\begin{align}
\label{:61B}&
\Pts = \Ps \circ \XB ^{-1}
.\end{align}
By definition $ \Pts $ is a probability measure on $ \WW $. 
We denote by $ (\mathbf{w},\mathbf{b}) $ generic elements in $ \WW $. 
Let $ \Ptsb $ be the probability measure on $ \WSN $ 
given by the regular conditional probability $ \Ptsb $ of $ \Pts $ such that 
\begin{align}\label{:61c}& 
 \Ptsb (\cdot ) = \Pts (\mathbf{w} \in \cdot \, | \bB ) 
%\equiv \Pts (\, \cdot \, | \mathbf{B}) \circ \Pitwo ^{-1}
.\end{align}
We introduce the conditions. 
% on $ \Ps $ and $ \XB $, which depend only on the distribution $ \Pts $: % 

\smallskip 
\noindent 
\ASTpath\ 
$\TpathTSN $ is $ \Ptsb $-trivial for $ \PBr $-a.s.\! $ \bB $. 
\\
\ASTpatH\ 
$ \TpathTSNone ;\Ptsb ) $ is independent of the distribution $ \Pts $ 
for $ \PBr $-a.s.\! $ \bB $ if $ \XB $ under $ \Ps $ is \7 and \ASTpath.

\smallskip
\noindent 
In other words, \ASTpatH\ means 
\smallskip 

\noindent 
\ASTpatH'\ 
 If $ \XB $ under $ \Ps $ and 
$ (\mathbf{X}',\mathbf{B}')$ under $ \Ps ' $ satisfy \iFcs\ and \ASTpath, 
then 
\begin{align*}&\quad \quad 
 \TpathTSNone ;\Ptsb ) = \TpathTSNone ;\Ptsb' ) 
\quad \text{ for $ \PBr $-a.s.\! $ \bB $}
.\end{align*}

\smallskip 

\noindent 
In this sense, \ASTpatH\ is a condition for ISDE \ISDEsixc rather than a specific solution $ \XB $ under $ \Ps $. 

\begin{remark}\label{r:61} \thetag{1} 
The conditions \ASTpath\ and \ASTpatH\ depend on $ \mathbf{s}$, which is the initial starting point of 
$ \mathbf{X}$. To indicate this we put the subscript $ \mathbf{s}$. 
\\\thetag{2} 
We emphasize that we fix $ \mathbf{s}$ throughout \sref{s:6}, 
while we do not fix $ \mathbf{s}$ in \sref{s:7}. 
We shall present a sufficient condition such that 
\ASTpath\ and \ASTpatH\ hold for a.s.\, $ \mathbf{s}$ with respect to the initial distribution of 
$ \mathbf{X}_0$. We remark that, unlike \sref{s:6}, 
$ \mathbf{X}$ does not necessarily start at a fixed single point in \sref{s:7}. 
\end{remark}

We note that \ASTpatH\ implies $ \TpathTSNone ;\Ptsb ) $ 
depends only on $ \bB $ for $ \PBr $-a.s.\! $ \bB $. 

We now state the main theorem of this section, which we shall prove in \sref{s:62}. 
In \tref{l:61} and \corref{l:61B} we consider ISDE \ISDEsixbb. 
Recall the notions of strong solutions starting at $ \mathbf{s} $ 
given by \dref{d:55} and 
a unique strong solution under the constraint of 
a condition \As{$ \bullet$} given by \dref{d:58}. 
\begin{theorem}[First tail theorem] \label{l:61} 
The following hold. 
\\\thetag{1} 
Assume that $ \XB $ under $ \Ps $ is \7 and \ASTpath. 
Then $ \XB $ under $ \Ps $ is a strong solution of \ISDEsixb. 
\\
\thetag{2} 
Make the same assumptions as \thetag{1}. 
We further assume that ISDE \ISDEsixc satisfies \ASTpatH. 
Then \ISDEsixb has a unique strong solution $ \Fs $ under the constraint of \iFcs\ and \ASTpath. 
\end{theorem}
\begin{corollary}\label{l:61B}
%Assume that $ \XB $ under $ \Ps $ is \7 and \ASTpath. Assume that ISDE \ISDEsixc satisfies \ASTpatH. 
Make the same assumptions as \tref{l:61} \thetag{2}. 
Then the following hold. 

\noindent \thetag{1}
For any weak solution $ (\mathbf{X}',\B ')$ of \ISDEsixb satisfying \iFcs\ and \ASTpath, 
it holds that $ \mathbf{X}'=\Fs (\B ')$. 
\\\thetag{2}
For any Brownian motion $ \mathbf{B}''$, $ \Fs (\mathbf{B}'') $ is a strong solution of \ISDEsixb satisfying 
\iFcs\ and \ASTpath. 
\\\thetag{3}
The uniqueness in law of weak solutions of \ISDEsixb holds under the constraint of \iFcs\ and \ASTpath. 
\end{corollary}

\Ssection{Infinite systems of finite-dimensional SDEs with consistency. }\label{s:62}
In this section we prove \tref{l:61}. 
Let $ \XB $ be a pair of an $ \Ft $-adapted continuous process $ \mathbf{X}$ and 
an $ \Ft $-Brownian motion $ \mathbf{B} $ defined on $ \2 $ satisfying 
\eqref{:62y} as before. 
We assume that $ \Ps $ satisfies: 
\begin{align}\notag &% \label{:62x}&
\Ps (\XB \in \Swc ) = 1 
.\end{align}
We denote by $ \PPPsb $ the regular conditional probability of $ \Ps $ 
conditioned by the random variable $ \B $ such that 
\begin{align}\label{:62z}& 
 \PPPsb (\cdot ) = \Ps ( (\mathbf{X},\B ) \in \cdot \, |\B = \bB ) 
.\end{align}
By construction $ \PPPsb (\Wsols \times \{ \bB \} ) = 1 $. 
This follows from the fact that the regular conditional probability in 
\eqref{:62z} is conditioned by the random variable $ \B $. 
We refer to \cite[(3.1) 15p]{IW} for proof. 
We set 
\begin{align}\label{:62c}&
\PPPs (\cdot ) = \int_{}\PPPsb (\cdot )\PBr (d\bB )
.\end{align}
From \eqref{:62y}, \eqref{:62z}, and \eqref{:62c} we have a representation of 
$ \Pts = \Ps \circ \XB ^{-1}$ such that for any 
$ C \in \mathcal{B}(\Wsols ) \times \mathcal{B}(\WRNz ) $ 
\begin{align}
\label{:62p}&
\Pts ( C ) = \PPPs ( C )
.\end{align}

\begin{remark}\label{r:62}
The probability measure $ \Ptsb $ in \eqref{:61c} resembles $ \PPPsb $ in \eqref{:62z}, 
and they are closely related to each other. The difference is that 
$ \Ptsb $ is a probability measure on $ \WSN $, whereas $ \PPPsb $ is on $ \WSN \ts \WRNz $. 
\end{remark}

In general, $ \overline{\mathcal{F} }^P$ denotes the completion of 
the $ \sigma $-field $ \mathcal{F} $ with respect to probability $ P $. 
Recall that $ \TpathTSN = \cap_{m\in\N }\sigma [\mathbf{w}^{m*}]$. 
Then it is easy to see that 
\begin{align}\label{:62s}&
\TpathTSN \times \mathcal{B}(\WRNz ) \subset 
\bigcap_{m\in \N }%^{\infty}
\Big\{ \sigma [\mathbf{w}^{m*}] \times \mathcal{B}(\WRNz ) \Big\} 
.\end{align}
Hence taking the closure of both sides in \eqref{:62s} with respect to $ \Pts $ we have 
\begin{align}\label{:62t} 
\overline{\TpathTSN \times \mathcal{B}(\WRNz )}^{\Pts } & \subset 
\overline{\bigcap_{m\in\N } \Big\{
\sigma [\mathbf{w}^{m*}] \times \mathcal{B}(\WRNz ) \Big\}
}^{\Pts }
\\ \notag & =
\bigcap_{m\in\N } 
\overline{\sigma [\mathbf{w}^{m*}] \times \mathcal{B}(\WRNz ) 
}^{\Pts }
.\end{align}

The equality in \eqref{:62t} follows from a general fact such that the completion of the countable intersection of 
$ \sigma $-fields with decreasing property coincides 
with the countable intersection of the completion of $ \sigma $-fields with decreasing property. 
 Indeed, if 
\begin{align}& \label{:62S}
K \in \overline{\bigcap_{m\in\N } \Big\{
\sigma [\mathbf{w}^{m*}] \times \mathcal{B}(\WRNz ) 
\Big\}
}^{\Pts }
,\end{align}
then there exist $ A , B \in \bigcap_{m\in\N } \sigma [\mathbf{w}^{m*}] \times \mathcal{B}(\WRNz ) $ such that 
$ A \subset K \subset B $ and that $ \Pts (B \backslash A ) = 0 $. 
Then $ A, B \in \sigma [\mathbf{w}^{m*}] \times \mathcal{B}(\WRNz )$ for each $ m $. 
This deduces 
\begin{align}&\label{:62T}
K \in \bigcap_{m\in\N } \overline{\sigma [\mathbf{w}^{m*}] \times \mathcal{B}(\WRNz ) }^{\Pts } 
.\end{align}
Conversely, if \eqref{:62T} holds, 
%$ K \in \bigcap_{m\in\N } \overline{\sigma [\mathbf{w}^{m*}] \times \mathcal{B}(\WRNz ) }^{\Pts }$, 
then there exist %$ A_m , B_m $ 
$ A_m , B_m \in \sigma [\mathbf{w}^{m*}] \times \mathcal{B}(\WRNz ) $ such that 
$ A_m \subset K \subset B_m $ and that $ \Pts (B_m \backslash A_m ) = 0 $ for each $ m \in \N $. 
It is clear that 
$ \limsup_m A_m $ and $ \liminf_m B_m $ are tail events such that 
\begin{align*}&
 \Pts ( \liminf_m B_m \backslash \limsup_m A_m ) \le \sum_{m=1}^{\infty}
\Pts (B_m \backslash A_m) = 0
.\end{align*}
Hence \eqref{:62S} holds. 
%$ K \in \overline{\bigcap_{m\in\N } \sigma [\mathbf{w}^{m*}] \times \mathcal{B}(\WRNz ) }^{\Pts }$. 

We need the inverse inclusion of \eqref{:62t}, which does not hold in general. 
Hence we introduce further completions of 
$ \TpathTSN \times \mathcal{B}(\WRNz )$ as follows. 
We set 
\begin{align}\label{:62Q}&
%\mathcal{K} = \overline{\TpathTSN \times \mathcal{B}(\WRNz )}^{\Pts } %[]変更注意
\mathcal{K} = {\TpathTSN \times \mathcal{B}(\WRNz )}%^{\Pts }
\end{align}
and $ \mathcal{I} $ intuitively given by 
\begin{align}
 \label{:62q}&
\mathcal{I} = \overline{
\bigcap_{\bB } %\big\{ 
\overline{\mathcal{K}}^{\, \PPPsb }
%\Big\}
}^{\, \PBr }
.\end{align}
Here the intersection is taken over $ \PBr $-a.s.\,$ \bB $. 
To be precise, we set for $ \U \subset \WRNz $ 
\begin{align}
\notag &% \label{:62q2}&
\mathcal{K}[\U ]= \bigcap_{\bB \in \U } \overline{\mathcal{K}}^{\PPPsb }
.\end{align}
Then the set $ \mathcal{I} $ is defined by 
\begin{align} \label{:62qq}&
\mathcal{I}= \{ \V \, ;\, 
\text{ there exist } \V_i \in\mathcal{K}[\U ] \, (i=1,2),\, 
 \U \in \mathcal{B}(\WRNz ) 
\\ \notag &\quad \quad \text{ such that }
\V_1 \subset \V \subset \V_2 , \ 
\PPPsb (\V_2\backslash \V_1) = 0 \text{ for all } \bB \in \U , \ 
 \PBr (\U )= 1 
\} 
.\end{align}

\begin{lemma} \label{l:69} 
Let 
$ \mathcal{W}_m = \sigma [\mathbf{w}^{m*}] \times \mathcal{B}(\WRNz ) $. 
Then the following holds. 
\begin{align}
 &\label{:62m}
\bigcap_{m=1}^{\infty} \ol{\mathcal{W}_m}^{\Pts } \subset \mathcal{I} 
.\end{align}
\end{lemma}
\begin{proof} 
We recall $ \Pts = \Ps \circ \XB ^{-1} $ by \eqref{:61B}. 
%Throughout the proof, let $ \ol{\ \cdot \ }$ be the completion of $ \cdot $ by $ \Pts $. 
%With the same reason as the equality in \eqref{:62t}, we easily see 
Similarly as the equality in \eqref{:62t}, %we see 
\begin{align}&\notag %\label{:69a}&
\cap_{m=1}^{\infty}\ol{\mathcal{W}_m }^{ \Pts } = 
\ol{\cap_{m=1}^{\infty}\mathcal{W}_m } ^{ \Pts } 
% ,\quad %\\label{:69aa}&
% \cap_{m=1}^{\infty}\ol{\sigma [\mathbf{w}^{m*}]}^{\Ptsb } = 
% \ol{\cap_{m=1}^{\infty}\sigma [\mathbf{w}^{m*}]}^{\Ptsb }
.\end{align}
From this we deduce that \eqref{:62m} is equivalent to 
\begin{align}\label{:69g}&
\ol{ \cap_{m=1}^{\infty}\mathcal{W}_m } ^{ \Pts } 
\subset 
\mathcal{I} 
.\end{align}

%
% Let $ \PPPsb (\cdot ) = \Ps ( (\mathbf{X},\B ) \in \cdot \, |\B = \bB ) $ be as \eqref{:62z}. 
% By definition we have 
% \begin{align}\label{:69b}
% \ol{\ \ol{\cap_{m=1}^{\infty}\mathcal{W}_m }^{ \Pts }  }^{\, \PPPsb } = 
% \{ \V \, ;\, 
% &
% \text{there exist $ \V _1 , \V _2 \in \ol{\cap_{m=1}^{\infty}\mathcal{W}_m} ^{ \Pts }  $ such that }
% \\ \notag 
% &
% \V _1 \subset \V \subset \V _2 , \quad 
%%\\ \notag &\quad \quad \quad \quad 
% \PPPsb (\V _2 \backslash \V _1 ) = 0 \} 
% .\end{align}
%%
% Set $ \Wb =(\WSN , \bB ) \equiv \{ (\mathbf{w},\bB )\, ;\, \mathbf{w}\in \WSN \} $ 
% for $ \bB \in \WRNz $. 
% Then \eqref{:69b} yields 
% \begin{align} \label{:69c}
% \ol{\ \ol{\cap_{m=1}^{\infty}\mathcal{W}_m }^{ \Pts }  }^{\, \PPPsb } = &
% \{ \V \, ;\, 
% \text{there exist $ \V _1 , \V _2 \in \ol{\cap_{m=1}^{\infty}\mathcal{W}_m} ^{ \Pts }  $ such that }
% \\ \notag &\quad \quad 
% \V _1 \cap \Wb \subset \V \cap \Wb \subset \V _2 \cap \Wb , \quad 
%%\\ \notag &\quad \quad \quad \quad 
% \PPPsb (\V _2 \backslash \V _1 ) = 0 \} 
% .\end{align}

Let $ V  \in  \ol{\bigcap_{m=1}^{\infty}\mathcal{W}_m} ^{ \Pts } $. 
Then there exist $ V_{i} \in \cap_{m=1}^{\infty}\mathcal{W}_m $ $ (i=1,2)$ such that 
\begin{align}& \label{:69C}
V_{1} \subset V \subset V_{2} , \quad \Pts (V_{2} \backslash V_{1}) = 0 
.\end{align}
For $ \bB \in \WRNz $, we set 
\begin{align*}&
\text{$ \Wb :=(\WSN , \bB )  = \{ (\mathbf{w},\bB )\, ;\, \mathbf{w}\in \WSN \} $}
.\end{align*}
Note that $ \V _{i} \cap \Wb = (A _{i}^{\mathbf{b}},\bB )$ for a unique 
$ A _{i}^{\mathbf{b}} \subset \WSN $. 
Clearly, $ \Wb \in \cap_{m=1}^{\infty}\mathcal{W}_m  $. 
From these and  $ V_{i} \in \cap_{m=1}^{\infty}\mathcal{W}_m $ we see 
\begin{align}\notag &% \label{:69d}&
(A _{i}^{\mathbf{b}},\bB ) = \V _{i} \cap \Wb \in {\cap_{m=1}^{\infty}\mathcal{W}_m } 
.\end{align}
Hence 
$ (A _{i}^{\mathbf{b}},\bB ) \in {\mathcal{W}_m } $ for all $ m \in \N $. 
Then $ A _{i}^{\mathbf{b}} \in \cap_{m=1}^{\infty}\sigma [\mathbf{w}^{m*}]  $.
This implies 
\begin{align} \label{:69e}&
A _{i}^{\mathbf{b}} \in {\TpathTSN } 
.\end{align}
Let $ \Ptsb $ be as in \eqref{:61c}. 
Let $ A ^{\mathbf{b}}$ be a unique subset of $  \WSN $ such that 
$  (A ^{\mathbf{b}},\bB ) =  \V \cap \Wb $. 
From  $ (A _{i}^{\mathbf{b}},\bB ) =  \V _{i} \cap \Wb $ combined with 
\eqref{:61c}, \eqref{:69C}, and \eqref{:69e},  there exists a set 
$ U \in \mathcal{B}(\WRNz ) $ such that $ \PBr ( U ) = 1 $ and that 
\begin{align} \label{:69E}&
A ^{\mathbf{b}} \in \ol{\TpathTSN } ^{\Ptsb } 
\quad \text{ for all $ \mathbf{b} \in U$}
%\quad \text{ for $ \PBr $-a.s.\,$ \mathbf{b}$}
.\end{align}
%
%Let $ \PPPsb (\cdot ) = \Ps ( (\mathbf{X},\B ) \in \cdot \, |\B = \bB ) $ be as in \eqref{:62z}. 
Recall the relation between $ \Ptsb $ and $  \PPPsb $ given by \rref{r:62}. 
Then \eqref{:69E} implies 
\begin{align}\notag%
 (A ^{\mathbf{b}},\bB ) \in \ol{\TpathTSN \ts \mathcal{B}(\WRNz )}^{\PPPsb }
\quad \text{ for all $ \mathbf{b} \in U$}
%\quad \text{ for $ \PBr $-a.s.\,$ \mathbf{b}$}
.\end{align}
This together with $ (A ^{\mathbf{b}},\bB ) = \V  \cap \Wb $ yields  
\begin{align}&\label{:69f}
 \V  \cap \Wb \in \ol{\TpathTSN \ts \mathcal{B}(\WRNz )}^{\PPPsb }
\quad \text{ for all $ \mathbf{b} \in U$}
%\quad \text{ for $ \PBr $-a.s.\,$ \mathbf{b}$}
.\end{align}
Recall that $ \mathcal{K} = {\TpathTSN \times \mathcal{B}(\WRNz )}$ by \eqref{:62Q}. 
Then we rewrite \eqref{:69f} as 
\begin{align}\notag &
%\label{:69h}&
 \V  \cap \Wb \in \ol{\mathcal{K}} ^{\PPPsb } 
\quad \text{ for all }\mathbf{b} \in U 
.\end{align}
Obviously, $ \Wb ^c \in \ol{\mathcal{K}} ^{\PPPsb } $. 
Hence we deduce %from \eqref{:69h} 
\begin{align}\label{:69i}&
\big( \V  \cap \Wb \big) \cup \Wb ^c 
\in \ol{\mathcal{K}} ^{\PPPsb } 
\quad \text{ for all }\mathbf{b} \in U 
.\end{align}
Because $ \Wb = \WSN \ts \{ \mathbf{b} \} $ and %$ \PPPsb ( \Wb ^c)= 0 $, 
$ \PPPsb $ is concentrated on $ \Wb $,  
we deduce from \eqref{:69i} 
\begin{align}\notag &%\label{:69j}&
V \in  \ol{\mathcal{K}} ^{\PPPsb } \quad \text{ for all $ {\mathbf{b} \in U }$}
.\end{align}
Then we obtain 
\begin{align}\label{:69k}&
V \in \bigcap_{\mathbf{b} \in U } \ol{\mathcal{K}} ^{\PPPsb } 
.\end{align}
From \eqref{:62qq}, \eqref{:69k}, and $ \PBr ( U ) = 1 $,  
we obtain \eqref{:69g}. This completes the proof.  
\qed 
\end{proof}

For $ \XB $ as above, we assume $ \XB $ satisfies \iFcs\ with $ \Fmss $. We set 
\begin{align} \label{:62e}&
\FopX = \FFF 
.\end{align}
Then we see $ \FopX \in \Wsols $. 
By construction, $ (\uL{\FopX }, \B )$ satisfies the SDE in integral form such that, 
for $ i = 1,\ldots,m $, 
\begin{align} \label{:62f} &
\barF _{\mathbf{s}}^{m,i} \XB _t= s_i + \int_0^t 
 \sigma ^i ( \FopX )_u dB_u^i + 
\int_0^t b ^i ( \FopX )_u du 
,\end{align}
where 
$ \Fsmbar = 
(\barF _{\mathbf{s}}^{m,1},\ldots,\barF _{\mathbf{s}}^{m,m},X^{m+1},X^{m+2}, \ldots ) $ 
by construction. 
We write 
\begin{align}\label{:62g}& 
 \Fsi \XB = \lim_{m\to\infty } \uL{\FopX }\text{ in $ \Wsol $ under $ \Ps $}
\end{align}
if $ \Fsi \XB \in \Wsols $ and limits \eqref{:62h}--\eqref{:62j} 
converge in $ \WT (\sS ) $ for $ \Ps $-a.s.\,for all $ i \in \N $: 
\begin{align}\label{:62h}
\lim_{m\to\infty } \barF _{\mathbf{s}}^{m,i} \XB = & \Fs ^{\infty , i} \XB 
,\\\label{:62i}
\limi{m} 
 \int_0^{\cdot } \sigma ^i ( \FopX )_u dB_u^i = &
 \int_0^{\cdot } \sigma ^i ( \FsiXB )_u dB_u^i 
,\\\label{:62j}
\limi{m} \int_0^{\cdot } b ^i ( \FopX )_u du = & 
\int_0^{\cdot } b ^i ( \FsiXB )_u du 
.\end{align}
\begin{lemma} \label{l:62} 
\6\ 
Then the following hold. 
\\\thetag{1} The sequence of maps $ \{ \Fsmbar \}_{m\in\mathbb{N}} $ 
 is consistent in the sense that, for $ \Ps $-a.s., 
\begin{align}\label{:62a}&
\barF _{\mathbf{s}}^{m,i} \XB = \barF _{\mathbf{s}}^{m+n,i} \XB 
\quad \text{ for all } 1\le i \le m ,\ m, n \in \mathbb{N}
.\end{align}
Furthermore, \eqref{:62g} holds and the map $ \Fsi $ is well defined.

\noindent \thetag{2} 
$ \XB $ is a fixed point of $ \Fsi $ in the sense that, for $ \Ps $-a.s., 
\begin{align}\label{:62b}&
( \mathbf{X}, \B ) = ( \FsiXB , \B ) 
.\end{align}

\noindent \thetag{3} 
$ \Fsi (\cdot ,\bB )$ is $ \TpSNsB $-measurable for $ \PBr $-a.s.\! $ \bB $, where 
\begin{align}
\label{:62bb}&
 \TpSNsB = \overline{\Tpath (\SN )}^{\Ptsb }
.\end{align}
$ \Fsi $ is $ \mathcal{I} $-measurable, where 
$ \mathcal{I} $ is given by \eqref{:62q} and \eqref{:62qq}. 
 \end{lemma}
\begin{proof}
The consistency \eqref{:62a} is clear because $ \XB $ under $ \Ps $ is \7. 
\eqref{:62b} is immediate from the consistency \eqref{:62a}. 
We have thus obtained \thetag{1} and \thetag{2}. 

Let $ \ol{\mathcal{W}_m }$ be the completion of $ \mathcal{W}_m $ 
with respect to $ \pPs $ as before. 
From 
\dref{d:41} and \eqref{:62e} we see $ \Fsmbar $ is $ \ol{\mathcal{W}_m }$-measurable. 
Clearly, 
\begin{align}\notag &% \label{:62nn}&
\text{$ \ol{\mathcal{W}_m }\supset \ol{\mathcal{W}_n }$ 
for $ m \le n $.}
\end{align}
Hence, $ \Fsn $ are $ \ol{\mathcal{W}_m }$-measurable 
for all $ n\ge m $. Combining this with \eqref{:62g}, we see that the limit function 
$ \Fsi $ is $ \ol{\mathcal{W}_m }$-measurable for each $ m \in \mathbb{N}$. 
Then $ \Fsi $ is $ \{\cap_{m=1}^{\infty}\ol{\mathcal{W}_m }\}$-measurable. %
Hence $ \Fsi $ is $ \mathcal{I} $-measurable from \lref{l:69}. 
The first claim in \thetag{3} follows from the second and the definition of $ \mathcal{I} $. 
We have thus obtained \thetag{3}. 
\qed \end{proof}

The next theorem reveals the relation between the existence and pathwise uniqueness 
of strong solutions and tail triviality of the labeled path space $ \WSN $. 
The definition of the strong solution staring at $ \mathbf{s}$ is given by 
\dref{d:55} with replacement of \eqref{:50a}--\eqref{:50c} by \eqref{:60a}--\eqref{:60c}. 

Recall that $ \XB $ is a continuous process defined on $ \2 $ satisfying \eqref{:62y} 
introduced at the beginning of this section, and $ \Pts $ is the distribution of $ \XB $ under $ \Ps $. 
From $ \Ps $, we set $\Ptsb $ as \eqref{:61c} and $ \TpathTSNone ; \cdot )$ as \eqref{:61b}. 
\begin{theorem} \label{l:64} 
\thetag{1} 
Assume that $ \XB $ under $ \Ps $ is \7. 
Then $ \XB $ under $ \Ps $ is a strong solution of 
\ISDEsixb if and only if \ASTpath\ holds. 

\noindent \thetag{2} 
Let $ \mathbf{X} $ and $ \mathbf{X}' $ be strong solutions of \ISDEsixb defined on $ \2 $ 
with the same $ \Ft $-Brownian motion $ \B $. 
Assume that $ \XB $ and $ (\mathbf{X}' ,\mathbf{B})$ under $ \Ps $ satisfy \iFcs. Then, 
\begin{align}\label{:64x}& 
 \Ps ( \mathbf{X} = \mathbf{X}' ) =1 
\end{align} 
if and only if for $ \PBr $-a.s.\! $ \bB $ 
\begin{align}\label{:64y}& \quad 
 \text{$ \TpathTSNone ;\Ptsb ) = \TpathTSNone ;\Ptsb ')$}
 .\end{align}
Here $ \Ptsb ' $ is defined by \eqref{:61B} and  \eqref{:61c} 
by replacing $ \mathbf{X}$ by $ \mathbf{X}'$. 
\end{theorem}
\begin{proof} 
We prove \thetag{1}. 
We note that by \lref{l:62} \thetag{2}
\begin{align}
\label{:64H}&
\text{$ ( \mathbf{X}, \B ) = ( \FsiXB , \B ) $ for $ \Ps $-a.s}
.\end{align}
The fixed point property \eqref{:64H} is a key of the proof. 
We shall utilize the structure $ \mathbf{X}= \FsiXB $. 
By assumption $ \XB $ under $ \Ps $ is a weak solution of \ISDEsixbb. 
Then $( \FsiXB , \B ) $ under $ \Ps $ is a weak solution by \eqref{:64H}.

Suppose \ASTpath. 
Then $ \Tpath (\SN ) $ is $ \Ptsb $-trivial for $ \PBr $-a.s.\! $ \bB $. 
By \lref{l:62} \thetag{3}, we see 
$ \Fsi (\cdot ,\bB )$ is $ \TpSNsB $-measurable for $ \PBr $-a.s.\! $ \bB $. 
From these, we see that $ \Fsi (\cdot ,\bB ) $ under $ \Ptsb $ is a constant. 
Hence, $ \Fsi \XB $ under $ \Ps $ becomes a function in $ \B $. 
So we write, under $ \Ps $, 
\begin{align} \label{:64a}&
 \Fsi \XB = \FFFs (\B ) 
.\end{align}

We next prove that $ \FFFs $ is a $ \mathcal{B}(\PBr ) $-measurable function.  
Let $ A \in \mathcal{B}(\WSN )$ and set $ A ' = (\Fsi )^{-1} ( A ) $. 
%Let $  \PPPsb $ be as \eqref{:62z}. 
From \eqref{:62z},  \eqref{:64H}, and \eqref{:64a},  we deduce that 
there exists $ \mathcal{N} $ such that $ \PBr (\mathcal{N} ) = 0 $ and that 
\begin{align}\label{:64b}
 \PPPsb ( A ' )  & = \Ps ((\Fsi \XB , \B )  \in  A '  | \B = \bB ) 
\\ \notag &
  = \Ps ((\FFFs (\B ) , \B )  \in  A '  | \B = \bB ) 
\\ \notag &= 
1_{ A ' } ((\FFFs (\bB ) , \bB ) )  \quad \text{ for all $ \bB \notin\mathcal{N}  $} 
\\ \notag &= 
1_{ A } (\FFFs (\bB )) = 1_{\FFFs ^{-1} ( A )} (\bB )
.\end{align}
Note that  $  \PPPsb ( A ' )  $ is $ \mathcal{B}(\WRNz ) $-measurable in $ \bB $ 
because $ \PPPsb $ is the regular conditional probability of $  \Ps $ induced by $ \B $. 
Hence from \eqref{:64b} we see $ \FFFs ^{-1} ( A ) \in  \mathcal{B}(\PBr )  $. 
This implies 
$ \FFFs $ is a $ \mathcal{B}(\PBr ) $-measurable function.

Similarly, we can prove $ \FFFs $ to be a $ \Bt (\PBr ) / \Bt $-measurable function for each $ t $. 
Indeed, let 
$ \Xt  = \{\mathbf{X}_s\}_{0\le s \le t } $ 
and 
$ \mathbf{B} ^{[0,t]}= \{ \mathbf{B}_s \}_{0\le s \le t }  $ 
be given by the restriction of the time parameter set $ [0,\infty )$ to $ [0,t]$. 
Then $ \XBt $ is a weak solution of \ISDEsixb on $ [0,t]$. 
The filtered space of $ \XBt $ is taken as 
$ ( \Omega , \mathcal{F}_t  , \Ps , \{ \mathcal{F}_s  \}_{0\le s \le t}  )$. 

We set $\Fst (\mathbf{B} )  =  \{\Fs (\mathbf{B}) _s \}_{0\le s \le t } $. 
Then from \eqref{:64H} and \eqref{:64a} we see  under $ \Ps $ 
\begin{align}\label{:64c}&
 \Fst (\mathbf{B} )  = %\Fsit \XB  = 
\Xt 
% \Fst (\mathbf{B} ) \quad \text{  under $ \Ps $}
.\end{align}
Applying the argument for $ \XB $ to $ \XBt $, we obtain 
the counterparts to \eqref{:64H} and \eqref{:64a} of $ \XBt $, that is, 
we have functions 
$ \Fs ^{\infty ,t} $ of $ \XBt $ and 
 $ \Fs ^t $ of $ \mathbf{B}^{[0,t]}$ satisfying   under $ \Ps $ 
\begin{align}
\label{:64e}&
\XBt  = ( \Fs ^{\infty ,t} \XBt , \mathbf{B}^{[0,t]} ) 
,\\
\label{:64C}& 
 \Fs ^{\infty ,t}  \XBt  = \Fs ^t (\mathbf{B}^{[0,t]} )
.\end{align}
Hence from \eqref{:64c}, \eqref{:64e}, and \eqref{:64C} we have   under $ \Ps $ 
\begin{align}\label{:64cc}&
 \Fst (\mathbf{B} )  =  \Fs ^t (\mathbf{B}^{[0,t]} ) 
.\end{align}

We set 
$ \WTzt = \{ \www \in C ([0,t]; \RdN ) ; \www _0  = \mathbf{0} \} $. 
We regard $ \PBr $ as a probability measure on $ (\WTzt , \mathcal{B}(\WTzt ) )$ and 
denote by $\overline{\mathcal{B}(\WTzt )}$ the completion of $ \mathcal{B}(\WTzt )$ 
with respect to $ \PBr $. 
We replace 
 $ \mathcal{B}(\PBr ) $ and $ \mathcal{B}(\WSN ) $  with  
 $\overline{\mathcal{B}(\WTzt )}$ and $ \Bt $, respectively. 
Then, applying the argument above to $ \XBt $, we deduce  $ \Fs ^t $ is 
$ \overline{\mathcal{B}(\WTzt )}/ \Bt  $-measurable. 
Because we can naturally identify $ \overline{\mathcal{B}(\WTzt )}$ with $ \Bt (\PBr ) $, 
we see $ \Fs ^t $ is $ \Bt (\PBr ) / \Bt $-measurable. 
This and \eqref{:64cc} imply $ \Fst $ is $ \Bt (\PBr ) / \Bt $-measurable. 
Hence, we conclude $ \FFFs $ is $ \Bt (\PBr ) / \Bt $-measurable immediately.

Collecting these we see $ \FFFs (\B ) = \FsiXB $ under $ \Ps $ is a strong solution. 
In particular, the function $ \FFFs $ is a strong solution in the sense of \dref{d:55}.

Suppose that $ \FsiXB $ under $ \Ps $ is a strong solution. 
Then by \dref{d:55} there exists a function $ \FFFs $ such that 
$ \Fsi \XB = \FFFs (\mathbf{B})$ for $ \Ps $-a.s.
%
%Hence we see that $ \Fsi \XB $ is constant for $ \Ps $-a.s. %なぜ必要かわからなかった[]
%
By \eqref{:61B} and \eqref{:61c} we have 
$ \Ptsb = \Ps (\mathbf{X}\in \cdot | \mathbf{B}= \bB ) $ for given $ \mathbf{b}$. 
Hence the distribution of $ \Fsi (\mathbf{w} , \bB ) $ under $ \Ptsb $ is the delta measure 
$ \delta _{\mathbf{z}}$ concentrated at a non-random path $ \mathbf{z}$, say. 
Therefore, we deduce that $ \TpSNsB $ is $ \Ptsb $-trivial for $ \PBr $-a.s.\! $ \bB $. 
 We have thus obtained \thetag{1}.

We proceed with the proof of \thetag{2}. 

We suppose \eqref{:64x}. Then the image measures $ \Ptsb $ and $ \Ptsb ' $ 
defined by \eqref{:61c} for $ \PBr $-a.s.\! $ \bB $ are the same. 
We thus obtain \eqref{:64y}. 

We next suppose \eqref{:64y}. 
By assumption $ \XB $ and $ (\mathbf{X}' ,\mathbf{B})$ under $ \Ps $ are strong solutions. 
Hence from \thetag{1} we see $ \Tpath (\SN ) $ is trivial with respect to $ \Ptsb $ and $ \Ptsb ' $ 
for $ \PBr $-a.s.\! $ \bB $. 
Combining this with \eqref{:64y} we obtain for $ \PBr $-a.s.\! $ \bB $ 
\begin{align}
\label{:64p}& 
\Ptsb |_{\Tpath (\SN ) } = \Ptsb '|_{\Tpath (\SN ) } 
\end{align}
and in particular 
\begin{align}
\label{:64P}&
 \TpSNsB = \TpSNsBB 
.\end{align}
Here we set $ \TpSNsBB $ by \eqref{:62bb} by replacing $ \Ptsb $ with $ \Ptsb '$.

Let $ \PPPsb $ and $ \PPPsb ' $ be defined by \eqref{:62z} for 
$ \XB $ and $ (\mathbf{X}' ,\mathbf{B})$ under $ \Ps $, respectively. 
Then from \eqref{:64p} we easily see for $ \PBr $-a.s.\! $ \bB $ 
\begin{align}
\label{:64q}& 
\PPPsb |_{\Tpath (\SN ) \ts \{ \bB \} } = \PPPsb '|_{\Tpath (\SN ) \ts \{ \bB \} } 
.\end{align}
Here we set $ \Tpath (\SN ) \ts \{ \bB \} = \Tpath (\SN ) \ts \sigma [\{\bB \}]$.

Let $ \mathcal{I} $ and $ \mathcal{I}' $ be the $ \sigma $-fields defined by \eqref{:62qq} for 
$ \XB $ and $ (\mathbf{X}' ,\mathbf{B})$ under $ \Ps $, respectively. 
Then we obtain from \eqref{:64p} and \eqref{:64q} combined with \eqref{:62c} and \eqref{:62p} 
\begin{align}
\label{:64r}&
\mathcal{I} = \mathcal{I}' 
\end{align}
and for $ \PBr $-a.s. $ \bB $
\begin{align}
\label{:64R}&
\overline{\mathcal{K} } ^{\, \PPPsb }=
\overline{\mathcal{K}} ^{\, \PPPsb '}
.\end{align}

Let $ \Fsi $ and $ \Fs '^{\infty}$ be as \lref{l:62} for 
$ \XB $ and $ (\mathbf{X}' ,\mathbf{B})$ under $ \Ps $, respectively. 
Then by \eqref{:64r} and \lref{l:62} \thetag{3} both $ \Fsi $ and $ \Fs '^{\infty}$ are 
$ \mathcal{I} $-measurable. 
%
%Furthermore, by \eqref{:64R} both $ \Fsi $ and $ \Fs '^{\infty} $ are 
%$ \overline{\mathcal{K} } ^{\, \PPPsb }$-measurable for $ \PBr $-a.s. $ \bB $. 
%
From \eqref{:64P} and \lref{l:62} \thetag{3} we have 
both $ \Fsi (\cdot ,\bB )$ and $ \Fs '^{\infty} (\cdot ,\bB )$ are 
$ \TpSNsB $-measurable for $ \PBr $-a.s. $ \bB $. 

Let $ \UU , \VV \in \mathcal{I} $ be subsets with $ \pPs (\UU )=\pPs ' (\VV )=1$ 
 such that $ \Fsi $ and $ \Fs '^{\infty}$ satisfy \eqref{:62b} on $ \UU $ and $ \VV $, respectively. 
From \eqref{:62b} we have $ \Fsi \vw = \mathbf{w} $ on $ \UU $ and $ \pPs (\UU )=1$. 
Then we take a version of $ \Fsi $ under $ \pPs $ such that $ \Fsi \vw = \mathbf{w}$ on 
$ \UU \cup \VV $. 
Similarly, we take a version of $ \Fs '^{\infty} $ under $ \pPs '$ such that 
$ \Fs '^{\infty} \vw = \mathbf{w}$ on $ \UU \cup \VV $. 
We thus have %Then from \eqref{:62b} 
\begin{align}\notag & 
\Fsi \vw = \mathbf{w} = \Fs '^{\infty} \vw 
\quad \text{ for $ \vw \in \UU \cup \VV $}
.\end{align}
Hence for $ \pPs $- and $ \pPs '$-a.s.\!% $ \vw $, 
\begin{align}
\label{:64d}&
\Fsi = \Fs '^{\infty} 
.\end{align}
We then deduce from 
 \eqref{:64d} for $ \PBr $-a.s. $ \bB $ 
\begin{align}
\label{:64D}&\quad \quad 
\Fsi (\cdot ,\bB )= \Fs '^{\infty} (\cdot ,\bB ) \quad \text{ for $ \Ptsb $- and $ \Ptsb ' $-a.s.}
\end{align}

We recall that $ \Fsi (\cdot ,\bB )$ and $ \Fs '^{\infty} (\cdot ,\bB ) $ are 
$ \TpSNsB $-measurable for $ \PBr $-a.s. $ \bB $, and that 
$ \TpSNsB $ is $ \Ptsb $- and $ \Ptsb ' $-trivial for $ \PBr $-a.s.\! $ \bB $. 
From these, we write 
\begin{align}\label{:64f}& 
\Fs (\bB ) = \Fsi (\cdot ,\bB ) ,\quad \Fs ' (\bB ) = \Fs '^{\infty} (\cdot ,\bB )
.\end{align}
Then the strong solutions $ \mathbf{X}$ and 
$ \mathbf{X}'$ are given by $ \mathbf{X} = \Fs (\mathbf{B})$ and 
$ \mathbf{X}'=\Fs ' (\mathbf{B})$. 
Putting \eqref{:64p}, \eqref{:64D}, and \eqref{:64f} 
together, we obtain $ \Fs (\bB ) = \Fs ' (\bB ) $ for $ \PBr $-a.s.\ $ \bB $. This deduces 
$ \mathbf{X} = \Fs (\mathbf{B}) = \Fs ' (\mathbf{B}) = \mathbf{X}'$ under $ \Ps $, which yields 
\eqref{:64x}. 
This concludes \thetag{2}. 
\qed \end{proof}

Let $ \XB $ be a continuous process defined on $ \2 $ satisfying \eqref{:62y} as before. 
Recall the notion of pathwise uniqueness starting at $ \mathbf{s}$ under the constraint of \iFcs\ 
given by \dref{d:54b} and \rref{r:55} \thetag{2}. 

\begin{theorem}	\label{l:65}
Assume \ASTpatH. 
Then the pathwise uniqueness of weak solutions starting at $ \mathbf{s}$ 
under the constraint of \iFcs\ and \ASTpath\ holds. 
\end{theorem}
\begin{proof}
Let $ \XB $ and $ (\mathbf{X}',\mathbf{B})$ be \8 and \ASTpath\ 
with the same Brownian motion $ \mathbf{B}$. 
Then $ (\mathbf{X},\mathbf{B})$ and $(\mathbf{X}',\mathbf{B})$ are 
strong solutions by \tref{l:64} \thetag{1}. 
By \ASTpatH\ we see \eqref{:64y} holds. 
Hence applying \tref{l:64} \thetag{2} we obtain \eqref{:64x}. 
From this we deduce the pathwise uniqueness of weak solutions starting at $ \mathbf{s}$ 
under the constraint of \iFcs\ and \ASTpath. 
\qed\end{proof}

\noindent 
{\it Proof of \tref{l:61}. } 
The claim \thetag{1} follows from \tref{l:64} \thetag{1} immediately. 

Let $ (\hat{\mathbf{X}},\hat{\mathbf{B}})$ be 
any weak solution of \ISDEsixb satisfying \iFcs\ and \ASTpath. 
% Assume that \ISDEsixb satisfies \ASTpatH. %
Then from \thetag{1} we deduce that 
$ (\hat{\mathbf{X}},\hat{\mathbf{B}})$ becomes a strong solution with $ \hat{\Fs } $
such that $ \hat{\mathbf{X}} = \hat{\Fs } (\hat{\mathbf{B}})$. 

Because the distribution of $ \mathbf{B}$ coincides with that of $ \hat{\mathbf{B}}$, 
 $ (\Fs (\mathbf{B}),\mathbf{B}) $ and $ (\Fs (\hat{\mathbf{B}}) ,\hat{\mathbf{B}})$ 
 have the same distribution. 
Hence $ (\Fs (\hat{\mathbf{B}}) ,\hat{\mathbf{B}})$ is \7. 
By assumption, \ASTpath\ and  \ASTpatH\ hold. Then \tref{l:65} yields 
\begin{align}\notag &% \label{:65d}&
\Ps ( \Fs (\hat{\mathbf{B}}) =\hat{\Fs } (\hat{\mathbf{B}}) ) = 1
.\end{align} 
Thus, we deduce $ \Fs = \hat{F}_{\mathbf{s}}$ for $ \PBr $-a.s. 
Hence we obtain 
$$ \hat{\mathbf{X}} = \hat{F}_{\mathbf{s}} (\hat{\mathbf{B}})= \Fs (\hat{\mathbf{B}})
.$$
This implies \eqref{:50s}. 

Let $ \mathbf{B}'$ be any $ \{ \mathcal{F}_t' \} $-Brownian motion defined on 
$ (\Omega ', \mathcal{F}',P', \{ \mathcal{F}_t' \} )$. 
Because $ \mathbf{B}'$ and $ \mathbf{B}$ have the same distribution, 
$ (\Fs (\mathbf{B}'), \mathbf{B}')$ and $ (\Fs (\mathbf{B}), \mathbf{B})$ have the same distribution. 
Hence $ (\Fs (\mathbf{B}'), \mathbf{B}')$ is \7\ and \ASTpath. 
This implies $ \Fs (\mathbf{B}')= \Fsi (\Fs (\mathbf{B}'), \mathbf{B}')$ 
similarly as the argument at the beginning of the proof. 
Then by \tref{l:64} \thetag{1} we deduce that $ \Fs (\mathbf{B}') $ is a strong solution. 
We have thus completed the proof of \thetag{2}. 
\qed 

\medskip 

We conclude this section with two notions of solutions. 
%The following notion is an analogy of the \IFC\ solution of \eqref{:50a}--\eqref{:50b}
%introduced in \dref{d:59}. 
%The IFC solution in the next definition treats very general ISDEs beyond those in \sref{s:5}. 
%
\begin{definition}[{IFC solution}]\label{dfn:61}
%A probability measure $ \Ps $ with 
A continuous process $ \XB $ under $ \Ps $ is called an \IFC\ solution of \ISDEsixb if 
$ \XB $ under $ \Ps $ is a weak solution of \ISDEsixb satisfying \iFcs. 
% A weak solution satisfying the \IFC\ condition \As{\iFc} is called an \IFC\ solution. 
Also, the distribution of an \IFC\ solution is called an \IFC\ solution. 
\end{definition}
We introduce a new notion of solutions of ISDEs. 
The notion is an equivalent notion of the weak solution of \ISDEsixb in terms of 
{\em an asymptotic infinite system of finite-dimensional SDEs with consistency} (A\IFC ).
We do not use the notion of A\IFC\ solutions in the present paper; still, 
we expect that it will be useful to solve ISDEs. 
\begin{definition}\label{d:61AIFC} 
%A probability measure $ \Ps $ with $ \XB $ defined on $ \2 $ %defined on $ \OF $ 
A continuous process $ \XB $ under $ \Ps $ is called an A\IFC\ solution of \ISDEsixb 
if $ \XB $ under $ \Ps $ satisfies \iFcs\ and \eqref{:62g}. 
Also, the distribution of an A\IFC\ solution is called an A\IFC\ solution. 
%Also, we call its distribution $ \Pts $ or $ \XB $ under $ \Ps $ an A\IFC\ solution. 
\end{definition}

\begin{remark}\label{r:63}
We note that \IFC\ solutions in \dref{dfn:61} are always A\IFC\ solutions. 
Conversely, we can construct an \IFC\ solution %a weak solution 
of \ISDEsixb from an A\IFC\ solution. 
Indeed, assuming that $ \XB $ under $ \Ps $ is an A\IFC\ solution of \ISDEsixbb\ and 
letting $ \Fsi $ be as \eqref{:62g}, 
we see $ (\FsiXB , \B )$ under $ \Ps $ is a weak solution of \ISDEsixbb. 
This is obvious from \eqref{:62f} and \eqref{:62h}--\eqref{:62j}. 
Then $ (\FsiXB , \B )$ under $ \Ps $ is an \IFC\ solution. 
Thus, constructing an \IFC\ solution and, in particular, a weak solution 
are reduced to constructing an A\IFC\ solution. 
\end{remark}

%\Section{Second tail theorem (\tref{l:70})}
\Section{Triviality of $ \TpathTSN $: Second tail theorem (\tref{l:70})} \label{s:7} 

Let $ \TpathTSN $ be the tail $ \sigma $-field of 
the labeled path space $ \WSN $ introduced in \eqref{:61a}. 
The purpose of this section is to prove triviality of $ \TpathTSN $ 
under distributions of IFC solutions, 
which is a crucial step in constructing a strong solution as we saw in \tref{l:61}. 
This step is very hard in general because $ \WSN $ is a huge space and 
its tail $ \sigma $-field $ \TpathTSN $ is topologically wild. 
To overcome the difficulty, 
we introduce a sequence of well-behaved tail $ \sigma $-fields, and deduce
 triviality of $ \TpathTSN $ from that of $ \TS $ under the stationary distribution of unlabeled dynamics $ \mathsf{X}$ step by step along this sequence of tail $ \sigma $-fields. 

The space $ \sSS $ is a {\em tiny} infinite-dimensional space compared with $ \SN $ and 
$ \WSN $. 
Hence, $ \sSS $ enjoys a nice probability measure $ \mu $ unlike $ \SN $ and $ \WSN $. 
We can take $ \mu $ to be associated with $ \sSS $-valued stochastic dynamics 
satisfying the $ \mu $-absolute continuity condition \As{\muAC}\ and 
the no big jump condition \As{\nbj}. 
This fact is important to the derivation.

Let $ \la $ be a probability measure on $ (\sSS , \mathcal{B}(\sSS ) )$. 
We write $ \ww =\{ \wwt \} \in \WS $. 
Let $ \Pl $ be a probability measure on $ (\WS ,\mathcal{B}(\WS ) )$ such that 
\begin{align} & \label{:70v}
 \Pl \circ \wW _0^{-1} = \la
.\end{align}
If necessary, we extend the domain of $ \Pl $ using completion of measures. 
We call $ \Pl $ a lift dynamics of $ \la $ if $ \Pl $ and $ \la $ satisfy \eqref{:70v}. 
For a given random point field $ \la $ there exist many lift dynamics of it. 
We shall take a specific lift dynamics given by \eqref{:70u}. 

Let $ \mathbf{X}$ be an $ \SN $-valued continuous process on $ \OFPF $. 
We set $ \upath (\mathbf{X})= \mathsf{X} = \sum_{i\in \N } \delta_{X^i}$ as before. 
We assume that the associated unlabeled process 
$ \upath (\mathbf{X}) $ is an $ \sSS $-valued continuous process 
such that $ P \circ \ulab (\mathbf{X}_0)^{-1} = \lambda $ and define $ \Pl $ by 
\begin{align}& \label{:70u}
\Pl := (P \circ \mathbf{X}^{-1}) \circ \upath ^{-1} = 
P \circ \upath (\mathbf{X})^{-1} = P \circ \mathsf{X}^{-1}
.\end{align}
Let $ \mu $ be a probability measure on $ (\sSS , \mathcal{B}(\sSS ) )$ and 
let $ \WSsiNE $ be as \eqref{:23gg}. 
We make the following assumptions originally introduced in \sref{s:5d}. 

\smallskip 

\noindent 
\As{\Cone} \quad \ \ $ \mu $ is tail trivial. 

\noindent 
\As{\Ctwo} 
\quad \ \ 
$ \Pl \circ \wW _t^{-1} \prec \mu $ for all $ \zzti $.

\noindent 
\As{\Cthree} \quad \, \!$ \Pl ( \WSsiNE ) = 1 $. 
%and $ \pP \circ \mathbf{X}^{-1} (\WT (\SSSsde )) = 1 $. 

\noindent 
\As{\Cfour} \quad $ P ( \{ \mrXX < \infty \} ) = 1 $ for all $ r, T \in \N $.

\medskip

We denote by $ \Ps $ the regular conditional probability such that 
\begin{align}
\label{:70p}&
\Ps = P (\cdot | \mathbf{X}_0=\mathbf{s})
.\end{align}
Let $ \mathbf{B} $ be a continuous process defined on the filtered space $ \OFPF $. 
Let $ \Pts $ denote the distribution of $ \XB $ under $ \Ps$: 
\begin{align}
 \label{:70r}&
\Pts = \Ps \circ \XB ^{-1}
.\end{align}
Let $ \PsB $ be the regular conditional probabilities of $ \Ps $ such that 
\begin{align}\label{:71c}&
\PsB = \Ps ( \mathbf{X} \in \cdot \, | \B = \bB )
.\end{align}
Then we easily see 
\begin{align} &\notag %\label{:71C}&
 \Pts = \Pll (\XB \in \cdot | \mathbf{X}_0= \mathbf{s})
,\\ \notag &
\PsB = \Pll ( \mathbf{X} \in \cdot \, | \mathbf{X} _0 =\mathbf{s},\, \B = \bB )
= \Pts (\mathbf{w} \in \cdot | \mathbf{b})
.\end{align}
We assume $ \XB $ under $ \OFPF $ is a weak solution of ISDE \eqref{:60a}--\eqref{:60b}. 
Then $ \XB $ under $ \2 $ is a weak solution of ISDE \ISDEsixc for 
$ P \circ \mathbf{X}_0 $-a.s.\! $ \mathbf{s}$. 
Hence we can apply the results in \sref{s:6} to $ \XB $ under $ \Ps $.

We now state the main theorem of this section. 
\begin{theorem}[Second tail theorem]\label{l:70}
Assume that $ \mu $ and $ \Pl $ satisfy \Conefour. 
Assume that there exists a label $ \lab $ such that 
$ \lab \circ \ulab (\mathbf{s}) = \mathbf{s} $ 
for $ \PXz $-a.s.\! $ \mathbf{s}$. 
\Dzero. 
Then $ \Pts $ satisfies \ASTpath\ for $ \PXz $-a.s.\! $ \mathbf{s}$. 
Furthermore, ISDE \ISDEsixc satisfies \ASTpatH\ for $ \PXz $-a.s.\! $ \mathbf{s}$. 
\end{theorem}

To explain the strategy of the proof, 
we introduce the notions of cylindrical tail $ \sigma $-fields on $ \WSN $ and $ \WS $. 
We set 
\begin{align}& \label{:74Y}
\mathbf{T} = \{ \mathbf{t}=(t_1,\ldots,t_m)\, ;\,
0 < t_i < t_{i+1} \, (1\le i < m),\, \ m \in \mathbb{N}\}
.\end{align}% 
We remark that we exclude $ t_1 = 0$ in the definition of $ \mathbf{T}$ in \eqref{:74Y}. 
%, that is, if $ \mathbf{t}=(t_1,\ldots,t_m)$, then $ t_i \not = 0 $ for each $ i=1,\ldots,m $. 

Let $\pirc $ be the projection $ \map{\pirc }{\sSS }{\sSS }$ such that 
$ \pirc (\mathsf{s}) = \mathsf{s} (\cdot \cap \Sr ^c)$. 
For $ \ww =\{ \wwt \} \in \WS $ and $ \mathbf{t}=(t_1,\ldots,t_m) \in \mathbf{T}$ we set 
$ \pirc (\wW _{\mathbf{t}} ) = (\pirc (\wW _{t_1}), 
\ldots, \pirc (\wW _{t_m} ) ) \in \sSS ^m $. 

Let $ \TpathS $ be the cylindrical tail $ \sigma $-field 
of $ \WS $ such that 
\begin{align} & \label{:74y} 
\TpathS = \bigvee_{\mathbf{t} \in \mathbf{T}} \bigcap_{r=1}^{\infty}
\sigma [\, \pirc (\wW _{\mathbf{t}} ) \ ]
.\end{align}
Let $ \tTpath (\SN ) $ be the cylindrical tail $ \sigma $-field of $ \WSN $ defined as 
\begin{align}\label{:70a}&
\tTpath (\SN ) = \bigvee_{\mathbf{t} \in \mathbf{T}} \bigcap_{n=1}^{\infty} 
\sigma [\wwww _{\mathbf{t}}^{n*}]
.\end{align}
Here $\wwww _{\mathbf{t}}^{n*} = (\wwww _{t_i}^{n*})_{i=1}^m $ 
for $ \mathbf{t} \in \mathbf{T}$, where 
$ \wwww _t^{n*}=(\w _t^k)_{k=n +1}^{\infty}$ for $ \wwww = (\w ^n )_{n\in\N }$.

\medskip 

We shall prove \tref{l:70} along with the following scheme: 
\begin{align}\notag %\label{:70f}
&
\TS 
\xrightarrow[ \text{\tref{l:74}}\atop \As{\Cone}, \ \As{\Ctwo}]{(\textrm{Step I})}
&& 
\TpathS 
\xrightarrow[\text{\tref{l:77}}
\atop \As{\Cone},\, \As{\Ctwo},\, \As{\Cthree},\, \As{\Cfour}]{(\textrm{Step II})}
&& 
 \tTpath (\SN ) 
\xrightarrow[\text{\tref{l:79}}\atop {\As{\iFc}_{\mathbf{s}}}]{(\textrm{Step III})}
&&
\TpathTSN 
\\ \notag 
& \ \mu && \ \Pl = P \circ \mathsf{X}^{-1}&& \ P \circ \mathbf{X}^{-1} 
&& \ \ \PsB 
.\end{align}
%For this purpose, w
We shall prove triviality of each tail $ \sigma $-field in the scheme under the distribution put under the tail $ \sigma $-field. The theorems under the arrows correspond to each step and the conditions there indicate what are used at each passage. 
Our goal is to obtain triviality of $ \TpathTSN $ under $ \PsB $ 
for $ \PBr $-a.s.\! $ \mathbf{b}$ and for $ \PXz $-a.s.\! $ \mathbf{s}$.

%Here $ \PsB $ is given by \eqref{:71c}. 
\begin{remark}\label{r:70a}
\tref{l:74} needs only \As{\Cone} and \As{\muAC} for $ \mu $. 
\tref{l:77} needs only \Conefour. 
In \tref{l:74} and \tref{l:77}, we do not use 
any properties of ISDE. That is, $ \Pts $ is not necessary \7. 
Such generality of these theorems would be interesting and useful in other aspects. 
\tref{l:79} requires the property of \8 unlike \tref{l:74} and \tref{l:77}. 
The map $ \Fsi $ in \eqref{:62g} given by \7 
plays an important role in the proof of \tref{l:79}. 
\end{remark}

\begin{remark}\label{r:C4} 
An example of unlabeled path $ \ww $ such that 
$ \mrw = \infty $ is given by \rref{r:35}. 
Such a large fluctuation $ \mrw = \infty $ of unlabeled path 
$ \ww $ yields difficulty to control $ \TpathTSN $ 
by the cylindrical tail $ \sigma $-field $\TpathS $ of unlabeled paths. 
Hence we assume \As{\Cfour}. 
\end{remark}

\Ssection{Step I: From $ \TS $ to $ \TpathS $.} \label{s:71} 

Let $ \TS $ be the tail $ \sigma $-field of $ \sSS $ and 
 $ \la $ a probability measure on $ \sSS $ with lift dynamics $ \Pl $ as before. 
We shall lift $ \mu $-triviality of $ \TS $ to $ \Pl $-triviality of 
the cylindrical tail $ \sigma $-field $\TpathS $ of $ \WS $. 
For a probability $ Q $ on $ \mathcal{B}( \WS ) $, we set 
\begin{align}&\notag 
 \TpathSone Q ) = 
 \{ \mathcal{X} \in \TpathS \, ;\, Q ( \mathcal{X} ) = 1 \}
.\end{align}
We state the main theorem of this subsection. 
\begin{theorem} \label{l:74} 
Assume \As{\Cone} for $ \mu $ and \As{\Ctwo} for $ \mu $ and $ \Pl $. 
The following then hold. 
\\\thetag{1} 
$ \TpathS $ is $ \Pl $-trivial. 
\\\thetag{2} $ \TpathSone \Pl ) $ \DOmu. 
\end{theorem}

\begin{remark}\label{r:74} 
 \tref{l:74} \thetag{2} means, if $ \mu $ is invariant under the dynamics, 
\begin{align}\notag &% & \label{:74z}
\TpathSone \Pl ) = \TpathSone \Pm ) 
.\end{align}
In case of \tref{l:5A}, we can take $ \Pm $ as the distribution of the diffusion in \lref{l:51} starting from the stationary measure $ \mu $. We shall identify the distribution $ \Pl |_{ \TpathS }$ in \pref{l:73}. 
From this we can specify $ \TpathSone \Pl ) $. 
\end{remark}

For a probability $ \nu $ on $ \mathcal{B}(\sSS ) $, we set 
$ \TSone \nu ) = \{ \mathsf{A} \in \TS \, ;\, \nu (\mathsf{A} ) = 1 \} $. 
\begin{lemma} \label{l:71} 
Under the same assumptions as \tref{l:74}, the following hold. 

\noindent 
\thetag{1} $ \TS $ is $ \Pl \circ \wW _t^{-1} $-trivial for each $ t > 0$. 
\\
\thetag{2} $ \TSone \Pl \circ \wW _t^{-1} ) = \TSone \mu ) $ for each $ t > 0$
.\\
%\thetag{3} $ \TSone \Pl \circ \wW _t^{-1} ) $ is independent of 
%the particular choice of $ \Pl $ in \As{\Ctwo}. 
\end{lemma}
\begin{proof} 
\thetag{1} is obvious. Indeed, let $ \mathsf{A}\in \TS $ and suppose 
$ \Pl \circ \wW _t^{-1} (\mathsf{A}) > 0 $. Then 
$ \mu (\mathsf{A}) > 0 $ by \As{\Ctwo}. 
Hence from \As{\Cone}, we deduce $ \mu (\mathsf{A}) = 1 $. 
This combined with \As{\Ctwo} implies $ \Pl \circ \wW _t^{-1} (\mathsf{A}) = 1 $.
% We next suppose that $ \Pl \circ \wW _t^{-1} (\mathsf{A}) = 0 $. 
% Then $ \Pl \circ \wW _t^{-1} (\mathsf{A}^c) > 0 $. 
% With the same argument as above we have $ \mu (\mathsf{A}^c) = 1 $. 
% Hence we have $ \mu (\mathsf{A}) = 0 $. 
We thus obtain \thetag{1}. 
We deduce from \As{\Cone}, \As{\Ctwo}, and \thetag{1} that 
$\mathsf{A} \in \TSone \Pl \circ \wW _t^{-1} ) $ if and only if 
$\mathsf{A} \in \TSone \mu ) $. This implies \thetag{2}. 
% \thetag{3} is clear from \thetag{2}. 
\qed \end{proof}
We extend \lref{l:71} for multi-time distributions. 
Let $ \Tpath ^{\mathbf{t}} (\mathsf{S})$ be the cylindrical tail $ \sigma $-field 
conditioned at $ \mathbf{t}=(t_1,\ldots,t_n) \in \mathbf{T}$: 
\begin{align} \label{:72x}
&
\Tpath ^{\mathbf{t}} (\mathsf{S}) = 
\bigcap_{r=1}^{\infty}
\sigma [\, \pirc (\wW _{\mathbf{t}} ) \, ]
.\end{align}
Using \eqref{:72x} we can rewrite $ \TpathS $ as 
$ \TpathS = \bigvee_{\mathbf{t} \in \mathbf{T}} \Tpath ^{\mathbf{t}} (\mathsf{S})$.

We set $ \sSS ^{\bfu }=\sSS ^{n}$, 
$ \mu ^{\ot \bfu } = \mu ^{\ot n}$, and $ |\bfu | = n $
if $ \bfu =(u_1,\ldots,u_n) \in \mathbf{T}$.

\begin{lemma} \label{l:72} \thetag{1} Let 
$ \bfu =(u_i)_{i=1}^p ,\bfv =(v_j)_{j=1}^q\in \mathbf{T} $ such that $ u_p<v_1$
and set $ (\bfu ,\bfv ) =(u_1,\ldots,u_p,v_1,\ldots,v_q) \in \mathbf{T} $. Then, 
\begin{align}\label{:72p}&
\Tpath ^{\bfu } (\mathsf{S}) \times \Tpath ^{\bfv } (\mathsf{S})
\subset \Tpath ^{\bfw } (\mathsf{S}) 
.\end{align}
\thetag{2} For each $ \mathsf{C} \in \Tpath ^{\bfw } (\mathsf{S})$ and 
$ \mathsf{y} \in \Sv $, the set 
$ \mathsf{C}^{[\sfy ]}$ is $ \Tpath ^{\bfu } (\mathsf{S})$-measurable. 
Here 
\begin{align} \notag %\label{:72q} 
&
\mathsf{C}^{[\sfy ]} = \{ \sfx \in \Su \, ;\, (\sfx , \sfy )\in \mathsf{C}\} 
.\end{align}
Furthermore, $ \mu ^{\ot \bfu }(\mathsf{C}^{[\sfy ]} )$ 
is a $ \Tpath ^{\bfv } (\mathsf{S})$-measurable function in $ \mathsf{y}$. 
\end{lemma}
\begin{proof}
We easily deduce for each $ r \in \N $ 
%from $ \pirc (\wW _{\mathbf{t}} ) = (\pirc (\wW _{t_1}), 
%\ldots, \pirc (\wW _{t_m} ) ) $ for $ \mathbf{t}=(t_1,\ldots,t_m)$ 
%that 
\begin{align}\label{:72a}& 
\sigma [\, \pirc (\wW _{\bfu } ) \, ] \ts 
\sigma [\, \pirc (\wW _{\bfv } ) \, ]
= 
\sigma [\, \pirc (\wW _{\bfw } ) \, ]
.\end{align}
Then, taking intersections in the left-hand side step by step, we obtain 
\begin{align}\label{:72c}& 
\bigcap_{s=1}^{\infty} 
\sigma [\, \pisc (\wW _{\bfu } ) \, ] \ts 
\bigcap_{t=1}^{\infty} 
\sigma [\, \pitc (\wW _{\bfv } ) \, ]
\subset 
\sigma [\, \pirc (\wW _{\bfw } ) \, ]
.\end{align}
Hence, taking the intersection in $ r $ on the right-hand side of \eqref{:72c}, 
we deduce \eqref{:72p}. 

Next, we prove \thetag{2}. By assumption and from the obvious inclusion, we see that 
\begin{align*}& \text{
$ \mathsf{C} \in \Tpath ^{\bfw } (\mathsf{S})
\subset \sigma [\, \pisc (\wW _{\bfw } ) \, ] $ for all $ s \in \N $. }
\end{align*}
Then, from \eqref{:72a}, we see that 
\begin{align}\notag %\label{:72d}
& \mathsf{C} \in 
\sigma [\, \pisc (\wW _{\bfu } ) \, ] \ts \sigma [\, \pisc (\wW _{\bfv } ) \, ]
\quad \text{ for all } s \in \N 
.\end{align}
Hence, $ \mathsf{C}^{[\sfy ]}$ is $ \sigma [\, \pisc (\wW _{\bfu } ) \, ] $-measurable 
for all $ s \in \N $. From this we deduce that 
$ \mathsf{C}^{[\sfy ]}$ is $ \Tpath ^{\bfu } (\mathsf{S})$-measurable. 
The second claim in \thetag{2} can be proved similarly. 
\qed \end{proof}

To simplify the notation, we set 
$ \PmT = \Pl \circ \wW _{\mathbf{t}}^{-1}$ for $ \mathbf{t} \in \mathbf{T}$. 
\begin{proposition}	\label{l:73} 
%Assume \As{\Cone} and \As{\Ctwo} for $ \mu $. 
Under the same assumptions as \tref{l:74}, the following hold. \\
For each $ \mathbf{t}=(t_1,\ldots,t_n) \in \mathbf{T}$, 
\begin{align}\label{:73a}&
\PmT | _{\Tpath ^{\mathbf{t}} (\mathsf{S}) } = 
\mu ^{\ot \bft } |_{ \Tpath ^{\mathbf{t}} (\mathsf{S}) } 
.\end{align}
In particular, for any $ \mathsf{C}\in \Tpath ^{\mathbf{t}} (\mathsf{S}) $, 
the following identity with dichotomy holds. 
\begin{align}\label{:73b}&
\PmT ( \mathsf{C}) = \mu ^{\ot \bft } (\mathsf{C}) \in \{ 0,1 \} 
.\end{align}
\end{proposition}

\begin{proof}
We prove \pref{l:73} by induction with respect to $ n = |\mathbf{t}| $. 

If $ n = 1$, then the claims follow from \lref{l:71}. 
Here we used \As{\Cone} and \As{\Ctwo} to apply \lref{l:71}. 

Next, suppose that the claims hold for $ n-1$. 
Let $ \bfu $ and $ \bfv $ be such that $ (\bfu , \bfv ) = \bft $ and that 
$ 1 \le |\bfu |, |\bfv | < n $. 
We then see that $ |\bfu |+|\bfv | = |\bft | = n $. 
We shall prove that the claims hold for $ \bft $ with $ |\bft | = n $ in the sequel. 

Assume that $ \mathsf{C}\in \Tpath ^{\mathbf{t}} (\mathsf{S}) $. 
Then we deduce from \lref{l:72} that 
$ \mathsf{C}^{[\sfy ]} \in \Tpath ^{\bfu } (\mathsf{S})$ and that 
$ \mu ^{\ot \bfu }( \mathsf{C}^{[\sfy ]} ) $ 
is a $ \Tpath ^{\bfv } (\mathsf{S})$-measurable function in $ \sfy $. 
By the induction hypothesis, we see that 
$ \Pl ( \wW _{\bfu } \in \mathsf{C}^{[\sfy ]} ) $ 
is $ \Tpath ^{\bfv } (\mathsf{S})$-measurable in $ \sfy $ 
and that the following identity with dichotomy holds. 
\begin{align}\label{:73c}&
\Pl ( \wW _{\bfu } \in \mathsf{C}^{[\sfy ]} ) 
= 
\mu ^{\ot \bfu } ( \mathsf{C}^{[\sfy ]} )
\in \{ 0,1 \} 
\quad \text{ for each } \sfy \in \Sv 
.\end{align}

By disintegration, we see that 
\begin{align}\label{:73d}& 
\Pl ( \wW _{\bfu } \in \cdot ) = 
\int _{\Sv }
\Pl ( \wW _{\bfu } \in \cdot | 
\wW _{\bfv } = \sfy ) 
\PmV (d\sfy )
.\end{align}
By the induction hypothesis, we see that 
$ \Pl ( \wW _{\bfu } \in \mathsf{A} ) \in \{ 0,1 \} $ 
for all $ \mathsf{A} \in \Tpath ^{\bfu } (\mathsf{S}) $. 
Hence let $\mathsf{A} \in \Tpath ^{\bfu } (\mathsf{S}) $ and 
suppose that 
\begin{align}\label{:73D}&
\text{$ \Pl ( \wW _{\bfu } \in \mathsf{A} ) = a $, % in \eqref{:73c}, 
}
\end{align}
where $ a \in \{ 0,1 \} $. We then obtain from \eqref{:73d} 
and $ \Pl ( \wW _{\bfu } \in \mathsf{A} ) \in \{ 0,1 \} $ that 
\begin{align} & \label{:73e}
\Pl ( \wW _{\bfu } \in \mathsf{A} | 
\wW _{\bfv } = \sfy ) = a 
\quad \text{ for $ \PmV $-a.s.\! 
$ \sfy $.}
\end{align}
From \eqref{:73d}, \eqref{:73D}, and \eqref{:73e}, 
we deduce that for each $ \mathsf{A} \in \Tpath ^{\bfu } (\mathsf{S}) $ 
\begin{align}& \label{:73P} 
\Pl ( \wW _{\bfu } \in \mathsf{A} | 
\wW _{\bfv } = \sfy ) 
=
\Pl ( \wW _{\bfu } \in \mathsf{A}) 
\quad \text{ for $ \PmV $-a.s.\! $ \sfy $} 
.\end{align}

We next remark that 
\begin{align}
\label{:73PP}&
\text{$ \Pl ( \wW _{\bfv } = \sfy | \wW _{\bfv } = \sfy ) = 1 $ 
for $ \PmV $-a.s.\! $\sfy $}
.\end{align}
We refer the reader to the corollary of Theorem 3.3 on page 15 of \cite{IW} 
for the general result from which \eqref{:73PP} is derived.

For all $ \mathsf{A} \in \Tpath ^{\bfu } (\mathsf{S}) $ and 
$ \mathsf{B} \in \mathcal{B}(\mathsf{S}^{\bfv }) $ we deduce from 
\eqref{:73P} and \eqref{:73PP} that 
\begin{align}\label{:73G}
\Pl ( 
 \wW _{\bfu } \in \mathsf{A}, \, 
 \wW _{\bfv } \in \mathsf{B}) 
 &= 
\int _{\Sv }
\Pl ( 
 \wW _{\bfu } \in \mathsf{A}, \, 
 \wW _{\bfv } \in \mathsf{B}\, 
 | 
 \wW _{\bfv } = \sfy ) \, 
\PmV (d\sfy )
% \\ \notag 
% &= 
% \int _{\Sv }
% \Pl (
% \wW _{\bfu } \in \mathsf{A}, \, 
% \wW _{\bfv } \in \mathsf{B}, \, 
% \wW _{\bfv } = \sfy \, 
% | \wW _{\bfv } = \sfy \, ) \, 
% \PmV (d\sfy )
\\ \notag 
&= 
\int _{\mathsf{B} }
\Pl ( \wW _{\bfu } \in \mathsf{A} | 
\wW _{\bfv } = \sfy ) \, 
\PmV (d\sfy )
\quad \text{ by \eqref{:73PP}}
\\ \notag &= 
\int _{\mathsf{B} }
\Pl ( \wW _{\bfu } \in \mathsf{A} ) \, 
\PmV (d\sfy ) 
\quad \quad \quad \quad \text{ by \eqref{:73P}}
\\ \notag 
&= 
\Pl ( \wW _{\bfu } \in \mathsf{A} ) \, 
\PmV (\mathsf{B} ) 
.\end{align}

From \eqref{:73G} and the monotone class theorem, we deduce that 
$ \PmT = \Pl ^{(\wW _{\bfu },\wW _{\bfv })}$ 
restricted on 
$ \Tpath ^{\bfu } (\mathsf{S}) \ts \mathcal{B}(\mathsf{S}^{\bfv })$ 
is a product measure. We thus obtain 
\begin{align}\label{:73H}&
(\PmT , \Tpath ^{\bfu } (\mathsf{S}) \ts \mathcal{B}(\mathsf{S}^{\bfv }) ) 
= 
(\PmU |_{\Tpath ^{\bfu } (\mathsf{S}) } \ts 
 \PmV |_{\mathcal{B}(\mathsf{S}^{\bfv }) } , 
\Tpath ^{\bfu } (\mathsf{S}) \ts \mathcal{B}(\mathsf{S}^{\bfv }) ) 
.\end{align}
That is, 
\begin{align*}&
\PmT |_{ \Tpath ^{\bfu } (\mathsf{S}) \ts \mathcal{B}(\mathsf{S}^{\bfv }) } 
= \PmU |_{\Tpath ^{\bfu } (\mathsf{S}) } \ts 
\PmV |_{\mathcal{B}(\mathsf{S}^{\bfv }) } 
.\end{align*}
In particular, from \eqref{:73H}, we deduce that 
for $ \PmV $-a.s.\! $\sfy $ 
\begin{align}\label{:73H6}&
\Pl ( 
 \wW _{\bfu } \in \mathsf{A} \, 
| \wW _{\bfv } = \sfy \, ) 
= 
\Pl ( \wW _{\bfu } \in \mathsf{A} \, ) 
\quad \text{ for all $ \mathsf{A} \in \Tpath ^{\bfu } (\mathsf{S}) $}
.\end{align}

For any $ \mathsf{C} \in \mathcal{B} (\mathsf{S}^{\bft }) $, 
we deduce that 
\begin{align}\label{:73GGG}
\Pl ( 
 \wW _{\bft } \in \mathsf{C}) 
 &= 
\int _{\Sv }
\Pl ( 
 \wW _{\bfu } \in \mathsf{C}^{[\sfy ]}, \, 
 \wW _{\bfv } = \mathsf{y}\, 
 | 
 \wW _{\bfv } = \sfy ) \, 
\PmV (d\sfy )
\\ \notag 
&= 
\int _{\Sv }
\Pl ( 
 \wW _{\bfu } \in \mathsf{C}^{[\sfy ]} \, 
% \wW _{\bfv } = \mathsf{y}\, 
| \wW _{\bfv } = \sfy \, ) \, 
\PmV (d\sfy )
.\end{align}
Here, we used 
$ \Pl ( \wW _{\bfv } = \sfy | \wW _{\bfv } = \sfy ) = 1 $ 
for $ \PmV $-a.s.\! $\sfy $, which follows from \eqref{:73PP}.

Assume $ \mathsf{C} \in \Tpath ^{\bft } (\mathsf{S}) $. 
Then from \lref{l:72} \thetag{2}, we obtain 
\begin{align}\label{:73H4z}&
\mathsf{C}^{[\sfz ]} \in \Tpath ^{\bfu } (\mathsf{S}) 
\quad \text{for all } \sfz \in \mathsf{S}^{\bfv }
.\end{align}
Hence from \eqref{:73H6} and \eqref{:73H4z}, we see that 
for $ \PmV $-a.s.\! $ \sfy $ 
\begin{align}\label{:73H4}&
\Pl ( 
 \wW _{\bfu } \in \mathsf{C}^{[\sfz ]} \, 
| \wW _{\bfv } = \sfy \, )
=
\Pl ( \wW _{\bfu } \in \mathsf{C}^{[\sfz ]} \, ) 
\quad \text{ for all $ \sfz \in \mathsf{S}^{\bfv}$}
.\end{align}
We emphasize that \eqref{:73H4} holds for {\em all }$ \sfz \in \mathsf{S}^{\bfv} $. 
Hence we can take $ \sfz = \sfy $ in \eqref{:73H4} for $ \PmV $-a.s.\! $ \sfy $. 
This yields for $ \PmV $-a.s.\! $ \sfy $ 
\begin{align}\label{:73H5}&
\Pl ( 
 \wW _{\bfu } \in \mathsf{C}^{[\sfy ]} \, | \wW _{\bfv } = \sfy \, ) =
\Pl ( \wW _{\bfu } \in \mathsf{C}^{[\sfy ]} \, ) 
.\end{align}

From \eqref{:73GGG}, \eqref{:73H5}, and \eqref{:73c} 
 we obtain 
\begin{align}\label{:73g}
\Pl ( \wW _{\bft } \in \mathsf{C}) 
&= 
\int _{\Sv }
\Pl ( \wW _{\bfu } \in \mathsf{C}^{[\sfy ]}) 
\PmV (d\sfy ) 
= 
\int _{\Sv }
\mu ^{\ot \bfu } (\mathsf{C}^{[\sfy ]})
\PmV (d\sfy ) 
.\end{align}
From \lref{l:72} \thetag{2}, we see that 
$ \mu ^{\ot \bfu } (\mathsf{C}^{[\sfy ]}) $ is a 
$ \Tpath ^{\bfv } (\mathsf{S}) $-measurable function in $ \sfy $. 
% and that $ \mathsf{C}^{[\sfy ]} \in \Tpath ^{\bfu } (\mathsf{S}) $ 
%for each $ \sfy \in \Sv $. 
Combining this with the induction hypothesis \eqref{:73a} for $ |\bfv |< n $, we obtain 
\begin{align} \label{:73i}
\int _{\Sv }
\mu ^{\ot \bfu } (\mathsf{C}^{[\sfy ]})
\PmV (d\sfy ) 
&= 
\int _{\Sv }
\mu ^{\ot \bfu } (\mathsf{C}^{[\sfy ]})
\PmV |_{ \Tpath ^{\bfv } (\mathsf{S}) } (d\sfy ) 
\\ \notag &
= 
\int _{\Sv }
\mu ^{\ot \bfu } (\mathsf{C}^{[\sfy ]}) 
\mu ^{\ot \bfv } |_{ \Tpath ^{\bfv } (\mathsf{S}) }
(d\sfy ) 
\\ \notag &
= 
\int _{\Sv }
\mu ^{\ot \bfu } (\mathsf{C}^{[\sfy ]})
\mu ^{\ot \bfv } (d\sfy ) 
% \\ \notag &
= 
\mu ^{\ot \bft } (\mathsf{C})
.\end{align}
From \eqref{:73g} and \eqref{:73i}, we obtain \eqref{:73a} and the equality in \eqref{:73b} 
for $ |\bft | = n $. 

We deduce $ \mu ^{\ot \bft } $-triviality of $ \Tpath ^{\mathbf{t}} (\mathsf{S})$
from \lref{l:72} \thetag{2} and the equality % in \eqref{:73i}
\begin{align}
\label{:73j}&
 \int _{\Sv }\mu ^{\ot \bfu } (\mathsf{C}^{[\sfy ]})\mu ^{\ot \bfv } (d\sfy ) = 
\mu ^{\ot \bft } (\mathsf{C})
\end{align}
by induction with respect to $ n=|\mathbf{t}|$. 
Indeed, because $\mu^{\otimes \mathbf{u}}(\mathsf{C}^{[\sfy ]}) \in \{0,1\}$ 
by the assumption of induction 
and $ \mu^{\otimes \mathbf{u}}(\mathsf{C}^{[\sfy ]})$ is 
$ \Tpath ^{\bfv } (\mathsf{S})$-measurable in $\mathsf{y}$ by \lref{l:72} \thetag{2}, 
we obtain $ \mu ^{\ot \bft } $-triviality of $ \Tpath ^{\mathbf{t}} (\mathsf{S})$ 
from \eqref{:73j}. 
Then from this we see that $ \mu ^{\ot \bft } (\mathsf{C}) \in \{ 0,1 \} $ holds. 
This completes the proof. 
\qed \end{proof}

\noindent 
{\em Proof of \tref{l:74}. } 
For $ \bft =(t_1,\ldots,t_n)\in \mathbf{T}$ and $ \A _i \in \TS $, let 
\begin{align}\label{:74a}&
\mathcal{X} = \{ \wW \in \WS \, ;\, \wW _{t_i}\in \A _{i}\ (i=1,\ldots,n) \} 
.\end{align}
Let $ \mathit{C}\TpathS $ denote the family of the elements 
of $ \TpathS $ of the form \eqref{:74a}. 
$ \mathit{C}\TpathS $ is then $ \Pl $-trivial by \pref{l:73}. 
The first claim \thetag{1} follows from this and the monotone class theorem. 
The second claim \thetag{2} immediately follows from the equality \eqref{:73b} 
in \pref{l:73}. 
\bbbbb 

\Ssection{Step II: From $\TpathS $ to $ \tTpath (\SN ) $.}\label{s:72} 

Let $ \TpathS $ and $ \tTpath (\SN ) $ be as 
\eqref{:74y} and \eqref{:70a}, respectively: 
\begin{align}\notag &
\TpathS = \bigvee_{\mathbf{t} \in \mathbf{T}} \bigcap_{r=1}^{\infty}
\sigma [\, \pirc (\wW _{\mathbf{t}} ) \ ]
,\quad %\\\notag &
 \tTpath (\SN ) = 
\bigvee_{\mathbf{t} \in \mathbf{T}} 
\bigcap_{n=1}^{\infty}\sigma [\wwww _{\mathbf{t}}^{n*}]
.\end{align}
This subsection proves the passage from $\TpathS $ to $ \tTpath (\SN ) $. 
 
For a given label $ \lab $, we can define the label map 
$ \map{ \lpath }{\WSsiNE }{\WSN } $ by \eqref{:23i}. 
Let $ \lab $ be a label such that 
$ \lab \circ \ulab (\mathbf{w}_0) = \mathbf{w}_0 $ for $ \PX $-a.s. 
By \As{\Cthree} the map $ \lpath $ for $ \lab $ is well defined for $ \Pl $-a.s. 
Then $ \lpath \circ \upath (\mathbf{w})= \mathbf{w} $ holds for $ \PX $-a.s. 
%Then the relation $ \lpath ^{-1}= \upath $ holds for $ \PX $-a.s. 
Hence from \eqref{:70u} we deduce 
\begin{align}\label{:70o}& 
\Pl \circ \lpath ^{-1} = \PX 
.\end{align}

Assumption \As{\Cfour} is a key to constructing the lift from the unlabeled path space 
to the labeled path space. 
We use \As{\Cfour} in \lref{l:76} to control the fluctuation of the trajectory of the 
{\em labeled} path $ \mathbf{X}$. Indeed, we see the following. 
\begin{lemma} \label{l:76} 
Assume \As{\Cthree} and \As{\Cfour}. Then 
\begin{align}\label{:76a}&
\lpath ^{-1}(\tTpath (\SN ) ) \subset \tTpath (\sSS ) \quad \text{ under $\Pl $} 
.\end{align}
\end{lemma}
Here in general for sub $ \sigma $-fields $ \mathcal{G} $ and $ \mathcal{H} $ 
on $ (\Omega , \mathcal{F}, P )$ we write $ \mathcal{G} \subset \mathcal{H} $ under $ P $ if 
$ \mathcal{G} \subset \mathcal{H} $ holds up to $ P $, that is, for each $ A \in \mathcal{G} $ 
there exists an $ A' \in \mathcal{H} $ such that $ P (A \ominus A' ) = 0 $, 
where $ A \ominus A' = \{A\cup A'\}\backslash \{ A \cap A' \} $ denotes the symmetric difference. 

\begin{proof} 
Let $ \mrT $ be as \eqref{:53n}. By \eqref{:70o} we can rewrite \As{\Cfour} as 
\begin{align}\label{:70n}&
\text{$ \Pl (\lpathw ) = 1 $ for all $ r, T \in \N $}
.\end{align}
% 
%Let $ \mathbf{N}_{r,T}= \lpathw $ be the image of no big jump paths. 
From \eqref{:70n} we easily see 
\begin{align}
\label{:76B}&
\text{$ \Pl (\bigcap_{r=1}^{\infty} \lpathw ) = 1 $}
\text{ for all $ T \in \N $}
.\end{align}

Let $ \mathbf{A} \in \bigcup_{\mathbf{t} \in \mathbf{T}} 
\bigcap_{n=1}^{\infty}\sigma [\wwww _{\mathbf{t}}^{n*}] $. 
%be an arbitrary element. 
Then there exists a $ \mathbf{t}\in \mathbf{T}$ such that 
\begin{align} & \notag %\label{:76z}
\mathbf{A} \in 
\bigcap_{n=1}^{\infty}\sigma [\wwww _{\mathbf{t}}^{n*}] 
.\end{align}
Let %$ T \in \N $ such that $ t_k < T $. %, where $ \mathbf{t}=(t_1,\ldots,t_k) $. 
$ T \in \N $ such that $ t_k < T $, where $ \mathbf{t}= (t_1,\ldots,t_k) $. 
Then we deduce for each $ r \in \mathbb{N} $
\begin{align}&\label{:76b}
 \lpathA \cap \lpathw 
\\ \notag &\quad \quad \quad \quad 
\in \ \sigma [\pirc ( \wW _{\mathbf{t}} ) ] \cap \lpathw 
,\end{align}
where $ \mathcal{F} \cap A = \{ F \cap A ; F \in \mathcal{F} \} $ 
for a $ \sigma $-field $ \mathcal{F} $ and a subset $ A $. 

Combining \eqref{:76B} and \eqref{:76b}, we obtain under $ \Pl $ 
\begin{align}\notag %\label{:76c}
&&
 \lpathA 
&= \lpathA \bigcap \Big\{ \bigcap _{r=1}^{\infty} \lpathw \Big\} 
%\quad \quad 
&&\text{by \eqref{:76B}}
\\ \notag 
&&&= 
\bigcap _{r=1}^{\infty}
\Big\{ \lpathA 
\bigcap 
\lpathw \Big\}
\\ \notag &&
&\in
\bigcap _{r=1}^{\infty} \Big\{ 
\sigma [\pirc ( \wW _{\mathbf{t}} ) ] 
\bigcap 
\lpathw 
\Big\} 
%
%\quad 
&&\text{by \eqref{:76b}}
\\ \notag &&
& = 
\bigcap_{r=1}^{\infty} \sigma [\pirc ( \wW _{\mathbf{t}} ) ] 
%\quad \quad 
&&\text{by \eqref{:76B}}
.\end{align}
By \eqref{:74y} we see 
$ \bigcap_{r=1}^{\infty} \sigma [\pirc ( \wW _{\mathbf{t}} ) ] \subset \tTpath (\sSS ) $. 
Thus, for arbitrary $ \mathbf{A} \in 
 \bigcup_{\mathbf{t} \in \mathbf{T}} 
 \bigcap_{n=1}^{\infty}\sigma [\wwww _{\mathbf{t}}^{n*}] $, 
 we see 
 $ \lpathA \in \tTpath (\sSS ) $ under $\Pl $ from the argument above. 
Hence we obtain 
\begin{align}
\label{:76C}& 
 \lpath ^{-1}\big(
\bigcup_{\mathbf{t} \in \mathbf{T}} 
\bigcap_{n=1}^{\infty}\sigma [\wwww _{\mathbf{t}}^{n*}] \big)
\subset \tTpath (\sSS ) \quad \text{ under $ \Pl $}
.\end{align}
Applying the monotone class theorem, we then deduce \eqref{:76a} from \eqref{:76C}. 
\qed \end{proof}

For a probability measure $ Q $ on $ \WSN $, we set 
\begin{align}\notag &%\label{:76d}& 
\5 (\SN ; Q ) = \{ \mathbf{A} \in \tTpath (\SN ); Q (\mathbf{A}) = 1 \} 
.\end{align}
We set two conditions: 

\smallskip

\noindent 
\As{\Done} 
$ \tTpath (\SN ) $ is $ P \circ \mathbf{X}^{-1} $-trivial. 
\\\As{\Dtwo} 
$ \tTpathTSNone P \circ \mathbf{X}^{-1} ) $ \DOMUL. 
\begin{theorem} \label{l:77} 
Assume \As{\Cone} for $ \mu $ and \As{\Ctwo} for $ \mu $ and $ \Pl $. 
Assume \As{\Cthree} and \As{\Cfour}. 
Then \AsDD\ hold. 
\end{theorem}
\begin{proof}
Let $ \mathbf{A} \in \tTpath (\SN )$. 
Then we see $ \lpathA \in \tTpath (\sSS )$ under $\Pl $ from \lref{l:76}. 
\tref{l:74} \thetag{1} deduces that $ \TpathS $ is $ \Pl $-trivial. 
Hence we have 
\begin{align}
\label{:77a}& 
 \Pl (\lpathA ) \in \{ 0,1 \} 
.\end{align}
We see from \eqref{:70o} 
\begin{align} \label{:77b}&
 \Pl (\lpathA ) = % (P \circ \mathbf{X}^{-1}) \circ \upath ^{-1} (\lpathA ) = 
P \circ \mathbf{X}^{-1} (\mathbf{A})
.\end{align}
Combining \eqref{:77a} and \eqref{:77b}, we obtain \As{\Done}. 

From \lref{l:76} we see that 
$ \tTpathTSNone \Pl \circ \lpath ^{-1}) $ depends only on 
$ \TpathSone \Pl ) $ and $ \lpath $. 
From \tref{l:74} \thetag{2} we see that $ \TpathSone \Pl ) $ depends only on $ \mu $. 
%
%Furthermore, the map $ \lpath $ depends only on $ \lab $. 
Recall that $ \lpath (\ww ) $ is defined for all $ \ww \in \WSsiNE $ and 
 $ \lpath $ is unique for a given $ \lab $. 
Collecting these, we have 
%\smallskip 

\noindent 
\As{\Dtwo}' \quad 
$ \tTpathTSNone \Pl \circ \lpath ^{-1}) $ 
depends only on $ \mu $ and $ \lab $. 

For given $ \mul $ and $ \mu $ we see the label $ \lab $ is uniquely determined 
for $ \mu $-a.s. 
That is, if $ \hat{\lab}$ is a label such that $\mu \circ \hat{\lab }^{-1} = \mul $, 
then $ \hat{\lab } (\mathsf{s} )= \lab (\mathsf{s} )$ for $ \mu $-a.s.\ $ \mathsf{s} $. 
It is clear that $ \mu $ is uniquely determined by $ \mul $ because 
$ (\mul )\circ \ulab ^{-1} = \mu $.

Combining this with \As{\Dtwo}' we deduce 

\noindent 
\As{\Dtwo}'' \quad 
$ \tTpathTSNone \Pl \circ \lpath ^{-1}) $ \DOMUL. 

\noindent 
From \eqref{:70o} we have $ \Pl \circ \lpath ^{-1} = \PX $. Then 
 \As{\Dtwo}'' is equivalent to \As{\Dtwo}. Thus we obtain \As{\Dtwo}. 
\qed \end{proof}

\Ssection{Step III: From $ \tTpath (\SN ) $ to $ \TpathTSN $: Proof of \tref{l:70}.}\label{s:73} 

Let $ \Pts = \Ps \circ \XB ^{-1} $ be as \eqref{:70r} and 
let \AsDD\ be as \tref{l:77}. 
We state the main result of this section. 
\begin{theorem} \label{l:79} \Dzero. 

\noindent \thetag{1} 
Assume \As{\Done}. 
Then $ \Pts $ satisfies \ASTpath\ for $ \PXz $-a.s.\! $ \mathbf{s}$. 
\\\thetag{2} 
Assume \AsDD. 
Then \ISDEsixb satisfies \ASTpatH\ for $ \PXz $-a.s.\! $ \mathbf{s}$. 
\end{theorem}

Triviality in \AsDD\ is with respect to 
the {\it anneal} probability measure $ P \circ \mathbf{X}^{-1} $. 
The pair of assumptions \ASTpath\ and \ASTpatH\ is its quenched version. 
So \tref{l:79} derives the quenched triviality from the annealed one. 
Another aspect of \tref{l:79} is the passage of triviality of 
the cylindrical tail $ \sigma $-field of the labeled path space 
to that of the full tail $ \sigma $-field of the labeled path space.

Let $ \PsB = \Ps ( \mathbf{X} \in \cdot \, | \B = \bB )$ be as \eqref{:71c}. 
To simplify the notation, we set 
\begin{align} & \label{:71g} 
\Upsilon = \Pll \circ (\mathbf{X}_0,\mathbf{B} )^{-1} 
.\end{align}
Let $ \Fsi $ be the map in \eqref{:62g}. 
Such a map $ \Fsi $ exists for $ \PXz $-a.\! s.\! $ \mathbf{s}$ 
because $ \Pts $ is \9. 
Then for \mB 
\begin{align}& \label{:79z} \quad \quad 
\Fsiwb = \wwww 
\quad \text{ for $ \PsB $-a.s.\! $ \wwww $}
\end{align}
and the map $ \Fs $ in \tref{l:61} is given by $ \Fs (\bB ) = \Fsiwb $. 
Recall that $ \Fsi $ is $ \mathcal{I} $-measurable by \lref{l:62} \thetag{3}. 
$\Fsi $ is however not necessary $ \TpathTSN \times \mathcal{B}(\WRNz )$-measurable. 
We define the $ \Wsols $-valued map $ \FisB $ on a subset of $ \Wsols $ as 
\begin{align} \label{:79x}
&
 \FisB (\mathbf{\cdot }) = \Fsi (\mathbf{\cdot },\bB )
.\end{align}
The map $ \FisB $ is defined for $ \PsB $-a.s. Let 
$ \WWfixsB = \{ \wwww \in \Wsols \, ;\, \FisB (\wwww ) = \wwww \} $. 
\begin{lemma} \label{l:7X} \MST. Then, $ \PsB (\WWfixsB ) = 1 $ for \mB. 
\end{lemma}
\begin{proof}
We deduce \lref{l:7X} from \eqref{:71c}, \eqref{:79z}, and \eqref{:79x} immediately. 
\qed \end{proof}

We next recall the notion of a countable determining class. 

Let $ (U,\mathcal{U} )$ be a measurable space and let 
$ \mathcal{P} $ be a family of probability measures on it. 
Let $ \mathcal{U}_0 $ be a subset of $ \cap_{P\in\mathcal{P}} \mathcal{U}^{P}$, 
where $ \mathcal{U}^{P}$ is the completion of the $ \sigma $-field $ \mathcal{U} $ 
with respect to $ P $. 
We call $ \mathcal{U}_0 $ % of elements of $ \mathcal{U} $ 
a determining class under $ \mathcal{P} $ if any 
two probability measures $ P $ and $ Q $ in $ \mathcal{P} $ are equal 
if and only if $ P (A) = Q (A) $ for all $ A \in \mathcal{U}_0 $. 
Here we extend the domains of $ P $ and $ Q $ to 
$ \cap_{P\in\mathcal{P}} \mathcal{U}^{P}$ in an obvious manner. 
Furthermore, $ \mathcal{U}_0 $ is said to be a determining class of 
the measurable space $ (U,\mathcal{U} )$ if $ \mathcal{P} $ can be taken 
as the set of all probability measures on $ (U,\mathcal{U} )$. 
A determining class $ \mathcal{U}_0 $ is said to be countable if 
its cardinality is countable.

It is known that a Polish space $ X $ equipped with the Borel $ \sigma $-field 
$ \mathcal{B}(X) $ has a countable determining class. 
If we replace $ \mathcal{B}(X) $ with a sub-$ \sigma $-field $ \mathcal{G}$, 
the measurable space $ (X,\mathcal{G} )$ does not necessarily have 
any countable determining class in general. 
One of the difficulties to carry out our scheme is that measurable spaces 
$ \WS $ and $ \WSN $ equipped with tail $ \sigma $-fields do not 
have any countable determining classes. 
In the sequel, we overcome this difficulty using $ \FisB $ and $ \WWfixsB $. 

As $ \SN $ is a Polish space with the product topology, 
$ \WSN $ becomes a Polish space. 
Hence, there exists a countable determining class $ \mathcal{V} $
of $ (\WSN , \mathcal{B}(\WSN ) )$. 
We can take such a class $ \mathcal{V} $ as follows. 
Let $ \mathbf{S}_1 $ be a countable dense subset of $ \SN $, and 
\begin{align*}&
 \mathcal{U} = \mathcal{A} [ 
\{ U_{r}(\mathbf{s}) ; \, 0 < r \in \mathbb{Q},\, \mathbf{s} \in \mathbf{S}_1 
\} ] 
.\end{align*}
Here $ \mathcal{A}[\cdot] $ denotes the algebra generated by $ \cdot $, and 
$ U_{r}(\mathbf{s}) $ is an open ball in $ \SN $ with center $ \mathbf{s}$ and
radius $ r $. We also take a suitable metric defining the same topology of 
the Polish space $ \SN $. 
We note that $ \mathcal{U} $ is countable because the subset 
$ \{ U_{r}(\mathbf{s}) ; \, 0 < r \in \mathbb{Q},\, 
\mathbf{s} \in \mathbf{S}_1 \} $ is countable.
 Let 
\begin{align}\label{:7Jz}
&
 \mathcal{V} = \bigcup_{{k}=1}^{\infty} 
\{ (\wwww _{\mathbf{t}})^{-1} (\mathbf{A})\, ;\, 
\mathbf{A} \in \mathcal{U}^{k}, 
\mathbf{t} \in (\mathbb{Q}\cap (0,\infty))^{k} \} 
.\end{align}
We then see that $ \mathcal{V} $ is a countable determining class of 
$ (\WSN , \mathcal{B}(\WSN ) )$. 

%\smallskip 

%
\begin{lemma} \label{l:7J} \MST. 
Then for each $ \mathbf{V} \in \mathcal{V} $ and for \mB 
\begin{align} \label{:7Jb} 
 %\mathbf{V} \bigcap \WWfixsB = 
 (\Fsbi ) ^{-1} (\mathbf{V}) \bigcap \WWfixsB 
& \in \tTpath (\SN ) \sBc %\bigcap \WWfixsB 
.\end{align}
Here $ \tTpath (\SN ) \sBc $ is the completion of the $ \sigma $-field 
$ \tTpath (\SN ) $ with respect to $ \PsB $. 
\end{lemma}
%\marginpar{[]証明変更\\From \lref{l:62} \thetag{3}, }
\begin{proof} 
Let $ \TpSNsB $ be as \eqref{:62bb}. 
Then from \lref{l:62} \thetag{3}, we deduce that 
$ \FisB $ is a $ \TpSNsB $-measurable function for \mB. 
Then for \mB 
\begin{align}
\notag %\label{:7Ja}
(\Fsbi ) ^{-1} (\mathbf{V}) & \in \TpSNsB 
.\end{align}
Hence for \mB 
\begin{align}
\label{:7Ja}&
(\Fsbi ) ^{-1} (\mathbf{V}) \bigcap \WWfixsB \in \TpSNsB \bigcap \WWfixsB 
.\end{align}
Here for a $ \sigma $-field $ \mathcal{F} $ and a subset $ A $, we set 
$ \mathcal{F} \cap A = \{ F \cap A ; F \in \mathcal{F} \} $ as before. 

Suppose $ \wwww \in (\Fsbi ) ^{-1} (\mathbf{V}) \bigcap \WWfixsB $. 
Then we see 
$ \Fsbi (\wwww ) \in \mathbf{V}$ from $ \wwww \in (\Fsbi ) ^{-1} (\mathbf{V}) $ 
and 
$\www = \Fsbi (\wwww ) $ from $ \wwww \in \WWfixsB $. 
Hence $\www \in \mathbf{V} \bigcap \WWfixsB $. 
We thus obtain 
$$ (\Fsbi ) ^{-1} (\mathbf{V}) \bigcap \WWfixsB \subset \mathbf{V} \bigcap \WWfixsB .$$
Suppose $\www \in \mathbf{V} \bigcap \WWfixsB $. 
Then $ \Fsbi (\wwww ) = \www \in \mathbf{V}$. 
Hence $ \wwww \in (\Fsbi ) ^{-1} (\mathbf{V}) \bigcap \WWfixsB $. 
Thus, %We thus obtain 
$$ \mathbf{V} \bigcap \WWfixsB \subset (\Fsbi ) ^{-1} (\mathbf{V}) \bigcap \WWfixsB 
.$$
Combining these, we see that, for \mB, 
\begin{align}\label{:7Jc}& 
 (\Fsbi ) ^{-1} (\mathbf{V}) \bigcap \WWfixsB = 
\mathbf{V} \bigcap \WWfixsB 
%\in \sigma [\wwww _{\mathbf{t}}] \bigcap \WWfixsB 
.\end{align}
Because $ \mathbf{V} \in \mathcal{V} $, there exists 
a $ \mathbf{t} \in (\mathbb{Q}\cap (0,\infty))^{k} $ such that 
$ \mathbf{V}= (\wwww _{\mathbf{t}})^{-1} (\mathbf{A})$ 
for some $ \mathbf{A}\in \mathcal{U}^{k}$. 
From $ \mathbf{V}= (\wwww _{\mathbf{t}})^{-1} (\mathbf{A})$ we have 
$ \mathbf{V} \in \sigma [\wwww _{\mathbf{t}}] $. Hence we obtain 
\begin{align}
\label{:7Jcc}&
\mathbf{V} \bigcap \WWfixsB 
\in \sigma [\wwww _{\mathbf{t}}] \bigcap \WWfixsB 
.\end{align}
 Combining \eqref{:7Jc} and \eqref{:7Jcc}, we obtain for \mB 
\begin{align}\label{:7Jp}
(\Fsbi ) ^{-1} (\mathbf{V}) \bigcap \WWfixsB 
\in & \sigma [\wwww _{\mathbf{t}}] \bigcap \WWfixsB 
.\end{align}

From \eqref{:7Ja} and \eqref{:7Jp} we deduce for \mB 
\begin{align} \label{:7Jq}
(\Fsbi ) ^{-1} (\mathbf{V}) \bigcap \WWfixsB 
\in & \big\{ \TpSNsB \bigcap \sigma [\wwww _{\mathbf{t}}] \big\} \bigcap \WWfixsB 
.\end{align}
We easily see for \mB 
$$ \TpSNsB \bigcap \sigma [\wwww _{\mathbf{t}}] % \bigcap \WWfixsB 
\subset 
\tTpath (\SN ) \sBc % \bigcap \WWfixsB 
.$$
This together with \eqref{:7Jq} yields 
\begin{align}
\label{:7Js}&
(\Fsbi ) ^{-1} (\mathbf{V}) \bigcap \WWfixsB 
 \in \tTpath (\SN ) \sBc \bigcap \WWfixsB 
.\end{align}
By \lref{l:7X} we have $ \PsB (\WWfixsB ) = 1 $. 
Then \eqref{:7Js} implies \eqref{:7Jb} immediately. 
\qed \end{proof}

\begin{lemma} \label{l:7Q} 
\thetag{1} 
Assume \As{\Done}. Then, for each 
$ \mathbf{A} \in \tTpath (\SN ) $, 
\begin{align}\label{:7Qa}&
 \PsB (\mathbf{A}) \in \{ 0,1 \} \quad \text{ for \mB}
.\end{align}
\thetag{2} Assume \AsDD. 
Then $ \tTTT _{\mathrm{path},\ \Upsilon}^{\{1\}}(\SN ) $ \DOMUL, 
where 
\begin{align*}&
\tTTT _{\mathrm{path},\ \Upsilon}^{\{1\}}(\SN ) = 
\big\{ \mathbf{A} \in \tTpath (\SN ) \, ;\, \PsB (\mathbf{A}) = 1 \text{ for \mB} \big\} 
.\end{align*}
\end{lemma}
\begin{proof} We first prove \thetag{1}. 
By construction we have a decomposition of $ \Pll $ such that 
\begin{align}\label{:7Qx}&
\PX ( \mathbf{A} ) = 
\int_{ \SN \ts \WRNz } \PsB ( \mathbf{A} ) 
\, \Upsilon \dsbb 
.\end{align}
Let 
$ B_{\mathbf{A}} = \{ (\mathbf{s},\bB ) \, ;\, 0 < \PsB ( \mathbf{A} ) < 1 \} $. 
Suppose that \eqref{:7Qa} is false. 
Then there exists an $ \mathbf{A} \in \tTpath (\SN ) $ 
such that $ \Upsilon ( B_{\mathbf{A}} ) > 0$. 
Hence we deduce that 
\begin{align}\label{:7Qb} 
0 < &
 \int_{ B_{\mathbf{A}} } \PsB ( \mathbf{A} ) \Upsilon 
(d\mathbf{s}d\bB ) 
 < 1 
.\end{align}
From \eqref{:7Qx} and \eqref{:7Qb} together with $ \Upsilon ( B_{\mathbf{A}} ) > 0$ 
we deduce $ 0 < \PX ( \mathbf{A} ) < 1$. 
This contradicts \As{\Done}. 
Hence we obtain \thetag{1}. 

We next prove \thetag{2}. 
Suppose that $ \mathbf{A} \in \tTpathTSNone \PX ) $. 
%$=\{ \mathbf{A} \in \tTpath (\SN ); \Pll (\mathbf{A}) = 1 \} $. 
%, where $ \mathbf{A} \in \tTpathTSNone \Pll ) $ is given by \eqref{:76d}. 
%
Then from \eqref{:7Qa} and \eqref{:7Qx}, we deduce that 
$ \PsB (\mathbf{A}) = 1 $ for \mB. %$ \Upsilon $-a.s.\! $ (\sB )$. 
Furthermore, $ \tTpathTSNone \PX )$ \DOMUL\ 
% is independent of $ \lambda $ and the lift dynamics $ \Pl $ 
by \As{\Dtwo}. 
Collecting these, we obtain \thetag{2}. 
\qed \end{proof}

We next prepare a general fact on countable determining classes. 
\begin{lemma} \label{l:7K} \thetag{1} 
Let $ (U,\mathcal{U} )$ be a measurable space 
with a countable determining class 
$ \mathcal{V} =\{ V_n \}_{n\in \mathbb{N}} $. 
Let $ \nu $ be a probability measure on $ (U,\mathcal{U} )$. 
Suppose that $ \nu ( V_n ) \in \{ 0,1 \} $ for all $ n \in \N $. 
Then, $ \nu ( A ) \in \{ 0,1 \} $ for all $ A \in \mathcal{U} $. 
Furthermore, there exists a unique $ \Vnu \in \mathcal{U} $ such that 
$ \Vnu \cap A \in \{ \emptyset , \Vnu \} $ and 
$ \nu ( A ) = \nu ( A \cap \Vnu ) $ for all $ A \in \mathcal{U} $. \\
\thetag{2} In addition to the assumptions in \thetag{1}, we assume 
$ \{ u \} \in \mathcal{U} $ for all $ u \in U $. 
Then there exists a unique $ a \in U $ such that $ \nu = \delta_a $. 
\end{lemma}
\begin{proof} We first prove \thetag{1}. 
Let $ N (1) = \{ n \in \mathbb{N}; \nu (V_n) = 1 \} $. Then we take 
\begin{align*}&
 \Vnu = \big( \bigcap_{n \in N (1)} V_n \big) 
\bigcap 
\big( \bigcap_{n \not\in N (1)} V_n^c \big) 
. \end{align*}
Clearly, we obtain $ \nu ( \Vnu ) = 1$. 

Let $ A \in \mathcal{U} $. 
Suppose that $ \Vnu \cap A \not\in \{ \emptyset , \Vnu \}$. 
We then cannot determine the value of $ \nu ( \Vnu \cap A )$ 
from the value of $ \nu (V_n)$ ($ n \in \mathbb{N}$). 
This yields a contradiction. Hence, 
$ \Vnu \cap A \in \{ \emptyset , \Vnu \}$. 
If $ \Vnu \cap A = \emptyset $, then $ \nu ( A ) = 0 $. 
If $ \Vnu \cap A = \Vnu $, then $ \nu ( A ) \ge \nu ( \Vnu ) = 1 $. 
We thus complete the proof of \thetag{1}. 

We next prove \thetag{2}. 
Since $ \nu (\Vnu ) = \nu (\Vnu \cap U ) = \nu (U) = 1$, we have $ \# \Vnu \ge 1 $. 
 
Suppose $ \# \Vnu = 1$. Then there exists a unique $ a \in U $ such that $ \Vnu = \{ a \} $. 
This combined with \thetag{1} yields that 
 $ \nu (A) = \nu (A \cap \Vnu ) 
= \nu (A \cap \{ a \} ) \in \{ 0,1 \} $ for all $ A \in \mathcal{U} $. 
Hence we see that $ \nu = \delta_a $. 

Next suppose $ \# \Vnu \ge 2 $. 
Then $ \Vnu $ can be decomposed into two non-empty measurable subsets 
$ \Vnu = \Vnu^1 + \Vnu^2$ because $ \{ u \}\in \mathcal{U} $ for all $ u \in U $. 
From \thetag{1}, we have proved that $ \nu ( \Vnu^1 ), \nu ( \Vnu^2 ) \in \{ 0,1 \} $ and that 
$ \Vnu $ is unique. 
Hence such a decomposition $ \Vnu = \Vnu^1 + \Vnu^2$ yields a contradiction. 
This completes the proof of \thetag{2}. 
\qed \end{proof}

Let $ \Pts $ and $ \PsB $ be as \eqref{:70r} and \eqref{:71c}, respectively. 
\begin{lemma} \label{l:79a}
\MST. Assume \As{\Done}. 
Then $ \PsB $ is concentrated at a single path 
$ \wwww = \wsb \in \WSN $, that is, 
\begin{align}\label{:7Aa}&
 \PsB = \delta _{\wwww (\sB )} 
\end{align}
In particular, $ \wwww $ is a function of $ \bB $ 
under $ \Pts $ for $ \PXz $-a.s. $ \mathbf{s}$. 
\end{lemma}

\begin{proof}
From \eqref{:79z} and \eqref{:79x}, we deduce $ \FisB (\mathbf{w}) = \mathbf{w}$ 
for $ \PsB $-a.s.\ $ \mathbf{w}$ for \mB. 
Hence by \eqref{:71c} we deduce for \mB 
\begin{align}\label{:7Ab}&
 \PsB \circ (\FisB )^{-1} = \PsB 
.\end{align}

Let $ \mathcal{V} $ be the countable determining class 
given by \eqref{:7Jz}. 
Then, we deduce from \lref{l:7X}, \lref{l:7J}, \lref{l:7Q}, 
and the definition of $ \PsB $ 
 that, for all $ \mathbf{V} \in \mathcal{V} $, 
\begin{align}\label{:7Ac}
\PsB \circ (\Fsbi ) ^{-1} (\mathbf{V}) = & 
\PsB ( (\Fsbi ) ^{-1} (\mathbf{V}) )
\\\notag = & 
\PsB \big((\Fsbi ) ^{-1} (\mathbf{V}) \bigcap \WWfixsB \big) 
\in 
\{ 0,1 \}\quad \text{ for \mB}
.\end{align}
% Hence from \eqref{:7AC} we see, for all $ \mathbf{V} \in \mathcal{V} $, 
% \begin{align}\label{:7Ac}
% \PsB \circ (\Fsbi ) ^{-1} (\mathbf{V}) & 
% \in \{ 0,1 \} \quad \text{ for \mB}
% .\end{align}
%
Because $ \mathcal{V} $ is countable, we deduce from \eqref{:7Ac} that, for \mB, 
\begin{align}& \label{:7Ad}
\PsB \circ (\Fsbi ) ^{-1} (\mathbf{V})
 \in \{ 0,1 \}
\quad \text{ for all } \mathbf{V} \in \mathcal{V} 
.\end{align}

We denote by $ \Bbar $ the completion of $ \mathcal{B}(\WSN ) $ with respect to $ \Ptsb $. 
From \eqref{:7Ad} and \lref{l:7K}, we obtain for \mB, 
\begin{align}& \notag %\label{:7Ae}
\PsB \circ (\Fsbi ) ^{-1} (\mathbf{A})
 \in \{ 0,1 \} \quad \text{ for all } \mathbf{A} \in \Bbar 
.\end{align}
Furthermore, for \mB, there exists a unique 
$ \mathbf{U}_{\sB } \in \Bbar $ 
such that 
\begin{align}
\notag &% \label{:7Ap}&
 \mathbf{U}_{\sB } \cap \mathbf{A} 
\in \{ \emptyset , \mathbf{U}_{\sB } \} 
,\quad \PsB (\mathbf{A} ) = \PsB (\mathbf{A} \cap \mathbf{U}_{\sB } )
\end{align}
for all $ \mathbf{A} \in \Bbar $. 
Hence we deduce that the set $ \mathbf{U}_{\sB } $ consists of 
a single point $ \{ \wsb \} $ for some $ \wwww = \wsb \in \WSN $. 
Then the probability measure $ \PsB \circ (\Fsbi ) ^{-1} $ is 
concentrated at the single path 
$ \wwww = \wsb \in \WSN $, that is, 
\begin{align*}& \quad \quad 
\PsB \circ (\Fsbi ) ^{-1} = 
 \delta _{\wsb } 
.\end{align*} 
In particular, $ \wwww $ is a function of $ \bB $ under $ \Pts $ 
because of \eqref{:70p}, \eqref{:70r}, and \eqref{:71c}. 
This combined with \eqref{:7Ab} implies \eqref{:7Aa}. 
\qed\end{proof}

\noindent {\em Proof of \tref{l:79}. }
From \eqref{:7Aa}, we immediately obtain \ASTpath. Hence we obtain \thetag{1}. 
Let $ \wsb $ and $ \mathbf{U}_{\sB } $ be as \lref{l:79a}.
From \lref{l:7Q} \thetag{2}, we deduce that the set $ \mathbf{U}_{\sB } $ \DOMUL. 
In particular, $ \wsb $ \DOMUL. 
This completes the proof of \thetag{2}. 
\qed

\ms

\noindent {\em Proof of \tref{l:70}. } 
By assumption, \Conefour\ are satisfied, and furthermore, $ \Pts $ is \9. 
Then we apply \tref{l:77} to obtain \AsDD. 
Hence we deduce from \tref{l:79} that \ASTpath\ and \ASTpatH\ hold 
for $ \PXz $-a.s.\! $ \mathbf{s}$. 
\qed

%%%%%%%%%%%%%%%%%%%%%%%%%%%%%%%%%%%%%%%%%%%%%%%
\section{Proof of Theorems \ref{l:5A}--\ref{l:5B}.}\label{s:Z}
This section is devoted to proving Theorems \ref{l:5A}--\ref{l:5B}. 
Throughout this section, $ \sS $ is $ \Rd $ or a closed set satisfying the assumption in \sref{s:5}, 
and $ \sSS $ is the configuration space over $ \sS $. 
We have established two tail theorems: \tref{l:61} and \tref{l:70}. 
Then \tref{l:5A} and \tref{l:5B} are immediate consequences of them. 

\ms 

\noindent {\em Proof of \tref{l:5A}.} 
%%%%%%%%%%%%%%%%%%%%%%%%%%%%%%
By assumption, $ \mu $ is tail trivial and $ \XB $ under $ P $ is a weak solution satisfying \IASN. 
Hence we deduce \ASTpath\ and \ASTpatH\ for $ \pP \circ \mathbf{X}_0^{-1} $-a.s.\ $ \mathbf{s}$ 
from the second tail theorem (\tref{l:70}). 
We therefore conclude from the first tail theorem (\tref{l:61}) that 
%\ISDEb has a family of unique strong solutions $ \{\Fs \} $ 
\ISDEb has a unique strong solution $ \Fs $ 
starting at $ \mathbf{s}$ for $ \pP \circ \mathbf{X}_0^{-1} $-a.s.\ $ \mathbf{s}$ 
under the constraints of \iFcs\ and \ASTpath. 

Because the family of strong solutions $ \{ \Fs  \} $ is given by a weak solution $ \XB $ under $ \pP $, 
$ \{ \Fs  \} $ satisfies \As{$\mathbf{MF}$}. 
Recall that $ \pPF  = \pP \circ \mathbf{X}^{-1}$ by \eqref{:5Ay}.  
Hence $ \pPF $ satisfies \IASN. 

We next check \thetag{i} and \thetag{ii} in \dref{dfn:43}.

Let $ (\hat{\mathbf{X}},\hat{\mathbf{B}})$ under $ \hat{P}$ be a weak solution in \thetag{i}. 
Then $ (\hat{\mathbf{X}},\hat{\mathbf{B}})$  under $ \hat{P}$ satisfies \IASN\ by assumption. 
This implies $ (\hat{\mathbf{X}},\hat{\mathbf{B}})$  under $ \hat{P}_{\mathbf{s}}$ satisfies 
\iFcs, \As{\sIn}$ _{\mathbf{s}}$, and \As{\nbj}$_{\mathbf{s}} $ 
for $ \hat{\pP }\circ \hat{\mathbf{X}}_0^{-1} $-a.s.\ $ \mathbf{s}$ 
by Fubini's theorem. 
Because $ \Fs $ is a unique strong solution 
starting at $ \mathbf{s}$ for $ \pP \circ \mathbf{X}_0^{-1} $-a.s.\ $ \mathbf{s}$ 
under the constraints of \iFcs\ and \ASTpath\ and 
$\hat{P}\circ \hat{\mathbf{X}}_0^{-1} \prec P \circ \mathbf{X}_0^{-1} $, we obtain \thetag{i}. 
Condition \thetag{ii} is clear because $ \Fs $ is a unique strong solution. 
This completes the proof. 
\bbbbb

\noindent {\em Proof of \tref{l:5B}.} 
By \eqref{:54c} $ \mut $ is tail trivial. Hence \As{\muTT} for $ \mut $ is satisfied. 
By assumption, we have $\XB $ under $ \pP ^{\aa } $ satisfies \As{\muAC} for $ \mut $ 
for $ \mu $-a.s.\ $ \aa $. 

We see that, for $ \mu $-a.s.\ $ \aa $, \As{\iFc}, \As{\sIn}, and \As{\nbj} 
for $ \XB $ under $ \pP ^{\aa } $ follow from Fubini's theorem and those for $ \XB $ under $ \pP $. 

Indeed, from \eqref{:54b}, \eqref{:54h}, and \eqref{:54j} we obtain 
\begin{align}
\label{:Z1a}& 
\pP = \int \pP ^{\aa } \mu (d\aa )
.\end{align}
Then from \As{\sIn} and \As{\nbj} for $ \XB $ under $ \pP $ and \eqref{:Z1a} we deduce 
\begin{align*}&
1=\pP (\mathsf{X} \in \WSsiNE ) = 
\int_{\sSS } \pP ^{\aa } (\mathsf{X} \in \WSsiNE ) \mu (d\aa ) 
,\\&
 1 = \pP ( \mrXX < \infty ) = \int_{\sSS } \pP ^{\aa }( \mrXX < \infty ) \mu (d\aa )
.\end{align*}
Hence $ \pP ^{\aa } (\mathsf{X} \in \WSsiNE ) =\pP ^{\aa }( \mrXX < \infty )=1 $ for $ \mu $-a.s.\ $ \aa $. 
This implies \As{\sIn} and \As{\nbj} for $ \XB $ under $ \pP ^{\aa } $ holds 
for $ \mu $-a.s.\ $ \aa $.

By disintegration of $ \pP \circ \mathbf{X}_0^{-1} $, 
\begin{align}\label{:Z1b}&
\pP \circ \mathbf{X}_0^{-1} = \int \pP ^{\aa }\circ \mathbf{X}_0^{-1} \mu (d\aa )
.\end{align}
We set $ \Ps = P (\cdot | \mathbf{X}_0 = \mathbf{s})$. 
By definition, \As{\iFc} for $ \XB $ under $ \pP $ implies 
\iFcs\ for $ \XB $ under $ \Ps $ for $ \pP \circ \mathbf{X}_0^{-1} $-a.s.\ $ \mathbf{s}$. 
Then by \eqref{:Z1b} this implies \iFcs\ for $ \XB $ under $ \Ps $ holds 
for $ \pP ^{\aa }\circ \mathbf{X}_0^{-1} $-a.s.\ $ \mathbf{s}$ and for $ \mu $-a.s. $ \aa $. 
We easily see for $ \mu $-a.s. $ \aa $ 
%From \eqref{:54i}, we easily see for $ \mu $-a.s. $ \aa $ 
\begin{align*}&\Ps = \Ps ^{\aa }
\quad \text{for $ \pP ^{\aa }\circ \mathbf{X}_0^{-1} $-a.s.\ 
 $ \mathbf{s}$. }
\end{align*}
Hence, for $ \mu $-a.s. $ \aa $, we have \iFcs\ for $ \XB $ under $ \Ps ^{\aa }$ 
for $ \pP ^{\aa }\circ \mathbf{X}_0^{-1} $-a.s.\ $ \mathbf{s}$. 
Then, for $ \mu $-a.s.\ $ \aa $, we deduce \As{\iFc} for $ \XB $ under $ \pP ^{\aa }$. 

We have thus seen that, for $ \mu $-a.s. $ \aa $, 
$ \XB $ under $ \pP ^{\aa } $ fulfills the assumptions of \tref{l:5A}. 
Hence \tref{l:5B} follows from \tref{l:5A}. 
\bbbbb

\Section{The Ginibre interacting Brownian motion}\label{s:30}

In this section, we apply our theory to the special example of the Ginibre interacting Brownian motion and prove \existence\ of strong solutions and \pathwiseuniqueness. Our proof is based on the idea explained in Introduction. Behind it there are two general theories called tail theorems. These two theories are robust and can be applied to various kinds of infinite-dimensional stochastic (differential) equations with symmetry beyond interacting Brownian motions.
The purpose of this section is to clarify the roles of these two theories by applying them to the Ginibre interacting Brownian motion.

\Ssection{ISDE related to the Ginibre random point field.}\label{s:3} 
In \sref{s:3}, we introduce the Ginibre interacting Brownian motion and prepare the result of the first step. 

The Ginibre random point field $ \mug $ is a random point field on $ \Rtwo $ 
whose $ n $-point correlation function $ \rhog ^n $ with respect to the Lebesgue measure is given by 
\begin{align}&\label{:30a}
\rhog ^{n}(x_1,\ldots,x_n) = \det [\kg (x_i,x_j)]\ijn 
,\end{align} 
where $ \map{\kg }{\Rtwo \ts \Rtwo }{\mathbb{C}}$ 
is the kernel defined by 
\begin{align} \label{:30b}
&
\kg (x,y) = \pi ^{-1} 
e^{-\frac{|x|^{2}}{2}-\frac{|y|^{2}}{2}}\cdot 
e^{x \bar{y}}
.\end{align}
Here we identify $ \Rtwo $ as $ \mathbb{C}$ by the obvious correspondence 
$ \Rtwo \ni x=(x_1,x_2)\mapsto x_1 + \sqrt{-1} x_2 \in \mathbb{C}$, 
and 
$ \bar{y}=y_1-\sqrt{-1} y_2 $ is the complex conjugate in this identification. 

It is known that $ \mug $ is translation and rotation invariant. 
Moreover, $ \mug (\Ssi )=1 $ and $ \mug $ is tail trivial \cite{bqs,ly.18,o-o.tail}. 

We next introduce the Dirichlet form associated with $ \mug $ and 
construct $ \sSS $-valued diffusion. 
Let $ (\mathcal{E}^{\mug } ,\di ^{\mug } )$ be a bilinear form 
on $ \Lmg $ defined by 
\begin{align} 
\notag &% \label{:30c} & 
\di ^{\mug } = \{ f \in \di \cap \Lmg \, ;\, 
\mathcal{E}^{\mug }(f,f) < \infty \} 
,\\\notag &% \label{:30d}& 
\mathcal{E}^{\mug } (f,g) = \int_{\sSS } \DDD [f,g] \, \mug (d\mathsf{s})
,\\ \notag &% \label{:30e}& 
\DDD [f,g] (\mathsf{s}) = \frac{1}{2} 
\sum_{i} (\partialsi \check{f} , \partialsi \check{g} )_{\R ^2 }
.\end{align}
Here $ \mathsf{s}=\sum_i \delta_{s_i}$, 
$ \partialsi = (\PD{}{s_{i1}},\PD{}{s_{i2}}) $, and 
$(\cdot , \cdot )_{\R ^2 }$ denotes the standard inner product in 
$ \Rtwo $. $ \check{f}$ is defined before \eqref{:20j}.

\begin{lemma} [{\cite[Theorem 2.3]{o.rm}}] \label{l:31} 
\thetag{1}
The Ginibre random point field $ \mug $ is 
a \\
$ (|x|^2,- 2\log |x-y|)$-quasi Gibbs measure. 
\\\thetag{2} 
$ (\mathcal{E}^{\mug } ,\di ^{\mug } )$ is closable on $ \Lmg $. 
\\\thetag{3} 
The closure $ (\mathcal{E}^{\mug } ,\dom ^{\mug } )$ of 
$ (\mathcal{E}^{\mug } ,\di ^{\mug } )$ on $ \Lmg $ is a quasi-regular Dirichlet form. 
\\\thetag{4}
There exists a diffusion $ \{ \PPs \}_{\mathsf{s}\in \sSS } $
 associated with $ (\mathcal{E}^{\mug } ,\dom ^{\mug } )$ on $ \Lmg $. 
\end{lemma}

A family of probability measures $ \{ \Pssf \}_{\mathsf{s}\in \sSS } $ 
on $ (\WS ,\mathcal{B}(\WS )) $ is called a diffusion if the canonical process 
$ \mathsf{X}=\{ \mathsf{X}_t \} $ under $ \Pssf $ 
is a continuous process %$ \mathsf{X}=\{ \mathsf{X}_t \} $ 
with the strong Markov property starting at $ \mathsf{s}$. 
Here $ \mathsf{X}_t(\ww ) = \wwt $ for $ \ww =\{ \wwt \} \in \WS $ by definition. 
$ \mathsf{X}$ is adapted to $ \Ft $, where 
$ \mathcal{F}_t = \cap_\nu \mathcal{F}_t^{\nu} $ and 
the intersection is taken over all Borel probability measures $ \nu $, 
$ \mathcal{F}_t^{\nu}$ is the completion of 
$ \mathcal{F}_t^+ = \cap_{\epsilon > 0 } \mathcal{B}_{t+\epsilon } (\sSS ) $ 
with respect to $ \PP _{\nu }=\int \PPs \nu (d\mathsf{s})$. 
The $ \sigma $-field $ \Bt (\sSS ) $ is defined by 
\begin{align}\label{:32a}&
\Bt (\sSS ) = \sigma [\ww _s; 0 \le s \le t ]
.\end{align}
Furthermore, $ \{ \Pssf \}_{\mathsf{s}\in \sSS } $ is called stationary 
if it has an invariant probability measure. 

We refer to \cite{mr} for the definition of quasi-regular Dirichlet forms. 
We also refer to \cite{fot.2} for the Dirichlet form theory. 

We recall the definition of capacity [Chapter 2.1 in \cite{fot.2}].
Denote by $\mathcal{O}$ the family of all open subsets of $\mathsf{S}$. 
Let $(\mathcal{E}, \mathcal{D})$ be a quasi-regular Dirichlet form on $L^2(\mathsf{S},\mu)$, and $\mathcal{E}_1(f,f)=\mathcal{E}(f,f)+ ( f,f )_{L^2(\mathsf{S},\mu)}$.
For $\mathsf{B}\in \mathcal{O}$ we define
\begin{align} \notag &
\mathrm{Cap}^\mu (\mathsf{B})=
\begin{cases}
\inf_{u\in\mathcal{L}_{\mathsf{B}}}\mathcal{E}_1(u,u),
\quad &\mathcal{L}_{\mathsf{B}}\not= \emptyset
\\ 
\infty \quad &\mathcal{L}_{\mathsf{B}}= \emptyset,
\end{cases}
\end{align}
where 
$ \mathcal{L}_{\mathsf{B}}=\{u \in \mathcal{D} : u\ge 1, \mbox{ $\mu$-a.e on $\mathsf{B}$} \}$, 
and we let for all set $\mathsf{A}\subset \mathsf{S}$ 
\begin{align} \notag 
\mathrm{Cap}^\mu (\mathsf{A})=
\inf_{\mathsf{A}\subset \mathsf{B} \in \mathcal{O} }\mathrm{Cap}^\mu(\mathsf{B})
.\end{align}
We call this one-capacity of $\mathsf{A}$ or simply the capacity of $\mathsf{A}$. 

We recall the notion of quasi-everywhere and quasi-continuity. 
Let $\mathsf{A}$ be a subset of $\mathsf{S}$.
A statement depending on $\mathsf{s}\in\mathsf{A}$ is said to hold 
quasi-everywhere (q.e.) on $\mathsf{A}$ if there exists a set $\mathcal{N}\subset \mathsf{A}$ of zero capacity such that the statement is true for every $\mathsf{s}\in \mathsf{A}\setminus \mathcal{N}$. 
When $\mathsf{A}=\mathsf{S}$, quasi-everywhere on $\mathsf{S}$ is simply said quasi-everywhere.
Let $u$ be an extended real valued function defined q.e. on $\mathsf{S}$. 
We call $u$ quasi-continuous if there exists for any $\epsilon>0$ 
an open set $\mathsf{G}\subset \mathsf{S}$ such that 
$\mathrm{Cap}^\mu(\mathsf{G})<\varepsilon$ and the restriction 
$u|_{\mathsf{S}\setminus \mathsf{G}}$ of $u$ on $\mathsf{S}\setminus \mathsf{G}$ is finite continuous.

We set the probability measure $ \PP _{\mug } $ by 
\begin{align} \label{:32d}& \quad \quad 
 \PP _{\mug } (A) = \int _{\sSS } \PPs (A) \mug (d\mathsf{s}) 
\quad \text{ for } A \in \mathcal{F} 
.\end{align}
Under $ \PP _{\mug } $, the unlabeled process $ \mathsf{X}$ is a $ \mug $-reversible diffusion. 
Let $ \mathrm{Cap}^{\mug } $ be the capacity on the Dirichlet space 
$ (\mathcal{E}^{\mug } ,\dom ^{\mug } , \Lmg ) $. 
Let $ \Ssi $ and $ \WSsiNE $ be as \eqref{:23g} and \eqref{:23gg}, respectively . 

%We recall: 
\begin{lemma}[{\cite{o.col,o.isde}}] \label{l:32} 
$ \PP _{\mug } (\WSsiNE ) = 1 $ and $ \mathrm{Cap}^{\mug }(\Ssi ^c) = 0 $. 
\end{lemma}
\begin{proof}
$ \PP _{\mug } (\WT _{\mathrm{NE}}(\sSS )) = 1 $ follows from \cite[\thetag{2.10}]{o.isde}. 
Let $ \mathsf{S}_{\mathrm{s}}$ be as \eqref{:23g}. 
Then by \cite[Theorem 2.1]{o.col}, we see 
$ \mathrm{Cap}^{\mug }(\mathsf{S}_{\mathrm{s}} ^c)= 0 $. 
Clearly, $ \mug (\Ssi )= 1 $. 
Combining these imply the unlabeled diffusion never hits the set consisting of the finite configurations. 
We thus see the capacity of this set is zero. 
Hence we obtain $ \mathrm{Cap}^{\mug }(\Ssi ^c) = 0 $. 
In particular, $ \PP _{\mug } (\WSsi ) = 1 $, which together with 
$ \PP _{\mug } (\WT _{\mathrm{NE}}(\sSS )) = 1 $ implies the first claim. 
\qed \end{proof}

\begin{lemma} [{\cite[Theorem 61, Lemma 72]{o.isde}}]\label{l:33} 
The Ginibre random point field $ \mug $ has 
a logarithmic derivative $ \dG $. 
Furthermore, $ \dG $ has plural expressions: 
\begin{align}\label{:33a}&
 \dG (x, \mathsf{s}) = 2 \limi{r} \sum_{|x-s_i| < r } \frac{x-s_i}{|x-s_i|^2} 
, \\ \label{:33b}&
 \dG (x, \mathsf{s}) = -2x + 
 2 \limi{r} \sum_{|s_i| < r } \frac{x-s_i}{|x-s_i|^2}
.\end{align}
Here the convergence takes place in 
$ \LlocpG $ for any $ 1 \le p < 2 $. 
\end{lemma}

From \lref{l:32} the process 
$ \mathbf{X}=(X^i)_{i\in\N } = \lpathX \in \WT (\RtwoN )$ 
is well defined, where $ \RtwoN = (\R ^2)^{\N }$. 
The unlabeled process $ \mathsf{X}$ is defined on the canonical filtered space 
as the evaluation map $ \mathsf{X}_t (\mathsf{w}) = \mathsf{w}_t$. 
That is, $\OFPsF $ is given by 
$ \Omega = \WSsi $, $ \mathcal{F}=\mathcal{B}(\WSsi ) $, 
$ \{\PPs \}$ is the family of diffusion measures given by \lref{l:31} \thetag{4}. 
$ \{\PPs \}$ can be regarded as diffusion on $ \Ssi $ by \lref{l:32}. 
The labeled process $ \mathbf{X}= \lpathX $ is thus defined on $ \OFPsF $. 
The ISDE satisfied by $ \mathbf{X}$ is as follows. 
\begin{lemma}[{\cite[Theorem 21, Theorem 22]{o.isde}}] \label{l:34}
Let $ \mathbf{X}= \lpathX $ be the stochastic process defined on 
$ \OFPsF $ as above. 
Then there exists a set $ \mathsf{H} $ such that 
$$ \mug (\mathsf{H} ) = 1 
,\quad 
\mathsf{H} \subset \Ssi 
$$
and that $ \mathbf{X}$ under $ \PPs $ satisfies both ISDEs on 
$ \RtwoN $ starting at each point $ \mathbf{s}=\lab (\mathsf{s}) \in \lab (\mathsf{H} ) $
\begin{align} \tag{\ref{:11g}} & 
dX_t^i = dB_t^i + \lim_{r\to\infty } 
\sum_{|X_t^i-X_t^j|<r ,\, j\not=i} 
\frac{X_t^i-X_t^j}{|X_t^i-X_t^j|^{2}} dt \quad (i \in\N )
,\\ \tag{\ref{:11h}} & 
dX_t^i = dB_t^i - X_t^i dt + \lim_{r\to\infty }
\sum_{|X_t^j|<r ,\, j\not=i} \frac{X_t^i-X_t^j}{|X_t^i-X_t^j|^{2}} dt 
\quad (i \in\N )
.\end{align}
Furthermore, %$ \mathbf{X}$ satisfies 
$ \mathbf{X} \in \WT (\ulab ^{-1}(\mathsf{H} ))$. The process $ \mathbf{B}= (B^i)_{i=1}^{\infty}$ 
in \eqref{:11g} and \eqref{:11h} is the same and is the $ \RtwoN $-valued, \FtB given by the formula 
\begin{align}&\label{:34c}
B_t^i = X_t^i -X_0^i -\int_0^t 
\lim_{r\to\infty } \sum_{|X_s^i-X_s^j|<r ,\, j\not=i} 
\frac{X_s^i-X_s^j}{|X_s^i-X_s^j|^{2}} ds 
,\\ &\label{:34d}
B_t^i = X_t^i -X_0^i +\int_0^t X_s^i ds - \int_0^t 
\lim_{r\to\infty } \sum_{|X_s^j|<r ,\, j\not=i} %{|X_s^i-X_s^j|<r ,\, j\not=i} 
\frac{X_s^i-X_s^j}{|X_s^i-X_s^j|^{2}} ds 
.\end{align}
\end{lemma}

We note that $ \XB $ above is a solution of 
the two ISDEs \eqref{:11g} and \eqref{:11h}. 
We refer to \dref{d:51} for the definition of the solution in \lref{l:34}, 
which is often called a weak solution. 
We also note that the identity between \eqref{:34c} and \eqref{:34d} follows from the plural expressions of $ \dG $ in \lref{l:33}. 

\begin{remark}\label{r:G3} \thetag{1} 
We cannot replace $ \WT (\ulab ^{-1}(\mathsf{H} ))$ by $ \WT (\lab (\mathsf{H} ))$ in 
\lref{l:34}. Indeed, $ \lab (\mathsf{H} ) \subset \ulab ^{-1}(\mathsf{H} ) $ and 
$ \mathbf{X}_t \not\in \lab (\mathsf{H} ) $ for some $0 < t < \infty $. 
\\
\thetag{2}
The Brownian motion $ \mathbf{B} = (B^i)_{i=1}^{\infty}$ 
in \lref{l:34} is given by a functional of $ \mathsf{X}$. Hence we write 
$ \mathbf{B}(\mathsf{X})=(B^i(\mathsf{X}))_{i=1}^{\infty}$. 
It is given by the martingale part of the Fukushima decomposition of the Dirichlet process 
$ X_t^i-X_0^i$. Indeed, for \eqref{:11g} and \eqref{:11h}, $ B^i$ is given 
by \eqref{:34c} and \eqref{:34d}, respectively. 

We note that the coordinate function $ x_i$ does not belong to 
the domain of the unlabeled Dirichlet form even if locally. 
Hence we introduced in \cite{o.isde} the Dirichlet space of the $ m $-labeled process $ (X^1,\ldots,X^m ,\sum_{j>m}\delta_{X^j})$ to apply the Fukushima decomposition to 
the coordinate function $ x_i$, where $ i \le m $. 
The consistency of the Dirichlet spaces plays a crucial role 
for this argument \cite{o.tp,o.isde}. 
See \sref{s:90} for the definition of the Dirichlet space of the $ m $-labeled process and their consistency, which is given under the general situation. 
\\
\thetag{3} The function in the coefficient in \eqref{:11h} belongs to the domain of the $ m $-labeled Dirichlet space locally. 
Indeed, regarded as a function on $ \Rtwo \ts \mathsf{S}$, we prove that it 
belongs to the domain of one-labeled Dirichlet form locally (see \lref{l:W7}). 
By the same argument we can prove that it is 
in the domain of the $ m $-labeled Dirichlet form 
if we regard as a function on 
$ (\Rtwo )^m \ts \mathsf{S}$ in an obvious manner. 
Hence, by taking a quasi-continuous version of the function, 
the drift term becomes a Dirichlet process. We thus see that the process 
\begin{align*}&
 - X_t^i + \lim_{r\to\infty }
\sum_{|X_t^j|<r ,\, j\not=i} \frac{X_t^i-X_t^j}{|X_t^i-X_t^j|^{2}} 
\end{align*}
in the drift term is a continuous process. This makes the meaning of the drift term more explicit because we usually take a predictable version of the coefficients (see pp 45--46 in \cite{IW}). The key point here is that the coefficient is in the domain of the Dirichlet space. 
All examples in the present paper enjoy this property. 
We shall assume this in \As{\CONE}--\As{\CTWO} and use this in the proof of \pref{l:U2}. 
\end{remark}

\Ssection{Main result for the Ginibre interacting Brownian motion.}
\label{s:32}
All above results in \sref{s:3}
belong to the first step explained in Introduction (\sref{s:1}). 
Our purpose in \sref{s:32} is to obtain the existence of strong solutions and 
the pathwise uniqueness of the solution, which is the main result for 
the Ginibre interacting Brownian motion. 

Let $ \OFPFs $ be as \lref{l:34}. 
Let $ \XB $ be the solution of ISDEs \eqref{:11g} and \eqref{:11h} given by \lref{l:34}. 
Recall that $ \mathbf{X}= \lpathX $, where 
$ \mathsf{X}$ is the canonical process such that 
$ \mathsf{X}_t (\mathsf{w}) = \mathsf{w}_t$ and $ \mathbf{B}$ 
is the $ \Ft $-Brownian motion given by \eqref{:34d}. 
Let $ \PP _{\mug } $ be as \eqref{:32d}. 

\begin{theorem}	\label{l:35} %Let $ \XB $ be the solution given by \lref{l:34}. \\
\thetag{1} 
 $ \XB $ under $ \PP _{\mug } $ 
satisfies the same conclusions as \tref{l:5A},  \corref{l:53c1}, and \corref{l:53c2}  
for ISDE \eqref{:11g}. 
\\\thetag{2} 
There exists a subset $ \HHz $ of $ \SSsdeg $ such that $ \mug (\HHz ) = 1 $ and that 
$ \XB $ is a strong solution of ISDE \eqref{:11g} 
 starting at $ \mathbf{s}=\lab (\mathsf{s}) \in \lH $ defined on $ \OFPFs $. 
% \\\thetag{2} 
% A solution $ \hat{\mathbf{X}} $ of ISDE \eqref{:11g} defined on $ (\Omega , \mathcal{F}, P ,\Ft )$ 
% satisfying \Inmg\ is pathwise unique in the sense that, if $ \mathbf{X}'$ and $ \mathbf{X''}$ are 
% such solutions with the same $ \Ft $-Brownian motion and $ \mathbf{X}_0'=\mathbf{X}_0''$ a.s., then 
% \begin{align}\label{:35p}& P (\mathbf{X}'=\mathbf{X}'') = 1 .\end{align}
% For $ \mug $-a.s.\,$ \mathsf{s}$, the distribution of $ \mathbf{X}'$ starting at 
% $ \lab (\mathsf{s})$ coincides with $ \mathbf{X}$. 
\\\thetag{3} 
The same statements as \thetag{1} hold for ISDE \eqref{:11h}. 
\\\thetag{4} 
The same statements as \thetag{2} hold for ISDE \eqref{:11h}. 
\\
\thetag{5} 
Let $( \mathbf{X}',\mathbf{B}')$ and $ (\mathbf{X}'',\mathbf{B}')$ under $ \pP '$ be weak solutions of 
\eqref{:11g} and \eqref{:11h} satisfying \Inmg, respectively. 
Suppose that both solutions are defined on the same filtered space $ \OFPFF $
with the same $ \{ \mathcal{F}_t' \} $-Brownian motion $ \mathbf{B}'$ 
and that $ \mathbf{X}_0' = \mathbf{X}_0''$ a.s.\ 
Then $ \pP ' (\mathbf{X}_t'= \mathbf{X}_t'' \text{ for all }t ) = 1 $. 
\end{theorem}

We remark that ISDEs \eqref{:11g} and \eqref{:11h} are different ISDEs in general. 
\tref{l:35} \thetag{5} asserts that, if the unlabeled particles start from a support 
$ \SSsdeg $ of $ \mug $ and if the label $ \lab $ is common, 
then these two labeled dynamics are equal all the time. 
The intuitive explanation of this fact is as follows. 
One may regard the set $ \ulab ^{-1} (\SSsdeg )$ 
as a sub-manifold of $ \RtwoN $ and the drift terms $ b_1 $ and $ b_2$ 
\begin{align}& \notag %\label{:24F}&
b_1(x_i, \sum_{j\not=i} \delta_{x_j}) = 
\limi{r} \sum_{|x_i-x_j|<r} \frac{x_i-x_j}{|x_i-x_j|^2}
,\\& \notag 
b_2(x_i,\sum_{j\not=i} \delta_{x_j}) = - x_i + 
\limi{r} \sum_{|x_j|<r} \frac{x_i-x_j}{|x_i-x_j|^2}
\end{align}
of each ISDE are regarded as \lq\lq tangential vectors on $ \ulab ^{-1} (\SSsdeg ) $''. 
In \cite{o.isde}, it was shown that 
both drifts are equal on $ \ulab ^{-1} (\SSsdeg ) $. This implies the coincidence of 
ISDEs \eqref{:11g} and \eqref{:11h} on $ \ulab ^{-1} (\SSsdeg ) $. 
Since the drift terms $ b_1$ and $ b_2 $ are tangential, the solutions stay in 
$ \ulab ^{-1} (\SSsdeg ) $ all the time, 
which combined with the pathwise uniqueness of the solutions of ISDEs yields \thetag{4}. 
% It would be mre precise that we interpret $ \ulab ^{-1} (\SSsdeg ) $ as a subset of 
% $ \RtwoN $, 
%and regarded as a universal covering of $ \SSsdeg $. 

% Note that for a given Dirichlet space, the associated diffusion is unique. 
% Hence a solution of ISDE associated with the Dirichlet space is unique. 
% This solution however depends on the Dirichlet space, 
% and thus we do not have uniqueness of a solution 
% without the uniqueness of the Dirichlet space. 
% The pathwise uniqueness of solution of ISDE was therefore left open in \cite{o.isde}. 

We note that the unlabeled dynamics $ \mathsf{X}$ are $ \mug $-reversible because $ \mathsf{X}$ is given by the symmetric Dirichlet form $ (\mathcal{E}^{\mug } ,\dom ^{\mug } )$ on $ \Lmg $ (see \lref{l:31} \thetag{4}). Hence, the distribution of $ \mathsf{X}_t$ with initial distribution $ \mug $ satisfies 
$ \mug \circ \mathsf{X}_t^{-1} = \mug $ for each $ \zzti $. 
In contrast, the labeled dynamics $ \mathbf{X}_t$ 
are trapped on a very thin subset of the huge space $ \RtwoN $. 
We conjecture that the distribution of $ \mathbf{X}_t$ is singular 
to the initial distribution $ \mug \circ \lab ^{-1} $ for some $ t > 0 $.

\section{Proof of \tref{l:35}.}\label{s:4}

In \sref{s:4}, we prove the main theorem (\tref{l:35}) using \tref{l:5A}.

\Ssection{Localization of coefficients and Lipschitz continuity. }\label{s:41}

Recall that the labeled process $ \mathbf{X}=(X^i)_{i \in \N }= \lpathX $ is obtained from 
the unlabeled process $ \mathsf{X}$ under $\PP _{\mug } $. % such that $ X^i = \lpathi (\mathsf{X}) $. 
Set the $ m $-labeled process $ (\mathbf{X}^m,\mathsf{X}^{m*})$ by 
\begin{align}\notag &%\label{:41a}&
(\mathbf{X}^m,\mathsf{X}^{m*})=
(X_1,\ldots,X_m,\sum_{i=m+1}^{\infty}\delta_{X^i})
.\end{align}
This correspondence is similar to the correspondence of 
$ m $-labeled path in the sense of \eqref{:23p} and \eqref{:23q}. 
This relationship between the unlabeled process and the $m$-labeled process is called consistency.

The $ m $-labeled process is associated with the $m$-labeled Dirichlet space given in \sref{s:E} for the $m$-Campbell measure $ \mug ^{[m]}$ of $ \mug $. Here 
\begin{align}\notag &% \label{:41b}&
\mug ^{[m]} (d\mathbf{x}d\mathsf{y}) = \rhog ^m (\mathbf{x}) \mugxx ( d\mathsf{y})d\mathbf{x} 
,\end{align}
where $ \rhog ^m $ is the $ m $-point correlation function of $ \mug $ 
with respect to the Lebesgue measure $ d\mathbf{x} $ on $ \R ^{2m} $, and 
$ \mugxx $ is the reduced Palm measure conditioned at $ \mathbf{x} \in \R ^{2m}$. 
The $ m $-labeled Dirichlet form is given by 
\begin{align} \label{:41c} &
\E ^{\mug ^{[m]}} (f,g) = \int_{\R ^{2m} \ts \sSS } \Big\{
\frac{1}{2}\sum_{i=1}^{m} 
( \PD{f}{x_i}, \PD{g}{x_i})_{\Rtwo } + 
\DDD [f,g] \Big\} \mug ^{[m]} (d\mathbf{x} d\mathsf{y})
,\end{align}
where $ \partial /\partial {x_i}$ is the nabla in $ \Rtwo $. 
This coincides with the Dirichlet form \eqref{:90b} with $ d=2$, $ \sS = \R ^2$, 
and $ \mum = \mug ^{[m]}$. Furthermore, $ a (\mathbf{x},\mathsf{y})$ in \eqref{:90b} is taken to be 
the $ 2\ts 2$ unit matrix. 
We denote by $ \mathrm{Cap}^{\mug ^{[m]}}$ the capacity 
given by the Dirichlet form $ \E ^{\mug ^{[m]}}$.

\smallskip 

Let $ \mathbf{a}=\{ \ak \}_{\qq \in\mathbb{N}} $ be an increasing sequence of increasing sequences 
$ \ak = \{ \ak (r) \}_{r\in\mathbb{N}} $ such that 
$ \ak (r) < \akk (r)$ and $ \ak (r) < \ak (r+1)$ for all $ \qq , r \in \N $ 
and that $ \limi{r} \ak (r) = \infty $ for all $ \qq \in \N $. 
We take $ \ak (r) \in \N $. 

Let $ \mathsf{K}[\ak ] =\{ \mathsf{s}\, ;\, \mathsf{s} (\Sr ) \le \ak (r) \text{ for all } r \in \mathbb{N} \} $. 
Then $ \mathsf{K}[\ak ] \subset \mathsf{K}[\akk ]$ for all $ \qq \in\mathbb{N}$. 
It is easy to see that $ \mathsf{K}[\ak ]$ is a compact set 
 in $ \sSS $ for each $ \qq \in \mathbb{N}$. 
Let 
\begin{align}\notag &% &\label{:40a}
\Ka = \bigcup_{\qq =1}^{\infty} \mathsf{K}[\ak ] 
.\end{align}
We take $ \ak (r) = \qq r^2 $. Then because $ \mug $ is translation invariant, we have 
\begin{align}\notag &% &\label{:40b}
\mug (\Ka ) = 1 
.\end{align}

We introduce an approximation of $ \R ^{2m} \ts \sSS $ consisting of compact sets. 
Let 
\begin{align}
\notag &% \label{:40x}&
\Ssi ^{[m]} = \{ (\mathbf{x},\mathsf{s})\in \R ^{2m} \ts \sSS \, ;\, 
\ulab (\mathbf{x}) + \mathsf{s} \in \Ssi \} 
,\end{align}
where $ \mathbf{x}=(x_1,\ldots,x_m)$ and $ \ulab (\mathbf{x}) = \sum_{i=1}^m \delta_{x_i}$. 
Let $ \ak ^+$ such that $ \ak ^+ (r) = \ak (r+1)$ and 
\begin{align}\label{:41Y}&
\RRprs = \big\{ \mathbf{x} \in \Stwoqm \, ;\,
\min_{j\not=k } |x_j-x_k | \ge 2^{-\pp } ,\ 
\inf_{l,i} |x_l-s_i| \ge 2^{-\pp } \big\} 
,\end{align}
where $ j,k,l=1,\ldots,m $, $\mathsf{s}=\sum_i \delta_{s_i}$, and 
$ \Srr ^m = \{ x \in \Rtwo ;\, |x| \le \rr \}^m $. 
We set 
\begin{align} \label{:41y} &
\Hb = \big\{ (\mathbf{x},\mathsf{s}) \in \Ssi ^{[m]} \, ;\, \ 
 \mathbf{x} \in \RRprs ,\ \mathsf{s} \in \Kakk \big\} 
,\\
\label{:41z}& 
\Ha = \bigcup_{\rr =1}^{\infty} \Ha _{\rr } 
, \quad 
 \Ha _{\rr } = \bigcup_{\qq =1}^{\infty} \Ha _{\qq ,\rr }, 
\quad 
\Ha _{\qq ,\rr } = \bigcup_{\pp =1}^{\infty} \Hb 
.\end{align}
Although $ \RRprs , \Hb , \Ha _{\qq ,\rr } , \Ha _{\qq } $, and $ \Ha $ 
depend on $ m \in \N $, we suppress $ m $ from the notation. 
To simplify the notation we set $ \NNN = \NNN _1 \cup \NNN _2 \cup \NNNthree $, where 
\begin{align}
\label{:42u}&
\NNN _1 = \{ \rr \in \N \}, \ 
\NNN _2 = \{ (\qq , \rr ) \ ;\, \qq , \rr \in \N \} , \ 
\NNNthree = \{ (\pqr ) \ ;\, \pqr \in \N \}
\end{align}
and for $ \nnNNN $ we define $ \nn + 1 \in \NNN $ such that 
\begin{align}
\label{:42v}&
 \nn + 1 = 
\begin{cases}
(\pp + 1, \qq , \rr ) &\text{ for $ \nn = (\pqr ) \in \NNNthree $, }\\
(\qq +1, \rr )&\text{ for $ \nn = (\qq , \rr ) \in \NNN _2$, }\\
 \rr +1 &\text{ for $ \nn = \rr \in \NNN _1$. }
\end{cases}
\end{align}
We shall take a limit in $ \nn $ along with the order $ \nn \mapsto \nn + 1$. 
We write $ \Han = \Ha _{\pqr }$ for $ \nn = (\pqr ) \in \NNNthree $. 
We set $ \Han $ for $ \nn = (\qq ,\rr ) \in \NNN _2$ and $ \nn = \rr \in \NNN _1$ similarly. 

We remark that $ \Han $ is compact for $ \nn \in \NNNthree $. 
This property has critical importance in the proof of \pref{l:42}.

We set 
$ \sS _{\rr }^{m,\circ} = \{ |x| < \rr ,\, x \in \Rtwo \}^m $. 
Let $ \RRprCs $ be the open kernel of $ \RRprs $: 
\begin{align}
\label{:43n}
\RRprCs = \big\{ \mathbf{x} \in \sS _{\rr }^{m,\circ} \, ;\, \ &
\inf_{j\not=k } |x_j-x_k | > 2^{-\pp } ,\ 
\inf_{l,i} |x_l-s_i| > 2^{-\pp } \big\} 
.\end{align}
For $ \nn = (\pqr ) \in \NNNthree $ we set %let $ \HanC $ be such that 
\begin{align}
\label{:43m}&
\HanC =\Ha _{\pqr } ^{\circ } := \big\{ (\mathbf{x},\mathsf{s}) \in \Ssi ^{[m]} \, ;\, \ 
%\inf_{j\not=k } |x_j-x_k | > 2^{-\pp } ,\ \inf_{l,i} |x_l-s_i| > 2^{-\pp }\\ \notag & 
%\ \mathsf{s} \in \Kakk , \ \mathbf{x} \in \sS _{\rr }^{m,\circ} 
\mathbf{x} \in \RRprCs ,\, \mathsf{s} \in \Kakk 
\big\} 
,\\ \label{:43M}%\notag 
&\Ha ^{\circ }= \bigcup_{\rr =1}^{\infty} \Ha _{\rr }^{\circ } 
, \quad 
 \Ha _{\rr }^{\circ } = \bigcup_{\qq =1}^{\infty} \Ha _{\qq ,\rr }^{\circ }, 
\quad 
\Ha _{\qq ,\rr }^{\circ } = \bigcup_{\pp =1}^{\infty} \Hb ^{\circ } 
.\end{align}
Then $ \HanC \subset \Han $ and $ \HanC \cup \partial \HanC = \Han $. 
Note that $ \HanC $ has a compact closure for each $ \nn \in \NNNthree $. 
The next lemma gives a localization of the $ m $-labeled process. 

\begin{lemma} \label{l:41} 
For each $ m \in \N $ the following holds: 
\begin{align}& \label{:42w}
\PPs ( 
\limi{\nn } 
\tau_{\HanC } 
(\mathbf{X}^m,\mathsf{X}^{m*}) = \infty ) = 1
\quad \text{ for $ \mug $-a.s.\! $ \mathsf{s}$}
.\end{align}
Here $ \tau _A $ denotes the exist time from a set $ A $, 
 $ (\mathbf{X}^m,\mathsf{X}^{m*})$ is the $ m $-labeled process 
given by $ \XB $ in \tref{l:35}, and 
\begin{align}& \label{:42z}
\limi{\nn } = \limi{\rr }\limi{\qq }\limi{\pp }
.\end{align}
\end{lemma}
\begin{proof}
Let 
$ \mathbf{a}^+ =\{ \ak ^+ \}_{\qq \in \mathbb{N}} $. 
Then from \cite[Lemma 2.5 (4)]{o.dfa}, we obtain 
\begin{align}\notag &%\label{:41j}&
\limi{\qq }\mathrm{Cap}^{\mug } ( \Kakk ^c) = 0 .
\end{align}
Then we see for $ \mug $-a.s.\! $ \mathsf{s}$ 
\begin{align}
\label{:41J}&
\PPs ( \limi{\qq }\tau _{\Kakk }(\mathsf{X})=\infty ) = 1 
.\end{align}
%where $ \mathsf{X}=\sum_{i=1}^{\infty} \delta_{X^i}$. 
%is the unlabeled process. % of $ \mathbf{X}$. 
%
By definition we deduce 
$ \tau _{\Kakk }(\mathsf{X})\le \tau _{\Kakk }(\mathsf{X}^{m*})$. 
This combined with \eqref{:41J} yields 
\begin{align}\label{:41E}&
\PPs ( \limi{\qq }\tau _{ \Kakk }(\mathsf{X}^{m*})=\infty ) = 1 
.\end{align}
From \lref{l:32} we have $ \PP _{\mug } (\WSsiNE ) = 1 $. 
Then tagged particles neither collide nor explode. Hence we deduce 
for $ \mug $-a.s.\! $ \mathsf{s}$ 
\begin{align}
\label{:41h}&
\PPs ( \tau _{\Ssi ^{[m]}} (\mathbf{X}^m,\mathsf{X}^{m*}) =\infty )= 1 
,\\
\label{:41f1}&
\PPs \big( \limi{\rr } 
\tau _{\sS _{\rr }^{m,\circ} }(\mathbf{X}^{m})
=\infty \big) = 1 
.\end{align}

From \eqref{:41h} we have for $ \mug $-a.s.\! $ \mathsf{s}$ 
\begin{align}\label{:41H}&
\PPs ( 
\limi{\pp } \tau_{\Ha _{\pp ,\qq ,\rr } ^{\circ }} (\mathbf{X}^m,\mathsf{X}^{m*})
 = 
 \tau_{\Ha _{\qq ,\rr } ^{\circ }} (\mathbf{X}^m,\mathsf{X}^{m*})
 ) = 1
.\end{align}
By \eqref{:41E} we see for $ \mug $-a.s.\! $ \mathsf{s}$ 
\begin{align}\label{:41l}&
\PPs ( 
\limi{\qq } \tau_{\Ha _{\qq ,\rr } ^{\circ }} (\mathbf{X}^m,\mathsf{X}^{m*})
 = 
 \tau_{\Ha _{\rr } ^{\circ }} (\mathbf{X}^m,\mathsf{X}^{m*})
 ) = 1
.\end{align}
From \eqref{:41h} and \eqref{:41f1} we see for $ \mug $-a.s.\! $ \mathsf{s}$ 
\begin{align}\label{:41m}&
\PPs ( 
\limi{\rr } \tau_{\Ha _{\rr } ^{\circ }} (\mathbf{X}^m,\mathsf{X}^{m*}) = \infty 
 ) = 1
\end{align}
Putting \eqref{:41H}, \eqref{:41l}, and \eqref{:41m} together 
we conclude \eqref{:42w}. 
\qed \end{proof}

Let $ b^m =(b^{m,i})_{i=1}^m$ be the drift coefficient of the SDE describing $ \mathbf{X}^m $ and let 
$ \mathbf{B}^m$ be the $ 2m $-dimensional Brownian motion. Then 
\begin{align}
\label{:41n}&
d\mathbf{X}_t^m = d\mathbf{B}_t^m + b^m (\mathbf{X}_t^m,\mathsf{X}_t^{m*}) dt 
\end{align}
and $ b^{m,i}$ is given by %be the $ i $-th component $ X^i $ of $ b^m $. Then 
\begin{align}&\notag %\label{:42x}
b^{m,i}(\mathbf{x},\mathsf{s})=
\half \mathsf{d}^{\mug } (x_i,\sum_{j=1,j\not=i}^m \delta_{x_j} + \mathsf{s})
,\quad \mathbf{x}=(x_1,\ldots,x_m)
.\end{align}
Let $ \fg $ be a version of $ b^m =(b^{m,i})_{i=1}^m$ with respect to $ \mug ^{[m]} $. 
We shall prove in \lref{l:W7} that $ b^{m,i} $ are locally in the domain 
of the $ m $-labeled Dirichlet form, and we shall take a (locally) quasi-continuous version of $ b^m $ later.

Let $ \map{\Pitwo }{\Rtwom \ts \sSS }{\sSS }$ be the projection such that 
$ (x,\mathsf{s}) \mapsto \mathsf{s}$. 
Let $ \{ \IIm \}_{\mm \in \N } $ be an increasing sequence of closed sets in $ \Rtwom \ts \sSS $. 
Then by definition 
\begin{align} \label{:Hmn}
\HmnPc = &\{ \mathsf{s}\in \sSS \, ; \, 
\Hmnc \cap \big(\Rtwom \ts \{ \mathsf{s}\} \big) \not=\emptyset \}
.\end{align}
We set 
\begin{align} \label{:42d}
\HImnC = &\bigcup_{\mathsf{s}\in \Pitwo (\Hmnc )} \RRprCs \ts \{ \mathsf{s} \} 
.\end{align}
For $ \nn =(\pqr ) \in \NNNthree $, let $ \Ct (\mm ,\nn )\label{;GY3a} $ 
be the constant such that $ 0\le \cref{;GY3a} \le \infty $ and that 
\begin{align} \label{:42y} 
\cref{;GY3a} = \sup \{
\frac{|\fg (\mathbf{x},\mathsf{s}) - \fg (\mathbf{y},\mathsf{s}) |}
{|\mathbf{x}-\mathbf{y}|} ; \, & \ 
 \mathbf{x}\not=\mathbf{y}
,\ 
\mathsf{s} \in \HmnPc ,\, 
\\ \notag & 
\ (\mathbf{x},\mathsf{s}) , (\mathbf{y},\mathsf{s}) \in \RRprCs ,\ 
(\mathbf{x},\mathsf{s}) \simpr (\mathbf{y},\mathsf{s}) 
\} 
.\end{align}
Here $ (\mathbf{x},\mathsf{s}) \simpr (\mathbf{y},\mathsf{s}) $ means 
$ \mathbf{x} $ and $ \mathbf{y} $ are in the same connected component of 
$ \RRprCs $. 

\begin{proposition} \label{l:42} 
There exist a $ \mug ^{[m]} $-version $ \fg $ of $ b^m $ and an increasing sequence of closed sets 
$ \{ \IIm \}_{\mm \in \N } $ 
such that for each $ \mm \in \N $ and $ \nn \in \NNNthree $ 
\begin{align}\label{:42a}&
\cref{;GY3a} (\mm ,\nn ) < \infty 
% ,\\ & \label{:42b}
% \limi{\mm }\mathrm{Cap}^{\mug ^{[m]}} (\Han \backslash \IIm ) = 0 
% \quad \text{ for each }
% \nn \in \NNNthree 
,\\& \label{:42c}
\PPs ( 
\limi{\nn } \limi{\mm }
\tau_{\HImnC } ( \mathbf{X}^m,\mathsf{X}^{m*}) = \infty ) = 1 
\quad \text{ for $ \mug $-a.s.\! $ \mathsf{s}$}
.\end{align}
Here $ \tau_{\HImnC }(\mathbf{X}^m,\mathsf{X}^{m*}) $ denotes the exit time 
of $ (\mathbf{X}^m,\mathsf{X}^{m*})$ from the set $\HImnC $ 
and $ (\mathbf{X}^m,\mathsf{X}^{m*})$ is given by $ \XB $ in \tref{l:35} 
starting at $ \lab (\mathsf{s})$. 
\end{proposition}

We shall prove \pref{l:42} in \sref{s:W2}. 

From $ \cref{;GY3a} (\mm ,\nn ) < \infty $ we see 
$ \fg (\mathbf{x},\mathsf{s}) $ restricted on each connected component of 
$ \RRprCs $ is Lipschitz continuous in $ \mathbf{x}$ for each fixed $ \mathsf{s}$ and that the Lipschitz constant is bounded by $ \cref{;GY3a} (\mm ,\nn )$. 
Thus, \pref{l:42} implies a local Lipschitz continuity of the coefficients of the $ m $-labeled SDE \eqref{:35c}. Using this we shall obtain the pathwise uniqueness and the existence of a strong solution of finite-dimensional SDEs. 

The ideas of the proof of \pref{l:42} are twofolds. 
One is the property that $ b^m $ are in the domain of Dirichlet forms, hence we can take a quasi-continuous version of them, which enable us to control the maximal norm with suitable cut off due to $ \IIm $. The second is the Taylor expansion of $ b^m $ using the logarithmic interaction potential. We note here that differential gains the integrability of coefficients at infinity, which is a key point of the proof of \pref{l:42}. 
We refer to \sref{s:9C} for Taylor expansion, and to \sref{s:W2} for a specific calculation in case of the Ginibre interacting Brownian motion.

\Ssection{Proof of local Lipschitz continuity of coefficients: \pref{l:42}}\label{s:W2}

This section proves \pref{l:42} to complete the proof of \tref{l:35}. 
For simplicity we prove only for $ m = 1 $. 
Let $ \NNNthree $ be as \eqref{:42u}. 
Let $ \chin $ ($ \nnNNN_3 $) be the cut-off function defined on $ \R ^2 \ts \sSS $ 
introduced by \eqref{:U2u} with $ m = 1 $. 
 Then by \lref{l:U1} the function $ \chin $ satisfies the following. 
 \begin{align} \label{:W5p} 
 &
 \chin (x,\mathsf{s}) = 
 \begin{cases}
 0 & \text{ for } (x,\mathsf{s}) \not\in \Hann 
 \\
 1& \text{ for } (x,\mathsf{s}) \in \Han 
 \end{cases}
 ,\quad 
 \chin \in \dom _{\mathrm{Gin}} ^{[1]} 
 ,\\ \notag %
 & 
 0 \le \chin (x,\mathsf{s}) \le 1 ,\quad |\nabla_x \chin (x,\mathsf{s}) |^2 \le \cref{;83a} 
 ,\quad 
 \DDD [\chin ,\chin ] (x,\mathsf{s}) \le \cref{;83b}
 .\end{align}
 Here 
 $ \cref{;83a}(\nn )$ and $ \cref{;83b}(\nn )$ are positive constants independent of $ (x,\mathsf{s}) $ 
 in \lref{l:U1}, and $\dom _{\mathrm{Gin}}^{[1]} $ is the domain of the Dirichlet form of 
 the $ 1 $-labeled process $ \XoneX $ given by \eqref{:41c}. 
 Moreover, $ \nabla_x = (\PD{}{x_1},\PD{}{x_2})$ for $ x =(x_1,x_2) \in \Rtwo $. 

We refine the result in \lref{l:33} from $ \LlocpG $ ($ 1\le p < 2$) to $ \LchiG $. 
\begin{lemma} \label{l:W5} 
 $ \dG \in \LchiG $ holds and the convergence \eqref{:33a} and \eqref{:33b} takes place in $ \LchiG $ 
for each $ \nnNNN_3 $. 
\end{lemma}
\begin{proof}
From \cite[Lemma 72]{o.isde}, we deduce the convergence in $ \LlocG $ of the series 
\begin{align}\notag %\label{:33c}
&
 \dG _{1+} (x, \mathsf{s}) := 
 2 \limi{\RR } \sum_{1 \le |x-s_i| < \RR } \frac{x-s_i}{|x-s_i|^2} 
, \\ \notag %\label{:33d}
&
 \dG _{2+} (x, \mathsf{s}) := -2x + 
 2 \limi{\RR } \sum_{1 \le |s_i| < \RR } \frac{x-s_i}{|x-s_i|^2}
.\end{align}
By the definition of $ \chin $, this yields the convergence in $ \LchiG $. 
Because the weight $ \chin $ cuts off the sum around $ x $, 
%By the definition of $ \chin $, 
we easily see that 
\begin{align}\notag %\label{:33e}
&
 \dG _{1-} (x, \mathsf{s}) := 
2 \sum_{|x-s_i| < 1 } \frac{x-s_i}{|x-s_i|^2} \in \LchiG 
, \\ \notag %\label{:33f}
&
 \dG _{2-} (x, \mathsf{s}) := -2x + 
 2 \sum_{ |s_i| < 1 } \frac{x-s_i}{|x-s_i|^2} \in \LchiG 
.\end{align}
As $ \dG = \dG _{1+} + \dG _{1-} = \dG _{2+} + \dG _{2-} $, we conclude \lref{l:W5}. 
\qed \end{proof}

Let $ \vR \in C_0^{\infty}(\Rtwo ) $ be a cut-off function such that 
\begin{align}& \label{:W6u}
\text{$ 0 \le \vR (x) \le 1 , \quad |\nabla \vR (x)| \le 2 , \quad \vR (x) = \vRt (|x|) $ 
\quad for all $ x \in \Rtwo $}
,\end{align}
where $ \vRt \in C_0^{\infty}(\R )$ is such that 
\begin{align}& \label{:W6z}
\vRt (t )=
\begin{cases}1 &\text{for $ |t| \le \RR $, }\\
0&\text{for $ |t| \ge \RR +1 $}
\end{cases} 
\quad \text{ for all $ t \in \R $}
.\end{align}
We set 
\begin{align}\label{:W6a}&
\dG _{\RR } (x,\mathsf{s}) = -2 x + 
 2 \sum_{i} \vR (s_i) \frac{x-s_i}{|x-s_i|^2} 
.\end{align}
We write $ \dG _{\RR } = {}^t (\dG _{{\RR },1},\dG _{{\RR },2})$ and $ \partial_p =\PD{}{x_p}$, 
where $ x= {}^t(x_1,x_2)\in \Rtwo $. 
Then a straightforward calculation shows for $ j,k,l \in \{ 1,2 \} $
\begin{align}\label{:W6y}& 
 \partial_j \dG _{{\RR },k} (x,\mathsf{s}) = -2 \delta_{jk} + 
2 \sum_{i} \vR (s_i) 
 \frac{A_{jk} (x-s_i)}{|x-s_i|^4} 
,\\\label{:W6qq}&
 \partial_j\partial_k \dG _{{\RR },l} (x,\mathsf{s})= 
2 \sum_{i} \vR (s_i) \partial_j \big\{ 
 \frac{A_{kl} (x-s_i)}{|x-s_i|^4} \big\}
,\end{align}
where $ \map{A}{\Rtwo }{\mathbb{R}^4} $ is the $ 2\ts 2$ matrix-valued function defined by 
\begin{align} & \label{:W6b}
A (x) = [A_{ij}(x)]_{i,j=1}^2= 
\begin{pmatrix}
-x_1^2+x_2^2 & \ -2x_1x_2 
\\
-2x_1x_2 & \ x_1^2-x_2^2 
\end{pmatrix}
\quad \text{ for } x = (x_1,x_2) \in \Rtwo 
.\end{align}
%By a straightforward calculation 
We easily see there exist constants 
$ \Ct \label{;W6a}$ and $ \Ct \label{;W6}$ such that for all $ x \in \Rtwo $ 
\begin{align}\label{:W6h}
%\Big|
\frac{|A_{kl} (x )|}{|x |^4} %\Big| 
\le &\frac{\cref{;W6a}}{|x |^{2}}
,\\ \label{:W6p}
\Big| \partial_j \big\{ \frac{A_{kl} (x )}{|x |^4} \big\} \Big| 
\le & \frac{\cref{;W6}}{|x |^{3}}
.\end{align}
We write $ \dG = {}^t (\dG _1,\dG _2)$ similarly as 
$ \dG _{\RR } = {}^t (\dG _{{\RR },1},\dG _{{\RR },2})$. 
\begin{lemma} \label{l:W6} 
For any $ \nn = (\pqr ) \in \NNNthree $ and $ j,k,l \in \{ 1,2 \} $
\begin{align} &\label{:W6f}
 \dG = \limi{{\RR }} \dG _{\RR } \text{ in $ \Lc $}
\\
\label{:W6q}&
\partial_j \dG _{k}= \limi{{\RR }} \partial_j \dG _{{\RR },k} 
 \quad \text{ weakly in $ \Lc $}
,\\\label{:W6qqq}&
\partial_j\partial_k \dG _l (x,\mathsf{s}) = 
\limi{{\RR }} \sum_{|s_i|\le {\RR }} \partial_j \big\{ \frac{A_{kl} (x-s_i)}{|x-s_i|^4} \big\}
\quad \text{ for each $ (x,\mathsf{s}) \in \Han $}
.\end{align}
Here the sum converges absolutely and uniformly in $ \Han $. 
In particular, 
\begin{align}
\label{:W6c}&
\sup \{ |\partial_j\partial_k \dG _l (x,\mathsf{s}) |\, ;\, (x,\mathsf{s}) \in \Han \} < \infty 
\end{align}
and $ \partial_j\partial_k \dG _l (x,\mathsf{s}) $ 
is continuous in $ x $ for each $ \mathsf{s}$ on $ \Han $. 
\end{lemma}
\begin{proof} 
We deduce \eqref{:W6f} from \lref{l:W5}, \eqref{:W6u}, \eqref{:W6z}, and \eqref{:W6a} immediately. 

We next prove \eqref{:W6q}. For this purpose it is enough to show the summation term 
 in \eqref{:W6y} 
is bounded in $ \Lc $ as $ {\RR } \to \infty $ for each $ \nn = (\pqr ) \in \NNNthree $. 
Indeed, we deduce from this that the sequence $ \{\partial_j \dG _{{\RR },k} \}$ is relatively compact 
under the weak convergence in $ \Lc $ and that the limit points are unique by \eqref{:W6f}.

Let $ \rhog ^1 $ and $ \rhogx ^1$ be the one-point correlation functions of 
$ \mug $ and $ \mugx $ with respect to the Lebesgue measure, respectively. Then 
\begin{align}&\label{:W6g}
\rhog ^1(x) = \frac{1}{\pi },\quad 
\rhogx ^1(s) = \frac{1}{\pi }- \frac{1}{\pi } e^{-|x-s|^2 }
.\end{align}
Here $ \rhog ^1(x) = \frac{1}{\pi }$ follows from \eqref{:30a} and \eqref{:30b}, and 
$ \rhogx ^1(s) = \frac{1}{\pi }- \frac{1}{\pi } e^{-|x-s|^2 }$ follows from the formula due to 
Shirai-Takahashi \cite{ST03} such that the determinantal kernel 
$ \mathsf{K}_{\mathrm{Gin},x}$ of Palm measure $ \mugx $ is given by 
\begin{align}\notag &%\label{:W6G}&
\mathsf{K}_{\mathrm{Gin},x}(y,z) = \{\kg (y,y) \kg (z,z) - \kg (y,z)\kg (z,y) \}/\kg (x,x) 
.\end{align}
Let $ H(x,\pp ) =\{ s \in \Rtwo ;|x-s|\ge 2^{-\pp }\}$, where $ \pp \in \N $ and $ x \in \Rtwo $. Then we see 
\begin{align}
\label{:W6i}& 
 \Big| 
E^{\mugx}
 [\langle 1_{H(x,\pp )} 
\vR \frac{A_{kl} (x- \cdot )}{|x-\cdot|^4}, \mathsf{s}\rangle] 
 \Big|
\\ \notag = & 
\Big|
\int _{\Rtwo }1_{H(x,\pp ) }(s) \vR (s) \frac{A_{kl} (x-s)}{|x-s|^4} \rhogx ^1(s) ds 
 \Big|
\\ \notag \le & 
\Big|
\int _{\Rtwo }1_{H(x,\pp ) }(s) \vR (s) \frac{A_{kl} (x-s)}{|x-s|^4} \frac{1}{\pi }ds 
 \Big|
+ 
\int _{\Rtwo }1_{H(x,\pp ) }(s) 
\frac{\cref{;W6a}}{|x-s|^{2}}
\frac{1}{\pi } e^{-|x-s|^2 }
ds 
.\end{align}
Here we used \eqref{:W6g} and \eqref{:W6h} for the last line. 
Note that by \eqref{:W6u} and \eqref{:W6b} 
\begin{align}
\label{:46I}&
\int _{\Rtwo }1_{H(0,\pp ) }(s) \vR (s) \frac{A_{kl} (0-s)}{|0-s|^4} \frac{1}{\pi }ds = 0 
.\end{align}
Then we easily see from \eqref{:46I} 
\begin{align}
\label{:W6ii}&
\Big|
\int _{\Rtwo }1_{H(x,\pp ) }(s) \vR (s) \frac{A_{kl} (x-s)}{|x-s|^4} \frac{1}{\pi }ds 
 \Big| 
\\ \notag =&
\Big|
\int _{\Rtwo }
\Big\{
1_{H(x,\pp ) }(s) \frac{A_{kl} (x-s)}{|x-s|^4} %\frac{1}{\pi }ds 
- 
1_{H(0,\pp ) }(s) \frac{A_{kl} (0-s)}{|0-s|^4} \Big\} \vR (s) \frac{1}{\pi }ds 
 \Big| 
\\ \notag \le &
\int _{\Rtwo }\Big|
1_{H(x,\pp ) }(s) \frac{A_{kl} (x-s)}{|x-s|^4} 
 - 
1_{H(0,\pp ) }(s) \frac{A_{kl} (0-s)}{|0-s|^4} \Big| \frac{1}{\pi }ds 
.\end{align}
We set $ I(x,\pp ) = \{ s \in \Rtwo ;|x-s| < 2^{-\pp }\}$. Then 
$ I(x,\pp )= \Rtwo \backslash H(x,\pp )$. Hence 
\begin{align}\label{:W6i1}&
%\Big|
1_{H(x,\pp ) }(s) \frac{A_{kl} (x-s)}{|x-s|^4} - 
1_{H(0,\pp ) }(s) \frac{A_{kl} (0-s)}{|0-s|^4} %\Big| 
\\ \notag &
=
1_{H(x,\pp ) }(s) 1_{I(0,\pp ) }(s) \frac{A_{kl} (x-s)}{|x-s|^4} -
1_{I(x,\pp ) }(s) 1_{H(0,\pp ) }(s) \frac{A_{kl} (0-s)}{|0-s|^4}
\\ \notag &
+
1_{H(x,\pp ) }(s) 1_{H(0,\pp ) }(s)
 \Big( \frac{A_{kl} (x-s)}{|x-s|^4}- \frac{A_{kl} (0-s)}{|0-s|^4}\Big)
.\end{align}
It is clear that 
\begin{align}\label{:W6i2}&
\sup_{|x| \le \rr +1 } 
\int _{\Rtwo }\Big|
1_{H(x,\pp ) }(s)1_{I(0,\pp ) }(s) \frac{A_{kl} (x-s)}{|x-s|^4} 
 \Big| \frac{1}{\pi } ds < \infty 
,\\ &\notag %\label{:W6i3}&
\sup_{|x| \le \rr +1 } 
\int _{\Rtwo }\Big|
1_{I(x,\pp ) }(s) 1_{H(0,\pp ) }(s)
 \frac{A_{kl} (0-s)}{|0-s|^4} 
 \Big| \frac{1}{\pi } ds < \infty 
,\\ \notag &
\sup_{|x| \le \rr +1 } 
\int _{|s| \le 2(\rr +1 )}
1_{H(x,\pp ) }(s) 1_{H(0,\pp ) }(s)
 \Big|
\frac{A_{kl} (x-s)}{|x-s|^4}- \frac{A_{kl} (0-s)}{|0-s|^4}
 \Big| \frac{1}{\pi } ds < \infty 
.\end{align}
We write $ x =(x_1,x_2) \in \Rtwo $. 
Using \eqref{:W6p}, 
%and taking Taylor expansion at $ 0-s$ in $ x $, 
% e see there exists a constant $ \Ct \label{;W6}=\cref{;W6}(\pp , \rr )$ such that 
we see for all $ |x| \le \rr +1$ and $| s | > 2( \rr + 1 )$ 
\begin{align}\notag %\label{:W6i4}&
\Big|\frac{A_{kl} (x-s)}{|x-s|^4} - \frac{A_{kl} (0-s)}{|0-s|^4}\Big| & = 
\Big|
\int _0^1 \sum_{j=1}^2 x_j 
\partial_j \big\{ \frac{A_{kl} (\cdot )}{|\cdot |^4} \big\} (tx -s) dt 
\Big|
\\ \notag & \le 
\int _0^1 \sum_{j=1}^2 |x_j| 
\sup_{|x| \le \rr +1 } \big\{ 
\frac{\cref{;W6}}{|tx -s |^3 } \big\} dt 
\quad \text{ by \eqref{:W6p}}
\\ \notag & \le 
2(\rr + 1 ) \big\{ 
\frac{\cref{;W6}2^3}{|s |^3 } \big\} 
%\quad \text{ by $ |x| \le \rr +1$, $| s | > 2( \rr + 1 )$} 
.\end{align}
Here we used $ | s | /2 < |tx -s |$ in the last line. This follows from 
$ 0 \le t \le 1 $, $ |x| \le \rr +1$, and $ | s | > 2( \rr + 1 ) $. 
Note that $ 1_{H(x,\pp ) }(s) 1_{H(0,\pp ) }(s) =1$ for $ |x| \le \rr +1$ and $| s | > 2( \rr + 1 )$. 
Hence we deduce 
\begin{align}\label{:W6i5}&
\sup_{|x| \le \rr +1 } 
\int _{ | s | > 2( \rr + 1 )}
1_{H(x,\pp ) }(s) 1_{H(0,\pp ) }(s)
 \Big|
\frac{A_{kl} (x-s)}{|x-s|^4}- \frac{A_{kl} (0-s)}{|0-s|^4}
 \Big| \frac{1}{\pi } ds 
\\ \notag &\quad \quad \quad \quad 
\le \cref{;W6}2^4 (\rr +1)\int _{ | s | > 2( \rr + 1 )} \frac{1}{|s|^3}
 \frac{1}{\pi }ds
 < \infty 
.\end{align}

Collecting \eqref{:W6ii}--\eqref{:W6i5} together, we obtain 
\begin{align}
\label{:W6I}&
\limsupi{{\RR }} \sup_{|x| \le \rr +1 } 
\Big|
\int _{\Rtwo }1_{H(x,\pp ) }(s) \vR (s) \frac{A_{kl} (x-s)}{|x-s|^4} \frac{1}{\pi }ds 
 \Big| 
< \infty 
.\end{align}
Clearly, we have 
\begin{align}
\label{:W6III}& 
\int _{\Rtwo }1_{H(x,\pp ) }(s) 
\frac{\cref{;W6a}}{|x-s|^{2}}
\frac{1}{\pi } e^{-|x-s|^2 }
ds < \infty 
.\end{align}
Putting \eqref{:W6I} and \eqref{:W6III} into \eqref{:W6i}, we obtain 
\begin{align}
\label{:W6II}&
\limsupi{{\RR }} \sup_{|x| \le \rr +1 } 
E^{\mugx}
 [\langle 1_{H(x,\pp )} 
\vR \frac{A_{kl} (x- \cdot )}{|x-\cdot|^4}, \mathsf{s}\rangle] 
< \infty 
.\end{align}
From the inequality $ \mathrm{Var}[\langle f,\mathsf{s}\rangle ] \le \int |f|^2 \rho^1(s)ds $ 
valid for determinantal random point fields with Hermitian symmetric kernel, we have 
\begin{align}
\label{:W6j}&
\limsupi{{\RR }} \sup_{|x| \le \rr +1 } 
\mathrm{Var}^{\mugx}
 [\langle 1_{H(x,\pp )} 
\vR \frac{A_{kl} (x- \cdot )}{|x-\cdot |^4}, \mathsf{s}\rangle] 
\\ \notag \le &
\limsupi{{\RR }} \sup_{|x| \le \rr +1 } 
\int _{\Rtwo }\Big| 1_{H(x,\pp )} (s)\vR (s) \frac{A_{kl} (x-s)}{|x-s|^4} \Big|^2\rhogx ^1(s)ds 
\\ \notag < & \infty 
\quad \text{by \eqref{:W6h} and \eqref{:W6g}}
.\end{align}

Putting \eqref{:W6II} and \eqref{:W6j} together we immediately deduce 
\begin{align}\label{:W6k}&
\limsupi{{\RR }} \sup_{|x| \le \rr +1 } 
E^{\mugx}
 [ \Big|\langle 1_{H(x,\pp )} 
\vR \frac{A_{kl} (x- \cdot )}{|x-\cdot|^4}, \mathsf{s}\rangle \Big|^2 ] 
< \infty 
.\end{align}
From \eqref{:W5p}, \eqref{:W6g}, \eqref{:W6k} and recalling $ \nn = (\pqr ) \in \NNNthree $ 
we easily obtain 
\begin{align}
 \label{:W6m}&
\limsupi{{\RR }}\int 
 \Big| 2 \sum_{i} \vR (s_i) \frac{A_{kl} (x-s_i)}{|x-s_i|^4} \Big|^2 
\cmu < \infty 
.\end{align}
Then from \eqref{:W6y} and \eqref{:W6m} with a simple calculation we see
\begin{align}
\notag &% \label{:W6m}&
\limsupi{{\RR }}\int | \partial_j \dG _{{\RR },k} |^2 \cmu < \infty 
.\end{align}
Hence $ \{\partial_j \dG _{{\RR },k} \}$ is relatively compact in the weak topology in $ \Lc $ 
for each $ \nn =(\pqr )\in \NNNthree $. This and \lref{l:W5} yield \eqref{:W6q}.

We next prove \eqref{:W6qqq} and \eqref{:W6c}. We see from \eqref{:W6p} 
\begin{align} \label{:W6pP}& 
\sum_i 
\Big| \partial_j \big\{ \frac{A_{kl} (x-s_i)}{|x-s_i|^4} \big\} \Big| 
\le \cref{;W6}\sum_i 
\frac{1}{|x-s_i|^{3}}
.\end{align}
We have for any $ (x,\mathsf{s}) \in \Han $
\begin{align}\label{:W6o} 
\sum_{i }\frac{1}{|x-s_i|^{3}} 
= & 
\sum_{|s_i| < \rr +1}\frac{1}{|x-s_i|^{3}}
 +
\sum_{ \ttt =\rr +2}^{\infty} \, \sum_{\ttt - 1 \le |s_i| < \ttt } 
\frac{1}{|x-s_i|^{3}}
\\ \le & \notag 
 2 ^{3\pp }\mathsf{s}(\sS _{\rr +1} )
 +
 \sum_{ \ttt =\rr +2 }^{\infty}\, \sum_{\ttt - 1 \le |s_i| < \ttt }
 \frac{1 }{ ( |s_i| - \rr )^{ 3 }}
.\end{align}
Here we used $ |x-s_i|\ge 2^{-\pp } $ and $ |x| \le \rr $. By a straightforward calculation we have 
\begin{align}
\label{:W6o1}&
 \sum_{ \ttt =\rr +2 }^{\infty}\, \sum_{\ttt - 1 \le |s_i| < \ttt }
 \frac{1 }{ ( |s_i| - \rr )^{ 3 }}
\\ \notag 
\le & 
\sum_{ \ttt =\rr +2 }^{\infty} 
 \frac{\mathsf{s}(\sS _{\ttt } ) - \mathsf{s} (\sS _{\ttt - 1 })}
{ (\ttt - 1 -\rr )^{3 }} 
 \\ \notag \le &
 \lim_{R\to\infty } \Big\{ 
 \frac{\mathsf{s} (\sS _R ) }{(R -1-\rr )^{3}}
+ 
 \sum_{ \ttt =\rr +3}^{R} \mathsf{s}(\sS _{\ttt - 1 })
\Big\{ \frac{1}{ ( \ttt -2-\rr )^{3}}- \frac{ 1 }{ ( \ttt -1-\rr )^{3}} 
\Big\} \Big\} 
.\end{align}
Recall that $ \mathsf{s} (\sS _{\ttt } ) \le \ak ( \ttt ) $ by $ (x,\mathsf{s}) \in \Han $. 
Then \eqref{:W6pP}--\eqref{:W6o1} yield 
\begin{align} \label{:W6o2} &
\sum_i \Big| \partial_j \big\{ \frac{A_{kl} (x-s_i)}{|x-s_i|^4} \big\} \Big| 
%&\sum_{i }\frac{1}{|x-s_i|^{3}} 
\le \cref{;W6}2 ^{3\pp }\ak (\rr +1) 
 \\ \notag
&\quad + \cref{;W6}
\lim_{R\to\infty } \Big\{ 
 \frac{\ak (R)}{(R -1-\rr )^{3}} %{R^{3}}
+ \
 \sum_{ \ttt =\rr +3}^{R} \ak (\ttt - 1 ) 
\Big\{ \frac{1}{ ( \ttt -2-\rr )^{3}}- \frac{ 1 }{ ( \ttt -1-\rr )^{3}} \Big\} \Big\} 
.\end{align}
Because $ \ak (\ttt ) = k \ttt ^2 $, the sum in \eqref{:W6o2} converges in $ \Han $ uniformly. 
Hence, we obtain \eqref{:W6qqq} and \eqref{:W6c} 
from \eqref{:W6u}, \eqref{:W6qq}, and \eqref{:W6o2} immediately. 
The last claim follows from \eqref{:W6qqq} and 
the uniform convergence of the series in \eqref{:W6qqq}. 
\qed \end{proof}

Recall that $ \sigma (x,\mathsf{s}) = E $ and $ b (x,\mathsf{s}) = \frac{1}{2} \dG (x,\mathsf{s}) $. 
Let $ \dom _{\mathrm{Gin}} ^{[1]} $ be the domain of the 1-labeled Dirichlet form of the Ginibre interacting Brownian motion. 
\begin{lemma} \label{l:W7} 
$ \chin \dG ,\, \chin \partial_j \dG ,\, \chin \partial_j\partial_k \dG \in \dom _{\mathrm{Gin}} ^{[1]} $ 
$ (j,k\in\{ 1,2 \})$ for all $ \nn \in \NNNthree $. 
\end{lemma}
\begin{proof}
We only prove $ \chin \dG \in \dom _{\mathrm{Gin}} ^{[1]} $ 
because the other cases can be proved in a similar fashion. 
We set $ \DDD [f] = \DDD [f,f]$. 
Recall that $ \nabla_x = (\PD{}{x_1},\PD{}{x_2})= (\partial_1,\partial_2)$. 
Then 
\begin{align}\notag &%\label{:W7a}& 
 \E ^{\muone _{\mathrm{Gin}}} 
 (\chin \dG ,\chin \dG ) = 
 \int_{\Rtwo \ts \sSS } \frac{1}{2}|\nabla_x ( \chin \dG ) |^2 + 
 \, \DDD [\chin \dG ] 
d \muone _{\mathrm{Gin}}
.\end{align}
From \eqref{:W5p}, \lref{l:W5}, and \lref{l:W6}, we deduce that 
\begin{align} \notag %\label{:W7b}
& 
 \int_{\Rtwo \ts \sSS } |\nabla_x (\chin \dG ) |^2 d \muone _{\mathrm{Gin}}
\le 2\int_{\Hann } \{ \chin ^2 |\nabla_x \dG |^2 + 
|\nabla_x \chin |^2 |\dG |^2 \}
d \muone _{\mathrm{Gin}} 
< \infty 
.\end{align}
We set $ \PD{}{s_{j}} =(\PD{}{s_{j,k}} )_{k=1,2}$. 
From \eqref{:W6a} and $ \dG _{\RR } = {}^t (\dG _{{\RR },1},\dG _{{\RR },2})$ we have 
\begin{align}
&\notag %\label{:}&
 \PD{}{s_{j,k}} \dG _{{\RR },l} (x,\mathsf{s}) = 
2 (\partial _k \vR ) (s_j)\frac{(x-s_j)_l}{|x-s_j|^2} - 
2 \vR (s_j) 
 \frac{A_{kl} (x-s_j)}{|x-s_j|^4} 
.\end{align}
Then from \eqref{:W5p} and \eqref{:W6h}, we similarly deduce that 
\begin{align}\notag %\label{:W7c} 
 \int_{\Rtwo \ts \sSS } \, \DDD & [\chin \dG ] 
d \muone _{\mathrm{Gin}} 
 \le \, 2 
 \int_{\Rtwo \ts \sSS } \, \big\{
 \chin ^2 \DDD [ \dG ] + \DDD [\chin ] \, |\dG |^2 \big\} \, 
 d \muone _{\mathrm{Gin}} 
\\\notag 
 \le \, & 
 2 \int_{\Hann } \, \big\{ 
 \DDD [ \dG ] + \cref{;83b}^2 |\dG |^2 \big\} 
 \, 
 d \muone _{\mathrm{Gin}} 
 \\ \notag 
\le \, & 2
 \int_{ \Hann } \big\{ \cref{;X4} \big(
 \sum_{i} \frac{1}{|x-s_i |^4} \big) 
 + \cref{;83b}^2 |\dG |^2 \big\} 
 \, 
 d \muone _{\mathrm{Gin}} 
 < \, \infty 
.\end{align}
Here $ \Ct \label{;X4}$ is a finite, positive constant and 
$ \cref{;83b}(\nn )$ is the positive constant in \lref{l:U1}. 
The finiteness of the integral in the last line 
follows from the translation invariance of $ \mug $ and \lref{l:W5}. 

Combining these, 
 we obtain $ \E ^{\muone _{\mathrm{Gin}}} (\chin \dG ,\chin \dG ) < \infty $. 
We proved $ \chin \dG \in L^2 (\muone _{\mathrm{Gin}} ) $ in \lref{l:W5}. 
Hence, we see $ \chin \dG \in \dom _{\mathrm{Gin}} ^{[1]} $. 
This completes the proof. 
\qed \end{proof}

\medskip

\noindent 
{\em Proof of \pref{l:42}. }
For simplicity we prove the case $ m =1 $; the general case follows from the same argument. %

By \lref{l:W6} we can take a version of 
$ \partial_j\partial_k \dG _l (x,\mathsf{s}) $ 
which is continuous in $ \mathbf{x}$ for each $ \mathsf{s}$ on 
$ \Han $. We always take this version and denote it by the same symbol $ \partial_j\partial_k \dG _l $. 
By \eqref{:W6c} we have a finite constant $ \Ct (\nn ) \label{;!2c}$ ($ \nn \in \NNNthree $) such that 
\begin{align}& \label{:!2M}
\cref{;!2c}=\sup 
\{ |\partial_j\partial_k \dG _l (x,\mathsf{s}) |\, ;\, (x,\mathsf{s}) \in \Han , \ 
j,k,l \in \{ 1,2 \} \} < \infty 
.\end{align}

We denote by $ \widetilde{f } $ a quasi-continuous version of $ f \in \dom _{\mathrm{Gin}} ^{[1]} $. 
From \lref{l:W7} we take a quasi-continuous version %$ \widetilde{(\chin \dG _l ) } $ 
of $ \chin \dG _l $ 
being commutative with $ \partial_k $: 
\begin{align} \label{:!2h}&
\partial_k \widetilde{(\chin \dG _l )} = \widetilde{(\partial_k \chin \dG _l ) } 
.\end{align}
Then by definition there exists an increasing sequence 
$ \{ \IIm \}_{\mm =1}^{\infty} $ of closed sets such that 
$ \widetilde{\chin \dG _l } $ and $ \widetilde{(\partial_k \chin \dG _l )} $
are continuous on $ \IIm $ for each $ \mm $ and that 
\begin{align*}&
\text{$\limi{\mm }\mathrm{Cap}^{\mugone } (\Hann \backslash \IIm ) = 0 $
 for each $\nn \in \NNNthree $}
.\end{align*}
We used here $ \chin = 0 $ on $ \Hann ^c $ by \eqref{:W5p}. 
Because $ \Han \subset \Hann $, we have 
\begin{align}& \label{:!2g}
\limi{\mm }\mathrm{Cap}^{\mugone } (\Han \backslash \IIm ) = 0 
\quad \text{ for each }
\nn \in \NNNthree 
.\end{align}
Note that $ \chin = 1$ on $ \Han $ by \eqref{:W5p}, 
$ \Han \uparrow \Ha $, and $ \mugone (\Ha ^c) = 0$. 
We set 
\begin{align*}&
\text{ 
$ \widetilde{\dG _l } (x,\mathsf{s})
=\widetilde{\chin \dG _l } (x,\mathsf{s}) $\quad for $ (x,\mathsf{s}) \in \Han $}
.\end{align*}
Thus $ \widetilde{\dG _l }$ is a version of $ \dG _l $ such that 
$ \widetilde{\dG _l }$ and $ \partial_k \widetilde{\dG _l } $ 
are continuous on $ \Han \cap \IIm $ for each $ \nn \in \NNNthree $. 
Let $ \nn \in \NNNthree $ and $ \mm \in \N $. 
Let $ \Ct (\nn ,\mm )\label{;!2K}$ be the constant such that 
\begin{align}%\notag &%
\label{:!2K}&
\cref{;!2K} = \sup\big\{ 
|\widetilde{\dG _l } (\xi ,\mathsf{s}) |,\, 
|\partial_k \widetilde{\dG _l } (\xi ,\mathsf{s}) | 
\, ; (\xi ,\mathsf{s}) \in \Hmn , \ k,l \in \{ 1,2 \} \big\} 
.\end{align}
Note that $ \Han $ is a compact set because $ \nn \in \NNNthree $. Recall that $ \IIm $ is a closed set. 
Then $ \Hmn $ is compact. 
Because $ \widetilde{\dG _l } $ and $ \partial_k \widetilde{\dG _l } $ 
are continuous on $ \Hmn $, these are bounded on $ \Hmn $. 
Thus we have $ \cref{;!2K} < \infty $. 

We suppose 
\begin{align} \label{:!2s}&
\mathsf{s} \in \HmnP 
.\end{align}
We write $ \nn = (\pqr )$ as before. Let $ \simpr $ be same as \eqref{:42y}. Fix $ \mm $ and $ \nn $. 
Then there exists a positive constant $ \Ct \label{;!2d} < \infty $ depending only on $ \mm $ and $ \nn $ 
such that the following holds: 
For any $ ( x ,\mathsf{s}) ,( \xi ,\mathsf{s}) \in \Han $ such that $ ( x ,\mathsf{s}) \simpr ( \xi ,\mathsf{s}) $, 
there exist $ \{ x_1,\ldots,x_{\q } \} $ with $ x_1 = \xi $ and $x_{\q } =x $ such that \eqref{:!2B} and \eqref{:!2d} hold. 
\begin{align}&\label{:!2B}
\text{$ [x_j,x_{j+1}]\ts \{ \mathsf{s} \} \subset \Hann $\quad for all $1\le j < \q $}
,\\
\intertext{where $ [x_j,x_{j+1}] \subset \Rtwo $ denotes 
the segment connecting $ x_j $ and $ x_{j+1} $, 
}
\label{:!2d}&
\sum_{j=1}^{\q -1} |x_j - x_{j+1}| \le \cref{;!2d}|\xi -x | 
.\end{align}
Using \eqref{:!2h} and Taylor expansion we deduce that for each $ 1 \le j < \q $ 
\begin{align}
\label{:!2i}&
\partial_k \widetilde{\dG _l }
 (x_{j+1},\mathsf{s}) - \partial_k \widetilde{\dG _l } (x_{j} ,\mathsf{s})
\\ \notag &
=\int_{0}^{1} \1 dt
.\end{align}
Here we set $ x_j = (x_{j,1},x_{j,2}) \in \Rtwo $. 
Taking the sum of both sides of \eqref{:!2i} over $ j=1,\ldots,\q -1 $ 
and recalling $x_{\q } =x $ and $ x_1 = \xi $ we see 
\begin{align}& \notag %\label{:!2j}&
\partial_k \widetilde{\dG _l } (x,\mathsf{s}) - \partial_k \widetilde{\dG _l } (\xi ,\mathsf{s})
%\\ \notag &
 = 
\sum_{j=1}^{\q -1}
\int_{0}^{1} 
\1 dt 
.\end{align}
This together with \eqref{:!2B} and \eqref{:!2d} yields 
\begin{align}
\label{:!2l}&
|\partial_k \widetilde{\dG _l } (x,\mathsf{s}) - \partial_k \widetilde{\dG _l } (\xi ,\mathsf{s}) | 
\\ \notag &
\le 
\sum_{j=1}^{\q -1}\int_0^1 | \1 | dt
\\ \notag &
\le \sqrt{2} \cref{;!2c} \cref{;!2d} |x-\xi |
\\ \notag &
\le \sqrt{2} \cref{;!2c} \cref{;!2d} 2\rr \quad \text{ by $ x,\xi \in \Srr $}
.\end{align}
Hence for each $ \mathsf{s}$ satisfying \eqref{:!2s}, 
$ \partial_k \widetilde{\dG _l } (x,\mathsf{s}) $ is bounded in $ x $ with bound 
\begin{align} \label{:!2r}&
|\partial_k \widetilde{\dG _l } (x,\mathsf{s})| \le 
|\partial_k \widetilde{\dG _l } (\xi ,\mathsf{s}) | 
+ \sqrt{2} \cref{;!2c} \cref{;!2d} 2\rr 
.\end{align}
Because of \eqref{:!2s}, we can take $ (\xi ,\mathsf{s}) \in \Hmn $. 
Then by \eqref{:!2K} 
\begin{align}\label{:!2rr}&
|\partial_k \widetilde{\dG _l } (\xi ,\mathsf{s}) | \le \cref{;!2K}
.\end{align}
From \eqref{:!2r} and \eqref{:!2rr} we obtain 
\begin{align}\notag &%\label{:2!R}&
|\partial_k \widetilde{\dG _l } (x,\mathsf{s})| 
\le 
\cref{;!2K}+ \sqrt{2} \cref{;!2c} \cref{;!2d} 2\rr =:\cref{;!2m} 
\end{align}
The bound $ \Ct \label{;!2m}$ 
depends only on $ \mm , \nn $ appearing in \eqref{:!2s}. Then 
\begin{align}
\label{:!2n}&
\sup \{ |\partial_k \widetilde{\dG _l } (x,\mathsf{s}) |\,; \, 
(x,\mathsf{s}) \in \Han ,\, \mathsf{s} \in \HmnP 
 \} 
\le \cref{;!2m}
.\end{align}

Using \eqref{:!2h} and Taylor expansion again we deduce that for each $ 1 \le j < \q $ 
\begin{align}\notag %\label{:!2e} &
\widetilde{\dG _l } (x_{j+1},\mathsf{s}) - \widetilde{\dG _l } (x_{j} ,\mathsf{s}) =& 
\int_0^1 \4 dt 
.\end{align}
Then from this and \eqref{:!2n} we deduce in a similar fashion as \eqref{:!2l} 
\begin{align}
\label{:!2o}
|\widetilde{\dG _l } (x,\mathsf{s}) - \widetilde{\dG _l } (\xi ,\mathsf{s}) |\le 
& \sqrt{2} \cref{;!2d}\cref{;!2m} |x-\xi | 
.\end{align}
Note that \eqref{:!2o} holds for all 
$ (x,\mathsf{s}) \simpr (\xi ,\mathsf{s}) \in \Han $ with $ \mathsf{s}$ satisfying \eqref{:!2s}. 
Here $ (x ,\mathsf{s}) $ is not necessary in $ \Hmn $, whereas $ (\xi ,\mathsf{s}) \in \Hmn $. 
The constant $\sqrt{2} \cref{;!2d}\cref{;!2m} $ depends only on $ \mm $ and $\nn $. 
Hence we deduce $ \cref{;GY3a} (\mm ,\nn ) < \infty $. 

By definition $ \HanC \supset \Ha _{\nn -1} $. 
Then we easily deduce from \eqref{:Hmn} and \eqref{:42d} 
\begin{align}
\label{:!2p}&
\LLL \HanC \cap \IIm \rRRR \supset 
\LLL \Ha _{\nn -1} \cap \IIm \rRRR 
\supset \Ha _{\nn -1} \cap \IIm 
.\end{align}
Let $ \XoneX =(X^1,\sum_{i=2}^{\infty} \delta_{X^i})$ be the one-labeled process. 
By \eqref{:!2g} we see %for each $ \nn \in \NNNthree $
\begin{align}&\notag % \label{:!2q}&
\PPs \big( 
\tau_{\Han } \XoneX = 
\limi{\mm }\tau_{\Han \cap \IIm } \XoneX \big) = 1
\end{align}
for all $ \nn \in \NNNthree $. Then this yields 
\begin{align}
\label{:!2q}&
\PPs \big( 
\limi{\nn }\tau_{\Han } \XoneX = \limi{\nn }
\limi{\mm }\tau_{\Han \cap \IIm } \XoneX \big) = 1
.\end{align}
Combining \eqref{:!2p}, \eqref{:!2q}, and \lref{l:41}, we obtain 
\begin{align}\notag %\label{:!2r}
\PPs \big( 
\limi{\nn } \limi{\mm }
\tau_{\HImnC } \XoneX = \infty \big) 
&
\ge \PPs \big( 
\limi{\nn } \limi{\mm }
\tau_{\Ha _{\nn -1} \cap \IIm } \XoneX = \infty \big) 
\\ \notag &
= \PPs \big( 
\limi{\nn } 
\tau_{\Ha _{\nn -1} } \XoneX = \infty \big) 
\\ \notag &
\ge \PPs \big( 
\limi{\nn } 
\tau_{\Ha _{\nn -1} ^{\circ}} \XoneX = \infty \big) 
\\ \notag &
= 1 
.\end{align}
Indeed, the first line follows from \eqref{:!2p}. The second line follows from \eqref{:!2q}. 
The inequality in the third line is clear by $ \Ha _{\nn -1} \supset \Ha _{\nn -1} ^{\circ } $. 
The last line is immediate from \lref{l:41}. We have thus proved \eqref{:42c}. 
\bbbbb

\Ssection{A unique strong solution of SDEs with random environment. } \label{s:4C}
Throughout \sref{s:4C}, $ \XB $ is the weak solution of \eqref{:11h} 
starting at $ \mathbf{s}=\lab (\mathsf{s})$ defined on $ \OFPFs $ obtained in \lref{l:34}. 
We write $ \mathbf{X}=(X^i)_{i=1}^{\infty}$ and 
denote by $ \mathbf{B}^m $ the first $ m$-components of the Brownian motion 
$ \mathbf{B}=(B^i)_{i=1}^{\infty}$. 

For each $ m \in \N $, we introduce the $ 2m $-dimensional SDE of $ \mathbf{Y}^m $ 
for $ \XB $ under $ \PPs $ 
such that $ \mathbf{Y}^m =(Y^{m,i})_{i=1}^m $ is 
 an $ \{ \mathcal{F}_t \} $-adapted, continuous process satisfying 
\begin{align}\label{:35c}&
dY_t^{m,i} = dB_t^i - Y_t^{m,i} dt + \sum_{j=1 \atop j\not=i}^m 
\frac{Y_t^{m,i}-Y_t^{m,j}}{|Y_t^{m,i}-Y_t^{m,j}|^{2}} dt + 
\lim_{r\to\infty } \sum_{j = m +1 ,\atop |X_t^j|<r } ^{\infty}
\frac{Y_t^{m,i}-X_t^j}{|Y_t^{m,i}-X_t^j|^{2}} dt 
,\\\label{:35d}& 
 \PPs 
(\mathbf{Y}_t^m \in \SSSsdemtg \, \text{ for all } t ) = 1 
,\\\label{:35e}& 
\mathbf{Y}_0^m = \labm (\mathsf{s}) 
.\end{align}
Here $ \SSSsdemtw $ is the set given by \eqref{:35a} with replacement of $ \sS $ by $ \Rtwo $. 
Furthermore, we set $ \labm (\mathsf{s})=(\lab _i(\mathsf{s}))_{i=1}^m$
for $ \lab (\mathsf{s})= (\labi (\mathsf{s}))_{i=1}^{\infty}$. 
For convenience we introduce a variant of notation of solution. 
Let $ \mathsf{X}^{m*}= \sum_{i=m+1}^{\infty} \delta_{X^i}$ as before. 
We say $ (\mathbf{Y}^m , \mathbf{B}^m , \mathsf{X}^{m*})$ is a solution of \eqref{:35c} 
if and only if 
$ (\mathbf{Y}^m , \mathbf{B}^m , \mathbf{X}^{m*})$ is a solution of \eqref{:35c}.

\begin{lemma}	\label{l:97G} 
\thetag{1} 
For each $ m \in \N $, 
the SDE \ginSDE\ has a pathwise unique, weak solution for $ \mug $-a.s.\,$ \mathsf{s}$. 
That is, for $ \mug $-a.s.\,$ \mathsf{s}$, arbitrary solutions 
$ (\mathbf{Y}^m , \mathbf{B}^m , \mathsf{X}^{m*})$ and 
$ (\hat{\mathbf{Y}}^m , \mathbf{B}^m , \mathsf{X}^{m*})$ 
of \ginSDE\ defined on $ \OFPFs $ satisfy 
\begin{align}
\label{:97oG}&
\PPs (\mathbf{Y}^m = \hat{\mathbf{Y}}^m) = 1 
.\end{align}
In particular, $ (\mathbf{Y}^m , \mathbf{B}^m , \mathsf{X}^{m*})$ coincides with 
$ (\mathbf{X}^m , \mathbf{B}^m , \mathsf{X}^{m*})$ under $ \PPs $. 
\\\thetag{2}
Let 
$ (\mathbf{Z}^m , \hat{\mathbf{B}}^m , \hat{\mathsf{X}}^{m*})$ and 
$ (\hat{\mathbf{Z}}^m , \hat{\mathbf{B}}^m , \hat{\mathsf{X}}^{m*})$ 
be weak solutions of the SDE \eqref{:35c}--\eqref{:35e} 
defined on a filtered space 
$(\Omega ',\mathcal{F}', P' ,\{ \mathcal{F}_t '\} )$ satisfying 
\begin{align} \notag &% \label{:97O}&
 (\hat{\mathbf{Z}}^m , \hat{\mathbf{B}}^m , \hat{\mathsf{X}}^{m*})
\elaw (\mathbf{X}^m , \mathbf{B}^m , \mathsf{X}^{m*})
.\end{align}
Then, 
for $ \mathbf{Z}^m $ and $\hat{\mathbf{Z}}^m$ satisfying $ \mathbf{Z}_0^m=\hat{\mathbf{Z}}_0^m$ a.s., 
 it holds that for $ \mug $-a.s.\,$ \mathsf{s}$ 
\begin{align}\label{:97G1}&
 P' (\mathbf{Z}^m = \hat{\mathbf{Z}}^m) = 1 
.\end{align}
\thetag{3} Make the same assumptions as \thetag{2} except that 
the filtrations of 
$ (\mathbf{Z}^m , \hat{\mathbf{B}}^m , \hat{\mathsf{X}}^{m*})$ and 
$ (\hat{\mathbf{Z}}^m , \hat{\mathbf{B}}^m , \hat{\mathsf{X}}^{m*})$, 
$ \{\mathcal{F}_t'\} $ and $ \{ \mathcal{F}_t'' \} $ say, are different 
(but on the same probability space $(\Omega ',\mathcal{F}', P' )$). 
Then \eqref{:97G1} holds for $ \mug $-a.s.\,$\mathsf{s}$. 
\end{lemma}

\begin{remark}\label{r:97G} 
In conventional situations, 
the pathwise uniqueness in \lref{l:97G} \thetag{3} is called 
the pathwise uniqueness in the strict sense \cite[Remark 1.3, 162p]{IW}.
\end{remark}

\begin{proof}
We first prove \thetag{1}. 
Recall that $ (\mathbf{X}^m , \mathbf{B}^m , \mathsf{X}^{m*})$ under $ \PPs $ 
is a weak solution of \ginSDE. 
So it only remains to prove \eqref{:97oG} for $ \hat{\mathbf{Y}}^m = \mathbf{X}^m $. 

Let \nNpqk. Let $ \RRprCs $ be as \eqref{:43n}. 

For $ \UV \in \WT (\Rtm \times \sSS )$ and 
$ (\mathbf{u},\mathsf{v},\mathbf{w}) \in \WT (\Rtm \times \sSS \ts \Rtm )$, let 
\begin{align}\label{:43o}&
\thetaone \UV = 
\inf \{ t > 0 \, ;\, \mathbf{u}_t \not\in \RRprC (\mathsf{v}_t) \} 
,\\ \label{:43r}&
\thetatwo (\mathbf{u},\mathsf{v},\mathbf{w}) = 
\inf \{ t>0 ; 
(\mathbf{u}_t,\mathsf{v}_t) 
\not\simpr 
(\mathbf{w}_t,\mathsf{v}_t) \} 
.\end{align}
Here $ (\mathbf{x},\mathsf{s}) \simpr (\mathbf{x}',\mathsf{s}) $ is same as \eqref{:42y}. 
$ (\mathbf{x},\mathsf{s}) \not\simpr (\mathbf{y},\mathsf{s}) $ means 
$ \mathbf{x} $ and $ \mathbf{y} $ are not in the same connected component of $ \RRprCs $. 
Let $ \IIm $ be as \pref{l:42}. Let $ \thetathree (\mathsf{v})$ 
be the exit time of $ \mathsf{v}$ from $\Pitwo (\Hmnc ) $: 
\begin{align}\notag &%\label{:43p}&
\thetathree (\mathsf{v}) = \inf \{ t > 0 \, ;\, \mathsf{v}_t \not\in \Pitwo (\Hmnc ) \} 
.\end{align}
Then we easily see for $ \nn = (\pqr ) \in \NNNthree $ and $ \mm \in \N $ 
\begin{align}
\label{:43Q}&
\min \{ \thetaone (\mathbf{u} ,\mathsf{v} ), \thetathree (\mathsf{v} ) \}
\ge 
\tau_{\HImnC } (\mathbf{u} ,\mathsf{v} ) 
.\end{align}
For $ \nn = (\pqr ) \in \NNNthree $ and $ \mm \in \N $ we set 
\begin{align}\notag &%\label{:43P}&
\vartheta_{\mm ,\nn }(\mathbf{u},\mathsf{v},\mathbf{w})= 
\min \{ 
\thetaone (\mathbf{u} ,\mathsf{v} ) , 
\thetaone (\mathbf{w} ,\mathsf{v} ) , 
\thetatwo (\mathbf{u} ,\mathsf{v} ,\mathbf{w} ), 
\thetathree (\mathsf{v} )
\}
.\end{align}

Let 
\begin{align}&\notag %\label{:43t}&
\tau (\mm ,\nn)= \vartheta_{\mm ,\nn } (\mathbf{X}^m,\mathsf{X}^{m*},\mathbf{Y}^m)
% \min \{ 
% \thetaone (\mathbf{X}^m,\mathsf{X}^{m*}) , 
% \thetaone (\mathbf{Y}^m,\mathsf{X}^{m*}) , 
% \thetatwo (\mathbf{X}^m,\mathsf{X}^{m*},\mathbf{Y}^m), 
% \thetathree (\mathsf{X}^{m*})
% \} 
.\end{align}
Then from \eqref{:41n} and \eqref{:42y} and a straightforward calculation we see 
\begin{align}\notag %\label{:43jj}
| \mathbf{X}_{t\wedge \tau (\mm ,\nn) }^m - \mathbf{Y}_{t\wedge \tau (\mm ,\nn) }^m |
 = & \big|
\int_0^{t\wedge \tau (\mm ,\nn) }
b^{m} (\mathbf{X}_u^m,\mathbf{B}_u^m,\mathsf{X}_u^{m*}) -
b^{m} (\mathbf{Y}_u^m,\mathbf{B}_u^m,\mathsf{X}_u^{m*}) du 
\big|
\\ \notag \le 
&\sqrt{m} \cref{;GY3a}
\int_0^{t\wedge \tau (\mm ,\nn) }
|\mathbf{X}_{u}^m -\mathbf{Y}_{u}^m |du
.\end{align}
From this we easily deduce 
\begin{align}& \label{:43j}
(\mathbf{X}_t^m,\mathbf{B}_t^m,\mathsf{X}_t^{m*}) = 
(\mathbf{Y}_t^m,\mathbf{B}_t^m,\mathsf{X}_t^{m*})
\quad \text{ for all $ t \le \tau (\mm ,\nn) $}
.\end{align}
%\thetathree (\mathsf{X}^{m*})
Suppose that 
$ \tau (\mm ,\nn) < \thetathree (\mathsf{X}^{m*})$. 
% $ \tau (\mm ,\nn) = \thetaone (\mathbf{Y}^m,\mathsf{X}^{m*}) $ or 
% $ \thetatwo (\mathbf{X}^m,\mathsf{X}^{m*},\mathbf{Y}^m) $. 
Note that $ \RRprCs $ is an open set. 
Combining this with \eqref{:43o}, \eqref{:43r}, and \eqref{:43j} we obtain 
\begin{align}\label{:43v}&
\thetaone (\mathbf{X}^m,\mathsf{X}^{m*})=
\thetaone (\mathbf{Y}^m,\mathsf{X}^{m*}) = 
\thetatwo (\mathbf{X}^m,\mathsf{X}^{m*},\mathbf{Y}^m) 
.\end{align}
Hence we deduce from the definition of $ \tau (\mm ,\nn) $ and \eqref{:43v}
\begin{align}\label{:43w}
\tau (\mm ,\nn) = &
\min \{ \thetaone (\mathbf{X}^m,\mathsf{X}^{m*}), \thetathree (\mathsf{X}^{m*}) \} 
.\end{align}

From \eqref{:43w}, \eqref{:43Q}, and \pref{l:42} we deduce for $ \mug $-a.s.\,$ \mathsf{s}$ 
\begin{align*}&
\PPs (\limi{\nn }\limi{\mm }\tau (\mm ,\nn) =\infty ) 
\\= &
\PPs (\limi{\nn }\limi{\mm }
 \min \{ \thetaone (\mathbf{X}^m,\mathsf{X}^{m*}), \thetathree (\mathsf{X}^{m*}) \}
=\infty ) 
%\quad \text{ by \eqref{:43w}}
\\\ge & 
\PPs ( \limi{\nn }\limi{\mm }
\tau_{\HImnC } (\mathbf{X}^m,\mathsf{X}^{m*}) = \infty ) 
%\quad \text{ by \eqref{:43Q} }
\\ \notag = &1
%\quad \quad \text{ by \pref{l:42}}
.\end{align*}
We therefore deduce that the equality in \eqref{:43j} holds for all $ \zti $. 
We thus see that the SDE \ginSDE\ has a pathwise unique weak solution.

By $ (\hat{\mathbf{Z}}^m , \hat{\mathbf{B}}^m , \hat{\mathsf{X}}^{m*})
\elaw (\mathbf{X}^m , \mathbf{B}^m , \mathsf{X}^{m*})$ we can prove \thetag{2} 
similarly as \thetag{1}. 

We note that the coefficients of the martingale terms of the SDE are constant. 
Hence we have not used any specific property of the filtrations in the proof of \thetag{1}. 
Then we can prove \thetag{3} similarly as \thetag{1}. 
\qed
\end{proof}

\begin{lemma} \label{l:43} \thetag{1} 
For each $ m \in \N $, the SDE \ginSDE\ has a unique strong solution $ \Fms $ for $ \XB $ 
under $ \PPs $ for $ \mug $-a.s.\ $ \mathsf{s} $, where $ \mathbf{s}= \lab (\mathsf{s})$. 

\noindent \thetag{2} For $ \mug $-a.s.\ $ \mathsf{s}$ 
\begin{align}\notag &% & \label{:43a} 
\Fms (\lBlmM ) = \mathbf{X}^{m} \quad \text{ under }\PPs 
.\end{align}
\end{lemma}
\begin{proof} 
In Yamada--Watanabe theory the existence of weak solutions and the pathwise uniqueness imply that 
the SDE \ginSDE\ has a unique strong solution (see Theorem 1.1 in \cite[163p]{IW}). 
We modify it to prove \thetag{1}. 
Indeed, we shall prove the general result \pref{l:98} including \thetag{1}. 
The appearance of new randomness $ \mathsf{X}^*$ requires 
a substantial modification of the theory. 

We next prove \thetag{2}. 
From \thetag{1}, the SDE \ginSDE\ has a unique strong solution $ \Fms $. %[] 
Because $ (\mathbf{X}^m , \mathbf{B}^m , \mathsf{X}^{m*})$ under $ \PPs $ 
is a weak solution, \thetag{2} follows from \thetag{1}. 
\qed \end{proof}

\Ssection{Completion of proof of \tref{l:35}. }
Let $ \XB $ and $ \Pmg = \int \PPs \mug (d\mathsf{s})$ be as \sref{s:32}. 
\begin{lemma} \label{l:44} 
$ \XB $ under $ \Pmg $ satisfies \As{\Cfour}. 
\end{lemma}
\begin{proof}
Let 
$ \ER (t) = \int_t^{\infty} (1/\sqrt{2\pi }) e^{-|x|^2/2} dx $ be a (scaled) complementary error function. 
Let $ \rho _{\mathrm{Gin}} ^1 (x) = 1/\pi $ be the one-point correlation function of $ \mug $. 
Then 
\begin{align}\label{:44a} & \quad 
 \int_{\Rtwo }\ER ( \frac{|x|-r}{T}) \rho _{\mathrm{Gin}} ^1 (x) dx < \infty 
\quad \text{ for each $ r,T \in \N $}
.\end{align}

Let $ (\E ^{\mug ^{[1]}}, \dom _{\mathrm{Gin}}^{[1]}) $ be the Dirichlet form 
of the $ 1 $-labeled process given by \eqref{:41c}. 
Let $ \PP_{(x,\mathsf{s})}^{[1]} $ be the associated diffusion measure 
starting at $ (x,\mathsf{s})$. We set 
\begin{align}&\label{:44z}
\mathsf{P}_{\mug ^{[1]} }^{[1]} = 
\int_{\Rtwo \ts \sSS } \PP_{(x,\mathsf{s})}^{[1]} \mug ^{[1]} (dx d\mathsf{s})
.\end{align}% 
Let $ \XoneX = (X^1 , \sum_{i=2}^{\infty} \delta_{X^i})$ 
denote the one-labeled process. 
Then applying the Lyons--Zheng decomposition (\cite[Theorem 5.7.1]{fot.2}) 
to the coordinate function $ x $, we have for each $ 0 \le t \le T $, 
\begin{align}\label{:44b}&
X_t^1 - X_0^1 = \frac{1}{2}B_t^1 + \frac{1}{2} \hat{B}_t^1 
\quad \text{ for $ \mathsf{P}_{\mug ^{[1]} }^{[1]} $-a.e.}
,\end{align}
where $ \hat{B}^1 $ is the time reversal of $ B^1 $ on $ [0,T]$ such that 
\begin{align}&\label{:44d}
\hat{B}_t^1 = B_{T-t}^1 (r_T)- B_{T}^1 (r_T)
.\end{align}
Here we set $ \map{r_T}{C([0,T]; \Rtwo \ts \sSS )}
{C([0,T]; \Rtwo \ts \sSS )}$ such that $ r_T(w )(t) = w (T-t)$. 

Because of the coupling in \lref{l:90} ({\cite[Theorem 2.4]{{o.tp}}}), 
the definition of the one-Campbell measure, and \eqref{:44z}, we see 
\begin{align}\label{:44c}
 \sum_{i=1}^{\infty} \Pmg ( \inf_{t\in[0,T]}|X_t^i| \le r ) 
= &
 \int_{S\times \sSS} \mug ^{[1]}
(dx d\mathsf{s}) \mathsf{P}_{(x,\mathsf{s})}^{[1]} ( \inf_{t\in[0,T]}|X_t^1| \le r )
\\ \notag = &
\PP _{\mug ^{[1]}} ^{[1]} 
( \inf_{t\in[0,T]}|X_t^1| \le r )
.\end{align}
Then we have, taking $ x = X_0^1$, 
\begin{align} \label{:44o}
 & \PP _{\mug ^{[1]}} ^{[1]} ( \inf_{t\in[0,T]}|X_t^1| \le r ) \le \, 
\PP _{\mug ^{[1]}} ^{[1]} ( \sup_{t\in[0,T]}|X_t^1- x | \ge | x | - r ) 
\\ \notag \le \, &
\PP _{\mug ^{[1]}} ^{[1]} ( \sup_{t\in[0,T]} |B_t^1| \ge | x | - r ) + 
\PP _{\mug ^{[1]}} ^{[1]} ( \sup_{t\in[0,T]} |\ \hat{B}_t^1| \ge | x | - r ) 
 \quad \quad \text{by \eqref{:44b}}
\\ \notag = \, & 2 
 \PP _{\mug ^{[1]}} ^{[1]} ( \sup_{t\in[0,T]} |B_t^1| \ge | x | - r ) 
\\ \notag \le &2
 \int_{\sSS } 
\ER ( \frac{| x |-r}{{\cref{;44}}T} ) \mug ^{[1]} (dxd\mathsf{s}) 
 \\ \notag =&2 
\int_{\Rtwo }\ER ( \frac{|x|-r}{{\cref{;44}}T}) \rho _{\mathrm{Gin}} ^1 (x) dx 
< \, \infty \quad \quad \text{by \eqref{:44a}}
,\end{align}
where $ \Ct \label{;44}$ is a positive constant. 
From \eqref{:44c} and \eqref{:44o} we deduce 
\begin{align*}&
 \sum_{i=1}^{\infty} \Pmg ( \inf_{t\in[0,T]}|X_t^i| \le r ) < \infty 
.\end{align*}
Hence from Borel--Cantelli's lemma we have for each $ r ,T \in \N $ 
\begin{align}&\label{:44p}
\Pmg ( \limsup_{i\to\infty }\{\inf_{t\in[0,T]}|X_t^i| \le r \}) =0 
.\end{align}
Then we deduce from \eqref{:53n} and \eqref{:44p} for each $ r ,T \in \N $ 
\begin{align*}&
\Pmg ( \mrXX < \infty ) = 1 - \Pmg ( \limsup_{i\to\infty }\{\inf_{t\in[0,T]}|X_t^i| \le r \}) = 1 
.\end{align*}
This completes the proof. 
\qed \end{proof}

\smallskip 
\begin{lemma} \label{l:46} \thetag{1}
$ \XB $ under $ \Pmg $ satisfies \As{\iFc} for ISDE \eqref{:11h}. 
\\\thetag{2} 
$ \XB $ under $ \Pmg $ satisfies \As{\Ctwo} for $ \mug $. 
\end{lemma}

\begin{proof}
By \lref{l:34}, $ \XB $ under $ \PPs $ is a weak solution of \eqref{:11h} 
starting at $ \mathbf{s} = \lab (\mathsf{s})$ for $ \mug $-a.s.\,$ \mathsf{s}$. 
By \lref{l:43} we see $ \XB $ under $ \PPs $ satisfies \iFcs\ for $ \mug $-a.s.\,$ \mathsf{s}$. 
Hence the first claim holds. 
The second claim is obvious because $ \{ \PPs \} $ is a $ \mug $-stationary diffusion 
and $ \Pmg = \int \PPs \mug (d\mathsf{s})$. 
\qed\end{proof}

 \smallskip 

\noindent {\em Proof of \tref{l:35}. } 
We use \tref{l:5A}. 
For this, we check that $ \mug $ is tail trivial and that 
$ \XB $ under $ \PP _{\mug } $ satisfies {\Inmg}. 

Because $ \mug $ is a determinantal random point field, 
$ \mug $ satisfies \As{\muTT}. 
\As{\iFc} for \eqref{:11h} and \As{\Ctwo} for $ \mug $ follow from \lref{l:46}. 
\As{\sIn} follows from \lref{l:32}. 
\As{\nbj} follows from \lref{l:44}. 
Thus all the assumptions of \tref{l:5A} are fulfilled. 
Hence we obtain \thetag{3} by \tref{l:5A}. 
\thetag{4} is clear. Indeed, 
we can take a subset $ \mathsf{H}$ such that $ \mug (\mathsf{H}) = 1 $ and that 
the same conclusions of \thetag{3} hold for all $ \mathsf{s} \in \mathsf{H}$.

If a weak solution $ (\hat{\mathbf{X}},\hat{\mathbf{B}})$ of \eqref{:11g} satisfies 
\As{\muAC} for $ \mug $. Then $ (\hat{\mathbf{X}},\hat{\mathbf{B}})$ satisfies the ISDE \eqref{:11h} 
because of the coincidence of the logarithmic derivatives \eqref{:33a} and \eqref{:33b} 
given by \lref{l:33}. Furthermore, the SDE \eqref{:35c} for 
$ (\hat{\mathbf{X}},\hat{\mathbf{B}}) $ coincides with 
\begin{align} \notag &% \label{:49a}&
dY_t^{m,i} = d\hat{B}_t^i + 
\sum_{j=1 \atop j\not=i}^m 
\frac{Y_t^{m,i}-Y_t^{m,j}}{|Y_t^{m,i}-Y_t^{m,j}|^{2}} dt + 
\lim_{r\to\infty } \sum_{j = m +1 ,\atop |Y_t^{m,i}-\hat{X}_t^j|<r } ^{\infty}
\frac{Y_t^{m,i}-\hat{X}_t^j}{|Y_t^{m,i}-\hat{X}_t^j|^{2}} dt 
.\end{align}
Thus $ (\hat{\mathbf{X}},\hat{\mathbf{B}})$ is a weak solution of \eqref{:11g} 
satisfying \As{\iFc } if and only if 
$ (\hat{\mathbf{X}},\hat{\mathbf{B}})$ is a weak solution of \eqref{:11h} satisfying \As{\iFc }. 
Hence \thetag{1} and \thetag{2} follow from \thetag{3} and \thetag{4}. 

We proceed with the proof of \thetag{5}. 
Recall that both solutions satisfy \As{\Ctwo} for $ \mug $. Then the drift coefficients are equal 
because of the coincidence of the logarithmic derivatives given by \lref{l:33}. 
Hence we deduce \thetag{5} from \thetag{1} and \thetag{3}. 
\qed

\section{Dirichlet forms and weak solutions} \label{s:D}

\subsection{Relation between ISDE and a random point field.}
\label{s:D1}

We shall deduce the existence of a strong solution and the pathwise uniqueness 
of solution of ISDE \ISDEb from the existence of a random point field $ \mu $ 
satisfying the assumptions in the sequel. 
Recall the notion of the logarithmic derivative $ \dmu $ 
given by \dref{d:21}. 
To relate the ISDE \eqref{:50a} with the random point field $ \mu $, 
we make the following assumption. 

\smallskip 
\noindent \As{A1} 
There exists a random point field $ \mu $ 
such that $ \sigma \in L_{\mathrm{loc}}^{\infty}(\muone ) $ and 
$ \bbb \in L_{\mathrm{loc}}^{1}(\muone )$ and that $ \mu $
has a logarithmic derivative 
 $ \dmu = \dmu (x,\mathsf{y}) $ satisfying the relation 
\begin{align}\label{:51a}&
\bbb (x,\mathsf{y}) = \frac{1}{2} \{\nabla_x \cdot \aaa (x,\mathsf{y}) + 
\aaa (x,\mathsf{y}) \dmu (x,\mathsf{y}) \} 
.\end{align}
Here $ \nabla_x \cdot \aaa (x,\mathsf{y})= 
(\sum_{q=1}^d \PD{\aaa _{pq}}{x_q}(x,\mathsf{y}))_{p=1}^d $ 
and $ \aaa (x,\mathsf{y})=\{ \aaa _{pq}(x,\mathsf{y}) \}_{p,q=1}^d $ 
is the $ d\ts d$-matrix-valued function defined by 
\begin{align}\label{:51b}& 
\aaa (x,\mathsf{y}) = \sigma (x,\mathsf{y}) ^t \sigma (x,\mathsf{y}) 
.\end{align}
Furthermore, we assume that there exists a bounded, lower semi-continuous, non-negative function 
$ \map{a_0}{\mathbb{R}^+}{\mathbb{R}^+}$ and positive constants 
$ \Ct \label{;A1a}$ and $ \Ct \label{;79} $ such that 
\begin{align}\label{:51c}& 
a_0 (|x|) |\xi |^2 \le (\aaa (x,\mathsf{y}) \xi , \xi )_{\Rd } 
 \le \cref{;A1a} a_0 (|x|) 
|\xi |^2 \quad \text{ for all $ \xi \in \Rd $}
,\\\label{:51cc}&
 \cref{;79}:= \sup_{t\in\R ^+} a_0 (t) < \infty 
\end{align}
and that $ \aaa (x,\mathsf{y}) $ is smooth in $ x $ for each $ \mathsf{y}$.

\begin{remark}\label{r:51a}
From \eqref{:20k}, we obtain an informal expression 
 $ \dmu = \nabla_x \log \muone $. We then interpret the relation \eqref{:51a} 
 as a differential equation of random point fields $ \mu $ as 
\begin{align}\label{:51d}& 
\bbb (x,\mathsf{y}) = \frac{1}{2} \{\nabla_x \cdot \aaa (x,\mathsf{y}) + 
\aaa (x,\mathsf{y}) \nabla_x \log \muone (x,\mathsf{y}) \} 
.\end{align}
That is, \eqref{:51d} is an equation of $ \mu $ for given coefficients 
$ \sigma $ and $ \bbb $. 
Theorems \ref{l:5A}--\ref{l:5B} deduce that, 
if the differential equation \eqref{:51d} of random point fields has 
a solution $ \mu $ satisfying the assumptions in these theorems, 
then the existence of a strong solution and the pathwise uniqueness 
of a solution of ISDE \ISDEb hold. 
\end{remark}

\subsection{A weak solution of ISDE (First step).} \label{s:5c}

In this subsection, we recall a construction of $ \mu $-reversible diffusion 
from \cite{o.dfa,o.rm}. 
Recall the notion of quasi-Gibbs property given by \dref{d:22}. 
We now make the following assumption. 

\smallskip 
\noindent 
\As{A2} $ \mu $ is a $ ( \Phi , \Psi )$-quasi Gibbs measure such that 
 there exist upper semi-continuous functions $ (\hat{\Phi }, \hat{\Psi })$ 
 and positive constants $ \Ct \label{;A21}$ and $ \Ct \label{;A22} $ satisfying 
\begin{align} & \notag %\label{:A2a}
\cref{;A21}^{-1} \hat{\Phi } (x)\le \Phi (x) \le \cref{;A21} \hat{\Phi }(x) ,\quad 
\cref{;A22}^{-1} \hat{\Psi } (x,y)\le \Psi (x,y) \le \cref{;A22} \hat{\Psi }(x,y) 
\end{align}
and that $ \Phi $ and $ \Psi $ are locally bounded from below.

\smallskip 
We refer to \cite{o.rm,o.rm2} for sufficient conditions of \As{A2}. 
These conditions give us the quasi-Gibbs property of the random point fields 
appearing in random matrix theory, such as 
sine$_{\beta} $, Airy$_{\beta} $ ($ \beta=1,2,4$), and Bessel$_{2, \alpha } $ 
($ 1\le \alpha $), and the Ginibre random point fields \cite{o.rm,o.rm2,o-t.airy,h-o.bes}.

Let 
$ \sigma _r^{m}$ be the $ m $-density function of $ \mu $ on $ \Sr $ 
with respect to the Lebesgue measure $ d\mathbf{x}^{m}$ on $ \Sr ^{m}$. 
By definition, $ \sigma _r^{m}$ is a non-negative symmetric function such that 
\begin{align}\notag % \label{:A3z}
& \frac{1}{m!}
\int_{\Sr ^m} \check{f} (\mathbf{x}^m)\sigma _r^{m} (\mathbf{x}^m) 
d\mathbf{x}^m = 
\int_{\sSS _r^{m}} f (\mathsf{x}) \mu (d\mathsf{x})
\quad \text{ for all } f \in C_b(\sSS )
.\end{align}
Here we set $ \mathbf{x}^m =(x_1,\ldots,x_m)$ and 
$ d\mathbf{x}^m = dx_1\cdots dx_m $. Furthermore, we denote by 
$ \sSS _r^{m} $ the subset of $ \sSS $ such that 
$ \sSS _r^{m} = \{ \mathsf{s} \in \sSS \, ;\, \mathsf{s} (\Sr ) = m \} $ as before. 
We make the following assumption. 

\medskip 

\noindent 
\As{A3} \ $ \mu $ satisfies for each $ m, r \in \N $ 
\begin{align}
\label{:51h}&
\sum_{k=m}^{\infty} \frac{k!}{(k-m)!} \mu (\sSS _r^k ) < \infty 
.\end{align}
Clearly, \eqref{:51h} is equivalent to 
$ \int_{\Srm }\rho ^m (\mathbf{x}^m ) d\mathbf{x}^m < \infty $ 
% for all $ r \in \mathbb{N}$ 
if the $ m $-point correlation function $ \rho ^m $ of $ \mu $ with respect to the Lebesgue measure exists. 
Under assumptions of \As{A2} and \As{A3} $ \mu $ has correlation functions and Campbell measures 
of any order.

Let $ (\Eamu ,\di ^{\amu } )$ be a bilinear form on $ \Lmu $ with domain $ \di ^{\amu }$ defined by 
\begin{align} 
\label{:51i} & 
\di ^{\amu } = \{ f \in \di \cap \Lmu \, ;\, \Eamu (f,f) < \infty \} 
,\\ \label{:51j}& 
\Eamu (f,g) = \int_{\sSS } \DDDa [f,g] \, \mu (d\mathsf{s}) 
,\\ \label{:51k}& 
\DDDa [f,g] (\mathsf{s}) = \frac{1}{2} 
\sum_{i} ( a (s_i , \mathsf{s}^{i\diamondsuit}) 
 \partialsi \check{f} , \partialsi \check{g} )_{\Rd }
.\end{align}
Here $ \mathsf{s}=\sum_i \delta_{s_i}$ and 
$ \mathsf{s}^{i\diamondsuit} = \sum_{j\not=i} \delta _{s_j}$, 
$ \partialsi = (\PD{}{s_{i1}},,\ldots,\PD{}{s_{id}})$ and 
$(\cdot , \cdot )_{\Rd }$ denotes the standard inner product in $ \mathbb{R}^d$. 
When $ a $ is the unit matrix, we often remove it from the notation; 
for example, $ \Eamu = \Emu $, $ \di ^{\amu } = \di ^{\mu } $, and 
$ \DDDa = \DDD $. 

\begin{lemma} \label{l:51} 
Assume \eqref{:51c}, \eqref{:51cc}, \As{A2}, and \As{A3}. 
Then $ (\Eamu ,\di ^{\amu } )$ is 
closable on $ \Lmu $, and its closure $ (\Eamu ,\Damu )$ 
 is a quasi-regular Dirichlet form on $ \Lmu $. 
 Moreover, the associated $ \mu $-reversible diffusion 
$ (\mathsf{X},\{\Pssf \}_{\mathsf{s}\in\mathsf{S} }) $ exists. 
\end{lemma}

\lref{l:51} is a refinement of \cite[119p.\! Corollary 1]{o.dfa} and can be proved similarly. 

In \tref{l:5B} we decompose $ \mu $ by taking the regular conditional probability 
with respect to the tail $ \sigma $-field. The refinement above is motivated by this decomposition. 
Indeed, \eqref{:51h} is stable under this decomposition (see \lref{l:Q1}). 

By construction, a diffusion measure $ \Pssf $ given by a quasi-regular Dirichlet form 
with quasi-continuity in $ \mathsf{s}$ is unique for quasi-everywhere starting point $ \mathsf{s}$. 
Equivalently, there exists a set $ \HH  $ such that the complement of $ \HH $ has capacity zero, 
and the diffusion measure $ \Pssf $ associated with the Dirichlet space above 
with quasi-continuity in $ \mathsf{s}$ is unique for all $ \mathsf{s} \in \HH  $ 
and $ \Pssf (\mathsf{X}_t \in \HH  \text{ for all }t ) = 1 $ 
for all $ \mathsf{s} \in \HH  $. The set $ \HH $ is unique up to capacity zero. 

Let $ \{\Pssf \}_{\mathsf{s}\in\mathsf{S} }$ be as in \lref{l:51}. 
%Note that $ \mu (\mathsf{H}) = 1 $. 
We set 
\begin{align}\label{:51u}&
\Pm = \int_{\mathsf{S}} \Pssf \, \mu (d\mathsf{s})
\end{align}
Let $ \WSsiNE $ be as in \eqref{:23gg}. 
We assume the following. 

\smallskip 

\noindent 
\As{\sIn}\quad $ \Pm $ satisfies $ \Pm (\WSsiNE ) = 1 $. 

\smallskip \noindent 
From \As{\sIn} we see $ \lpathX $ is well defined under $ \Pm $-a.s. 
We present an ISDE that $ \mathbf{X}= \lpathX $ satisfies. 
\begin{lemma}[{\cite[Theorem 26]{o.isde}}] \label{l:52} 
Assume \As{A1}--\As{A3}. 
Assume that $ \Pm $ satisfies \As{\sIn}. 
Then there exists an $ \mathsf{H} \subset \SSsde $ satisfying \eqref{:50x} and \eqref{:50y} such that 
$ (\lpathX ,\mathbf{B}) $ defined on $ \OFPFs $ 
is a weak solution of ISDE \eqref{:50a}--\eqref{:50c} for each $ \mathsf{s}\in \mathsf{H} $ and that 
$ \mu (\mathsf{H} ) = 1 $ and 
$ \Pssf (\mathsf{X}_t \in \mathsf{H} \text{ for }0 \le \forall t < \infty ) = 1 $. 
\end{lemma}

\begin{remark}\label{r:52} 
\thetag{1} The solution $ \lpathX $ of the ISDE \eqref{:50a}--\eqref{:50b} 
is defined on the space $ (\WS ,\mathcal{B}(\WS )) $ with filtering $ \Bt (\sSS ) $ given by \eqref{:32a}. 
\\\thetag{2} 
In \lref{l:52}, Brownian motion $ \mathbf{B}=(B^i)_{ i \in \N }$ 
is given by the additive functional of the diffusion $ (\mathsf{X},\{\Pssf \})$. 
In particular, $ \mathbf{B} $ is a functional of $ \mathsf{X}$. 
Hence we write $ \mathbf{B}(\mathsf{X})$ when we emphasize this. 
Summing up, the weak solution in \lref{l:52} is given by 
$ (\lpathX ,\mathbf{B}(\mathsf{X}))$ and 
defined on the filtered space $ (\WS ,\mathcal{B}(\WS ))$ with $ \{ \Bt (\sSS )\} $. 
\\\thetag{3} 
A sufficient condition of \As{\sIn} will be given in \sref{s:E}. 
\end{remark}

\Ssection{The Dirichlet forms of the $ m $-labeled processes and the coupling.}\label{s:90}
Let $ (\Eamu ,\Damu )$ be the Dirichlet form on $ \Lmu $ in \As{A3}, and let 
$ (\mathsf{X},\{ \PP _{\mathsf{s}} \}_{\mathsf{s}\in \sSS } )$ be the associated diffusion. 

Let $ \lab $ be a label.
Assume that \As{\sIn} holds. 
Then $ \lpath $ makes sense and we construct the fully labeled process 
 $ \mathbf{X} =(X^i)_{i\in\mathbb{N}}$ 
with $ \mathbf{X}_0 = \lab (\mathsf{X}_0)$ 
associated with the unlabeled process $ \mathsf{X} $
 by taking $ \mathbf{X} = \lpathX $, where $ \mathsf{X} = \sum_{i\in\mathbb{N}} \delta_{X^i} $. 

Let $\mathbf{X}^m = (X^i)_{i=1}^m $ and $ \mathsf{X}^{m*} = \sum_{m<i } \delta_{X^i} $. 
We call the pair $ (\mathbf{X}^m,\mathsf{X}^{m*} )$ the $ m $-labeled process. 
We shall present the Dirichlet form associated with the $ m $-labeled process. 

We write $ \xxxm =(x_i)_{i=1}^m \in (\Rd )^m$. 
Let $ \mu^{[m]} $ be the (reduced) $ m $-Campbell measure of $ \mu $ defined as 
\begin{align}
\label{:90A}&
\mu^{[m]} (d\xxxm d\mathsf{y}) = \rho^m (\xxxm ) \mu _{\xxxm } ( d\mathsf{y})d\xxxm 
,\end{align}
where $ \rho^m $ is the $ m $-point correlation function of $ \mu $ 
with respect to the Lebesgue measure $ dx $ on $ \Sm $, and 
$ \mu _{\xxxm } $ is the reduced Palm measure conditioned at $ \xxxm \in \Sm $. Let 
\begin{align} \label{:90b} &
\E ^{\amu^{[m]}} (f,g) = \int_{\sS ^m \times \mathsf{S}} \Big\{
\frac{1}{2}\sum_{i=1}^{m} 
( a (x_i , \mathsf{x}^{i\diamondsuit } + \mathsf{y}) \PD{f}{x_i}, \PD{g}{x_i})_{\Rd } + 
\DDDa [f,g] \Big\} \mu^{[m]} (d\xxxm d\mathsf{y})
,\end{align}
where $ \PD{}{x_i}$ is the nabla in $ \Rd $, 
$ \mathsf{x}^{i\diamondsuit } = \sum_{j\not=i}^m \delta_{x_j}$, and 
$ a (x,\mathsf{y}) $ is given by \eqref{:51b}. 
Moreover, 
$ \DDDa $ is defined by \eqref{:51k} naturally 
regarded as the carr\'{e} du champ on $ \sS ^m \times \mathsf{S}$, and 
\begin{align}
\notag &% \label{:90c}&
\di ^{\amu ^{[m]}}=\{ f \in C_0^{\infty}(\sS ^m)\otimes \di \, ;\, 
\E ^{\amu^{[m]}} (f,f) < \infty ,\ f \in L^2(\sS ^m \times \mathsf{S}, \mum ) \} 
.\end{align}
When $ m=0 $, we interpret $ \E ^{\amu^{[0]}} $ and $ \mu^{[0]} $ 
as the unlabeled Dirichlet form $ \Eamu $ and $ \mu $, respectively. 

The closablity of the bilinear form $ (\E ^{\amu^{[m]}} , C_0^{\infty}(\sS ^m)\otimes \di ) $ 
on $ L^2(\sS ^m \times \mathsf{S}, \mum ) $ follows from 
\eqref{:51c}, \eqref{:51cc}, \As{A2}, and \As{A3}. 
We can prove this in a similar fashion as the case $ m =0 $ as \lref{l:51} (\cite{o.dfa}). 
We then denote by $ (\E ^{\amu ^{[m]}} ,\dom ^{\amu ^{[m]}} )$ its closure. 
% of $ (\E ^{\amu^{[m]}} , C_0^{\infty}(\sS ^m)\otimes \di ) $ on on $ L^2(\sS ^m \times \mathsf{S}, \mum ) $. 
%
The quasi-regularity of $ (\E ^{\amu ^{[m]}} , \dom ^{\amu ^{[m]}} )$ 
is proved by \cite{o.tp} for $ \mu $ with bounded correlation functions. 
The generalization to \As{A3} is easy, and we left its proof. 

Let $\PP _{(\mathbf{s}^m,\mathsf{s}^{m*})}^{[m]}$ denote the diffusion measures associated with the 
$ m $-labeled Dirichlet form $ (\E ^{\amu^{[m]}} , \dom ^{\amu ^{[m]}} ) $ 
on $ L^2(\sS ^m \times \mathsf{S}, \mum ) $. 
(see \cite{o.tp}). 
We quote: 
\begin{lemma}[{\cite[Theorem 2.4]{{o.tp}}}] \label{l:90}
Let $ \mathsf{s}=\ulab ((\mathbf{s}^m,\mathsf{s}^{m*})) = 
\sum_{i=1}^m\delta_{s_i}+\mathsf{s}^{m*}$. Then 
\begin{align}
\notag &% \label{:90a}&
\PP _{(\mathbf{s}^{m},\mathsf{s}^{m*})}^{[m]} = 
\PP _{\mathsf{s}} \circ (\mathbf{X}^m, \mathsf{X}^{m*})^{-1} 
.\end{align}
\end{lemma}

Note that $ \PPs $ in the right-hand side is independent of $ m \in \N $. 
Hence this gives a sequence of coupled $ \Sm \ts \sSS $-valued continuous processes 
with distributions $\PP _{(\mathbf{s}^m,\mathsf{s}^{m*})}^{[m]}$. 
In this sense, there exists a natural coupling among the $ m $-labeled Dirichlet forms 
$ (\E ^{\amu ^{[m]}} ,\dom ^{\amu ^{[m]}}) $ on $ L^2(\sS ^m\ts \sSS ,\mum ) $. 
This coupling is a crucial point of the construction of weak solutions of ISDE in \cite{o.isde}.

Introducing the $ m $-labeled processes, 
we can regard $ \mathbf{X}^m$ as a Dirichlet process of the diffusion $ (\mathbf{X}^m, \mathsf{X}^{m*})$ associated with the $ m $-labeled Dirichlet space. 
%$ (\E ^{\amu ^{[m]}} ,\dom ^{\amu ^{[m]}}) $ on $ L^2(\sS ^m\ts \sSS ,\mum ) $. 
That is, one can regard 
$ A_t^{[\mathbf{x_m}]}:= \mathbf{X}_t^m - \mathbf{X}_0^m $ 
as a $ dm$-dimensional additive functional given by 
the composition of $ (\mathbf{X}^m, \mathsf{X}^{m*})$ with 
the coordinate function $ \xxxm =(x_1,\ldots,x_m) \in (\Rd )^m$. 
Although $ \mathbf{X}^m$ can be regarded as an additive functional of the unlabeled process 
$ \mathsf{X}=\sum_i \delta _{X^i}$, $ \mathbf{X}^m$ is 
no longer a Dirichlet process in this case. 
Indeed, as a function of $ \mathsf{X} $, we cannot identify 
$ \mathbf{X}_t^m$ without tracing the trajectory of 
$ \mathsf{X}_s = \sum_i \delta_{X_s^i}$ for $ s \in [0,t]$. 

Once $ \mathbf{X}^m$ can be regarded as a Dirichlet process, 
we can apply the It$ \hat{\mathrm{o}}$ formula (Fukushima decomposition), and 
Lyons--Zheng decomposition to $ \mathbf{X}^m$, 
which is important in proving the results in the subsequent sections.

\section{Sufficient conditions of \As{\sIn} and \As{\nbj}.}\label{s:E}

The purpose of this section is to give sufficient conditions of \As{\sIn} and \As{\nbj} 
for \XlP. 
The following assumptions are related to Dirichlet forms introduced in \sref{s:5c}. 
So we named the series of them as \As{A4} followed by \As{A3} in \sref{s:5c}. 
Let $ \ER = \ER (t) $ be a (scaled) complementary error function: 
\begin{align*}&
 \ER (t) = \int_t^{\infty} (1/\sqrt{2\pi }) e^{-|x|^2/2} dx
.\end{align*}
We set $ \langle f, \mathsf{s}\rangle = \sum_{i} f (s_i)$ for 
$ \mathsf{s}= \sum_i \delta_{s_i}$ as before.

\medskip 

\noindent \As{A4} \ \thetag{1} 
For each $ r,T \in \N $ 
\begin{align}\label{:54k}&
E^{\mu }\Big[\langle \ER ( \frac{|\cdot |-r}{T}) , \mathsf{s}\rangle \Big] < \infty 
,\\\intertext{and there exists a $ T>0 $ such that for each $ R>0 $ }
\label{:54l}&
\liminfi{r} \ER (\frac{r}{T}) \, 
E^{\mu }[\langle 1_{\sS _{r+R}}, \mathsf{s}\rangle] = 0 
.\end{align}
\thetag{2} Each $ X^i $ neither hits the boundary $ \partial \sS $ of $ S $ nor collides each other. 
That is, 
\begin{align}\label{:92z}&
\mathrm{Cap}^{\amu }(\{ \mathsf{s} ; \mathsf{s} (\partial \sS ) \ge 1 \})
 = 0
,\\\label{:91a}&
\mathrm{Cap}^{\amu } (\mathsf{S}_{\mathrm{s}}^c) = 0
.\end{align}
Here $ \mathrm{Cap}^{\amu } $ is the capacity of the Dirichlet form 
$ (\Eamu ,\Damu ) $ on $ \Lmu $. 
% Furthermore, 
% $ \mathsf{S}_{\mathrm{s}} = 
% \{ \mathsf{s} \in \sSS \, ;\, \mathsf{s}(\{ x \} ) \le 1 \text{ for all } x \in \sS \} $ 
% is the set of all configurations consisting of single-point masses. 

\medskip 

The condition \eqref{:54l} is an analogy of \cite[\thetag{5.7.14}]{fot.2}. 
Unlike this, the carr\'{e} du champ has uniform upper bounds $ \cref{;A1a}\cref{;79}$ 
by \eqref{:51b}--\eqref{:51cc}. 

We remark that \eqref{:54l} is easy to check. \eqref{:54l} holds if the 1-correlation function of $ \mu $ 
has at most polynomial growth at infinity. Obviously, condition 
\eqref{:92z} is always satisfied if $ \partial \sS = \emptyset $. 
We state sufficient conditions of \eqref{:91a}. 
\begin{lemma}[{\cite[Theorem 2.1, Proposition 7.1]{o.col}}] \label{l:91}
Assume \As{A1}--\As{A3}. Assume that $ \Phi $ and $ \Psi $ 
are locally bounded from below. We then obtain the following. 
\\\thetag{1} 
Assume that $ \mu $ is a determinantal random point field with a locally Lipschitz continuous kernel 
with respect to the Lebesgue measure. \eqref{:91a} then holds. 
\\\thetag{2} 
Assume that $ d \ge 2 $. \eqref{:91a} then holds. 
\end{lemma}
\begin{proof}
Claims \thetag{1} and \thetag{2} follow from 
\cite[Theorem 2.1, Proposition 7.1]{o.col}.
\qed \end{proof}

We now deduce \As{\sIn} from \As{A1}--\As{A4}. 
\begin{lemma} \label{l:92} 
Assume \As{A1}--\As{A4}. Then \XlP\ satisfies \As{\sIn} and 
\begin{align}\label{:92c}&
 \mathrm{Cap}^{\amu }(\SSsde ^c) = 0 ,\quad 
\Pm \circ \lpath ^{-1}(\WT (\SSSsde )) = 1 
.\end{align}
\end{lemma}%
\begin{proof} 
Applying the argument in \cite[Theorem 5.7.3]{fot.2} 
to $ \mathbf{X}^m $ of this Dirichlet form, 
we see that the diffusion $(\mathbf{X}^m, \mathsf{X}^{m*}) $ is conservative. 
Because this holds for all $ m \in \mathbb{N}$, 
\begin{align}
\label{:92a}&
\Pm ( \, \sup \{|X _t^i |;\, t \le T \} < \infty \, \text{ for all } T, i \in \mathbb{N} \, ) = 1 
.\end{align}

Because $ \mathsf{X}\in \WT (\sSS _{\mathrm{s}} )$ $ \Pm $-a.s. by \eqref{:91a}, 
we write $ \mathsf{X}_t = \sum_i \delta_{X_t^i} $ such that 
$ X^i \in C(I^i;\sS )$, where 
$ I^i = [0,b)$ or $ I^i = (a,b)$ for some $ a, b \in (0,\infty ]$. 
We shall prove $ I^i = [0,\infty )$ $ \Pm $-a.s.. 
Suppose that $ I^i = [0,b)$. 
Then, from \eqref{:92a}, we deduce that $ b = \infty $. 
Next suppose that $ I^i = (a,b)$. 
Then, applying the strong Markov property of the diffusion 
$ \{ \PP _{\mathsf{s}} \} $ at any $ a' \in (a,b)$ and using the preceding argument, 
we deduce that $ b = \infty $. As a result, we have $ I^i = (a,\infty )$. 
Because of reversibility, we see that such open intervals do not exist. 
Hence, we obtain $ I^i = [0,\infty )$ for all $ i $. 
From this, \eqref{:91a}, and $ \mu (\{ \mathsf{s}(\sS )=\infty \} ) = 1 $, 
we obtain $ \mathrm{Cap}^{\amu }(\{\mathsf{s} \in \sSS \, ; \mathsf{s}(\sS ) < \infty \} ) = 0 $. 
From this and \As{A4} \thetag{2} we have 
\begin{align}\label{:92b}&
\mathrm{Cap}^{a,\mu }(\Ssi ^c) = 0 
.\end{align}
From \eqref{:92a}, \eqref{:92b}, and \eqref{:92z} we obtain $ \Pm (\WSsiNE ) = 1 $. 
By \lref{l:52} we see the first claim in \eqref{:92c}. 
The second claim is clear from the first. We have thus completed the proof. 
\qed \end{proof}

We next deduce \As{\nbj} from \As{A1}--\As{A4}.

\begin{lemma} \label{l:94}
 Assume \As{A1}--\As{A4}. Then \XlP\ satisfies \As{\nbj}. 
\end{lemma}
\begin{proof}
Let $ \rho ^1 $ be the one-point correlation function of $ \mu $ with respect to the Lebesgue measure. 
Then \eqref{:54k} implies 
\begin{align}\label{:93c}
\int_{\sS } \ER (\frac{|x| - r}{\sqrt{\cref{;A1a}\cref{;79}}T} ) \rho ^1 (x) dx < \infty 
.\end{align}
Here $ \cref{;79}= \sup_{t\in\R ^+} a_0 (t) $ is finite by \eqref{:51cc}.

Let $ (\E ^{\amu^{[1]}} , \dom ^{\amu ^{[1]}} ) $ be the Dirichlet form 
of the $ 1 $-labeled process. 
Let $ \PP_{(x,\mathsf{s})}^{[1]} $ be the associated diffusion measure starting at 
$ (x,\mathsf{s}) \in \sS \ts \sSS $. 
We set 
\begin{align}&\label{:93d}
\mathsf{P}_{\muone }^{[1]} = 
\int_{\sS \ts \sSS } \PP_{(x,\mathsf{s})}^{[1]} \muone (dx d\mathsf{s})
.\end{align}
Let $ \XoneX = (X^1 , \sum_{i=2}^{\infty} \delta_{X^i}) $ denote the one-labeled process. 
Then applying the Lyons--Zheng decomposition (\cite[Theorem 5.7.1]{fot.2}) 
to the coordinate function $ x $, we have a continuous local martingale $ M^1 $ 
such that for each $ 0 \le t \le T $, 
\begin{align}\label{:93e}&
X_t^1 - X_0^1 = \frac{1}{2}M_t^1 + \frac{1}{2} \hat{M}_t^1 
\quad \text{ for $ \mathsf{P}_{\muone }^{[1]} $-a.e.}
,\end{align}
where $ \hat{M}^1 $ is the time reversal of $ M^1 $ on $ [0,T]$ such that 
\begin{align}&\label{:93f}
\hat{M}_t^1 = M_{T-t}^1 (r_T)- M_{T}^1 (r_T)
.\end{align}
Here we set $ \map{r_T}{C([0,T]; \Rd \ts \sSS )}
{C([0,T]; \Rd \ts \sSS )}$ such that $ r_T(w )(t) = w (T-t)$. 
By \eqref{:51a}, \eqref{:51b}, and \eqref{:90b} the quadratic variation of 
$ M^1=(M_1^1,\ldots,M_d^1) $ is given by 
\begin{align}
\label{:93g}&
\langle M_{p}^1 ,M_{q}^1 \rangle_t = 
\int_0^t \aaa _{pq}(X_u^1 , \sum_{i=2}^{\infty} \delta_{X_u^i})du 
.\end{align}

We can prove \lref{l:94} in a similar fashion as \lref{l:44}. 
Indeed, \eqref{:44a}, \eqref{:44z}, \eqref{:44b}, and \eqref{:44d} correspond to 
\eqref{:93c}, \eqref{:93d}, \eqref{:93e}, and \eqref{:93f}, respectively. 
Although the continuous local martingale $ M^1 $ in \eqref{:93e} is not Brownian motion, 
the quadratic variation process in \eqref{:93g} is controlled by \eqref{:51c} and \eqref{:51cc}. 
So the proof of \lref{l:44} is still valid for \lref{l:94}. 
Hence we omit the detail. 
\qed \end{proof}

\Section{Sufficient conditions of \As{\iFc}. } \label{s:9}

\Ssection{Localization of coefficients. }\label{s:9A}
Let $ \mathbf{a}=\{ \ak \}_{\qq \in\mathbb{N}} $ be a sequence of 
increasing sequences $ \ak = \{ \ak (r) \}_{r\in\mathbb{N}} $ 
 of natural numbers such that $ \ak (r) < \akk (r)$ for all $ r , \qq \in \mathbb{N}$. 
 We set for $ \mathbf{a} = \{\ak \}_{\qq \in\mathbb{N}} $ 
\begin{align}\label{:95w}&
\Ka = \bigcup_{\qq =1}^{\infty} \mathsf{K}[\ak ], \quad 
\mathsf{K}[\ak ] =\{ \mathsf{s}\, ;\, 
\mathsf{s} (\Sr ) \le \ak (r) \text{ for all } r \in \mathbb{N} \} 
.\end{align}
By construction, $ \mathsf{K}[\ak ] \subset \mathsf{K}[\akk ]$ 
for all $ \qq \in\mathbb{N}$. It is well known that $ \mathsf{K}[\ak ]$ is a compact set 
 in $ \sSS $ for each $ \qq \in \mathbb{N}$. 
To quantify $ \mu $ by $ \mathbf{a}$ we assume

\smallskip 
\noindent \As{\Bone} \quad 
 $ \mu (\Ka ) = 1 $. 

\smallskip 

We note that a sequence $ \mathbf{a}$ with $ \mu (\Ka ) = 1 $ always exists 
for any random point field $ \mu $ (see \cite[Lemma 2.6]{o.dfa}). 
We set $ \ak ^+ (r) = \ak (r+1)$ and $ \mathbf{a}^+ =\{ \ak ^+ \}_{\qq \in \mathbb{N}} $. Let 
\begin{align}
\notag &% \label{:95u}&
\Ssi ^{[m]} = \{ (\xxx ,\mathsf{s})\in \SmSS \, ;\, 
\ulab (\xxx ) + \mathsf{s} \in \Ssi \} 
,\end{align}
where $ \xxx =(x_1,\ldots,x_m)$ and $ \ulab (\xxx ) = \sum_{i=1}^m \delta_{x_i}$. 
Similarly as \eqref{:41Y}--\eqref{:41z} let 
%$ \ak ^+$ such that $ \ak ^+ (r) = \ak (r+1)$ and 
\begin{align} \notag &%\label{:95x}&
\RRprs = \big\{ \mathbf{x}\in \Srm \, ;\,
\min_{j\not=k } |x_j-x_k | \ge 2^{-\pp }
 ,\ 
\inf_{l,i} |x_l-s_i| \ge 2^{-\pp } \big\} 
,\end{align}
where $ j,k,l=1,\ldots,m $, 
$\mathsf{s}=\sum_i \delta_{s_i}$, and $ \Srr ^m = \{ x \in \sS ;\, |x| \le \rr \}^m $, and 
for $ \nn = (\pqr ) $
%Similarly as As \eqref{:41y} and \eqref{:41z}, set 
% 
\begin{align} \notag &%\label{:95y} &
\Han = \Ha _{\pqr} :=\big\{ (\mathbf{x},\mathsf{s}) \in \Ssi ^{[m]} \, ;\, \ 
 \mathbf{x} \in \RRprs ,\ \mathsf{s} \in \Kakk \big\} 
.\end{align}
We set $ \sS _{\rr }^{m,\circ} = \{ |x| < \rr ,\, x \in \sS \}^m $. 
Let $ \RRprCs $ be the open kernel of $ \RRprs $: 
\begin{align} \notag %\label{:97j}&
\RRprCs = \big\{ \mathbf{x} \in \sS _{\rr }^{m,\circ} \, ;\, \ &
\inf_{j\not=k } |x_j-x_k | > 2^{-\pp } ,\ 
\inf_{l,i} |x_l-s_i| > 2^{-\pp } \big\} 
.\end{align}
For $ \nn = (\pqr ) \in \NNNthree $ we set %let $ \HanC $ be such that 
\begin{align} \notag %\label{:95a}
\HanC = \Ha _{\pqr}^{\circ } := \big\{ (\mathbf{x},\mathsf{s}) \in \Ssi ^{[m]} \, ;\, \ &
\mathbf{x} \in \RRprCs , \ \mathsf{s} \in \Kakk \big\} 
.\end{align}
Similarly as \eqref{:41z} and \eqref{:43M}, we set 
$ \Ha $, $ \Ha _{\qq ,\rr } $, $ \Ha ^{\circ } $, and $ \Ha _{\qq ,\rr } ^{\circ } $. 
% we set 
% \begin{align*}&
% \Ha = \bigcup_{\rr =1}^{\infty} \Ha _{\rr } , \quad 
% \Ha _{\rr } = \bigcup_{\qq =1}^{\infty} \Ha _{\qq ,\rr }, \quad 
% \Ha _{\qq ,\rr } = \bigcup_{\pp =1}^{\infty} \Hb \\&
% \Ha ^{\circ } = \bigcup_{\rr =1}^{\infty} \Ha _{\rr } ^{\circ } , \quad 
% \Ha _{\rr } ^{\circ } = \bigcup_{\qq =1}^{\infty} \Ha _{\qq ,\rr } ^{\circ } , \quad 
% \Ha _{\qq ,\rr } ^{\circ } = \bigcup_{\pp =1}^{\infty} \Hb ^{\circ } 
% \end{align*}

\begin{lemma} \label{l:96} 
Assume \As{A1}--\As{A4} and \As{\Bone}. 
For each $ m \in \N $ the following then holds: 
\begin{align} & \label{:97w}
\PPs ( 
\limi{\nn } 
\tau_{\HanC } 
(\mathbf{X}^m,\mathsf{X}^{m*}) = \infty ) = 1 
\quad \text{ for $ \mu $-a.s.\, $ \mathsf{s}$}
.\end{align}
Here $ (\mathbf{X}^m,\mathsf{X}^{m*})$ is the $ m $-labeled process given by 
$ \XB $ in \lref{l:52}, 
$ \tau_{\HanC } $ is the exit time from $ \HanC $, and 
$ \limi{\nn } $ is same as \eqref{:42z}. 
\end{lemma}
\begin{proof} 
The proof is same as \lref{l:41}. Hence we omit it. 
\qed \end{proof}

\Ssection{A sufficient condition of \IFC\ and Yamada-Watanabe theory 
for SDE of random environment type.}\label{s:9B}

We set $ \mathbf{x} = (x_1,\ldots,x_m) \in \Sm $ and 
$ \mathsf{x}^{i\diamondsuit}=\sum_{j\not=i}^m \delta_{x_j}$. 
Let $ (\sigma , b ) $ be as \eqref{:50a}. We set 
 \begin{align} \notag &
 \sigma ^m (\mathbf{x},\mathsf{s}) = (\sigma (x_i, \mathsf{x}^{i\diamondsuit}+\mathsf{s}))_{i=1}^m ,\quad b ^m (\mathbf{x},\mathsf{s}) = ( b (x_i,\mathsf{x}^{i\diamondsuit}+\mathsf{s}))_{i=1}^m
 .\end{align}
Then the time-inhomogeneous coefficients of SDE \eqref{:53m} are given by 
$ \SB |_{\mathsf{s}=\mathsf{X}_t^{m*}} $, where 
$\SB $ is a version of $ (\sigma ^m, b ^m )$ with respect to $ \mum $. 
%We prove the local Lipschitz continuity in $ \mathbf{x} $ for fixed $ \mathsf{s}$. 

Let $ \NNN $ and $ \nn $ be as \eqref{:42u} and \eqref{:42v}. 
Let $ \{ \IIm \}_{\mm \in \N } $ be an increasing sequence of closed sets in $ \SmSS $. 
Let $ \map{\Pitwo }{\SmSS }{\sSS }$ be the projection such that 
$ (\mathbf{x},\mathsf{s}) \mapsto \mathsf{s}$. 
Then 
\begin{align}\notag &% \label{:97i}&
\HmnPc = \{ \mathsf{s}\in \sSS \, ; \, 
\Hmnc \cap \big(\Sm \ts \{ \mathsf{s}\} \big) \not=\emptyset \}
.\end{align}
For $ \nn =(\pqr ) \in \NNNthree $, let $ \Ct \label{;Y3a} (\mm ,\nn ) $ 
be the constant such that $ 0\le \cref{;Y3a} \le \infty $ and that 
\begin{align} \label{:97l} 
\cref{;Y3a} = \sup \{&
\frac{|\fgsigma (\mathbf{x},\mathsf{s}) - \fgsigma (\mathbf{y},\mathsf{s}) |}
{|\mathbf{x}-\mathbf{y}|} 
,\ 
\frac{|\fg (\mathbf{x},\mathsf{s}) - \fg (\mathbf{y},\mathsf{s}) |}
{|\mathbf{x}-\mathbf{y}|}
 ; \, \ 
 \mathbf{x}\not=\mathbf{y}
,\ \\ \notag & 
\mathsf{s} \in \HmnPc ,\, 
\ (\mathbf{x},\mathsf{s}) , (\mathbf{y},\mathsf{s}) \in \RRprCs ,\ 
(\mathbf{x},\mathsf{s}) \simpr (\mathbf{y},\mathsf{s}) 
\} 
.\end{align}
Here $ (\mathbf{x},\mathsf{s}) \simpr (\mathbf{y},\mathsf{s}) $ means 
$ \mathbf{x} $ and $ \mathbf{y} $ are in the same connected component of $ \RRprCs $. 

Let $ \XB = (\lpathX , \mathbf{B})$ be the weak solution of ISDE \ISDE\ given by \lref{l:52} 
defined on $ \OFPFs $. 
%
%For $ \mathsf{A} \subset \SmSS $ we set $ \LLL \mathsf{A} \rRRR := \Pitwo ^{-1} (\Pitwo (\mathsf{A}))$ similarly as \eqref{:42r}. 
%
We set similarly as \eqref{:42d} 
\begin{align}\notag &%\label{:97m}&
\HImnC = \bigcup_{\mathsf{s}\in \Pitwo (\Hmnc )} \RRprCs \ts \{ \mathsf{s} \} 
.\end{align}
We assume the following. 

\medskip 

\noindent 
\As{\Btwo} For each $ m \in \N $, there exist a $ \mum $-version $ \SB $ of $ (\sigma ^m, b ^m )$ 
and an increasing sequence of closed sets $ \{ \IIm \}_{\mm \in \N } $ 
such that for each $ \mm \in \N $ and $ \nn = (\pqr ) \in \NNNthree $
\begin{align}\label{:97n}&
\cref{;Y3a} (\mm ,\nn ) < \infty 
,\\ \label{:97k}&
\limi{\mm } \mathrm{Cap}^{a,\mum } \Big( \HIiC \Big) = 0 
.\end{align}

 \begin{remark}\label{r:97} 
We shall later assume the coefficients are in the domain of the $ m $-labeled Dirichlet form, 
and take $ \SB $ as a quasi-continuous version and 
$\{ \IIm \}_{\mm = 1}^{\infty}$ as a {\em nest}. 
We refer to \cite[pp 67-69]{fot.2} for quasi-continuous version and nest. 
The sequence of sets $ \{ \IIm \} $ plays a crucial role in the proof of \lref{l:U3}. 
\end{remark}
\begin{lemma}	\label{l:97}
Assume \AAEE. Then the following hold. 
\\\thetag{1} For each $ m \in \N $, 
the SDE \eqref{:53f}--\eqref{:53h} has a pathwise unique, weak solution starting at 
$ \mathbf{s}^m=\labm (\mathsf{s}) $ for $ \mu $-a.s.\ $ \mathsf{s}$ 
in the sense that arbitrary solutions $ (\mathbf{Y}^m , \mathbf{B}^m , \mathsf{X}^{m*})$ and 
$ (\hat{\mathbf{Y}}^m , \mathbf{B}^m , \mathsf{X}^{m*})$ 
of \eqref{:53f}--\eqref{:53h} defined on $ \OFPFs $ satisfy 
\begin{align}\notag &%\label{:97o}&
\PPs (\mathbf{Y}^m = \hat{\mathbf{Y}}^m) = 1 
.\end{align}
In particular, $ (\mathbf{Y}^m , \mathbf{B}^m , \mathsf{X}^{m*})$ coincides with 
$ (\mathbf{X}^m , \mathbf{B}^m , \mathsf{X}^{m*})$ under $ \PPs $. 
\\\thetag{2}
Let 
$ (\mathbf{Z}^m , \hat{\mathbf{B}}^m , \hat{\mathsf{X}}^{m*})$ and 
$ (\hat{\mathbf{Z}}^m , \hat{\mathbf{B}}^m , \hat{\mathsf{X}}^{m*})$ 
be weak solutions of the SDE \eqref{:53f}--\eqref{:53h} defined on a filtered space 
 $(\Omega ',\mathcal{F}', P' ,\{ \mathcal{F}_t '\} )$ %$ \OFPF $ 
% starting at 
% $\mathbf{s}^m$ for $ \mu \circ \lab ^{-1}$-a.s.\ $ \mathbf{s}$
% satisfying $ \mathbf{Z}_0^m=\hat{\mathbf{Z}}_0^m$ a.s. and 
satisfying 
\begin{align}
\notag &% \label{:97O}&
 (\hat{\mathbf{Z}}^m , \hat{\mathbf{B}}^m , \hat{\mathsf{X}}^{m*})
\elaw (\mathbf{X}^m , \mathbf{B}^m , \mathsf{X}^{m*})
%(\hat{\mathbf{B}}^m , \hat{\mathsf{X}}^{m*}) \elaw (\mathbf{B}^m , \mathsf{X}^{m*})
.\end{align}
Then it holds that for $ \mu $-a.s.\,$ \mathsf{s}$
\begin{align}
\label{:97oo}&
 P' (\mathbf{Z}^m = \hat{\mathbf{Z}}^m) = 1 
.\end{align}
\thetag{3} Make the same assumptions as \thetag{2} except that the filtrations of 
$ (\mathbf{Z}^m , \hat{\mathbf{B}}^m , \hat{\mathsf{X}}^{m*})$ and 
$ (\hat{\mathbf{Z}}^m , \hat{\mathbf{B}}^m , \hat{\mathsf{X}}^{m*})$, 
$ \{\mathcal{F}_t'\} $ and $ \{ \mathcal{F}_t'' \} $ say, are different 
(but on the same probability space $(\Omega ',\mathcal{F}', P' )$). 
\sigmaconst. 
Then \eqref{:97oo} holds for $ \mu $-a.s.\,$ \mathsf{s}$. 
\end{lemma}

\begin{proof} 
%Other properties of Ginibre interacting Brownian motion used in the proof of \lref{l:97G} 
%follow from \AAE. 
From \As{A1}--\As{A4} we can construct a weak solution 
 $ (\mathbf{X}^m , \mathbf{B}^m , \mathsf{X}^{m*})$ under $ \PPs $ 
of \eqref{:53f}--\eqref{:53h}. 
We remark that \AAEE\ yield the related claims of \pref{l:42}. 
Indeed, \eqref{:97n} corresponds to \eqref{:42a}. 
From \AAEE\ we apply \lref{l:96} to obtain \eqref{:97w}. 
Then from \eqref{:97w} and \eqref{:97k} we easily see 
\begin{align} &\label{:97N}
\PPs ( \limi{\nn } \limi{\mm }\tau_{\HImnC } ( \mathbf{X}^m,\mathsf{X}^{m*}) = \infty ) = 1 
\quad \text{ for $ \mu $-a.s.\, $ \mathsf{s}$}
.\end{align}
Here $ \tau_{\HImnC }(\mathbf{X}^m,\mathsf{X}^{m*}) $ is the exit time 
of $ (\mathbf{X}^m,\mathsf{X}^{m*})$ from the set $ \HImnC $. 
\eqref{:97N} corresponds to \eqref{:42c}. 
The claims in \pref{l:42} are used in the proof of \lref{l:97G}. 
We do not need any other specific properties of the Ginibre interacting Brownian motion 
in the proof of \lref{l:97G}. 
Hence the proof of \lref{l:97} is same as \lref{l:97G}. 
\qed \end{proof}

\begin{remark}\label{r:97b}
In \thetag{3} the reference families $ \{\mathcal{F}_t'\} $ and $ \{ \mathcal{F}_t'' \} $ 
of SDEs are different although the probability space and the Brownian motion is the same. 
The pathwise uniqueness in \thetag{3} is called the pathwise uniqueness in the strict sense in 
\cite[162p]{IW}. 
We shall use this refinement in the proof of \pref{l:98}. 
In the classical situation of Yamada-Watanabe theory, the pathwise uniqueness in the strict sense follows from the pathwise uniqueness and the existence of weak solutions as a corollary of their main result. 
In the current case, it has been not yet succeeded to generalize this part of the Yamada-Watanabe theory 
to SDEs of random environment type. So we add the additional assumption in \thetag{3}. 
\end{remark}

Recall that $ (\mathbf{X},\mathbf{B})$ is the weak solution of \eqref{:50a}--\eqref{:50c} 
 defined on $ \OFPFs $ given by \lref{l:52}. 
Let $ (\mathbf{X}^m , \mathbf{B}^m , \mathsf{X}^{m*})$ be the weak solution 
of \eqref{:53f} made of $ (\mathbf{X},\mathbf{B})$. 
To simplify the notation, we set 
$$ w=(\bx ) \in \WRdzm \ts \WS .$$
Let $\Pw $ be the regular conditional probability such that 
\begin{align}
\notag &% \label{:98a}&
\Pw = \PPs (\mathbf{X}^m \in \cdot \,| \, (\mathbf{B}^m , \mathsf{X}^{m*}) =w) 
.\end{align}

Let $ (\mathbf{Y}^m , \hat{\mathbf{B}}^m , \hat{\mathsf{X}}^{m*})$ be 
an independent copy of $ (\mathbf{X}^m , \mathbf{B}^m , \mathsf{X}^{m*})$. 
Let $ \hat{P}$ be the distributions of $ (\mathbf{Y}^m , \hat{\mathbf{B}}^m , \hat{\mathsf{X}}^{m*})$. 
Let $ \hat{P}_{w}= \hat{P} (\mathbf{Y}^m \in \cdot \,| \, (\hat{\mathbf{B}}^m , \hat{\mathsf{X}}^{m*})=w)$. 
Let $ Q $ be the distribution of $ (\mathbf{B}^m , \mathsf{X}^{m*})$. % under $ \PPs $. 
We set the probability measure $ R $ on 
\begin{align}\notag &%\label{:98b}&
W^{\bullet}:= \WSm \ts \WSm \ts \WRdzm \ts \WS
\end{align}
by
\begin{align}\label{:98c}
&R (d\uuu d\vvv dw) = \Pw ( d\uuu ) \hat{P}_{w}(d\vvv ) Q(dw) 
.\end{align}
We set $ \mathcal{G} $ to be the completion of 
the topological $ \sigma $-field $ \mathcal{B}(W^{\bullet}) $ by $ R $, and 
$ \mathcal{G}_t = \cap_{\varepsilon}( \mathcal{B}_{t+\varepsilon}(W^{\bullet}) \vee \mathcal{N} ) $, 
where 
$ \mathcal{B}_t = \sigma [(\mathbf{u}(u),\mathbf{v}(u),w(u))\, ;\, 0 \le u \le t ] $ and 
$ \mathcal{N} $ is the set of all $ R $-null sets.

\begin{proposition}	\label{l:98}
Assume \AAEE. 
\sigmaconst. 
Then \XlP\ satisfies \As{\iFc}. 
\end{proposition}

\begin{proof}
From \lref{l:97} \thetag{3} we have a pathwise unique, weak solution 
$ (\mathbf{X}^m , \mathbf{B}^m , \mathsf{X}^{m*}) $. 
Then from a generalization of the Yamada--Watanabe theory 
(see Theorem 1.1 in \cite[163p]{IW}) to SDE with random environment, 
we shall construct a strong solution. %We have thus obtained \As{\iFc}. 

Under $ R $, both $\ubx $ and $ \vbx $ are weak solutions of 
\eqref{:53f}--\eqref{:53h}. These solutions are defined on 
$ \Xi = (W^{\bullet},\mathcal{G} , R ,\{ \mathcal{G}_t \} )$ and 
%where the coefficients of \eqref{:53f}--\eqref{:53h} are regarded to be defined on 
%$ \Xi $ in an obvious fashion. 
 SDEs \eqref{:53f} for $ \ubx $ and $\vbx $ become as follows. 
\begin{align}\label{:98f}&
d u_t^{i} = 
 \sigma ^m d\mathbf{b}_t^i + 
\bbbxms (t, (u_t^{i},\mathbf{u}_t^{i\diamondsuit})) dt
,\\\label{:98g}&
d v_t^{i} = \sigma ^m d\mathbf{b}_t^i + 
\bbbxms (t, (v_t^{i},\mathbf{v}_t^{i\diamondsuit})) dt
.\end{align}
Here $ \bbbxms $ is defined by \eqref{:53b} with replacement of 
$ \mathsf{X}_t^{m*}$ by $ \mathsf{x}_t^{m*}=\sum_{j=m+1}^{\infty} \delta_{x_t^j}$. 

Solutions of SDEs \eqref{:98f} and \eqref{:98g} are defined on 
 $ \Xi = (W^{\bullet},\mathcal{G} , R ,\{ \mathcal{G}_t \} )$. 
Although we do not need the fact that $ \mathbf{b} $ 
under $ R $ is a $ \{ \mathcal{G}_t \} $-Brownian motion in the present proof, 
we shall prove this here combined with \lref{l:99} below. 
By construction $ \mathbf{b} $ under $ R $ is a Brownian motion. 
So for this, it only remains to prove $ \mathbf{b}(u)- \mathbf{b}(t)$ 
is independent of $ \mathcal{G}_t $ for all $ t < u $ under $ R $, which we prove in \lref{l:99}. 

Note that the distributions of both $\ubx $ and $ \vbx $ under $ R $ coincide with that of 
$ (\mathbf{X}^m , \mathbf{B}^m , \mathsf{X}^{m*}) $. 
Hence we obtain 
\begin{align}
 \label{:99p}&
 R (\uuu = \vvv ) = 1
\end{align}
from the pathwise uniqueness of weak solutions given by \lref{l:97} \thetag{3}. Let 
\begin{align}\label{:98e}&
 Q _w=
R ((\mathbf{X}^m ,\mathbf{Y}^m )\in \cdot | (\mathbf{B}^m , \mathsf{X}^{m*})=w )
.\end{align}
Then we deduce from \eqref{:98c} 
\begin{align}
\label{:98d}&
Q _w (d\uuu d\vvv )= \Pw ( d\uuu ) \hat{P}_{w}(d\vvv )
.\end{align}
The identity $ R (\uuu = \vvv ) = 1$ in \eqref{:99p} together with \eqref{:98e} 
implies $ Q _w (\uuu = \vvv ) = 1 $ for $ Q $-a.s.\,$ w $. 
Meanwhile, from \eqref{:98d} we deduce that $ \uuu $ and $ \vvv $ under $ Q _w $ are mutually independent. 
Hence the distribution of $ (\uuu ,\vvv )$ under $ Q _w$ is $ \delta _{(F(w),F(w))}$, where 
$ F (w)$ is a nonrandom element of $ \WSm $ depending on $ w $. 
Thus $ F $ is regarded as a function from $ \WRdzm \ts \WS $ to $ \WSm $ by 
$ w \longmapsto F(w)$. 
The distributions of $ \Pw ( d\uuu ) $ and $\hat{P}_{w}(d\vvv ) $ coincide with $ \delta_{F(w)}$. 
We therefore obtain 
$ (\mathbf{X}^m , \mathbf{B}^m , \mathsf{X}^{m*}) = 
(F(\mathbf{B}^m , \mathsf{X}^{m*}),\mathbf{B}^m , \mathsf{X}^{m*})$. 

We easily see that $ F $ is 
$ \overline{\mathcal{B}( \WRdzm \ts \WS )}^{Q} /\mathcal{B}(\WSm ) $-measurable. 
Indeed, letting $ \iota $ and $ \kappa $ be the projections such that 
$ \iota (\uuu ,\vvv ,w) = \uuu $ and $ \kappa (\uuu ,\vvv ,w)= w $, we see 
$ F =\iota \circ \kappa ^{-1}$ $ R $-a.s. 
Then $ \kappa ^{-1} (F^{-1}(A)) = \iota ^{-1}(A) $ $ R $-a.s.\ 
 for each $ A \in \mathcal{B}(\WSm ) $, that is, 
$$ R ( \kappa ^{-1}(F^{-1}(A)) \ominus \iota ^{-1}(A)) = 0 
,$$
where $ \ominus $ denotes the symmetric difference of sets. 
Note that 
$$ \iota ^{-1}(A) \in \mathcal{G}= 
\overline{\mathcal{B}( \WSm \ts \WSm \ts \WRdzm \ts \WS ) }^{R}
.$$
%where the completion is taken over $ R $. 
Then we see $ \kappa (\iota ^{-1}(A)) \in 
\overline{\mathcal{B}(\WRdzm \ts \WS ) }^{Q} $ by $ Q = R \circ \kappa ^{-1}$. 
This implies 
$$ F^{-1}(A) \in \overline{\mathcal{B}( \WRdzm \ts \WS )}^{Q} .$$
Hence $ F $ is 
$ \overline{\mathcal{B}( \WRdzm \ts \WS )}^{Q} $/$\mathcal{B}(\WSm ) $-measurable. 

We can prove that $ F $ is 
$ \overline{\mathcal{B}_t (\WRdzm \ts \WS )}^{Q} $/$\mathcal{B}_t(\WSm ) $-measurable 
for each $ t $ in a similar fashion. 
Here subscript $ t $ denotes the $ \sigma $-field being generated until time $ t $. 
Indeed, we can localize the ISDE and the solution in time to $ [0,T]$ for all $ 0 < T $. 
The restriction of the original weak solution $ (\mathbf{X}^m , \mathbf{B}^m , \mathsf{X}^{m*}) $ 
on the time interval $ [0,T]$ is also a solution of the ISDE. 
So with the same argument we can construct a strong solution $ F_T $ defined on
$ W_{\mathbf{0}}([0,T];\Rdm ) \ts W([0,T];\sSS )$. The solution $ F_T $ is 
$ \overline{\mathcal{B}(W_{\mathbf{0}}([0,T];\Rdm ) \ts W([0,T];\sSS ))}^{Q} $/
$\mathcal{B}(\WT ([0,T]; \Sm )) $-measurable. 
Because of the pathwise uniqueness of weak solutions, 
$ F_T(\mathbf{b},\mathsf{x}) (t)= F (\mathbf{b},\mathsf{x})(t) $ for all $ 0 \le t \le T $. 
We have natural identities 
\begin{align*}
 \mathcal{B}(W_{\mathbf{0}}([0,T];\Rdm ) \ts W([0,T];\sSS ))=
&\mathcal{B}_T (\WRdzm \ts \WS )
\\ 
\mathcal{B}(\WT ([0,T]; \Sm )) = &\mathcal{B}_T(\WSm ) 
.\end{align*}
From this $ F $ is 
$ \overline{\mathcal{B}_T (\WRdzm \ts \WS )}^{Q} $/$\mathcal{B}_T(\WSm ) $-measurable 
for each $ T $. 

Recall that $ \mathbf{s}$ is suppressed from the notation of $ F $. We set 
\begin{align*}&
\Fms (\lBlhatm ) = F (\mathbf{B}^m, \upath (\mathbf{X}^{m*}))
.\end{align*}
Because $ \lpath \circ \upath = \mathrm{id.}$\ for $ \Pm \circ \lpath ^{-1}$-a.s.\ 
and 
$ \upath \circ \lpath = \mathrm{id.}$\ for $ \Pm $-a.s., 
the function $ \Fms $ inherits necessary measurabilities in \dref{d:41} from those of $ F $. 
Hence we see $ \Fms $ is a strong solution of \eqref{:53f}--\eqref{:53h} for $ \XB $ \uPs\ 
for $ \Pm \circ \lab ^{-1}$-a.s.\ $ \mathbf{s}$. 
We have already proved the pathwise uniqueness. 
Thus, $ \Fms $ is a unique strong solution for $ \XB $ starting at $ \mathbf{s}^m$ 
for $ \Pm \circ \lab ^{-1}$-a.s.\ $ \mathbf{s}$. Hence \XlP\ satisfies \As{\iFc}. 
\qed\end{proof}

\begin{lemma} \label{l:99}
For each $ t $, 
%$ \sigma [\mathbf{b}]_t^*=\{ \mathbf{b}(\cdot ) -\mathbf{b}(t)\}_{t < u < \infty } $ 
$ \{\mathbf{b}_u-\mathbf{b}_t \}$ $ (t < u < \infty )$ 
are independent of $ \mathcal{G}_t $ under $ R $. 
\end{lemma}
\begin{proof}
Note that $ \Pw =R(\uuu \in \cdot | w) $ and $\Pw '=R(\vvv \in \cdot | w) $. 
Let $ F $ be the strong solution obtained in the proof of \pref{l:98}. 
Note that $ F $ is $ \overline{\mathcal{B}_t (\WRdzm \ts \WS )}^{Q} $-adapted. 
Then for any $ A_1\ts A_2\ts A_3 \in \mathcal{G}_t $ and $ \theta \in \Rdm $
\begin{align}
\notag %\label{:99c}
E^{R}[e^{\sqrt{-1}\langle \theta , \mathbf{b}_u-\mathbf{b}_t\rangle}1_{A_1\ts A_2 \ts A_3 }]
= &\int_{A_3} 
e^{\sqrt{-1}\langle \theta , \mathbf{b}_u-\mathbf{b}_t\rangle}
\Pw (A_1)\Pw '(A_2) Q (dw)
\\ \notag 
= &\int_{A_3} 
e^{\sqrt{-1}\langle \theta , \mathbf{b}_u-\mathbf{b}_t\rangle}
1_{A_1}(F(w))1_{A_2}(F(w))Q(dw)
\\ \notag =&
e^{-|\theta|^2/2(u-t)} 
\int_{A_3} 1_{A_1}(F(w))1_{A_2}(F(w))Q(dw) 
\\ \notag =&
e^{-|\theta|^2/2(u-t)} R (A_1\ts A_2\ts A_3 )
%\int_{A_3} \Pwt (A_1)\Pwt '(A_2) Q (d\mathbf{b})
.\end{align}
This implies the claim. 
\qed\end{proof}

Recall that $ \Pm $ is given by \eqref{:51u}. 
\begin{theorem} \label{l:UA} 
Assume that $ \mu $ and $ \lpathX $ under $ \Pm $ satisfy \As{\muTT}, \AAcEE. \sigmaconst. 
Then ISDE \ISDEb has a family of unique strong solutions $ \{ \Fs \} $ 
starting at $ \mathbf{s}=\lab (\mathsf{s})$ for $ \mu $-a.s. $ \mathsf{s}$ 
%under the constraints such that $ \{ \Fs  \} $ satisfies \As{$\mathbf{MF}$} and that $ \pPF $ satisfies \IASN. 
under the constraints of \As{$\mathbf{MF}$}, \IASN. 
\end{theorem}

\begin{proof} 
We take $ \XB $ and $ P $ in \tref{l:5A} as $ \mathbf{X}:=\lpathX $ and $ P := \Pm $. 
We check this satisfies all the assumptions of \tref{l:5A}. 
That is, we have to check \As{\muTT} for $ \mu $ and that 
$ \XB $ under $ P $ is a weak solution satisfying \IASN. 

We see 
\As{\muTT} for $ \mu $ follows by assumption. 
By \lref{l:92} $ \lpathX $ under $ \Pm $ satisfies \As{\sIn}. 
Hence by \lref{l:52} $ \lpathX $ under $ \Pm $ is a weak solution. 
Note that the coefficient $ \sigma^m $ is constant and that \AAEE\ hold by assumption. 
Then the assumptions of \pref{l:98} are fulfilled. 
Hence we obtain \As{\iFc} from \pref{l:98}. 
Using \lref{l:94} we obtain \As{\nbj}. 
Because $ \Pm $ is $ \mu $-reversible, \As{\muAC} for $ \mu $ is obvious. 
Thus, all the assumptions of \tref{l:5A} are fulfilled. 
Hence the claim follows from \tref{l:5A}. 
\qed \end{proof}

\begin{corollary}\label{l:UAA}
Under the same assumptions of \tref{l:UA}, 
$ \lpathX $ under $ \Pm $ is a weak solution of \ISDEb satisfying \IASN. 
Furthermore, $ \lpathX $ under $ \{\PPs \} $ is a family of unique strong solutions starting at 
$ \mathbf{s}=\lab (\mathsf{s})$ for $ \mu $-a.s.,$ \mathsf{s}$ under the constraints of \As{$\mathbf{MF}$},  \IASN. 
\end{corollary}

\Ssection{A sufficient condition of \As{\Btwo} and Taylor expansion of coefficients.}
 \label{s:9C}

In this section we give a sufficient condition of \As{\Btwo}. 
We begin by introducing the cut-off functions on $ \SmSS $. 
Let $ \vq \in C_0^{\infty}(\Sm ) $ be the cut-off function such that 
\begin{align*}&
\text{$ 0 \le \vq (\xxx ) \le 1 , \quad |\nabla \vq (\xxx )| \le 2 , \quad 
\vq (\xxx ) = \tilde{\varphi }_{\rr }(|\xxx |) $ \quad for all $ \xxx \in \Sm $}
,\end{align*}
where $ \tilde{\varphi }_{\rr }\in C_0^{\infty}(\R )$ is given by \eqref{:W6z}. 
Let 
$ \hp \in C_0^{\infty}(\R )$ $ (\pp \in \{ 0 \}\cup \N )$ such that 
\begin{align} \notag 
&
 \hp (t) = 
\begin{cases} 
 1 & (t \le 2^{-\pp -2}) ,\\
 %1 - 2^{\pp +1}t & (0 \le t \le 2^{-\pp -1}) , \\
 0 & (2^{-\pp -1} \le t )
\end{cases}
,\quad 
0 \le \hp (t) \le 1 , \quad |\hp '(t)| \le 2^{\pp + 3} 
\quad \text{ for all }t
.\end{align}
We write $ \xxx =(x_k)_{k=1}^m \in \Sm $ and $ \mathsf{s} = \sum_i \delta_{s_i}\in \sSS $. 
Let $ \map{h_{\pp }^{\dagger}}{\Sm \ts \sSS }{\mathbb{R}}$ such that 
\begin{align} \notag 
&\quad \quad \quad 
 h_{\pp }^{\dagger} ( \xxx ,\mathsf{s}) = 
\prod_{k=1}^{m}
\Big\{\prod_{j\not=k}^m \{ 1- \hp (|x_k-x_j|)\}\Big\} 
\Big\{\prod_{i}\{ 1- \hp (|x_k-s_i|)\}\Big\} 
.\end{align}
We label $ \mathsf{s}=\sum_i \delta_{s_i}$ in such a way that $ |s_i| \le |s_{{i+1}}|$ for all $ i $. 
We set 
$$ I(\qq ,r ,\mathsf{s}) = \{ i \, ;\, i > \ak (r) ,\ s_i \in \Sr \} .$$
Here $ \ak = \{ \ak (r) \}_{r\in\mathbb{N}} $ are the increasing sequences in \eqref{:95w}. 
We set 
\begin{align}\notag &
 d_{\ak } (\mathsf{s}) = 
 \{ \sum_{r=1}^{\infty} \sum_{i\in I(\qq ,r ,\mathsf{s}) } (r-|s_i|)^2 \}^{1/2}
,\\ \notag &
 \chi _{\ak } (\mathsf{s}) = h _{0} \circ d_{\ak } (\mathsf{s}) 
. \end{align}

We introduce the cut-off functions defined as 
\begin{align} \label{:U2u}
\begin{cases}
\chi _{\rr } (\xxx ,\mathsf{s}) = \vq (\xxx ) 
,\\ 
 \chi _{\qq ,\rr } (\xxx ,\mathsf{s}) = \vq (\xxx ) \chi _{\ak } (\mathsf{s}) 
,\\ 
\chi _{\pqr } (\xxx ,\mathsf{s}) = \vq (\xxx ) \chi _{\ak } (\mathsf{s}) 
 h_{\pp }^{\dagger} (\xxx ,\mathsf{s}) 
\end{cases}
.\end{align}
We easily see that 
\begin{align}\notag &% \label{:U2v} &
\limi{\rr } \limi{\qq } \limi{\pp } 
\chi _{\pqr } (\xxx ,\mathsf{s}) = 1 
\quad \text{ for all $ (\xxx ,\mathsf{s})\in \Ha $}
.\end{align}
Let $ \NNN = \{ (\pqr ), (\rr ,\qq ), \qq \ ;\, \pqr \in \N \} $, 
$ \NNN _1$, $ \NNN _2 $, and $ \NNNthree $ be as \eqref{:42u}. 
Recall that for $ \nn \in \NNN $ we define $ \nn + 1 \in \NNN $ such that 
\begin{align}
\notag &% \label{:42vU}&
 \nn + 1 = 
\begin{cases}
(\pp + 1, \qq , \rr ) &\text{ for $ \nn = (\pqr )$, }\\
(\qq +1, \rr )&\text{ for $ \nn = (\qq , \rr ) $, }\\
 \rr +1 &\text{ for $ \nn = \rr $. }
\end{cases}
\end{align}
We set $ \chin = \chi _{\pqr } $. Then $ \{ \chin \} $ are consistent in the sense that 
$ \chin (\xxx ,\mathsf{s}) = \chinn (\xxx ,\mathsf{s}) $ for $ (\xxx ,\mathsf{s}) \in \Han $. %
We suppress $ m $ from the notation of $ \chin $ although $ \chin$ depends on $ m \in \N $. 
By a direct calculation similar to that in \cite[Lemma 2.5]{o.dfa}, we obtain the following. 
\begin{lemma} \label{l:U1} 
For each $ m \in \N $, the functions $ \chin $ ($ \nnNNN $) satisfy the following. 
\begin{align}
\notag &% \label{:U1c} &
\chin (\xxx ,\mathsf{s}) = 
\begin{cases}
0 & \text{ for } (\xxx ,\mathsf{s}) \not\in \Hann 
\\
1& \text{ for } (\xxx ,\mathsf{s}) \in \Han 
\end{cases}
,\quad 
 \chin \in \dom ^{a ,\mum } 
,\\ \notag %
& 
0 \le \chin (\xxx ,\mathsf{s}) \le 1 ,\quad |\nabla_\xxx \chin (\xxx ,\mathsf{s}) |^2 \le \cref{;83a} 
,\quad 
\DDD [\chin ,\chin ] (\xxx ,\mathsf{s}) \le \cref{;83b}
.\end{align}
Here 
$ \Ct =\cref{;83a}(\nn ) \label{;83a}$ and 
$ \Ct =\cref{;83b}(\nn )\label{;83b}$ 
are positive finite constants independent of $ (\xxx ,\mathsf{s}) $, and 
$ \dom ^{a ,\mum }$ is the domain of the Dirichlet form of 
the $ m $-labeled process $ (\mathbf{X}^m, \mathsf{X}^{m*})$. 
\end{lemma}

We shall give a sufficient condition of $ \As{\Btwo}$ using Taylor expansion. Let 
\begin{align}
%\notag &% 
\label{:U2z}&
\mathbf{J}^{[l]} = \{ \mathbf{j}= (j_{k,i})_{k=1,\ldots,m,\, i=1,\ldots,d} 
;\ j_{k,i} \in \{ 0 \} \cup \N ,\, \sum_{k=1}^m \sum_{i=1}^{d} j_{k,i}=l \} 
.
\end{align}
We set 
%$ \partial _{\mathbf{j}} = \mathrm{id.} $ for $ \mathbf{j}=(0,\ldots,0) \in \mathbf{J}^{[0]}$ 
%and 
$ \partial _{\mathbf{j}} = \prod_{k,i} (\partial /\partial x_{k,i} )^{j_{k,i}} $ 
%$ ({\partial^{j_1}}/{\partial x_{1}},\ldots,{\partial^{j_m}}/{\partial x_{m}})$ 
for 
$ \mathbf{j}=( j_{k,i}) \in \mathbf{J}^{[l]} $, where $ ((x_{k,i})_{i=1}^d)_{k=1}^m \in \Rdm $ and 
$ ({\partial }/{\partial x_{k,i}})^{j_{k,i}}$ denotes the identity if $ j_{k,i}=0 $. 
For $ \ellell \in \N $, we introduce the following: 

\smallskip
\noindent 
\As{\CONE} 
For each $ \mathbf{j} \in \cup_{l=0}^{\ellell } \mathbf{J}^{[l]} $ and $ \nn \in \NNNthree $, 
\begin{align}
\notag &%\label{:}&
\text{$ \chin \partial _{\mathbf{j}} \sigma ^m, \chin \partial _{\mathbf{j}} b ^m \in \dom ^{a ,\mum } $}
.\end{align}
%\medskip 
\noindent 
\As{\CTWO} 
There exists a $ \mum $-version $ \SBhat $ of $ \{ \sigma ^m ,b ^m \} $ such that 
\begin{align}
\notag &% \label{:U2c} & 
\sup \big\{ |\partial _{\mathbf{j}}\ff (\xxx ,\mathsf{s})| ; \, 
(\xxx ,\mathsf{s}) \in \Han ,\ \ff \in \SBhat \big\} < \infty 
\end{align}
for each $ \mathbf{j} \in \mathbf{J}^{[\ellell ]}$ and $ \nn \in \NNNthree $.

\begin{remark}\label{r:U2} 
Note that $ \vq (\xxx ) $ and $ h_{\pp }^{\dagger} (\xxx ,\mathsf{s}) $ in \eqref{:U2u} 
are smooth in $ \mathbf{x}$, and that $\partial _{\mathbf{i}} \vq (\xxx ) $ and 
$ \partial _{\mathbf{i}} h_{\pp }^{\dagger} (\xxx ,\mathsf{s}) $ are bounded. 
Hence \eqref{:U2u} and \As{\CONE} with a straightforward calculation show 
$ (\partial _{\mathbf{i}}\chin)( \partial _{\mathbf{j}} \sigma ^m ) $ and 
$(\partial _{\mathbf{i}}\chin )(\partial _{\mathbf{j}} b ^m ) $ belong to $ \dom ^{a ,\mum } $. 
\end{remark}
\begin{proposition} \label{l:U2}
Assume that \As{\CONE} and \As{\CTWO} hold for some $ \ellell \in \N $. 
Then \As{\Btwo} holds. 
\end{proposition}
\begin{proof} 
We omit the proof of \pref{l:U2} because it is same as \pref{l:42}. 
Indeed, in the proof of \pref{l:42}, we need \lref{l:41}, \lref{l:W7}, and \eqref{:!2M}. 
These correspond to \lref{l:96}, \As{\CONE}, and \As{\CTWO}, respectively. 
\qed \end{proof}

\begin{theorem} \label{l:UB} 
Under the same assumptions as \tref{l:UA} with replacement \As{\Btwo} by \As{\CONE}--\As{\CTWO}, 
the same conclusions as \tref{l:UA} holds. 
\end{theorem}
\begin{proof} 
From \pref{l:U2} and \tref{l:UA} we obtain the claim. 
\qed \end{proof}

\begin{corollary}\label{l:UBB}
Under the same assumptions of \tref{l:UB}, the same conclusions as \corref{l:UAA} hold. 
% 
% $ \lpathX $ under $ \Pm $ is a strong solution of \ISDEb satisfying \IASN. 
% Furthermore, 
% $ \lpathX $ under $ \{\PPs \} $ is a family of unique strong solutions starting at 
% $ \mathbf{s}=\lab (\mathsf{s})$ for $ \mu $-a.s.\,$ \mathsf{s}$ under the constraints such that $ \Pm \circ  \lpathX ^{-1}$ satisfies \IASN. 
\end{corollary}

\begin{remark}\label{r:53q}
For sine$_{2} $, Airy$ _2$, and Bessel$ _{2,\alpha }$ random point fields, there is another construction of stochastic dynamics based on space-time correlation functions \cite{KT07b,p-spohn}. 
\tref{l:UB} combined with tail triviality obtained in \cite{bqs,ly.18,o-o.tail} proves that these two dynamics are the same \cite{o-t.airy,o-t.core,o-t.sm}. 
\end{remark}

%%%%%%%%%%%%%%%%%%%%%%%%%%%%%%%%%%%%%%%%%%%%%%%%%%%%%%

\Section{Sufficient conditions of \NaN\ for $ \mu $ with non-trivial tails.} \label{s:Q}

In this section, we deduce \NaN\ from the assumptions of $ \mu $. 
We shall show the stability of \NaN\ under the operation of 
conditioning with respect to the tail $ \sigma $-field $ \Tail (\sSS ) $. 
Recalling the decomposition \eqref{:54b} of $ \mu $ we have 
\begin{align}\label{:Q1a} &
\mu (\A ) = \int _{\sSS } \mut (\A ) \mu (d \aa )
. \end{align}

\begin{lemma}\label{l:Q1}
Assume \As{A1}--\As{A3} for $ \mu $. 
Then $ \mut $ satisfies \As{A1}--\As{A3} for $\mu $-a.s.\! $\aa $. 
\end{lemma}
\begin{proof}
\As{A1}--\As{A2} for $ \mut $ 
are clear from the definitions of logarithmic derivatives and quasi-Gibbs measures 
combined with Fubini's theorem, respectively. 
 \As{A3} for $ \mut $ follows from \As{A3} for $ \mu $ and Fubini's theorem. 
Indeed, we have 
\begin{align} \notag 
\int_{\sSS } \Big[ \sum_{k=m}^{\infty} \frac{k!}{(k-m)!}
 \, \mut (\sSS _r^k) \Big] \mu (d \aa )& = 
\sum_{k=m}^{\infty} \frac{k!}{(k-m)!} \int_{\sSS } \mut (\sSS _r^k) \mu (d \aa ) 
\\ \notag &
= 
\sum_{k=m}^{\infty} \frac{k!}{(k-m)!} \, \mu (\sSS _r^k)< \infty 
.\end{align}
Hence we see that $ \mut $ satisfies \eqref{:51h} for $ \mu $-a.s.\! $ \aa $, 
which implies \As{A3}. 
\qed\end{proof}

Let $ (\E ^{\amut }, \dom ^{\amut })$ be the Dirichlet form given by \eqref{:51i}--\eqref{:51j} 
with replacement of $ \mu $ by $ \mut $. 
By \lref{l:51} and \lref{l:Q1} $ (\E ^{\amut }, \dom ^{\amut })$ 
are quasi-regular Dirichlet forms on $ L^2(\sSS ,\mut )$ for $ \mu $-a.s.\ $ \aa $. 
We easily see that 
\begin{align}
\label{:Q1c}&
\text{$ \dom ^{\amut } \supset \Damu $ for $ \mu $-a.s.\! $ \aa $}
.\end{align}
Let $ \mathrm{Cap}^{a,\mu } $ and $ \mathrm{Cap}^{a,\mut } $ 
be the capacities associated with $ (\Eamu ,\Damu ) $ on $ \Lmu $ and 
$ (\E ^{\amut }, \dom ^{\amut })$ on $ L^2(\sSS ,\mut )$, respectively. 
Then by the variational formula of capacity and \eqref{:Q1c}, 
we easily deduce that for each $ A $ 
\begin{align*}&
\mathrm{Cap}^{a,\mut } ( A ) \le \mathrm{Cap}^{a,\mu } ( A ) 
\quad \text{ for $ \mu $-a.s.\,$ \mathsf{a}$}
.\end{align*}
This implies 
\begin{align}\label{:Q1d}&
\int _{\sSS } \mathrm{Cap}^{a,\mut } ( A ) \mu (d\aa ) \le \mathrm{Cap}^{a,\mu } ( A )
.\end{align}

\begin{lemma}\label{l:Q2}
Assume \As{A4} for $ \mu $. Then $ \mut $ satisfies \As{A4} for $\mu $-a.s.\! $\aa $. 
\end{lemma}
\begin{proof}
We begin by proving \As{A4} \thetag{1} for $ \mut $. 
By \eqref{:54k} for $ \mu $ and \eqref{:Q1a} we have 
\begin{align}
\notag &% \label{:Q2b}&
\int_{\sSS } E^{\mut }\Big[\langle \ER ( \frac{|\cdot |-r}{T}) , \mathsf{s}\rangle \Big] \mu (d\mathsf{a}) = 
E^{\mu }\Big[\langle \ER ( \frac{|\cdot |-r}{T}) , \mathsf{s}\rangle \Big] < \infty 
.\end{align}
Then $ \mut $ satisfies \eqref{:54k} by Fubini's theorem. 
By Fatou's lemma, \eqref{:54l}, and \eqref{:Q1a}% we have
\begin{align}\label{:Q2c}&
\int_{\sSS } \liminfi{r} \ER (\frac{r}{ T }) \, 
E^{\mut }[\langle 1_{\sS _{r+R}}, \mathsf{s}\rangle] \mu (d\aa )
\\ \notag \le & 
 \liminfi{r} \int_{\sSS } \ER (\frac{r}{ T }) \, 
E^{\mut }[\langle 1_{\sS _{r+R}}, \mathsf{s}\rangle] \mu (d\aa )
&& \text{ by Fatou's lemma}
\\ \notag = & 
 \liminfi{r} \ER (\frac{r}{ T }) \, 
 \int_{\sSS } 
E^{\mut }[\langle 1_{\sS _{r+R}}, \mathsf{s}\rangle] \mu (d\aa )
\\ \notag = &
\liminf_{r\to \infty} \ER (\frac{r}{ T }) \, 
E^{\mu }[\langle 1_{\sS _{r+R}}, \mathsf{s}\rangle] =0 
&& \text{ by \eqref{:Q1a}, \eqref{:54l}}
.\end{align} 
Hence from \eqref{:Q2c} we see 
$ \liminfi{r} \ER (\frac{r}{ T }) \, 
E^{\mut }[\langle 1_{\sS _{r+R}}, \mathsf{s}\rangle] = 0 $ for $ \mu $-a.s.\ $ \aa $. 
% \begin{align}
% \notag &\liminfi{r} \ER (\frac{r}{ T }) \, E^{\mut }[\langle 1_{\sS _{r+R}}, \mathsf{s}\rangle] = 0
% .\end{align}
This implies \eqref{:54l} for $ \mut $. 
We have thus obtained \As{A4} \thetag{1} for $ \mut $.

From \eqref{:Q1d} and \As{A4} \thetag{2} for $ \mu $, we deduce 
\begin{align}\notag &%\label{:Q2e}&
\int _{\sS } \mathrm{Cap}^{a,\mut } (\{ \mathsf{s} ; \mathsf{s} (\partial \sS ) \ge 1 \}) \mu (d\aa ) 
\le \mathrm{Cap}^{a,\mu }(\{ \mathsf{s} ; \mathsf{s} (\partial \sS ) \ge 1 \}) 
= 0 
\\ \notag &
\int _{\sS } \mathrm{Cap}^{a,\mut } (\mathsf{S}_{\mathrm{s}}^c) \mu (d\aa ) 
\le 
\mathrm{Cap}^{\amu } (\mathsf{S}_{\mathrm{s}}^c)
 = 0
.\end{align}
Hence we obtain \As{A4} \thetag{2} for $ \mut $ for $ \mu $-a.s. $ \aa $. This completes the proof. 
\qed\end{proof}

We can regard $ (\E ^{\amut }, \mathcal{D}^{\amut }) $ as a Dirichlet form on $ L^2( \HHa ,\mut ) $ 
from \eqref{:54d}--\eqref{:54g}. 
Let $ \PPsa $ be the distribution of the unlabeled diffusion starting at $ \mathsf{s}$ 
associated with the Dirichlet form 
$ (\E ^{\amut }, \mathcal{D}^{\amut }) $ on $ L^2( \HHa ,\mut ) $ given by \lref{l:51}. Let 
\begin{align}
\label{:Q2a}&
 \PPa (\cdot ) = \int _{\HHa } \PPsa (\cdot ) \mut (d\mathsf{s})
.\end{align}
Then $ \PPa $ is a $ \mut $-stationary diffusion on $ \HHa $.

\begin{lemma}\label{l:Q3}
Assume \As{A1}--\As{A4} for $ \mu $. 
Then the diffusion $ \PPa $ satisfies \As{\sIn} and \As{\nbj} for $\mu $-a.s.\! $\aa $. 
Furthermore, for $\mu $-a.s.\! $\aa $ 
\begin{align}\label{:Q3a}&
\mathrm{Cap}^{a,\mut } (\SSsde ^c) = 0 ,
\quad \PPa \circ \lpath ^{-1}(\WT (\SSSsde )) = 1 
.\end{align}
\end{lemma}
\begin{proof} 
By \lref{l:Q1} and \lref{l:Q2} $ \mut $ satisfy \As{A1}--\As{A4} for $\mu $-a.s.\! $\aa $, which are 
the assumptions of \lref{l:92} and \lref{l:94}. 
Then $ \mut $ and $ \PPa $ satisfy \As{\sIn}, \eqref{:Q3a}, and \As{\nbj} 
for $\mu $-a.s.\! $\aa $ by \lref{l:92} and \lref{l:94}. 
\qed \end{proof}

\begin{lemma}\label{l:Q4}
Assume that $ \mu $ and $ \{ \PPs \} $ satisfy \AAEE. \sigmaconst. 
Then $ \mathbf{X}=\lpathX $ under $ \PPa $ satisfies \As{\iFc} for $\mu $-a.s.\! $\aa $. 
 \end{lemma}
\begin{proof}
We use \pref{l:98} to prove \lref{l:Q4}. 
Our task is to check $ \mut $ and $ \{\PPa _{\mathsf{s}} \} $ fulfill 
the assumptions of \pref{l:98}: namely, \AAEE. 

We see $ \mut $ satisfies \As{A1}--\As{A4} for $\mu $-a.s.\! $\aa $ by \lref{l:Q1} and \lref{l:Q2}. 
From \As{\Bone} for $ \mu $ and the tail decomposition \eqref{:54b} of $ \mu $ 
combined with Fubini's theorem, we obtain \As{\Bone} for $ \mut $ for $\mu $-a.s.\! $\aa $. 

We next proceed with \As{\Btwo} for $ \mut $ and $ \PPa $. 
Let $ \IIm $ be the increasing sequence of closed sets in \As{\Btwo} for $ \mu $ and $ \{ \PPs \} $. 
Then \eqref{:97n} is satisfied. 
Indeed, the condition \eqref{:97n} is independent of the measure, so 
\eqref{:97n} for $ \mu $ implies \eqref{:97n} for $ \mut $.

Let $ \mu _{\mathrm{Tail}}^{\aa ,[m] }$ be the $ m $-Campbell measure of $ \mut $. 
%We set $ \IIi = \cup_{\mm = 1}^{\infty} \IIm $. 
%
From the analogy of \eqref{:Q1d} for the $ m $-labeled Dirichlet form, we deduce 
for all $ \nn \in \NNNthree $ and $ \mm \in \N $ 
\begin{align}\label{:83a}& 
\int \mathrm{Cap}^{a,\mu _{\mathrm{Tail}}^{\aa ,[m] } } 
\Big( \HIiC \Big) \mu (d\aa )
\le 
\mathrm{Cap}^{a,\mum } \Big( \HIiC \Big) 
.\end{align}
By \eqref{:97k} for all $ \nn \in \NNNthree $ 
\begin{align} \label{:83A}&
\limi{\mm } \mathrm{Cap}^{a,\mum } \Big( \HIiC \Big) = 0
.\end{align}
From \eqref{:83a}, \eqref{:83A}, and Fatou's lemma 
we obtain for $ \mu $-a.s.\! $ \mathsf{a}$ and for all $ \nn \in \NNNthree $
\begin{align} \notag &%\label{:83b}& 
\limi{\mm } \mathrm{Cap}^{a,\mu _{\mathrm{Tail}}^{\aa ,[m] } } 
 \Big( \HIiC \Big) = 0
.\end{align}
We hence obtain \eqref{:97k} for $ \mut $. 
Then we have verified \As{\Btwo} for $ \mut $ for $ \mu $-a.s.\! $ \mathsf{a}$. 

We thus see that all the assumptions of \pref{l:98} are fulfilled for $ \mut $ and $ \PPa $, 
which deduces that $ \PPa $ satisfies \As{\iFc} for $\mu $-a.s.\! $\aa $. 
\qed
\end{proof}

The next theorem claims quenched results follows from anneal assumptions. 
\begin{theorem} \label{l:QA} 
Assume that $ \mu $ and $ \{ \PPs \} $ satisfy \AAEE. 
Let $ \PPa $ as \eqref{:Q2a}. 
\sigmaconst. %[]
Then, for $ \mu $-a.s.\! $ \aa $, \ISDEb has a unique strong solution 
starting at $ \mathbf{s}=\lab (\mathsf{s})$ for $ \mut $-a.s.\! $ \mathsf{s}$ 
under the constraints of \As{$\mathbf{MF}$}, \IASNaa. 
\end{theorem}
\begin{proof} 
We use \tref{l:5A} for the proof. 
We take $ \XB $ and $ P $ in \tref{l:5A} as $ \mathbf{X}:=\lpathX $ and 
$ P :=\PPa $ for $ \mu $-a.s.\ $ \aa $. 
We check this choice satisfies all the assumptions of \tref{l:5A}, namely, 
$ \lpathX $ under $ \PPa $ is a weak solution satisfying {\IASNaa}, 
and that $ \mut $ satisfies \As{\muTT} for $ \mu $-a.s.\ $ \aa $. 

Below we omit \lq\lq for $ \mu $-a.s.\ $ \aa $''. 
From \lref{l:52}, \lref{l:Q1}, and \lref{l:Q3}, 
we see that $ \lpathX $ under $ \PPa $ is a weak solution of \ISDEb. 

 \As{\muTT} for $ \mut $ follows from \eqref{:54c}. 
Obviously, $ \lpathX $ under $ \PPa $ satisfies 
\As{\muAC} for $ \mut $ because $ \PPa $ is $ \mut $-reversible. 
\As{\sIn}, \As{\nbj}, and \As{\iFc} for $ \lpathX $ under $ \PPa $ follow from \lref{l:Q3} and \lref{l:Q4}. 
%
%\As{\iFc} for $ \lpathX $ under $ \PPa $ follows from \lref{l:Q4}. 

Thus, all the assumptions of \tref{l:5A} are fulfilled. %, which completes the proof. 
\qed \end{proof}

\begin{theorem}	\label{l:QB}
Under the same assumptions of \tref{l:QA} with replacement \As{\Btwo} by \As{\CONE}--\As{\CTWO} 
the same conclusions as \tref{l:QA} hold. 
\end{theorem}
\begin{proof} 
From \pref{l:U2}, \As{\Btwo} holds. Then from \tref{l:QA} we deduce the claim. 
\qed \end{proof}

\begin{corollary} \label{l:QC} 
Make the same assumptions of \tref{l:QA} or \tref{l:QB}. 
Then $ \lpathX $ under $ \PPsa $ is a unique strong solution of 
 \ISDEb starting at $ \mathbf{s}=\lab (\mathsf{s})$ 
under the same constraints as \tref{l:QA}. 
%of \As{\iFc}, \nNm\ for $ \mut $. 
\end{corollary}

\begin{remark}\label{r:54} \thetag{1}
We have two diffusions $ \{ (\mathsf{X}, \Pssf )\}_{\mathsf{s} \in \mathsf{S} } $ and 
$ \{\{ (\mathsf{X}, \Pssf ^{\aa })\}_{\mathsf{s} \in \HHa } \}_{\aa \in \mathsf{H}}$. 
The former is deduced from \As{A2} and \As{A3} for $ \mu $, and 
the associated Dirichlet form is given by 
$ (\E ^{\mu } , \mathcal{D} ^{\mu } ) $ 
on $ L^2( \sSS ,\mu ) $. 
The latter is deduced from \As{A2} and \As{A3} for $ \mut $, 
and is a collection of diffusions whose Dirichlet forms are 
$ (\E ^{\amut }, \mathcal{D}^{\amut }) $ on $ L^2( \HHa ,\mut ) $. %
These two diffusions are the same (up to quasi-everywhere starting points) 
when $ \mu $ is tail trivial. 
We emphasize that \tref{l:QA} does not follow from \tref{l:5B} because 
we do not know how to prove $ \{ (\mathsf{X}, \Pssf )\}_{\mathsf{s} \in \mathsf{S} } $ 
satisfies $ \As{\muAC}$ for $ \mut $. 
\\\thetag{2} In \tref{l:QA}, we have constructed the tail preserving solution. 
Its uniqueness however is under the constraint of \As{\muAC} for $ \mut $. 
Hence it does not exclude the possibility that 
there exists a family of solutions $ \mathbf{X}$ not satisfying this condition, 
and, in particular, a family of solutions whose distributions 
on the tail $ \sigma $-field $ \Tail (\sSS ) $ changing as $ t $ varies. 
We refer to \cite{kawa} for a tail preserving property of the Dyson model. 
\end{remark}

\section{Examples.}\label{s:V}
We devote to examples and applications of Theorems \ref{l:UA}, \ref{l:QA}, and \ref{l:QB}. 
Throughout this section, $ b (x,\mathsf{y}) = \frac{1}{2} \dmu (x,\mathsf{y})$, 
where $ \dmu $ is the logarithmic derivative of random point field $ \mu $ associated with ISDE. 

We present a simple sufficient condition of \As{\CTWO} introduced before \pref{l:U2}. 
We write $ \mathbf{x}=(x_1,\ldots,x_m) \in \Rdm $ and $ \mathsf{s}=\sum_i\delta_{s_i}$ as before. 
Let $ \mathbf{J}^{[\ellell ]} $ be as \eqref{:U2z} with $ l = \ell $. 
Recall that we suppress $ m $ from the notation. 
\begin{lemma} \label{l:U3} 
Assume that for each $ \ff \in \SBhat $, $ \mathbf{j} \in \mathbf{J}^{[\ellell ]}$, and 
 $ k \in \{ 1,\ldots,m \} $, there exists $ \g _{\mathbf{j},k} $ such that 
$ \g _{\mathbf{j},k} \in C (\sS ^2 \backslash \{ x = s \}) $ and that 
\begin{align} 
\notag &
\partial _{\mathbf{j}} \ff (\xxx ,\mathsf{s}) = 
\sum_{k=1}^m \big\{
 \sum_{p\not= k}^m \g _{\mathbf{j},k} (x_k , x_p ) 
+
 \sum_{i} \g _{\mathbf{j},k} (x_k , s_i ) 
\Big\}
\ \text{ for $ (\xxx ,\mathsf{s}) \in \Ha $}
.\end{align}
Assume \As{\Bone}. Assume that there exist positive constants 
$ \Ct \label{;88d} $ and $ \Ct \label{;88c} $ such that 
\begin{align} \label{:U3b}& 
\limsup_{r\to\infty} \frac{a_{\qq }(r)}{r^{\cref{;88d}}} < \infty ,\quad 
\sum_{r=1}^{\infty} \frac{a_{\qq }(r)}{r^{ \cref{;88d}+1 }} < \infty 
,\\ & \label{:U3c}
 | \g _{\mathbf{j},k} (x ,s ) | \le \frac{\cref{;88c} }{ (1+|s|)^{ \cref{;88d} }} 
\quad \text{ for all } (x ,s) \in H_{\pp ,\rr } 
.\end{align}
Here $ H_{\pp ,\rr } = \{ (x ,s)\in \sS ^2 ;\, 2^{-\pp } \le |x -s| ,\, x \in \sS _{\rr }\} $ for $ \pp ,\rr \in \N $. 
We then obtain \As{\CTWO}. 
\end{lemma}
\begin{proof} 
From \eqref{:U3c}, we deduce that for $ (x,\mathsf{s}) \in \Han $ 
\begin{align} \notag 
 \sum_{ i } | \g _{\mathbf{j},k} (x_k,s_i )| 
\le & \cref{;88c} \Big\{
\Big( 
 \sum_{r=1}^{\infty} \sum_{ s_i\in \Sr \backslash \sS _{r-1}}
 \frac{1 }{ (1+|s_i|)^{ \cref{;88d} }}\Big) 
+ \mathsf{s}(\sS _0 )
\Big\}
\\\notag \le & 
\cref{;88c} \Big\{ \Big(
\sum_{r=1}^{\infty} \sum_{ s_i\in \Sr \backslash \sS _{r-1}}
 \frac{1}{ (1+r-1)^{ \cref{;88d} }}\Big) 
+ \mathsf{s}(\sS _0 )\Big\}
 \\ \notag 
 = &\cref{;88c} \limsup_{R\to\infty } \Big\{ 
 \frac{\mathsf{s} (\sS _R ) }{R^{\cref{;88d}}}
+ 
 \sum_{r=2}^{R} \mathsf{s}(\sS _{r-1})
\Big\{ \frac{1}{(r-1)^{\cref{;88d}}}- 
 \frac{ 1 }{ r^{ \cref{;88d}}} \Big\} + \mathsf{s}(\sS _0 ) \Big\} 
.\end{align}
The last line is finite by \eqref{:95w} and \eqref{:U3b}. This yields \As{\CTWO}. 
\qed \end{proof}

In the rest of this section, $ \sigma $ is the unit matrix.

\Ssection{sine$_{\beta }$ random point fields/Dyson model in infinite dimensions.} \label{s:V1}

Let $ d = 1 $ and $ \sS = \mathbb{R}$. Recall ISDE \eqref{:11b} and take $ \beta = 1,2,4$. 
\begin{align}& \label{:V1a} 
 dX_t^i = dB_t^i + \frac{\beta }{2} 
 \lim_{r\to\infty }\sum_{|X_t^i - X_t^j |< r, \ j\not= i}^{\infty} 
 \frac{1}{X_t^i - X_t^j } dt \quad (i\in\mathbb{Z})
.\end{align}
Let $ \mu _{\sinbeta }$ be the sine$_{\beta }$ random point field 
\cite{Meh04,forrester} with $ \beta = 1,2,4$. 
$ \mu _{\sintwo } $ is the random point field on 
$ \mathbb{R}$ whose $ n $-point correlation function $ \rho_{\sintwo }^{n} $ 
%with respect to the Lebesgue measure 
is given by 
\begin{align}\label{:V1b}&
 \rho_{\sintwo }^{n} (\mathbf{x} ) = \det [\mathsf{K}_{\sintwo } (x_i,x_j)]\ijn 
. \end{align}
Here $ \mathsf{K}_{\sintwo } (x,y) = \sin \pi (x-y)/\pi (x-y)$ is the sine kernel. 
$ \musinone $ and $ \musinfour $ are also defined by correlation functions given by quaternion determinants \cite{Meh04}. 
% 
%The random point fields 
$ \mu _{\sinbeta }$ are solutions of the geometric differential equations \eqref{:51d} with 
$ a (x,\mathsf{y}) = 1 $ and 
\begin{align}\label{:V1c}&
 b (x,\mathsf{y}) = \frac{\beta}{2} \limi{r} \sum_{|x-y_i|<r} 
\frac{1}{x-y_i}
\quad \text{ in }L_{\mathrm{loc}}^1 (\mu _{\sinbeta }^{[1]})
.\end{align}
Here \lq\lq in $ L_{\mathsf{loc}}^1 (\mu _{\sinbeta }^{[1]})$'' means the convergence in 
 $ L^1(\Sr \ts \mathsf{S},\, \mu _{\sinbeta }^{[1]} )$ for all $ r \in \N $. 
Unlike the Ginibre random point field, \eqref{:V1c} is equivalent to 
\begin{align} \notag %\label{:V1d}
& b (x,\mathsf{y}) = \frac{\beta}{2} \limi{r} \sum_{|y_i|<r} \frac{1}{x-y_i}
\quad \text{ in }L_{\mathrm{loc}}^1 (\mu _{\sinbeta }^{[1]})
.\end{align}
We obtain the following. 
\begin{theorem}\label{l:V1} 
\thetag{1} 
The conclusions of \tref{l:UA} hold for $ \mu _{\sintwo } $. 
\\\thetag{2}
Let $ \beta = 1,4$. 
Let $ \mu _{\sinbeta ,\, \mathrm{Tail}}^{\aa }$ be defined as \eqref{:54a} for 
$ \mu _{\sinbeta }$. 
Then, for $ \mu _{\sinbeta }$-a.s.\! $ \aa $, 
the conclusions of \tref{l:QA} hold for 
$ \mu _{\sinbeta ,\, \mathrm{Tail}}^{\aa }$. 
\end{theorem}
\begin{remark}\label{r:V1}
When $ \beta = 2 $, the solution of ISDE \eqref{:V1a} 
 is called the Dyson model in infinite dimensions \cite{Spo87}. 
 The random point fields $ \mu _{\sinbeta }$ 
are constructed for all $ \beta >0 $ \cite{virag.1}. 
It is plausible that our method can be applied to this case. 
We also remark that Tsai \cite{tsai.14} solved 
ISDE \eqref{:V1a} for all $ \beta \ge 1 $ 
employing a different method. 
\end{remark}

To prove \tref{l:V1} we shall check the assumptions of \tref{l:UA} $ (\beta=2)$ 
and \tref{l:QA} $ (\beta=1,4)$ for $ \mu _{\sinbeta } $. Let $ \chin $ be as in \lref{l:U1}. 

\begin{lemma} \label{l:W2} Let $ \beta=1,2,4$. 
For each $ \nn \in \NNNthree $ the following hold. 
\\\thetag{1} 
The logarithmic derivative 
$ \dSine $ of $ \mu _{\sinbeta } $ exists in 
$ \LchiSine $. 
\\
\thetag{2} The logarithmic derivative $ \dSine $ has the expressions. 
\begin{align}\label{:W2a}&
 \dSine (x, \mathsf{s}) = \beta 
 \limi{r} \sum_{|x-s_i| < r } \frac{1}{x-s_i} \quad \text{ in }
 \LchiSine 
, \\ \label{:W2b}&
 \dSine (x, \mathsf{s}) = 
 \beta \limi{r} \sum_{|s_i| < r } \frac{1}{x-s_i}
\quad \text{ in } \LchiSine 
.\end{align}
\end{lemma}
\begin{proof}
\thetag{1} and \eqref{:W2a} follow from \cite[Theorem 82]{o.isde}. 
Set $ \Sr ^x = \{ s;|x-s| < r \} $ and let $ \Sr \ominus \Sr ^x $ be the symmetric difference of 
$ \Sr $ and $ \Sr ^x $. Then we have 
\begin{align}\label{:W2c}
& 
\limi{r} \sum_{s_i\in \Sr \ominus \Sr ^x } 
\frac{1}{x-s_i} = 0 \quad \text{ in } 
 \LlocSine 
\end{align}
because $ d = 1$ and one- and two-point correlation functions of $ \mu _{\sinbeta } $ are bounded. 
Hence, \eqref{:W2b} follows from \eqref{:W2a} and \eqref{:W2c}. 
\qed \end{proof}

Take $ \ellell = 1 $ in \tref{l:UA} and 
note that $ \sigma (x,\mathsf{s})= 1 $ and $ b(x,\mathsf{s}) = \frac{1}{2}\dSine (x,\mathsf{s})$. 
Let $ \dom _{\sinbeta } ^{[1]} $ be the domain of 
the Dirichlet form associated with the 1-labeled process. 
\begin{lemma} \label{l:W3} 
Let $ \beta=1,2,4$. Then $ \chin \dSine ,\, \chin \nabla _x \dSine 
 \in \dom _{\sinbeta } ^{[1]} $ for all $ \nnNNN_3 $. 
In particular, \As{\CONE} holds for $ \ellell = 1 $. 
\end{lemma}
\begin{proof} 
We only prove $ \chin \dSine 
\in \dom _{\sinbeta } ^{[1]} $ because 
$ \chin \nabla _x \dSine \in \dom _{\sinbeta } ^{[1]} $ 
can be proved similarly. 
We set $ \DDD [f] = \DDD [f,f]$ for simplicity. By definition, 
\begin{align}\label{:W3a}& 
 \E ^{\muSinb } (\chin \dSine ,\chin \dSine ) =
 \int_{\sS \ts \sSS } \frac{1}{2}|\nabla_x \chin \dSine |^2 + 
 \, \DDD [\chin \dSine ] d \muSinb 
.\end{align}
From \lref{l:U1} and \lref{l:W2} \thetag{1}, 
we deduce that 
\begin{align}\notag 
 \int_{\sS \ts \sSS } |\nabla_x \chin \dSine |^2 d \muSinb \le & 
2\int_{\Hann } \{ \chin ^2 |\nabla_x \dSine |^2 + |\nabla_x \chin |^2 |\dSine |^2 \}
d \muSinb 
\\< & \infty 
\notag %label{:W3b}
.\end{align}
From the Schwarz inequality and \lref{l:U1}, we deduce that 
\begin{align}\notag 
 \int_{\sS \ts \sSS } \, & \DDD [\chin \dSine ] 
d \muSinb 
% \\\notag 
 \le \, 2
 \int_{\sS \ts \sSS } \, 
 \chin ^2 \DDD [ \dSine ] + \DDD [\chin ] \, |\dSine |^2 \, 
 d \muSinb 
\\\notag 
\le \, & 2 \int_{\Hann } \, 
 \DDD [ \dSine ] + \cref{;83b} |\dSine |^2
 \, 
 d \muSinb 
 \\ \notag %\label{:W3c} 
\le \, & 2 \int_{ \Hann } 
 \frac{\beta^2}{2}\big( \sum_{i}%%{\frac{1}{\pnn } \le |x-s_i| } %
 \frac{1}{|x-s_i |^4} \big) + \cref{;83b} |\dSine |^2 \, d \muSinb 
 < \, \infty \ 
.\end{align}
Here the last line follows from a direct calculation and \lref{l:W2}. 
%We note that $ \chin \dSine \in L^2(\muSinb ) $ is obvious from \lref{l:W2}. 
Putting these inequalities into \eqref{:W3a}, 
 we obtain $ \chin \dSine \in \dom _{\sinbeta } ^{[1]} $ for all $ \nnNNN_3 $. 
\qed \end{proof}

\begin{lemma} \label{l:W1} 
$ \mu _{\sinbeta } $ satisfies \As{A1}--\As{A4} for $ \beta =1,2,4$. 
\end{lemma}
\begin{proof} 
\eqref{:51a} in \As{A1} follows from \cite[Theorem 82]{o.isde}. 
In {\cite[Theorem 2.2]{o.rm}}, it was proved that $ \mu _{\sinbeta } $ is 
a $ (0, -\beta \log |x-y|)$-quasi-Gibbs measure for $ \beta = 1,2,4$. 
This yields \As{A2}. 
 \As{A3} is immediate from \eqref{:V1b} for $ \beta=2$, and 
a similar determinantal expression of correlation functions in \cite{Meh04} 
for $ \beta=1,4 $. 
We next verify \As{A4}. 
\eqref{:92z} holds obviously. \eqref{:91a} is satisfied by \lref{l:91}. 
Hence we have \As{A4} \thetag{2}. 
Because $ \mu _{\sinbeta } $ is translation invariant, \As{A4} \thetag{1} holds 
We thus see \As{A4} holds. 
\qed \end{proof}

\medskip

\noindent 
{\em Proof of \tref{l:V1}. } We check the assumptions of \corref{l:QC}. 
\As{A1}--\As{A4} follow from \lref{l:W1}. 
Let $ \mathbf{a} = \{ \ak \} $ be as in \eqref{:95w}. Take $ \ak (r) = \qq r $. 
Then \As{\Bone} holds because $ \mu _{\sinbeta } $ is translation invariant. 
\As{\CONE} follows from \lref{l:W2} and \lref{l:W3}. 
From the Lebesgue convergence theorem, we obtain 
\begin{align}\notag %\label{:W4q}
&
\nabla_x b(x,\mathsf{s}) =
 \frac{1}{2}\nabla_x \dSine (x,\mathsf{s}) = 
 - \frac{\beta }{2} \sum_{i} \frac{1}{(x-s_i)^2} 
 \in L^{\infty } (\Han , \muSinb ) 
.\end{align}
Hence we can apply \lref{l:U3} to obtain \As{\CTWO}. 

Assume $ \beta = 2$. 
Then $ \mu _{\sintwo } $ satisfies \As{\muTT} because $ \mu _{\sintwo } $ 
is a determinantal random point field. Hence applying \tref{l:UB} we obtain \thetag{1}. 
Assume $ \beta = 1,4$. Then applying \tref{l:QB} we obtain \thetag{2}. 
\bbbbb

\Ssection{Ruelle's class potentials}\label{e:V2} 
Let $ \sS = \Rd $ and $ \Phi = 0$. 
Let $ \Psi $ be translation invariant, that is, $ \Psi (x,y) = \Psi (x-y,0)$ for all $ x,y \in \Rd $. 
We set $ \Psi (x) = \Psi (x,0)$. Then \eqref{:11a} becomes 
\begin{align}
\notag &% \label{:V2a} &
dX^i_t = dB^i_t - 
\frac{\beta }{2} \sum ^{\infty}_{j=1,j\ne i} 
\nabla \Psiz (X_t^i - X_t^j ) dt 
\quad (i\in \N )
. \end{align}

Assume that $ \Psiz $ is a Ruelle's class potential, smooth outside the origin. 
That is, $ \Psiz $ is super-stable and regular in the sense of Ruelle \cite{ruelle.2}. Here we say $ \Psiz $ is regular if there exists a non-negative decreasing function 
$ \map{\psi }{[0,\infty)}{[0,\infty)}$ and $ R_0 $ such that 
\begin{align}
\notag &% \label{:V2b}&
\Psiz (x) \ge - \psi (|x|) \quad \text{ for all } x ,
%\in \mathbb{R}^d , 
 \quad 
%\\ \notag &
\Psiz (x) \le \psi (|x|) \quad \text{ for all } 
 |x| \ge R_0 ,% \in \mathbb{R}^d ,
\\ \notag &
\int_0^{\infty} \psi (t)\, t^{d-1}dt < \infty 
.\end{align}
Let $ \mupsi $ be a canonical Gibbs measure with interaction $ \Psiz $. 
We do not a priori assume the translation invariance of $ \mupsi $. 
Instead, we assume a quantitative condition in \eqref{:V2c} below, 
which is obviously satisfied by the translation invariant canonical Gibbs measures. 
Let $ \rho ^m $ be the $ m $-point correlation function of $ \mupsi $. 
\begin{theorem} \label{l:V2} 
Let $ \sS = \Rd $ and $ \beta > 0$. 
Let $ \Psiz $ be an interaction potential of Ruelle's class smooth outside the origin. 
Assume that, for each $ \pp \in \N $, there exist positive constants 
$ \Ct \label{;29c} $ and $ \Ct \label{;29d} $ satisfying 
\begin{align} \label{:V2c}&
\sum_{r=1}^{\infty} 
\frac{\int_{\Sr }\rho ^1 (x)dx } {r^{ \cref{;29c} +1 }} < \infty , \quad 
\limsup_{r\to\infty } 
\frac{\int_{\Sr }\rho ^1 (x)dx } {r^{ \cref{;29c} }} < \infty 
,\\
\notag 
& 
 | \nabla \Psiz ( x )|,\ | \nabla ^2 \Psiz ( x )| \le
\frac{\cref{;29d} }{ (1+|x|)^{ \cref{;29c} }}
\quad \text{ for all $ x $ such that $ |x|\ge 1/\pp $}
.\end{align}
Assume either $ d \ge 2 $ or $ d=1$ with $ \mupsi $ such that 
there exist a positive constant $\Ct \label{;27E} $ and 
a positive function $ \map{h}{[0,\infty)}{[0,\infty]}$ 
depending on $ m , R \in \N $ such that 
\begin{align}\label{:V2d}& 
\int_{0\le t \le \cref{;27E}} \frac{1}{h(t)} dt = \infty 
,\\\notag & 
\rho ^m(x_1,\ldots,x_m) \le h (|x_i-x_j|) 
\text{ for all } x_i\not=x_j \in \sS _R 
.\end{align}
Then the conclusions of \tref{l:QA} hold for $ \mupsi $. 
\end{theorem}

We begin with the calculation of the logarithmic derivative. 
\begin{lemma} \label{l:W8} 
Assumption \As{A1} holds and the logarithmic derivative $ \dpsi $ is given by 
\begin{align}\label{:W8z}&\quad \quad 
\dpsi (x,\mathsf{y}) = - \beta \sum_i \nabla \Psiz (x-y_i)
\quad (\mathsf{y}=\sum_{i} \delta_{y_i})
.\end{align}
\end{lemma}
\begin{proof} 
This lemma is clear from the DLR equation. For the sake of completeness we give a proof. 
We suppose $ \mupsi (\mathsf{s}(\sS )=\infty ) = 1 $. 

Let $ \sSS _r^{[1],m} = \Sr \ts \SSrm $ for $ m \in \{ 0 \} \cup \N $, 
where $ \SSrm = \{ \mathsf{s}\in \sSS ; \mathsf{s}(\Sr ) = m \} $. 
Let $ \mu _{r,\eta }^{[1]} $ be a conditional measure of 
the 1-Campbell measure $ \mupsi ^{[1]}$ 
conditioned at $ \pirc (\mathsf{y}) = \pirc (\eta )$ for $ (x,\mathsf{y}) \in \sS \ts \sSS $. 
We normalize $ \mupsi ^{[1]}( \cdot \cap { \sSS _r^{[1],m} }) $ and we denote by 
 $ \sigma _{r,\eta }^{[1],m} $ the density function of 
$ \mu _{r,\eta }^{[1]} $ on $ \sSS _r^{[1],m} $. 
%We remark here that the logarithmic derivative is invariant under constant multiplication of the original measure. 
Then, by the DLR equation and the definitions of reduced Palm and Campbell measures we obtain 
%with a direct calculation, %we see that 
\begin{align} & \notag %\label{:W8a}
 \sigma _{r,\eta }^{[1],m} (x,\y ) = 
 \frac{1}{\mathcal{Z}_{r,\eta}^{m} } 
e^{- \beta \big\{\sum_{i=1}^m \Psiz (x-y_i) + \sum_{i<j}^m \Psiz (y_i-y_j) 
 + \sum_{\eta_k \in \Sr ^c} 
 \{ \Psiz (x-\eta_k ) + \sum_{i=1}^m \Psiz (y_i-\eta _k) \} \big\} }
,\end{align}
where $ \y =(y_1,\ldots,y_m) \in \Srm $, $ \eta = \sum_k \delta_{\eta _k}$, 
 and $ \mathcal{Z}_{r,\eta}^{m} $ is the normalizing constant. 
Then we see that 
\begin{align}\label{:W8b}&
\nabla_x \log \sigma _{r,\eta }^{[1],m} (x,\y ) = 
- \beta \Big\{\sum_{i=1}^{m} \nabla \Psiz (x-y_i) + 
 \sum_{\eta_k \in \Sr ^c} \nabla \Psiz (x-\eta _k) \Big\} 
.\end{align}
For $ \varphi \in C_0(\sS )\ot \dom _{\circ} $ such that 
$ \varphi (x,\mathsf{y})= 0 $ on $ \Sr ^c \ts \sSS $, we have that 
\begin{align}&\notag 
- \int _{ \sS \ts \sSS } \nabla_x \varphi (x,\mathsf{y}) 
\mupsi ^{[1]}(dxd\mathsf{y})
 \\ \notag & =
 - \sum_{m=0}^{\infty} \int _{ \Sr \ts \Sr ^{m} \ts \sSS } 
\{ \nabla_x \varphi (x,\sum_{i=1}^m \delta_{y_i}) \} \, 
 \sigma _{r,\eta }^{[1],m}(x,\y )\, dxd\y \, 
 \mupsi ^{[1]} \circ (\pirc )^{-1} (d\eta ) 
\\ \notag &
= 
 \sum_{m=0}^{\infty} \int _{ \Sr \ts \Sr ^{m} \ts \sSS } 
\varphi (x,\sum_{i=1}^m \delta_{y_i}) \, \{ 
\nabla_x \log \sigma _{r,\eta }^{[1],m}(x,\y ) \} \, 
 \sigma _{r,\eta }^{[1],m}(x,\y )\, dxd\y \, 
 \mupsi ^{[1]} \circ (\pirc ) ^{-1} (d\eta )\\ \notag &
= 
 \int _{ \sS \ts \sSS } 
\varphi (x,\mathsf{y}) 
\{ - \beta \sum_{i=1}^{\infty} \nabla \Psiz (x-y_i) \} \mupsi ^{[1]}(dxd\mathsf{y}) 
\quad \text{ by \eqref{:W8b}}
.\end{align}
From this, we obtain \eqref{:W8z}. 
\qed \end{proof}

\begin{lemma} \label{l:W9}
Assume that $ d=1 $ and also \eqref{:V2d}. Then \eqref{:91a} holds. 
\end{lemma}
\begin{proof}
We can prove \lref{l:W9} in a similar fashion as the proof of \cite[Theorem 2.1]{o.col}. 
We easily deduce from the argument of \cite{inu} that \cite[\thetag{4.5}]{o.col} holds under the assumptions 
\eqref{:V2d}, and the rest of the proof is the same as that of \cite[Theorem 2.1]{o.col}. 
\qed\end{proof}

\noindent {\em Proof of \tref{l:V2}.} 
We verify that $ \mupsi $ satisfies the assumptions 
\As{A1}--\As{A4}, \As{\Bone}, and \As{\CONE}--\As{\CTWO} with $ \ellell =1 $. 

\As{A1} follows from \lref{l:W8}. 
We obtain \As{A2} from the DLR equation and the assumption that 
$ \Psiz $ is smooth outside the origin. 
We deduce \As{A3} and \As{A4} from \eqref{:V2c}, \lref{l:W9}, 
and $ \partial \sS = \emptyset $. 
Taking $ a_{\qq }(r) = \qq \int_{\Sr }\rho^1(x)dx $, we obtain \As{\Bone}. 
From \eqref{:V2c} we can apply the Lebesgue convergence theorem to 
differentiate \eqref{:W8z} to obtain \As{\CONE}. 
\As{\CTWO} follow from \lref{l:U3} and \eqref{:V2c}. 
\bbbbb

\begin{remark}\label{r:V2} \thetag{1} 
Inukai \cite{inu} proved that the assumption \eqref{:V2d} is 
a necessary and sufficient conditions of the particles 
never to collide for finite particle systems. 
\\\thetag{2} One can easily generalize \tref{l:V2} even if a free potential $ \Phi $ exists. 
\\\thetag{3} 
For given potentials of Ruelle's class $ \Psiz $, 
there exist translation invariant grand canonical Gibbs measures 
associated with $ \Psiz $ such that 
the $ m $-point correlation function $ \rho ^m $ 
with respect to the Lebesgue measure satisfies $ \rho ^m (x_1,\ldots,x_m) \le \cref{;27Z}^m $
for all $ (x_1,\ldots,x_m) \in (\Rd )^m $ and $ m \in \mathbb{N}$ 
(see \cite{ruelle.2}). 
Here $ \Ct \label{;27Z}$ is a positive constant. 
\end{remark}

\section{Appendix: Tail decomposition.} \label{s:J}
%\subsection{}\label{s:X} 
This section proves the tail decomposition of $ \mu $ using \cite{geo}. %
The notation in this section are adjusted as much as possible with \cite{geo}. 
We begin by preparing several notions. 

\begin{definition}[probability kernels {\cite[13p]{geo}}]\label{dfn:Xa} 
Let $ (X,\mathcal{X} )$ and $ (Y,\mathcal{Y} )$ be measurable spaces. 
A function $ \map{\pi }{\mathcal{X}\ts Y }{[0,1]}$ is called a probability kernel from $ \mathcal{Y} $ to $ \mathcal{X} $ if 
\\
\thetag{i} $ \pi (\cdot | y)$ is a probability on 
$ (X,\mathcal{X} )$ for all $ y \in Y$. 
\\
\thetag{ii} $ \pi (A|\cdot )$ is $ \mathcal{Y} $-measurable for each $ A \in \mathcal{X} $. 
\end{definition}

For a measurable space $ *$, we denote by $ \mathcal{P}(*) $ 
the set of all probabilities on $ *$. 

Let $ (\Omega,\mathcal{F} )$ be a measurable space. 
Let $ \mathcal{A} \subset \mathcal{F} $ be a sub-$ \sigma$-field. 
Let $ \mathcal{P} $ be a non-empty subset of $ \mathcal{P}(\Omega,\mathcal{F} ) $. 
We set %$ \pi ^{\omega }(\cdot ) =\pi (\cdot | \omega ) $ and 
$ \Omega_{\mathcal{P} } = 
\{ \omega \in \Omega \, ;\, \pi (\cdot | \omega ) \in \mathcal{P} \}$. 
\begin{definition}
[$ (\mathcal{P},\mathcal{A})$-kernel {\cite[(7.21) Definition, 130p]{geo}}] 
\label{dfn:Xb}
A probability kernel $ \map{\pi}{\mathcal{F}\ts \Omega }{[0,1]}$ is called 
$ (\mathcal{P},\mathcal{A} )$-kernel if it satisfies the following: 
\\
\thetag{i} 
$ \pi (A|\cdot ) = \mu (A|\mathcal{A} ) $ 
$ \mu $-a.s.\! for all $ \mu \in \mathcal{P} $ and $ A \in \mathcal{F} $. \\
\thetag{ii} $ \Omega_{\mathcal{P} }\in \mathcal{A} $. \\
\thetag{iii} $ \mu ( \Omega_{\mathcal{P} }) =1 $ for all $ \mu \in \mathcal{P} $. 
\end{definition}
Let $ (\Omega,\mathcal{F} )$, $ \mathcal{A} \subset \mathcal{F} $, and 
$ \mathcal{P} $ be as in \dref{dfn:Xb}. 
Let 
\begin{align}
\notag &% \label{:X1z}&
\mathcal{P}_{\mathcal{A} } = \{ \mu \in \mathcal{P} \, ;\, 
\mu (A) \in \{ 0,1 \} \text{ for all } A \in \mathcal{A} 
\} 
.\end{align}
We set $ \mathfrak{e}({\mathcal{P}_{\mathcal{A} }})= 
\sigma [{e}_A; A \in \mathcal{F} ]$, where 
$ \map{{e}_A}{\mathcal{P}_{\mathcal{A} } }{[0,1]}$ 
such that $ {e}_A(\mu) = \mu (A)$. 
\begin{lemma}[{\cite[(7.22) Proposition, 130p]{geo}}]\label{l:Xc}
Assume that $ (\Omega,\mathcal{F} )$ has a countable determining class 
(also called a countable core in \cite{geo}). 
Suppose that there exists a $ (\mathcal{P},\mathcal{A} )$-kernel $ \pi $. 
Then the following holds:\\
\thetag{1} $ \mathcal{P}_{\mathcal{A} } \not= \emptyset $.\\
\thetag{2} For each $ \mu \in \mathcal{P} $ there exists a unique $ w \in 
\mathcal{P}(\mathcal{P}_{\mathcal{A} } ,
\mathfrak{e}(\mathcal{P}_{\mathcal{A} } ) ) $ such that 
\begin{align}
\notag &% \label{:X1a}&
\mu = \int_{\mathcal{P}_{\mathcal{A} } } \nu w (d\nu )
.\end{align}
Furthermore, $ w $ is given by 
$ w (M) = \mu (\{ \omega \in \Omega ; \pi (\cdot | \omega ) \in M \} ) $ for 
$ M \in \mathfrak{e}(\mathcal{P}_{\mathcal{A} } )$. 
\\\thetag{3} $ \mathcal{P}_{\mathcal{A} } $ satisfies the following:
\begin{align}\label{:X1b}&
\mathcal{P}_{\mathcal{A} } = \{ \mu \in \mathcal{P} \, ;\, 
\mu (\{\omega ;\pi (\cdot | \omega ) = \mu \} ) = 1 
\} 
.\end{align}
\thetag{4} 
$ \mu (\{\omega ;\pi (\cdot | \omega ) \in \mathcal{P}_{\mathcal{A} } \} ) = 1 $ 
for all $ \mu \in \mathcal{P} $. 
\end{lemma}
\begin{proof}
\thetag{1} and \thetag{2} follow from \cite[(7.22) Proposition]{geo}. 
\thetag{3} follows from \cite[(7.23)]{geo}. \thetag{4} follows from 2) of the proof of \cite[(7.22) Proposition]{geo} (see 12-17 lines, 131p in \cite{geo}). 
\qed \end{proof}

We introduce the notion of specifications $ \gamma $. We slightly modify it 
according to the present situation. 
Let $ X=Y$ and $ \mathcal{Y}\subset \mathcal{X} $ in \dref{dfn:Xa}.
Then a probability kernel $ \pi $ is called proper if $ \pi (A\cap B| \cdot ) = \pi(A | \cdot ) 1_{B}$ 
for all $ A \in \mathcal{X} $ and $ B \in \mathcal{Y} $. 

\begin{definition}[specification {\cite[16p]{geo}}]\label{dfn:Xc} 
A family $ \gamma =\{ \gamma_r \}_{r\in\N } $ of proper probability kernels $ \gamma_r $ 
from $ \sigma [\pi _r^c]$ to $ \mathcal{B}(\sSS ) $ is called a specification if 
it satisfies the consistency condition $ \gamma_s\gamma_r = \gamma_s $ 
when $ r\le s \in \N $, where for $ A \in \mathcal{B}(\sSS ) $ 
\begin{align} \notag %\label{:10z}&
\gamma_s\gamma_r (A,\mathsf{s}) = 
\int _{\sSS } \gamma_s (d\mathsf{t},\mathsf{s}) \gamma_r (A,\mathsf{t}) 
.\end{align}
The random point fields in the set 
$$ \mathcal{G}(\gamma ) = 
\{ \upsilon \in \mathcal{P}(\sSS ,\mathcal{B}(\sSS ) ) \, ;\, 
\upsilon (A | \sigma [\pi _r^c]) = \gamma_r (A ,\, \cdot ) \ 
\text{$ \upsilon $-a.s.\! for all $ A \in \mathcal{B}(\sSS ) $ and $ r \in \N $}
\} 
$$
are said to be specified by $ \gamma $. 
\end{definition}

With these preparations we recall the tail decomposition given by Georgii \cite{geo}. 
Let $ \mu $ and $ \mut = \mu (\, \cdot \, | \TS )(\mathsf{a})$ be as in 
Theorems \ref{l:5A}--\ref{l:5B}. 

\begin{lemma} \label{l:X2}
Assume \As{A2}. 
Let $ \mathsf{H}\in \mathcal{B}(\sSS ) $ such that $ \mu (\mathsf{H}) = 1 $. 
%be a subset of $ \SSsde $ satisfying \As{\sIn} and \As{\iFc}. 
There then exists a subset of $ \mathsf{H} $
denoted by the same symbol $ \HHz $, 
 such that $ \mu (\HHz ) = 1 $ and that for all $ \mathsf{a} \in \HHz $ 
\begin{align}\label{:X2a}&
 \mut (\A ) \in \{ 0,1 \} \text{ for all }\A \in \TS 
,\\\label{:X2b}&
 \mut (\{ \bb \in \sSS ; \mut = \mu _{\mathrm{Tail}}^{\bb } \} ) = 1 
,\\\label{:X2d}&
\text{$ \mut $ and $ \mu_{\mathrm{Tail}}^{\mathsf{b}}$ 
are mutually singular on $ \TS $ 
if $ \mut \not=\mu_{\mathrm{Tail}}^{\mathsf{b}}$}
.\end{align}
\end{lemma}

\begin{proof}
We use \lref{l:Xc} for the proof of \lref{l:X2}. 
Let $ (\Omega ,\mathcal{F})=(\sSS ,\mathcal{B}(\sSS ) ) $. 
Then 
$ (\Omega ,\mathcal{F}) $ 
% is countably determined and also has 
% a countable-determinating class \cite{bill} 
% (also called a countable core in \cite{geo}). 
has a countable determining class because $ \sSS $ is a Polish space. 

Let $ \gamma_r (\cdot , \mathsf{a}) = \mu (\cdot | \sigma [\pi_r^c])(\mathsf{a})$, 
%% for each $ r \in \N$, 
where $ \mu (\cdot | \sigma [\pi_r^c])(\mathsf{a})$ 
is the regular conditional probability of $ \mu $. 
By \As{A2}, we can take a version of $ \mu (\cdot | \sigma [\pi_r^c])(\mathsf{a}) $ in such a way that 
$ \gamma =\{ \gamma_r \}_{r\in\N} $ 
becomes a specification and $ \mu \in \mathcal{G}(\gamma ) $. 
We set 
\begin{align}
\notag &% \label{:X2c}&
 \pi (\cdot | \mathsf{a}) = \mut (\cdot ), \quad 
\mathcal{P} = \mathcal{G}(\gamma ) ,\quad \mathcal{A}= \TS 
.\end{align}
From \cite[(7.25) Proposition, 132p]{geo} 
we see that $ \pi (\cdot | \mathsf{a}) = \mut (\cdot ) $ 
becomes a $ (\mathcal{P},\mathcal{A} )$-kernel. 
Let $ \Omega_1 = \{ \mathsf{a}; \pi (\cdot | \mathsf{a}) \in \mathcal{G}(\gamma ) \} $. 
It is also proved in the proof of \cite[(7.25) Proposition, 132p]{geo} that 
$ \mu (\Omega_1) = 1 $ for all $ \mu \in \mathcal{G} (\gamma )$. 
Then \eqref{:X2a} follows from \thetag{4} of \lref{l:Xc}, and 
\eqref{:X2b} follows from \eqref{:X1b}; moreover, 
\eqref{:X2d} follows from \cite[(7.7) Theorem (d), 118p]{geo}. 
\qed \end{proof}

\begin{remark}\label{r:X2} 
In \cite{geo}, a parameter is taken to be a countable infinite set and spin is a general measurable set. 
The argument in \cite{geo} is valid in the present paper by considering parameters consisting of a countable partition of the continuous set $ \sS $ and taking a spin as a configuration space on each element of the partition. %
%We do not know whether the regular conditional probability $ \gamma_r (\cdot , \mathsf{a}) = \mu (\cdot | \sigma [\pi_r^c])(\mathsf{a})$ becomes a specification without \As{A2}. 
\end{remark}

\noindent 
\textit{Acknowledgements: }
H.O. is supported in part by a Grant-in-Aid for Scientific Research
(Grant Nos. 16K13764, 16H02149, 16H06338, and KIBAN-A, No. 24244010)
from the Japan Society for the Promotion of Science. 
H.T. is supported in part by a Grant-in-Aid for Scientific Research
(KIBAN-C, No. 15K04916, Scientific Research (B), No. 19H01793)
from the Japan Society for the Promotion of Science. 

%\bibliographystyle{spbasic} % basic style, author-year citations
%\bibliographystyle{spmpsci} % mathematics and physical sciences
%\bibliographystyle{spphys} % APS-like style for physics
%\bibliography{} % name your BibTeX data base


% ///
\end{document}